\documentclass[leqno,12pt]{article}
\usepackage{amsmath, amsthm, mathrsfs}
\usepackage{amsfonts}
\usepackage{amssymb}
\usepackage{graphicx}
\usepackage{esint,verbatim}
\usepackage[monochrome]{color}
\usepackage{hyperref}
\usepackage{empheq}
\usepackage{ifpdf} 
\usepackage{geometry}
\newtheorem{theorem}{Theorem}[section]
\newtheorem{corollary}[theorem]{Corollary}
\newtheorem{lemma}[theorem]{Lemma}
\newtheorem{proposition}{Proposition}[section]

\newtheorem{open}{Open Problem}

\newcommand{\re}{\mathbb{R}}
\newcommand{\ren}{{\mathbb{R}^N}}

\newcommand{\Hn}{\mathbb{H}^N}
\newcommand{\SP}{\mathbb{S}}
\newcommand{\ve}{\varepsilon}

\newcommand{\var}{\varphi}

\newcommand{\dx}{\,{\rm d}x}

\newcommand{\Sh}{\mathrm{Sh}}
\newcommand{\vp}{\varphi}
\newcommand{\ifl}{\Delta_\infty^s}

\numberwithin{equation}{section}
\def\qed{\,\unskip\kern 6pt \penalty 500
\raise -2pt\hbox{\vrule \vbox to8pt{\hrule width 6pt
\vfill\hrule}\vrule}\par}
\newcommand{\bc}{\color{blue}}

\definecolor{ao}{rgb}{0.0, 0.0, 1.0}

\newcommand{\nc}{\normalcolor}

\definecolor{darkblue}{rgb}{0.05, .05, .65}
\definecolor{darkgreen}{rgb}{0.1, .65, .1}
\definecolor{darkred}{rgb}{0.8,0,0}

\begin{document}

\title{\bf A survey on mass conservation, self-similarity and related topics \\ in nonlinear diffusion
}

\author{ \sc Juan Luis V\'azquez
\footnote{Depto. de Matem\'{a}ticas. Univ. Aut\'onoma de Madrid;  Real Academia de Ciencias, Spain}
}


\vskip -3cm

\maketitle

\begin{abstract}

We examine the validity of the principle of mass conservation for solutions of some typical equations in the theory of nonlinear diffusion, including equations in standard differential form and also  their fractional counterparts.

In Part 1,  consisting of the first 9 sections, we use as main examples the heat equation, the porous medium equation and the $p$-Laplacian equation. Though these equations have the form of conservation laws, it happens that in some ranges of exponents  the solutions posed in the whole Euclidean space actually lose mass in time. From the start we pay attention to the close connection between the validity of mass conservation and the existence of finite-mass self-similar solutions, as well as their role in the asymptotic behaviour of more general classes of solutions.  Describing the extent of this connection is the common thread throughout the manuscript.  When mass conservation does not hold, we are led to examine the situation when it \nc is replaced by its extreme alternative, extinction in finite time,  a very surprising fact.
The equations offer difficult borderline cases in the presence of  some critical parameters that we identify. New results are proved to complete the picture; we establish mass conservation or its failure in those pending cases. We also explain the occasional disappearance of the fundamental solutions in a very graphical way, as ``escape to infinity''.
In this part of the manuscript, we aim to provide a comprehensive overview of the types of behaviour to be expected in these basic models. Additionally, Section \ref{sec.dnle} on the Doubly Nonlinear Equation complements this foundational block.\nc

The next sections extend the detailed study to other models that occupy a relevant role in the current literature. Thus, the 3 sections of Part 2 are devoted to the discussion of mass conservation for some fractional nonlinear diffusion equations, where the situation is surveyed and a number of  new \nc theorems are proved.  This is followed by the analogous study of Scalar Conservation Laws in Part 3, with surprising similarities. We conclude with a long review of related equations and topics in Part 4. The survey contains an extensive bibliography.

Summing up, the document aims at surveying important topics in nonlinear diffusion related to mass conservation. At the time we solve a number of open problems on key issues and point out interesting connections and new directions.

\end{abstract}

\vskip .2cm

\vskip 0.2cm


\


\medskip

\vskip 1cm

\



\setcounter{tocdepth}{2}
\tableofcontents

\medskip

\newpage

\

\setcounter{section}{1}

\section*{Introduction}\label{se.intro}
\addcontentsline{toc}{section}{Introduction}

In physics and chemistry, the law of conservation of mass states that, for any system closed to all transfers of matter and energy, the amount of mass is conserved over time. The law implies that mass can neither be created nor destroyed in those systems, although it may be rearranged in space, or the substances associated with it may be changed in form. Indeed, the law is credited to the French chemist Antoine de Lavoisier who in 1789 wrote a treaty \cite{Lavois1798} explaining that mass is neither created nor destroyed in chemical reactions. He did that by carefully weighing the system before and after the reaction, and checking that the total mass had not changed. The concept of mass conservation is widely used in many fields such as chemistry, mechanics, and fluid dynamics, and of course in several areas of mathematics,  specially in the analysis of PDEs. Notice that we talk about closed systems, mass may escape or enter through the borders of an open system.

In the classical fields of fluid mechanics and general continuum mechanics, the mass conservation law is usually written in the form
	\begin{equation}\label{cont.eq}
 {\displaystyle \partial_t  \rho +\mbox{\rm div}_x\, (\rho \mathbf {v} )=0,} \tag{CE}
	\end{equation}
which is known as the {\sl continuity equation}, \cite{ChorinM2013, Gurtin, LandL}. Here,  ${\textstyle \rho } $ is the density (mass per unit volume) of the substance under study, the symbol $\mbox{div}_x$ denotes the divergence operator acting on the space variables, and $ {\textstyle \mathbf {v} }$ is the velocity field. This is the differential formulation. The law admits a corresponding integral version that reflects the following interpretation: for a given closed surface in the system, the change in time of the mass enclosed by the surface is equal to the mass that traverses the surface, positive if matter goes in and negative if matter goes out (this is usually called the mass flux through the boundary). For the whole isolated system, the conservation law implies that the total mass $ {\textstyle M}$ -- the sum of the masses of all components in the system -- does not change in time, i.e.,
	\begin{equation}\label{mss.conserv}\tag{IMC}
{\displaystyle {\frac {dM}{dt}}={\frac {d}{dt}}\int \rho \, dV=0},
\end{equation}
where $dV$ represents the volume element within the domain occupied by  such system. This short introduction follows the standard literature, like the Wikipedia article \cite{CM Wikipedia} to which we refer for further description and historical information. Note that the conservation of mass does not hold in processes where mass is transformed into energy, as occurs in various areas of physics, such as nuclear and relativistic processes. This possibility will not be considered here\nc.

The text we present below contains a survey of mass conservation and its relation to self-similarity, as found in the mathematical treatment of nonlinear diffusion by PDE methods, together with number of new insights and results. It builds on  material contained in our monographs \cite{Vazquez2006, VazquezPME2007}, the lecture notes \cite{VazNotesHE} and \cite{VazTNLD2016}, as well as the literature that is mentioned in those references. It contains many new developments and results.

A number of monographs have dealt since the 1970s with general aspects of nonlinear diffusion equations. In the parabolic setting  we  may mention \cite{Crank75,  DaskaBk,  DiBenedetto93, Ni1993, WuZY01}. The authors of \cite{WuZY01} say: ``The appearance of degeneracy or singularity makes the study of nonlinear diffusion equations more involved and challenging. Many new ideas and methods have been developed to overcome the special difficulties caused by the degeneracy and singularity, which enrich the theory of partial differential equations''.

\subsection{Conservation of mass in Diffusion}\label{sec.cm.diff}

 In this memoir \nc we have made an effort to provide an updated perspective on key topics in nonlinear diffusion,
 with the aim of helping interested readers navigate the expanding body of related literature.
The notes are united by mass conservation as leitmotif.  The existence and role of special self-similar solutions will be a second focus of attention. \nc We present numerous new results essential for a comprehensive and systematic exposition, along with diverse ideas, connections, and related directions.

To start with, the property of conservation of the total mass  (often abridged here as MC, mass conservation) seems quite natural in diffusion processes posed in the whole space. Actually, seen in probabilistic terms in the discrete case, such processes  describe how particles migrate from site to site  to fill available space, without creation of new particles or disappearance of some of them (the paradigm being random walk models). Conservation of mass amounts then to counting the total number of particles. In the continuous setting,  diffusion processes describe likewise how a density function spreads to fill the available continuous space without creation or disappearance of the underlying substance; the basic idea is the same but the counting process may get complicated for reasons we show below. Posed in analytical terms, if \ $u(x,t)$ \ is the probability density or concentration of a continuous medium, we want to prove that the quantity $M(t)=\int u(x,t)\,dx$ is an invariant of the evolution,  i.e., \nc $M(t)=$ constant in time.

The property also applies to particles or densities confined in bounded domains when  reflecting conditions are imposed on the boundary (Neumann conditions), but it fails when homogeneous Dirichlet conditions are imposed, since this amounts in probabilistic terms to ``killing particles upon exiting the domain'' (in the homogeneous case). All of this is known to apply to the classical heat equation \ $\partial_t u=\Delta u$ (HE for short), and also to uniformly parabolic equations in divergence form, like those written  as $\partial_t u=\mbox{div\,}(D\,\nabla u)$, where  $D=D(x,t)>0$ \nc is the diffusion coefficient.
In good situations this coefficient is bounded from above  and also from below by  a positive constant (i.e., we have uniformly parabolic equations), but most of our main models do not accept these severe restrictions. The equation becomes degenerate when the coefficient $D$ may approach zero somewhere, or else it becomes singular when $D$ grows to infinity somewhere. Degeneracy and singularity will be a common occurrence in our models and the origin of significantly new concepts and results\nc.

The primary focus of the work is devoted to examine the question of MC in the context of the typical equations of nonlinear diffusion. These equations form a topic that has been extensively studied since the last part of the 20th century, and it is known that, surprisingly, MC does not always hold. Moreover, we will show that its validity is exactly tied to very relevant facts of the nonlinear diffusion theory, namely the existence and properties of special self-similar solutions.  We will devote a first part of the work to so-called standard models of nonlinear diffusion (down to Section \ref{sec.dnle}), and a second part to the models called nonlinear fractional diffusion  in  Sections \ref{sec.nonlinfrac} to \ref{sec.loss.vf.fpme}\nc.
This is followed by a study of 1D scalar conservation laws in Part 3, following the same  lines of what was shown to work in nonlinear diffusion. We conclude with a long review of related equations and topics in  Part 4. It includes reviews of many diffusion problems, with suggestions for further work.

 The study of MC in this work puts much emphasis on the \nc close connection between self-similarity and mass conservation that holds in our basic diffusion models. This serves as a unifying theme throughout a large part of the manuscript. Indeed, one could say that for our equations, ``mass conservation and self-similarity go hand in hand.'' Together, they give rise to fundamental solutions with distinct patterns - beautiful mathematical objects with far-reaching consequences and applications.
This still vague statement  will be carefully documented in our three main cases of linear and nonlinear diffusion of differential type, as well as in the case of first-order scalar conservation laws. The topic  offers a number of difficulties and changes of behaviour that have been carefully investigated over  the last decades. \bc We provide here a fairly comprehensive account from the intended perspective, and provide detailed concepts, results and technical proofs as needed\nc.

A  further issue that turns out to be connected with MC and self-similarity is the asymptotic behaviour of finite-mass solutions for large times. The full connection we show  is a quite general phenomenon that can have application in other areas, as clearly shown in Parts 3 and 4 of this work.

As a conclusion, the different models we have considered produce a collection of  self-similar exponents and profiles (patterns). The way they operate  offers the reader of the present text very illustrative information on the nonlinear theories we discuss, and that information is  both qualitative and quantitative. This applies to pattern formation and asymptotic convergence, as well as a number of interesting functional inequalities that appear repeatedly in the PDE theory. At the time, the text offers a tour of themes, concepts and techniques that are very representative of present-day nonlinear diffusive PDEs. \bc In our effort to clarify  and complete the theory many new results are established and a number of open problems are clearly stated\nc.

\newpage

\subsection{General concepts, notations and comments}\label{ss.gencomm}

 \noindent {\bf Equations.} The main equations we consider have parameters that are also easily identifiable. Thus, in Section \ref{sec.nld} we introduce the two main models of Nonlinear Diffusion, the porous medium equation (PME): \ $\partial_t u=\Delta (|u|^{m-1}u)$,
and the $p$-Laplacian equation  (PLE): \ $\partial_t u =\Delta_p(u):= \mbox{div} (|\nabla u|^{p-2}\nabla u)$.
Unless explicit mention to the contrary, we assume the values of the exponents $0<m<\infty$ in porous medium flows, and  $1<p<\infty$ in $p$-Laplacian flows. Then, the limit cases $m=0$ and $p=1$ are called limit cases and will be treated with due care. \bc Also $m=\infty$ and $p= \infty$ are  mentioned an the salient features explained\nc. In the fractional models  of Sections \ref{sec.nonlinfrac} to \ref{sec.loss.vf.fpme} we consider fractional exponents in the range $0<s<1$.
The particular subrange of $m$, $p$, or $s$ discussed in any  given section will be always carefully announced.  By a universal constant we mean a constant $c>0$ that may depend on the parameters of the theories, typically $m,p,s$ and the space dimension $N$. If we need to stress that it depends on $m$ we write $c(m)$, if there are several constants in a row we write $c_i(m)$. When the parameters are not mentioned, they  must have  been identifiable from the context or previous statements\nc .

 As a general rule, our equations will be posed on the Euclidean space, $x\in \ren$, and $N\ge 1$ denotes the space dimension.
This is a suitable context for Mass Conservation and self-similarity. There is a huge literature for problems posed in bounded domains but the flavor and results are quite different and will not be tackled here\nc. The same  applies to parabolic equations \nc posed on Riemannian manifolds. There, a control of the curvature is also needed and plays a role in the theory. For this topic \nc  cf. \cite{Gri2009}. We briefly touch this subject in Subsection \ref{ssec.manif}.

\medskip

\noindent {\bf Classes of Solutions.}
The nonlinear models we deal with do not have classical solutions as a general rule, due to the degenerate or singular parabolic character induced by the nonlinearity/-ies. Hence. we  need to select a suitable concept of generalized solution that guarantees existence and uniqueness of a large enough class of the solutions of the evolution problems (and  sometimes stationary ones), allowing at the same time  for useful functional and computational properties. The concept of distributional solution is very useful in the analysis of PDEs, but it is often the case that it does not guarantee by itself the uniqueness of solutions of the nonlinear problems we treat. The same happens with the closely related concept of weak solution.

There are several approaches in the literature to address this problem, notably those based on re-formulating the solution of the PDEs under consideration as the solution of an abstract ODE problem of the form
	\begin{equation}\label{abs.ode}
 u_t+A(u)=0,
	\end{equation}
 with initial value $u(0)=u_0$. We then look for solutions $u=u(t)$ that take values in a Hilbert, Banach or metric space $X$, and $A$ is a possibly nonlinear operator $A$ defined in $ D(A)\subset X $ with values in $X$ or the dual  $X'$ and having  the required properties.  Typically, for every admissible initial function $u_0$  the result is a curve $u(t)\in C([0,T): X)$ and the set of such functions generates a semigroup $S_t: u_0\mapsto u(t)$, $t>0$.
 There are several popular ways of implementing  this program, notably

 (i) the theory of maximal monotone operators which generate semigroups of contractions in Hilbert Spaces, cf. Brezis \cite{BrezisBkOMM73} and Komura \cite{Komura}.

 (ii) the Crandall-Liggett Semigroup Generation Theorem (CL SGT) that works in Banach spaces, see the original \cite{CrandallLiggett71}; see also  \cite{BeTh72, Barbu1976, Ev78, Veron79}, and the discussion in Chapter 10 of the monograph \cite{VazquezPME2007},  and

 (iii) the theory of Gradient Flows in metric spaces for which the standard reference is the book by Ambrosio et al. \cite{AmGiSav05}, see also \cite{Clem09, Santamb2017}.

Notice that the theory of abstract semigroups  was first established for linear operators, cf. \cite{HillePh57, EN95, Gol85, Paz83}, but it is now firmly established for nonlinear operators as well.
 In what follows, we will generally adopt option (ii). This choice is motivated by the properties of the nonlinear operators in our model equations and by the fact that mass conservation naturally takes place in the space of integrable functions $L^1(\ren)$, which is not a Hilbert space. The change of emphasis from $L^2$ to $L^1$ as a preferred functional space emerged in the 1970s as a natural setting for some nonlinear elliptic  and parabolic problems, see also \cite{BS73, BBC75}. Actually,  the Crandall-Liggett SGT works perfectly on Banach spaces $X$ on the condition that the (possibly nonlinear) operator $A$ is $m$-accretive in $X$.

 By means of the SGT we establish the existence and uniqueness of suitable abstract solutions for the evolution problems we consider.
 Indeed, such solutions form a contraction semigroup, a property that plays a big role in our studies. Indeed, the key idea is to solve the abstract ODE \eqref{abs.ode} via the method of Implicit Time Discretization. The simplest example of such  $m$-accretive operators is the signed Laplacian, $A=-\Delta$ acting in $L^1(\ren)$. Briefly, the SGT  produces a concept of solution, $u(t)\in C([0,T): X)$  so-called \sl mild solution\rm, which is in principle a very abstract concept of solution. However, in every case we study it has been proved that most mild solutions are distributional solutions with enough regularity so that the calculus operations we perform are justified after suitable approximations. Exceptions in some limit cases are important for the sake of completeness of the theory and they are duly explained.

We remind the reader that the problems of existence, uniqueness, and basic theory are not our main concern and will be explained only as far as they are important to lay a firm foundation for the main  theory and to understand the arguments of the text. \bc A high level of technical detail is enforced in the main part of the work, while in  the last parts the interest is mainly to show examples of this theory working in new fields and a less comprehensive style is adopted\nc.
We have tried to give appropriate references  at all times\nc.

\medskip

 \noindent {\bc \bf Nonnegative and signed solutions.} In light of the most common applications of these PDEs, we confine the majority of our analysis to the class of nonnegative solutions. This is quite natural for instance when we talk about mass, probability, or density. However, the existence and uniqueness theory of these PDEs can be done in the setting of solutions with both signs.  We will deal many times with the version of our results for two-signed solutions as an extension of the results in the nonnegative setting, in some cases we do it only  briefly to avoid diverting from the main topic. There is an important  remark about terminology: for  two-signed functions  $f(x)$  it does not make sense to give  the name ``mass'' to the first integral of $f(x)$ over the domain of reference, $\int f(x)\,dx$. So the results about conservation of mass for nonnegative solutions become when dealing with signed solutions results about ``conservation of the first spatial integral'', never about  the $L^1$-norm $\int |f(x)|\,dx$. When we do not enter into details for signed solutions, the interested reader is referred to the relevant literature, or asked to fill the gaps if they are easy\nc.

\medskip

\noindent {\bf Notations.} We sometimes write $u(t)$ instead of $u(x,t)$ for convenience, when we want to stress that $u$ is a function of $x$ with parameter $t$.
In the sequel we often use the notation $a^m$ to mean the signed power $|a|^{m-1}a$ of a quantity $a\in \re$, $a\ne 0$.  All of this is done for brevity and only when no confusion is to be feared.

\

\noindent {\bf  \large Some final remarks}

There is a close relation between nonlocal evolution problems and probability theory, as we have mentioned above. Indeed, in that view of mathematics, diffusion processes are just a distinguished class of stochastic process, and the heat equation is another aspect of the most popular stochastic process, the Brownian motion or Wiener process, \cite{KarShreve1988}.  When one looks at a L\'evy process,  the nonlocal operator that appears naturally is a fractional power of the Laplacian, see \cite{Apple2009, Bert96}. Such a  reference to the fruitful interplay between nonlocal partial differential equations and probability is a big topic and must be covered by other experts. The PDE approach is visited below in Sections \ref{sec.nonlinfrac} and  following.

We point out that mass conservation is closely related to the very active field of optimal transportation, that deals with evolution of constant-mass functions, \nc for which we refer to the monographs \cite{Villani1, Villani2}, and recently \cite{FigGlaudo, Santamb15}. Prominent tools of this trade, like the Wasserstein distances, defined between nonnegative solutions of the same total mass, appear in the analysis of the topics we will study, cf. \cite{Otto2001, AmGiSav05, CarrMcV06}. \bc We will not cover such topic  as it is has an extensive  development of its own, is not central for the  results we have in mind, and we wish to remain focused on the main scope of this work\nc.

\bc Finally, the study of special solutions of self-similar form is closely tied to the study of travelling waves in different forms of nonlinear diffusion. We will also avoid this important study, we refer to a well-known book \cite{GilKer04}\nc.

\vskip 1cm

\newpage


{\bf \huge Part 1}


\section{Models of nonlinear diffusion (NLD)}\label{sec.nld}

As we have said, the law of MC has been tested in many  models, in particular in the diffusion models having finite propagation and  passive free boundaries, a general reference book can be \cite{FdmnBk}. A rather thorough examination of the issue has been done in the case of two of the most popular models, the  Porous Medium Equation
	\begin{equation}\label{pme1}
\partial_t u=\Delta (|u|^{m-1}u), \tag{PME}
	\end{equation}
and the $p$-Laplacian Evolution Equation:
	\begin{equation}\label{pleq1}
\partial_t u =\Delta_p(u):= \mbox{div} (|\nabla u|^{p-2}\nabla u). \tag{PLE}
	\end{equation}
 Both are  outstanding in the studies of nonlinear diffusion. The relative simplicity of the equations combined with the rich behaviour that  has been discovered, makes them a convenient foundation to advance in the more complex models currently found in theory and many areas of application. The simple structure allowed  to obtain a rather complete knowledge of the main concepts, their properties and their links. \nc We will devote our attention in this work to a parallel study of both of them  due to their strong self-similarity properties that make the whole theory quite clear and computable. We take also into account historical reasons as documented in \cite{Vazquez2006, VazquezPME2007}. Note that the PME becomes the HE if the exponent is $m=1$, and so does the PLE when $p=2$.  As we have mentioned above, our study is concerned with  real-valued solutions $ u(x, t) $ defined for  $x\in \ren$\nc.

The PME is very easily written in the conservation form \eqref{cont.eq} (see Introduction) by introducing the velocity vector ${\bf v}=-mu^{m-2}\nabla u$, thus the velocity is density dependent. In turn, the velocity  can be written as (minus) the gradient of a potential, also called pressure,
	\begin{equation}\label{pme.vel}
{\bf v}=-\nabla U, \quad \mbox{with } \ U=(m/(m-1)) \, u^{m-1}.
	\end{equation}
The formula is a version of Darcy's law, well-known in Fluid Mechanics\nc. This idea played a prominent role in the theory of the PME, mainly in the study of free boundaries, cf. \cite{Aron1986, CaffFried80, Pel81} and many later references, see extensive details in the monograph \cite{VazquezPME2007}. A similar approach for the PLE is valid but not so immediate, the pressure formulation was studied in \cite{EV88, EV90}, see Section \ref{sec.mc.ple}.

For both equations there is  a physically motivated class of finite-mass solutions, and it is relatively easy to prove that  mass is conserved for superlinear exponents, $m>1$ and $p>2$ respectively. Such ranges of exponents correspond to what is called \sl slow diffusion\rm, which gives rise to the existence of free boundaries, that confine the mass. This leads in an natural way to the validity of MC. For $m=1$ in PME or $p=2$ in PLE the equations reduce to the linear case, in which case MC is  a classical result, see \cite{EvansBk, Vazquez2006, VazquezPME2007}. So far so good.

However, the MC question becomes delicate (and quite interesting) for nonlinearities with smaller powers, i.e., $0<m<1$ and $1<p<2$ respectively, that correspond to a relevant physical phenomenon called \sl fast diffusion\rm. Such diffusion has received a plethora of studies and applications, \cite{BeHoll78, BeHoll90, DaskaBk, Diaz2, HerrPierre1984, Ki93, LionsTos97, Vazquez2006, VazTNLD2016} and the relevant theory grows by the moment. The intensity of the fast process depends on the exponents $m$ and $p$, and conservation of mass may break down if $m$ and $p$ are small enough, in a range called  \sl very fast diffusion \rm \bc (VFDE)\nc. The precise borderline case is well known: it is $m_c=(N-2)/N$ for the PME (we then assume that $N\ge 3$) and $p_c=2N/(N+1)$ for the PLE (with $N\ge 2$). The conservation of mass in those critical cases posed some nontrivial open problems that will be solved  here in the way that \bc  we carefully explain in Section \ref{sec.crit.exp}\nc. Extra attention is given to the special cases of  planar logarithmic diffusion \ (PME, case $m=0$, $N=2$, see Subsection \ref{sec.ricci}) \nc and one-dimensional total variation flow  \ (PLE with $p=1$, $N=1$, see Subsection \ref{sec.n1pc1})\nc. The detailed analysis of the interplay of fast diffusion with mass conservation is carefully explained in  this work, with special attention to the phenomenon of gradual loss of mass.

To start the study, we will introduce next a number of concepts and connections regarding the equations we have just presented, and we will try to give rigorous qualitative and quantitative information about them in order to form a coherent panoramic view around the leitmotif of mass conservation, with strong focus on its companion concepts: self-similarity and existence of self-similar fundamental solutions. This will lead us in the end to some strange situations that are like doors into new, in principle unexpected qualitative behaviour of the solutions.


\subsection{Outline and main results in standard nonlinear diffusion}\label{sec.main}

In Sections \ref{sec.cm.conn} and \ref{sec.mc.ple} we introduce the connection between  self-similarity and mass conservation  for the most typical nonlinear diffusion equations. The known theory is called upon, and the connection is closely examined in the cases of heat equation, porous medium equation, and then for the $p$-Laplacian equation. The corresponding self-similar patterns are identified and examined. Other related issues like $L^1$ to $L^\infty$ estimates and asymptotic behaviour are considered. A general paradigm is found, and it is proved to be valid until we reach a bifurcation of behaviour at the so-called critical exponents, $m_c$ for the PME and $p_c$ for the PLE. Mass conservation is closely related to the $L^1$ contraction property, a very well-known property of this equation, \cite{BeTh72, BCP84, CrandallLiggett71, VazquezPME2007}, that reads: if $u$ and $v$ are two solutions of the HE, the PME or the PLE with $u_0,v_0\in L^1(\ren) $, then for $t>0$ we have $u(\cdot,t), v(\cdot,t)\in L^1(\re)$ and
\begin{equation}\label{Licont}\tag{L1C}
{ \|u(x,t)- v(x,t)\|_1\le \|u_0(x)- v_0(x)\|_1. }
	\end{equation}
 We remind the reader that our definition of $L^1$ contraction  accepts the  symbol $\le$ in formula  \eqref{Licont}, so it need not be a strict contraction. Some literature would call it {non-expansion } instead of contraction.  In this work we will follow the above definition in view of the abundant literature in the field. We will pay attention to possible confusions\nc.

 In Section \ref{ssec.emc} we introduce the concept of relative mass, prove the corresponding relative mass conservation. This is a  new result that is shown to be valid until we reach a bifurcation of behaviour at the so-called critical exponents, $m_c$ for the PME and $p_c$ for the PLE. To recall, this is a  stronger form of mass conservation\nc.

In Section \ref{sec.break} we show the quite novel situation that happens when we go below the critical exponents, to the so-called \sl very fast diffusion regime\rm. Then all finite-mass solutions lose mass with time, in fact we show that the mass goes to zero as $t$ goes to infinity. In particular, many of such solutions vanish identically after a finite time, a phenomenon called \sl finite-time extinction\rm. We also recall the dichotomy between everywhere positivity (in space and time) versus extinction in fast diffusion. An explanation of the loss of mass is given in Section 5.5 of \cite{Vazquez2006} in terms of particle dynamics and described as ``particles that have reached spatial infinity in finite time''. We refer the reader to the mentioned text which is meaningful for our present purposes. We review the fact at the end of Section \ref{ss.vfpme}.

Section \ref{sec.crit.exp} addresses the issue of mass conservation for the critical exponents, a main topic of the paper because it settles two important questions. We prove that MC holds in both cases, the result for the PLE, Theorem \ref{thm.mc.cple}, is new to our knowledge and far from evident. We also provide an alternative proof for the known result for the PME.

The critical exponents are very special in the sense that they are unique to combine mass conservation and no self-similar fundamental solution. We devote some effort to show how the fundamental solutions of the PME for $m>m_c$ manage to disappear in the limit, see Theorem \ref{thm.dirac.pme}, a phenomenon that we explain as a kind of frozen diffusion. See also Theorem \ref{delta.pc} for the PLE.

The study of the critical exponents is completed in Section \ref{sec.ncbl} where the following two exceptional cases are treated as limits: \ the PME with $m=0$, $N=2$, \ and the PLE with $p=1$, $N=1$. The situation that arises combines a frozen mass at the origin while the rest of the mass is lost by an amount that is exactly calculated, $4\pi t$ in the first case, and \, $2t$ in the second one. The combination leads to an object called  \sl measure-valued solution \rm that has finite time extinction.
Let us point out that both exceptional cases are related to interesting applications;   the first case represents planar logarithmic diffusion which models 2D Ricci flow, while the PLE case represents 1D Total Variation Flow.

Useful reminder: Following our book \cite{Vazquez2006} \nc we call \sl good fast diffusion range \rm the set of parameters where $m_c<m<1$ for the PME, resp. $p_c<p<2$ for the PLE, while \sl very fast diffusion \rm refers to the exponents  $0<m<m_c$ for the PME, resp. $1<p<p_c$ for the PLE. Both $m_c$ and $p_c$ will also \nc depend on the fractional exponent $s$ in the fractional models, see Sections \ref{sec.nonlinfrac} to \ref{sec.loss.vf.fpme}\nc.

In this part of the manuscript, we also aim to provide a comprehensive overview of the types of behaviour to be expected in these basic models. As a complement to this detailed study, we devote Section \ref{sec.dnle} to present the Doubly Nonlinear Equation, that combines the two nonlinearities that appear  in both  models, the PME and PLE\nc. We outline the main features that are found and state the main results, thus observing a common trend of development. But we will not delve on the fine proofs of the theory because we want  to avoid the \nc  danger of deviating from our main goal, which is presenting important features of nonlinear diffusion from the vantage point of view of self-similarity and mass conservation without taking up excessive space in the exposition of the less basic sections.

\newpage
\section{First study of  mass conservation and self-similarity}\label{sec.cm.conn}

The first subsections will cover all needed material to set the basics and then proceed into new results. We deal here with the Heat Equation and the Porous Medium Equation. In both cases we work with solutions $u(x,t)$ defined in the whole space $x\in \ren$ for $0<t<\infty$, unless mention to the contrary\nc.

\subsection{Mass conservation and the Heat Equation}\label{ssec.cmheq}
One of simplest ways of deriving the fundamental solution of the HE, $\partial_t u =\Delta u$, consists in assuming that the fundamental solution $U(x,t)$ starts from a point mass $U(\cdot,0)= \delta$ at $t=0$,  and then evolves by the heat flow as a nonnegative density function preserving the same unit mass at all times, $\int u(x,t)\,dx=1$. In other words, it is  an evolving probability density.

Let us recall that there is a well-established theory for the HE. It  shows that for every initial datum at $t=0$ in the form of bounded measure $\mu_0$ (in particular, an integrable function or Dirac mass) there exists a unique weak solution $u(x,t)$ defined for $x\in\ren$ and $t>0$ that is a smooth and bounded solution of the HE for all strictly positive times. Moreover, nonnegative measures produce nonnegative solutions. If the initial data $u_0\in L^1(\ren)$  then the solution $u(t): C([0,T): L^1(\ren))$ and it coincides with the semigroup construction found by the Crandall-Liggett approach. It takes the initial data in the sense that  $u(t)\to u_0$ if $t\to 0$ in the sense of $L^1$ convergence.

\noindent {\bf Remark}.  The minus the Laplacian operator is an $m$-accretive operator in all $L^q(\ren)$ spaces for all $q>1$;  therefore the HE generates an  $L^q$ contraction semigroup in all these norms. For $q=2$ this is the well known property of being a maximal monotone operator, see \cite{BrezisBkOMM73}. Because of mass conservation we are interested in $q=1$. We mention in passing that the HE generates a gradient flow in $L^2(\ren)$ and in $H^{-1}(\ren)$. In the common domain of definition all these concepts of solution agree.

\medskip

\noindent {\bf Similarity analysis.} We show next how to construct the fundamental solution by assuming their self-similarity and radial symmetry. \nc The invariance of the equation under the mass-conservation scaling group implies that  the fundamental solution $U(x,t)$ has
the typical form
\begin{equation}\label{form.sssol}
U(x,t)=t^{-\alpha}F(x\,t^{-\beta}),
 \end{equation}
with $\alpha$, $\beta$ and $F$ to be determined, cf. the specialized monograph \cite{BarentSSSIA}\nc.
Upon substituting the self-similar formula into the equation $u_t=\Delta u$, we easily find
that this substitution is compatible with the equation if $\beta= 1/2$, and this leaves freedom to choose $\alpha$. The choice of the class of solutions with finite mass that conserve this quantity
implies that $\alpha=N\beta $ since
$$
\int U(x,t)\,dx=t^{-\alpha}\int F(x\,t^{-\beta})\,dx=t^{N\beta-\alpha}\int F(y)\,dy\,,
$$
with integrals over $\ren$. The consequence is that we have  exact values of the self-similar exponents
$$
\alpha=N/2, \qquad \beta= 1/2.
$$
Then, the profile $F(y)$ has to satisfy the elliptic equation
$$
\Delta F + (1/2) \,\mbox{div}(yF)=0.
$$
i.e., $\mbox{div}(\nabla F + (1/2) yF)=0$. In order to find a solution, we assume that the expression in parentheses vanishes so that  $\nabla F + (1/2) yF$=0. We now  make the assumption of radial symmetry (a natural assumption since the solution comes from initial delta mass);  we get by simple integration  the explicit formula
 $$
 F(y)=C\,e^{-|y|^2/4}.
 $$
 The constant $C$ is adjusted so that the solution has constant unit mass and we end up with the Gaussian evolution kernel $U(x,t)=G(x,t)$, where
\begin{equation}\label{form.gauss}
G(x,t)=\frac{1}{(4\pi t)^{N/2}}e^{-|x|^2/4t}.
\end{equation}

\noindent {\bf The Gaussian pattern.} This is how conservation of mass and self-similarity allow for a very easy derivation of the fundamental solution $G$. This is the first diffusive pattern found in the theories of the present paper. It is probably the most influential pattern in higher mathematics, applied mathematics and mathematical physics.
 The procedure we used is not the only way to derive the Gaussian kernel, but it is very illustrative of the way to proceed in order to obtain the fundamental solutions of other nonlinear models, like the PME and the PLE.

It  is well known that $G $ is the unique nonnegative weak solution starting from a unit Dirac mass (this fact  is another necessary step of the theory). A solution starting from a Dirac mass is sometimes called \sl source-type  solution\rm.
The relation $|x|^2/t\sim $ constant is important in many propagation calculations and is called the Brownian scaling, because of its use in the description of Brownian motion, cf. \cite{KLChung88, KarShreve1988}. We will return to the subject of natural scalings later.

\begin{figure}[ht!]
\centerline{\includegraphics[width=0.6\textwidth]{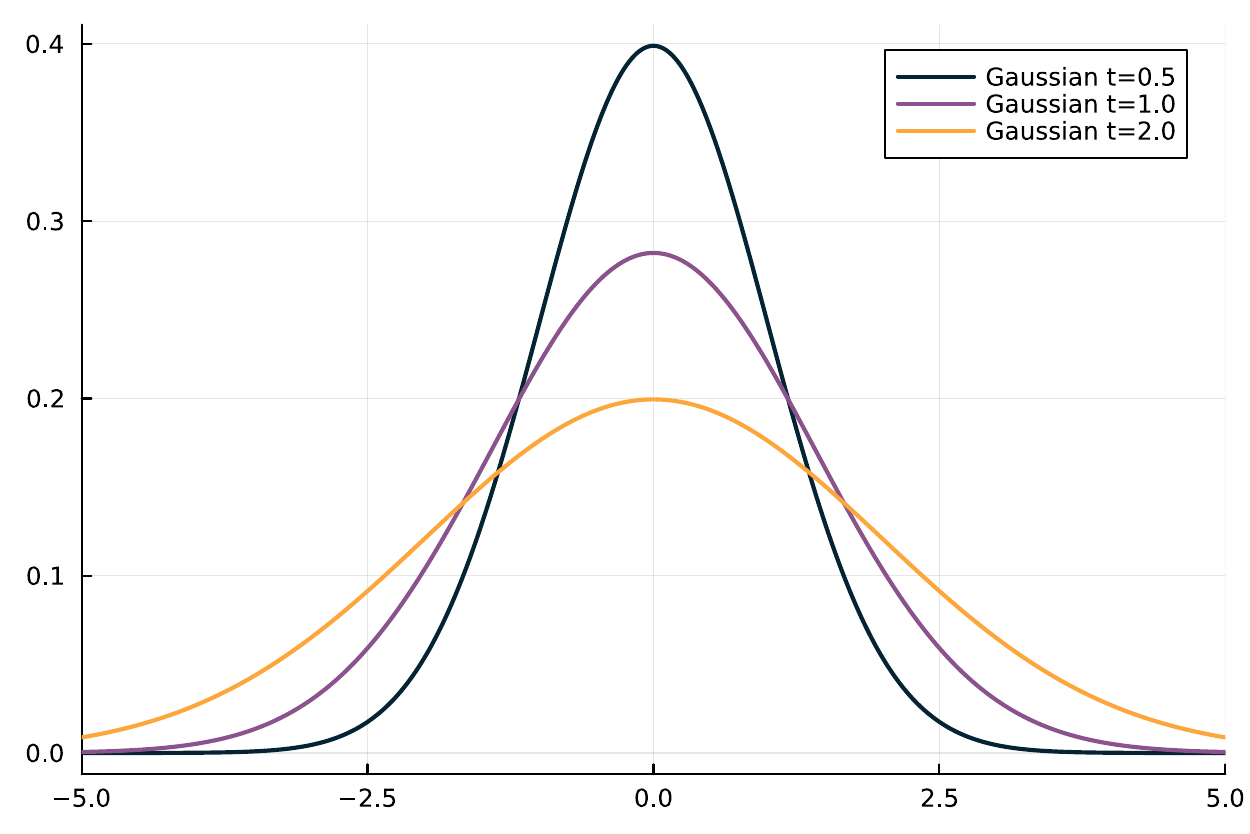}}
\vspace{-0.5cm}
\caption{Gaussian function evolution by the HE}
\label{fig:Gaussian}
\end{figure}

\noindent {\bf Gaussian representation.} Once $G$ is found, the linear representation formula
$$
u(x,t)= (G_t\ast u_0)(x):=\int_\ren G(x-y,t)\,u_0(y)\,dy
$$
solves the heat equation for a large class of initial data, in particular for all integrable data. It is  then easily proved that all solutions starting from initial data $u_0\in L^1(\ren)$ have a constant first integral in time:
\begin{equation}\label{CM.sec2}
\int_\ren u(x,t)\,dx=\int_\ren u_0(x)\,dx.
\end{equation}
For nonnegative data $u_0\ge 0$, this is just what we call  conservation of the total mass, MC.
The conservation result follows in a very elementary way from the convolution form of the representation and the MC that holds for $G$.  Therefore, the circle is closed:

Fundamental Solution - Self-similar solution - Gaussian function - MC.

\noindent We will see in moment a  consequence, attraction to the Gaussian.

Since in the linear setting the difference of solutions is a solution, the mass conservation law applies to the difference of two solutions, this will be called later relative mass conservation. Warning: no absolute value must be inserted in the integrands. Thus, in the following well-known law of $L^1$ contraction for any two solutions, see \eqref{Licont},
 we must point out that the equal sign almost never holds unless $u$ and  $v$ are ordered.
\medskip

\noindent {\bf Attraction to the Gaussian. Asymptotic behaviour.}
  In a further  step we can prove that for any $u_0\in L^1(\ren)$ with $\int u_0\,dx\ne 0$ we have the asymptotic behaviour  as $t\to \infty$ of the informal form $u(\cdot,t)\sim MG(\cdot,t)$, thus showing that the Barenblatt solution has a universal attraction power inside the set of solutions with the same mass. This takes very precise formulations, the simplest being the $L^1$ version
\begin{equation}\label{he.asymp}
\lim_{t\to \infty}\|u(\cdot,t)- MG(\cdot,t)\|_1=0.
\end{equation}
\bc Another asymptotic convergence result happens in $L^\infty$ norm, and then a normalization factor is needed for optimal formulation:
\begin{equation}\label{he.asymp.Li}
\lim_{t\to \infty} t^{N/2}\|u(\cdot,t)- MG(\cdot,t)\|_{\infty}=0.
\end{equation}
Convergence in all $L^p$ norms follows by interpolation for all $1\le p\le \infty$ with suitable time factors, $t^{N(p-1)/2p}$. \nc
Thus, not only every integrable solution conserves mass, indeed the mass $M$ (a number) and the heat profile $G$ (a function) are the only details the solution remembers from its long evolution when we look at very large times in first approximation. This was called by Barenblatt ``intermediate asymptotics'',   \cite{BarentSSSIA}. Of course, further details of the asymptotics are known for the heat equation with finer approximation, but this is another story, see \cite{VazNotesHE}.

These results are well-known in the literature, see e.g. \cite{CourHilb, EvansBk}. A very detailed account of this issue with description of  different methods of analysis of the behaviour of the solutions of the HE and its transformations is contained in the Notes \cite{VazNotesHE}. Note in particular that the $L^p(\ren)$ norms, $p>1$, of any finite mass solution are not conserved but go to zero with a certain rate that has some freedom.

\subsection{Mass conservation and its connections  for the PME}\label{ssec.pme}

We begin the analysis of these issues in the nonlinear setting by a typical example, the Porous Medium Equation, $\partial_t u=\Delta (|u|^{m-1}u)$, where the same approach works in many cases but offers an interesting twist. The analysis gets more complicated, so it is best addressed via studies in different ranges; we do first the ones that work as predicted by the current programme. We will often write the equation in the simpler form $\partial_t u=\Delta (u^m)$ that is completely correct for $u\ge 0$.

\subsubsection{Slow diffusion range of the  PME} We consider the PME  in principle for nonnegative solutions. We begin with the range of exponents $m>1$, the superlinear nonlinear range that leads to a regime called {\sl slow diffusion} for reasons that are well explained in the literature \cite{Aron1986, VazquezPME2007}. Indeed, the finite speed of propagation is already clear from the Barenblatt profile to be constructed below\nc.

\noindent $\bullet$ {\bf Barenblatt solution.} Arguing as before, we start our enquiry with the existence of a fundamental solution, assuming as before that self-similarity and mass conservation hold. This leads to a solution of the form
$$
U(x,t)=t^{-\alpha}F(x\,t^{-\beta}).
$$
Again, we perform the substitution of the self-similar formula into the equation, this time $u_t=\Delta u^m $, and we easily find that the cancelation of the powers of time that appear after substitution leads to the necessary condition \ $\alpha\,(m-1)+ 2\beta=1$. This leaves still one freedom to choose the parameters. We use  this freedom of  choice to impose conservation of mass in the class of solutions with finite mass, that implies  $\alpha=N\beta $ as we have seen in the case of the heat equation. An immediate algebraic calculation leads to  the exact values of the self-similar exponents
\begin{equation}\label{for.simexp.pme}
\alpha=\frac{N}{2+ N(m-1)}, \qquad \beta= \frac{1}{2+ N(m-1)}\,,
\end{equation}
while the profile $F(y)$ has to satisfy the nonlinear elliptic equation
\begin{equation}\label{profile.pme}
\Delta F^m + \beta \,\mbox{div}\,(yF)=0.
 \end{equation}
Eliminating (informally for the moment)  the divergence from the last equation we get $\nabla F^m = -\beta yF,$ which for $F\ne0$ reduces to $\frac{m}{m-1} \nabla F^{m-1} = -\beta y$. \nc This can be easily integrated for radial functions to get  a {\sl weak solution} of \eqref{profile.pme},
$$
\frac{m}{m-1}F(y)^{m-1}=\left(C_o- \frac{\beta}{2}y^2\right)_+=\max\left\{C_o- \frac{\beta}{2}y^2, 0\right\},
$$
that combines the analysis for $F>0$ and $F=0\nc$. Summing up, we obtain the so-called the Barenblatt solution $U=B(x,t)=B_m(x,t)$, with formula
\begin{equation}\label{for.Bar1}
B(x,t)=t^{-\alpha}\left(C-k\,\frac{|x|^2}{t^{2\beta}}\right)_+^{1/(m-1)} \quad \mbox{with} \ k=\frac{\beta (m-1)}{2m}.
\end{equation}
We will write $B_M(x,t)$ or $B(x,t;M)$ when we want to stress the dependence of the Barenblatt solution on the mass $M$.

$\blacktriangleright$
 The constant $C>0$ is free and can be used to fix the total mass $M$. It is very easy to see that this mass  is conserved in time. Indeed, a direct computation of the mass of $B_m (x,t)$  gives
$$
M=N \omega_N \int_0^{\infty} (C-kr^2)_+^{1/((m-1)}r^{N-1}dr= d(m,N)\,C^{1/(2\beta(m-1))},
$$
 Using known expressions Euler's Beta and Gamma functions  the value of $d(m,N)$ is calculated in pages 25, 26  of \cite{Vazquez2006} as:
\begin{equation}\label{for.d_m.Bar}
d(m,N) =(1/2) N \omega_N  \frac{\Gamma(\frac{N}2) \Gamma( \frac{m}{m-1})}{\Gamma(\frac{N}2 +\frac{m}{m-1})} k^{-N/2},
\end{equation}
where $\Gamma(\cdot)$ is Euler's Gamma function. From this we get the expression
\begin{equation}\label{for.C.Bar}
C=a(m,N)\,M^{2(m-1)\beta}, \qquad a(m,N)=d(m,N)^{-2\beta(m-1)},
\end{equation}
with $\beta$ as in \eqref{for.simexp.pme}.

$\blacktriangleright$ In this way we control the decay of the solution in $L^\infty$ norm
$$
\|B_m(\cdot,t)\|_\infty = B_m(0,t)=C^{1/(m-1)}t^{-\alpha}= \mathfrak{c}(m,N)\,M^{2\beta}\,t^{-\alpha},
$$
with  $\mathfrak{c}(m,N)=a(m,N)^{1/(m-1)}=d(m,N)^{-2\beta}.$

$\blacktriangleright$  Notice that the positivity set of $B_M(x,t)$ is an expanding  ball surrounded by the hypersurface (called \sl the free boundary\rm) of exact size
\begin{equation}\label{form,fbpme}
|x(t)|={L}_1\, ( M^{m-1}t )^{\beta},        \qquad L_1^2=\frac{2m}{\beta (m-1)} a(m,N).
\end{equation}
This formula is a quantitative version of the property of finite propagation of the PME.
Note that  ${L}_1 = {L}_1(m,N)$ is the normalized {\sl expansion radius\rm} of the signal emanating from $x=0$ for $t=1$ and $M=1$.\nc

$\blacktriangleright$ To make things easy we may often consider mass $M=1$. If we need to take a different mass into account we will write $B(x,t;M)$ or $B_M(x,t)$. There is a scaling formula to relate them
$$
B(x,t;M)=M\, B_1(x, M^{m-1} t).
$$
The details are carefully worked out in \cite{Vazquez2006}, Chapter 2.

We will keep track of these constants when we need their behaviour later on, see also Theorem \ref{limitminfty} below.\nc

\begin{figure}[ht!]
\centerline{\includegraphics[width=0.6\textwidth]{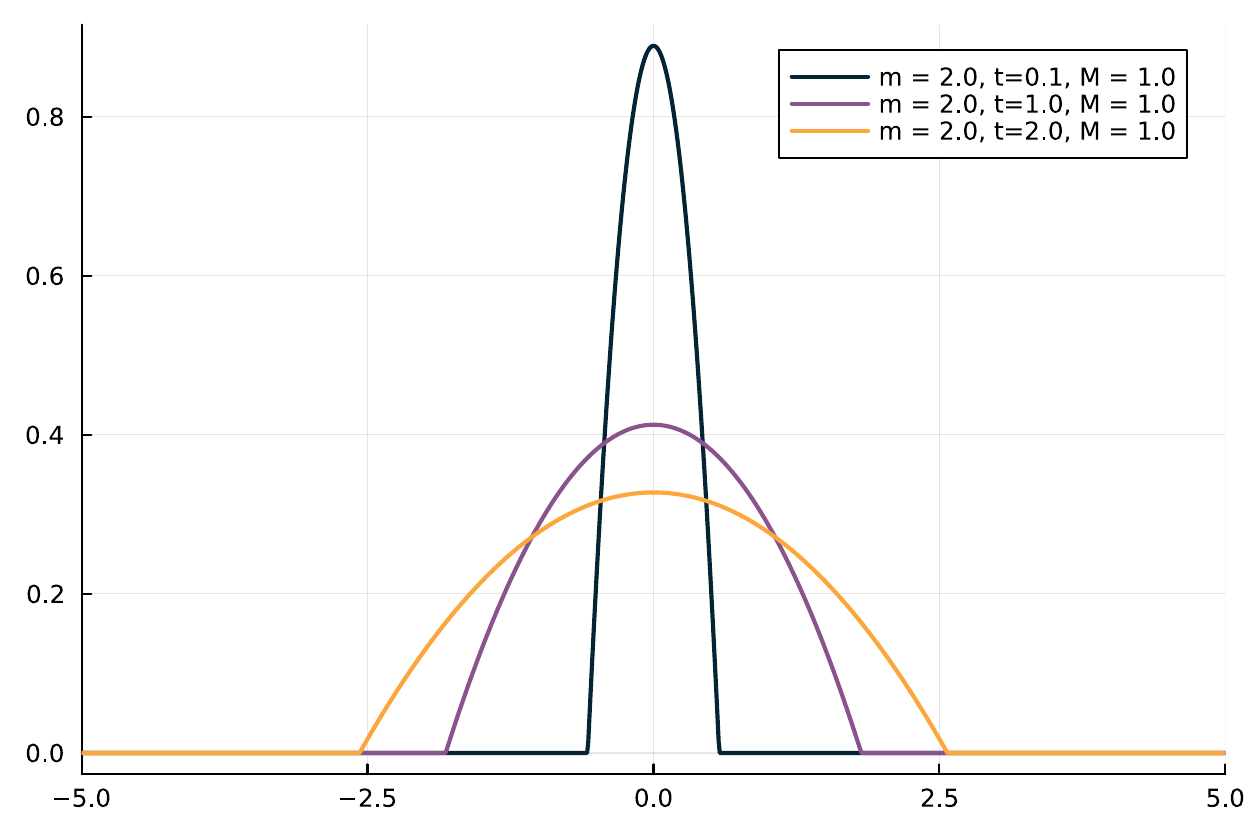}}
\vspace{-0.5cm}
\caption{Evolving  Barenblatt Solution represented at different times, starting from a Dirac delta at $t=0$. The free boundary $x(t)$ lies at a finite distance, moving away with time as indicated by the obtained formulas}
\label{fig:PME}
\end{figure}

 Summing up, we have found de typical Barenblatt patterns for diffusion with finite propagation speed and free boundaries, an investigation started in the early 1950s by Zeldovich and his group \cite{ZK50, Bar52}, see also \cite{Pattle59}. The authors saw that $B_M(x,t)$ is not a classical solution of the PME due to lack of regularity on the free boundary (look at the first derivative), but they accepted it as true solution on  physical grounds.
 This intuition was validated with the introduction of weak solutions by Oleinik et al. \cite{OKC58} some years later. In the class of weak solutions existence and uniqueness was proved in 1958 in dimension $N=1$, a landmark in the theory\nc. See \cite{Sb61} for the $N$-dimensional version.

\medskip

\noindent $\bullet$  {\bf About the good class of  integrable solutions of the PME.}  In this nonlinear equation, there is no reasonable way of representing a general solution with integrable initial data in terms of the fundamental solution, a problem that happens with almost all nonlinear equations. Moreover, the degeneracy of the equation complicates things, and in particular it produces the Barenblatt solution of Figure \ref{fig:PME} that we need to include into the set of admissible class of solutions because of physical evidence.
This is a standard situation in the analysis of PDEs; it is solved by introducing a suitable class of generalized solutions,  as explained in the Introduction.

In the case of the PME the existence and uniqueness of finite mass solutions is well-known. Thus, for every initial datum $u_0\in L^1(\ren)$
there exists a unique distributional solution $u(x,t)$ defined for $x\in\ren$ and $t>0$ such that $u(t): C([0,T): L^1(\ren))$ with  $u(0)=u_0$, and $u$ is  a continuous  and uniformly bounded function for all strictly positive times $t\ge \tau>0$.  This solution can be obtained via the semigroup construction obtained by  the  Semigroup Generation Theorem acting in the space $X=L^1(\ren)$. The reader can find extensive details in  our monograph \cite{VazquezPME2007}, where well-posedness and the available regularity are proved. See also the general references \cite{DaskaBk, DiBenedetto93, WuZY01} and the original works of Bénilan and Crandall et al, \cite{BeTh72, CrandallLiggett71, BC81}. We call them mild solutions, semigroup solutions, or simply solutions  if the context is clear. Moreover, nonnegative initial data produce nonnegative solutions.

 There  are many partial regularity  results that make computations with the PME easier. Thus, all bounded solutions are H\"older continuous for $t>0$. Let us mention a  set of data that is dense in $X=L^1_+(\ren)$ for which solutions are classical: if $u_0$ is continuous and positive then the solution $u(x,t)$ is $C^{\infty}$ smooth and positive for all positive times $t>0$.

\noindent {\bf The $L^1$ contraction property.} The semigroup map $S_t: u_0\mapsto u(t)$ is a contraction in the space norm $L^1(\ren)$. By this we mean that for any two $L^1$-semigroup solutions $u$ and $v$ we have
$$
\int_\ren |u(x,t)-v(x,t)|\,dx\le \int_\ren |u_0(x)-v_0(x)|\,dx
$$
as in the case of the HE. More precisely, we get the stronger form of contraction called {\sl ordered  $L^1$-contraction} or {\sl $L^1$ comparison principle},
$$
\int_\ren (u(x,t)-v(x,t))_+\,dx\le \int_\ren (u_0(x)-v_0(x))_+\,dx,
$$
see \cite{BeTh72}, which implies in an easy way comparison of solutions: $u_0\le v_0$  implies for every $t>0$ that $u(t)\le v(t)$. In particular, since $u=0$ is a solution, nonnegative initial $v_0\ge 0$ produces  a nonnegative solution $v(t)\ge 0$. Similarly,  for nonpositive data $u_0\le 0$\nc.

\noindent {\bf Measure data.} This theory does not directly incorporate the Barenblatt solution because it takes as initial datum a Dirac delta and not a function. The existence and unique theory has been extended to cover initial data in the form of any bounded Radon measure $\mu_0$. Existence is obtained by approximating  $\mu_0$ with a sequence of functions $u_{0n}$ that converge to $\mu_0$ in the measure norm, and passing to the limit in the corresponding semigroup solutions $u_n$. In the limit we find a function $u(x,t)$ such that $u(t): C([\tau,T): X)$   for every $ \tau>0$  and $u$ is a  semigroup solution for $t\ge \tau>0$ with $L^1$ data at $t=\tau$; finally, $u(\tau)\to \mu$ as $\tau\to 0$ as bounded measure. This is the concept of solution with measure-valued data. Uniqueness of such solution was a difficult problem. See details and references in Chapter 13 of \cite{VazquezPME2007}. Note that this approximation process applies to the Barenblatt solution by just putting
$$
u_n(x,t)=B(x,t+(1/n)), \qquad n\ge 1,
$$
as an approximation  sequence\nc.

\noindent {\bf Remark.} { The operator $Au= -\Delta(u^m)$ that appears in the PME is not maximal monotone and the theory of contraction semigroups in Hilbert spaces could not be applied, cf. \cite{BrezisBkOMM73}. Moreover, it is folklore knowledge that  the operator does not generate an $L^q(\ren)$ contraction semigroup for any exponent $q>1$. We will find a concrete negative example to confirm this folk  assertion  in Subsection \ref{ssec.lack}.}

\medskip

 \noindent $\bullet$ {\bf Properties of nonnegative solutions. } If $u\ge 0$ the \nc theory is significantly simplified by the existence of two a priori inequalities that are valid for all nonnegative solutions defined in the whole space. The first inequality  the Bénilan-Crandall homogeneity inequality   \cite{BC81b} that reads  for $m>1$\nc
\begin{equation}\label{BCest}
(m-1)t\,u_t\ge -u.
\end{equation}
The other one is the Aronson-Bénilan inequality \cite{Aron1979} that reads
\begin{equation}\label{ABest}
\Delta( u^{m-1})\ge -\frac{K}{t}, \quad K\nc=\frac{\alpha(m-1)}{m}.
\end{equation}
This estimate is exact for all the Barenblatt solutions. It also  implies after an easy computation that $t\,u_t\ge -\alpha u$, which is finer than  Bénilan-Crandall's. The estimate is valid also for the Fast Diffusion case $m<1$ on the condition that
$\alpha>0$, which means $m>m_c$, see that subsection below. Both estimates are valid for the class of nonnegative mild solutions with $L^1$ initial data that we have mentioned. By taking limits we may apply them to the solutions with optimal classes of data to be discussed below\nc.

\medskip

\noindent $\bullet$ Once this is settled, the indirect way to prove MC  is easy.  We observe that the set of bounded initial data with compact support forms a dense set in $L^1_{loc}(\ren)$. By plain comparison with a large Barenblatt solution $B_M(x,t+T)$,  we show that the corresponding \nc solutions $u(x,t)$ may expand in time but always keeping the property of compact support.  Mass conservation is then easily proved. By  density and $L^1$ contraction, MC follows for all solutions with $L^1$ data.

Since the assumption of mass conservation led to a case of self-similarity (the Barenblatt solution) with finite propagation, and MC can be proved  for all $L^1$ solutions by comparison with the Barenblatt solution, the circle we mentioned in the HE is closed:

   Fundamental Solution - Self-similar solution - Barenblatt Solution - MC, with its consequence: attraction to the Barenblatt.

\noindent {\bf  Asymptotic behaviour. Attraction to the Barenblatt.} \nc
  Once this is settled, we may prove that for any $u_0\in L^1(\ren)$ with $\int u_0\,dx\ne 0$ we have the asymptotic behaviour  as $t\to \infty$ of the informal form $u(\cdot,t)\sim B(\cdot,t;M)$, thus showing that the Barenblatt solution has a universal attraction power inside the solutions with the same mass. This allows for very precise formulations, the simplest being the $L^1$ version
\begin{equation}\label{form.asb.bar}
\lim_{t\to \infty}\|u(\cdot,t)- B(\cdot,t;M)\|_1=0.
\end{equation}
\bc Another version of the asymptotic convergence result happens in $L^\infty$ norm and then a normalization factor is needed for optimal formulation:
\begin{equation}\label{asymp.pme.Li}
\lim_{t\to \infty} t^{\alpha}\|u(\cdot,t)- B(\cdot,t;M))\|_{\infty}=0.
\end{equation}
with $\alpha$ the similarity exponent in \eqref{for.simexp.pme}.
Convergence in all $L^p$ norms follows by interpolation for all $1\le p\le \infty$ with suitable time factors, $t^{\alpha (p-1)/p}$\nc.

What is surprising is that the results are (almost) completely formally analogous to the ones for the heat equation \eqref{form.asb.bar}.
However, the different available proofs have nothing to do, there is no trace of any representation formula\nc. The first proof is due to Kamin \cite{Kamin73}, see then \cite{FmKam80}. For a detailed account of this theory see Chapter 18 of the book \cite{VazquezPME2007}. \bc  For solutions with compact support a very precise regularity of the solutions up to the free boundary holds for large times, and moreover the asymptotic convergence  can be improved into an infinite expansion, see \cite{KKochV18, Seis14}\nc.

\medskip

\noindent {\bf Optimal boundedness estimates.}
There is another important aspect in which the finite-mass self-similar solution has
an optimal performance, we are referring to the so-called {\sl $L^1$-$L^\infty$ smoothing effect}. Indeed, when we look at the $N$-dimensional Barenblatt solution \eqref{for.Bar1} we  have observed that
\begin{equation*}
\sup_x B_M(x,t)=  \mathfrak{c}(m,N)\,M^{2\beta}\,t^{-\alpha}, \quad  \mathfrak{c}(m,N)=a(m,N)^{1/(m-1)},
\end{equation*}  
where $a(m,N)$ is the constant that appears  in \eqref{for.Bar1}, \eqref{for.C.Bar}.  So the first optimality question is to decide if this inequality is satisfied by a large class of solutions. The positive answer is as follows.

\begin{proposition}\label{prop.2.1}  For every nonnegative solution $u(x,t)$ of the PME with initial datum $u_0\in L^1(\ren)$ and every $t>0$ we have
\begin{equation}\label{form.L1Linf}
\sup_x u(x,t)\le B_M(0,t)= \mathfrak{c}(m,N)\, M^{2\beta}\,t^{-\alpha},
\end{equation}
where $M=\int_{\ren} u_0(x)\,dx$, and $\alpha$ and $\beta$ are the similarity exponents found in \eqref{for.simexp.pme}\nc.
Equality holds at a certain time $t=t_1$ if $u(x,t_1)$ is the Barenblatt solution with mass $M$ or a space translation thereof.
\end{proposition}
\nc

\noindent {\sl Proof.}  The general inequality \eqref{form.L1Linf} comes  from  an argument of  Schwarz symmetrization \nc and concentration comparison that was in introduced by the author in \cite{V82sym} and is explained in Chapter 17 of \cite{VazquezPME2007}. We will refrain from explaining in detail these two techniques in this memoir since they are well documented, see [19, 265, 266, 271], they take
some space and depart from the main topics we have in mind. However, a brief idea of the proof could be useful: a nonnegative and integrable function real  $f(x)\ge0$ is called rearranged if it is radially symmetric and nonincreasing as a function of $|x|$,  see \cite{Talenti76}. It follows that a rearranged function  attains its supremum at the origin $x=0$.  A well-known result of the $L^1$ theory of the PME asserts that of rearranged initial data $u_0$  produce a solution $u(x,t)$ that is rearranged as a function of $x$ for every fixed $t>0$. Another main theorem the theory allows to compare the $L^P$ norms of a solution $u$ a time  with the radially symmetric solution $v$ starting from $v_0$, the rearranged function obtained from $u_0$. Comparing the radially symmetric $v$ with the  Barenblatt solution uses the concentration comparison theorem that applies since the Barenblatt has the most concentrated initial data, i.e., a Dirac delta.\qed

A direct consequence of the asymptotic result \eqref{form.asb.bar} is that $\sup_x u(x,t)$ approaches $B(0,t;M)$ in relative size as $t\to \infty$:
$$
\lim_{t\to \infty}\frac{\sup_x u(x,t)}{M^{2\beta}\,t^{-\alpha}}=\mathfrak{c}.
$$

We wonder what happens with the optimality of the inequality  at earlier times. The first question is to determine whether the Barenblatt solutions are the only solutions that satisfy the inequality with equality sign.
In the sequel we will try to quantify the error in the boundedness estimate,
$$
E(u,t)= B(0,t;M)-\sup_x u(x,t).
$$

\begin{proposition}\label{prop.2.2}  \bc For every nonnegative solution $u(x,t)$ of the PME with initial datum $u_0\in L^1(\ren)$ and every $t>0$, the error $E(u,t)$
 can be controlled from above by the $L^1$ norm of the difference,  $\int_{\ren} |u(t)-B(t)|\, dx$,  at the same time, or any previous time, see formula
\eqref{Error.above}
\end{proposition}
\nc

\noindent {\sl Proof.} \bc Notation: Here below the notation $c_i(m)$, $i=1,\cdots$, denotes universal positive constants whose exact value we do not need to chase\nc. Assume that optimality is not obtained. In order to relate the error $E(u,t,M)$ to the $L^1$ norm of the difference from $u$ to $B$ we may assume that $M=t=1$ without loss of generality (by using the scaling properties of the equation).

Let $E(u,1)=h>0$. Then
$u(x,1)\le \mathfrak{c}(m,N)-h$, so that the quantity
$$
B(x,1)-u(x,1)\ge B(x,1)+h -\mathfrak{c}
$$
 is positive in a ball or radius $R(h)$ such that, recalling that $\mathfrak{c}=C^{1/(m-1)}$:
 $$
(\mathfrak{c}^{m-1}-k(m)\,R^2)^{1/(m-1)}=\mathfrak{c}-h = \mathfrak{c}(1-h/\mathfrak{c})
$$
For $h$ small this gives $R\sim c_1(m)h^{1/2}$, so that
\begin{equation}\label{Error.above}
\|(B(x,1)-u(x,1))_+\|_1\ge \int_{|x|\le R(h)} |B(x,1)-u(x,1)|\,\dx\ge c_2(m) \,h^{1+(N/2)}.
\end{equation}
By conservation of  mass we have $\|(B(x,1)-u(x,1))_-\|_1=\|(B(x,1)-u(x,1))_+\|_1$. Hence,
$$
\|B(x,1)-u(x,1)\|_1\ge 2 \int_{|x|\le R(h)} |B(x,1)-u(x,1)|\,\dx.
$$
The $L^1$ contraction formula  \eqref{Licont} \nc implies that \nc this estimate is also true for $\|B(x,t)-u(x,t)\|_1$ for all $t<1$. This implies that a small error in  $\|B(x,t)-u(x,t)\|_1$ for some $t<1$ implies a small error $E=h$ in the functional inequality \eqref{form.L1Linf}. \qed

\medskip

To complete the control of $E(u,t)$ we would like to have a control in the opposite direction. Here is a partial answer.

\begin{proposition} \label{prop.2.3} \bc Assume that  we work in dimension  $N=1$. Then the optimal bound is attained at a certain time $t_1$ if and only if $u(x,t_1)$ is the Barenblatt solution with same mass $M$ or a space translation thereof \footnote{For all later times $t>t_1$ the equality $u(x,t)=B_M(x,t)$ would follow.}. In several dimensions, $N>1$, the result is true  if  the solution is radially symmetric and rearranged w.r.t. the space variable.

 In both cases, $E(u,t)$ is bounded from below  by a function of $\|B(x,t)-u(x,t)\|_1$, which is linear for small errors,
\begin{equation}\label{Error.below}
\|B(r,1)-u(r,1)\|_1 \le c_3(m)\,h\,,
\end{equation}
for small $h$. We have taken $M=1$ without loss of generality.\nc
\end{proposition}

\noindent {\sl Proof.}  We want to find conditions under which $E(u)=h\sim 0$ implies that $\|B(x,1)-u(x,1)\|_1\sim 0$. We will show that this happens in 1 space dimension for general integrable  data $u_0\ge 0$, and also in several space dimensions when the solution is radially symmetric and rearranged, $u=u(r,t)$ with $r=|x|$.

For the proof in several dimensions we argue as follows: Assume that we are in the situation $u(0,1)=\mathfrak{c}-h< \mathfrak{c}$. We choose $M^*=1-\ve<1$ the mass for which the Barenblatt $B(0,1;M^*)=\mathfrak{c}-h$.  Observe that since $\mathfrak{c}-h= \mathfrak{c} \,(M^*)^{2\beta}$, for small $h$ we have $\ve=1-M^*(h)\sim c_4(m)h$.

Then the Aronson-B\'enilan estimate says that
$$
\Delta u^{m-1}\ge \Delta B^{m-1}(r,1;M^*)=-\frac{K}t\,.
$$
It is crucial in what follows that the constant $K$ applies to both $u$ and $B$ (it is universal). Integrating twice this inequality from $r=0$ we get
$$
u(r,1)^{m-1
}\ge B^{m-1}(r,1;M^*)
$$
i.e, $u(r,1)\ge B(r,1;M^*)$.
Since $ B(r,1;1)$ and $u(r,1)$ have the same mass, we get
$$
\|(B(r,1;1)-u(r,1))\|_1 =2 \|(B(r,1;1)-u(r,1))_+\|_1\le
$$
$$
2\|(B(r,1;1)-B(r,1;M^*))_+\|_1 \le c_5(m)\,h\,.
$$

 This last estimate implies that when there is contact, $h=0$, then  $u(r,1)= B(r,1;1)$.

When we are in $N=1$ dimension, we first move the origin of coordinates to the point were the supremum
is attained and so the same argument applies with respect to the translated Barenblatt solution. No rearrangement is needed\nc.
\qed

\noindent {\bf Note.} The results of Propositions \ref{prop.2.2} and  \ref{prop.2.3} are new.
\medskip

\noindent {\bf Remark.} \bc Check that  estimate \eqref{Error.below} agrees with the example of a delayed Barenblatt solution $u(x,t)=B(x,t+\tau)$ with $\tau>0$\nc.

\begin{open}
We wonder if  a version of the last implication holds true in several space dimensions  when $u$ is  a nonnegative, finite-mass solution without symmetry assumptions on the data or the solution.  \nc
\end{open}



\subsection{The limit case $m\to \infty$ of the PME. A break of the theory}\label{ssec.pme.minfty}

\noindent The PME with large values of $m$ appeared in models to  describe heat propagation by radiation occurring in plasmas (ionized gases) at very high temperatures, see \cite{ZR}, among other applications. This case attracted attention when numerical calculations showed that for (not so) large values of  exponent $m$  the nonnegative solutions of the  PME with finite mass tended to take a shape that was called a \sl mesa \rm (the name taken from the Mesa plateaus of different arid regions of the world). And this happens for all $t>0$. \nc The natural reaction in mathematics is examining the limit $m=\infty$ with the tools of analysis.

\noindent $\bullet$ We will do that by looking at the limit of the Barenblatt solution  $B_m (x,t;M)$ with mass $M$ for exponent $m>1$. Indeed, in order to  pass to the limit $m\to\infty$ we first see that both similarity exponents  $\alpha, \beta\to 0$;  then we observe the optimal upper bound  \eqref{form.L1Linf} and observe that all three factors in the right-hand side tend to 1 in the limit, so that we get
$$
\lim_{m\to \infty} B_m (x,t;M)\le  1 \quad \forall x\in\ren, \ t>0.
$$
The collapse  of the limit solution into a flat pile of height 1 is shown next. The proof is inspired on Theorem 1.2 of \cite{FdmnHoll87}.

\begin{theorem}\label{limitminfty}
As $m\to\infty$  there is convergence  of $B_m (x,t;M)$ to the mesa profile
\begin{equation}\label{form.mesa}
B_\infty (x,t)= \mathcal{M}(x),
\end{equation}
where $\mathcal{M}(x)$  is the characteristic function of the ball $B_R(0)$ with volume $M$. In other words, $\mathcal{M}(x)=\chi_{B_R(0)}(x)$ with $M=\omega_N R^N$\nc.
The convergence cannot be uniform everywhere since $\mathcal{M}(x)$ is not continuous, it happens in the $L^\infty$-weak star topology for finite time intervals $0<a<t<b$.
\end{theorem}

\begin{figure}[ht!]
\centerline{\includegraphics[width=0.6\textwidth]{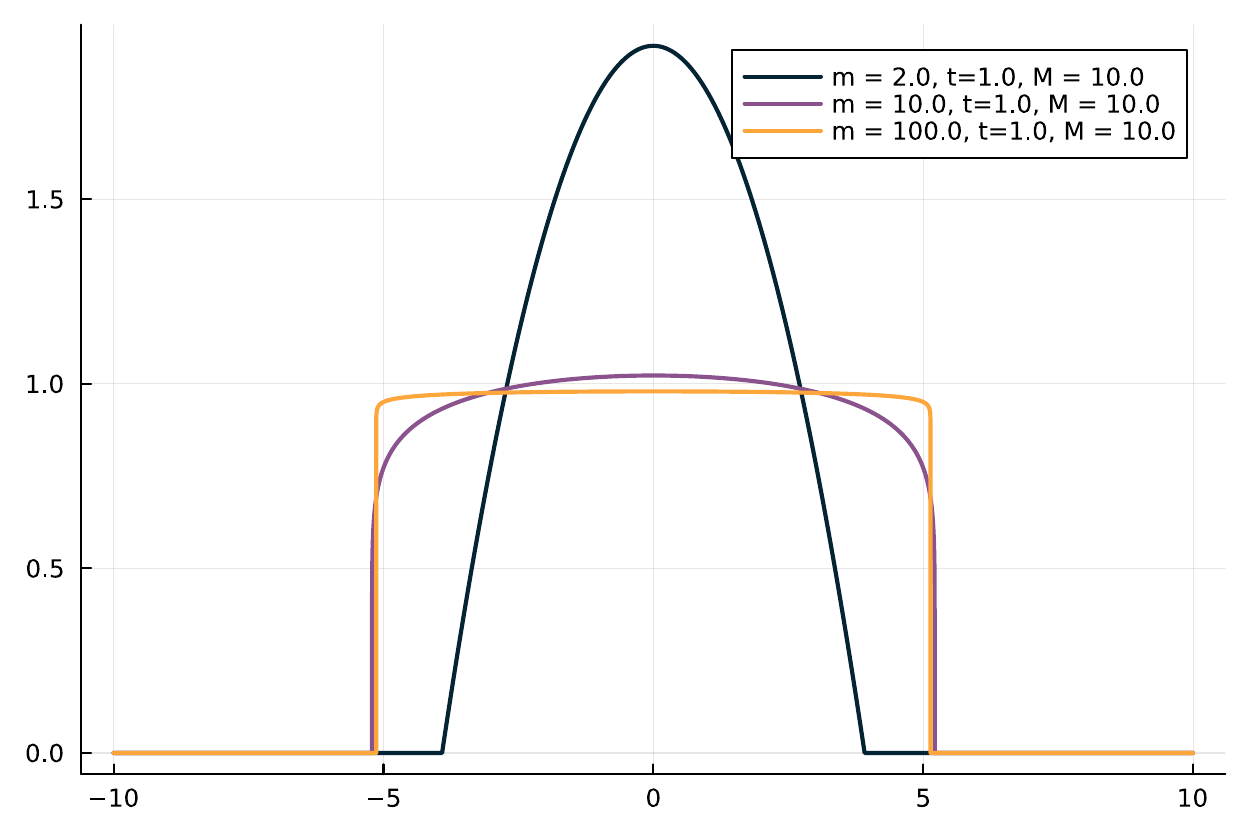}}
\vspace{-0.5cm}
\caption{ Barenblatt solutions for large $m$. Mesa formation }
\label{fig:PMEmesa}
\end{figure}

\noindent {\sl Proof.}
We need some precise calculations.
(i) If we recall the previous computation of the mass of $B_m (x,t;M)$ to get $M= d(m,N)\,C^{1/(2\beta(m-1))}$,
$$
M=N \omega_N \int_0^{\infty} (C-kr^2)_+^{1/((m-1)}r^{N-1}dr= d(m,N)\, C^{1/(2\beta(m-1))}.
$$
Besides, the formula \eqref{for.d_m.Bar} for $d(m,N)$ was
$$
d(m,N) =\omega_N  \frac{N}2  \frac{\Gamma(N/2) \Gamma( m/(m-1))}{\Gamma((N/2)+( m/(m-1))} k^{-N/2},
$$
 We can write it as
$$
d(m,N)=\omega_N\,\frac{N }{N + \frac{2}{(m-1)}}\, \frac{\Gamma(\frac{N}2)\,\Gamma( \frac{m}{(m-1)})}{\Gamma(\frac{N}2 + \frac{1}{(m-1)})} \,k^{-N/2}\,.
$$
For $m\sim \infty$ we get the estimate
$$
d(m,N)\sim \omega_N \,k^{-N/2}.
$$
Moreover, since \ $2\beta(m-1)\sim 2/N $ \ and \ $k=\beta(m-1)/2m\sim 1/(2mN)$,  we see that $d_m$ diverges like
$$
 d(m,N)\sim \omega_N (2Nm)^{N/2}.
$$

\noindent (ii) We move forward to calculate other constants. Since $\beta\sim 1/(Nm)$ we have
$$
\mathfrak{c}(m,N)=d(m,N)^{-2\beta}\sim (c m^{N/2})^{2\beta}\sim m^{1/m} \to 1.
$$
This settles 1 as supremum of the limit solution $B_\infty$ for all times $t>0$.

\noindent (iii)
On the other hand, the limit of the free boundaries is located at a distance $R_\infty$ such that $R_\infty^2$  is limit of
$$
\frac{C(M)t^{2\beta}}{k}= M^{2\beta(m-1)}\, d^{-2\beta(m-1)}\, t^{2\beta}k^{-1}.
$$
As $m\to \infty$, we note that $t^{2\beta}\to 1$ uniformly for $0<a<t<b$, so  that the expression tends to
$$
 M^{2/N}\omega_N^{-2/N}k^{\beta(m-1)N-1}.
$$
Since $k\sim (Nm)^{-1}$ and $\beta(m-1)N-1=2\beta\sim 2/(Nm)$, the last factor tends to 1 and we get
$$
R_\infty(M,t)=  M^{1/N}\omega_N^{-1/N}
$$
which is the radius of the ball with volume equal to $M$. It does not depend on time. For $M=1$, $L_\infty(M,t)=\omega_N^{-1/N}$.

(iii) It is  now easy task to conclude the convergence result. Assume to simplify that $M=1$ (no loss of generality). We can write the formula for $B_m(x,t)$ as
\begin{equation}\label{for.Bar,L}
B_m(x,t)=k^{1/(m-1)} t^{-\alpha}\left( L_1(m)^2 - |x|^2t^{-2\beta} \right)_+^{1/(m-1)}
\end{equation}
where $L_1(m)=L_1(m,N)$ is the length  defined  in \eqref{form,fbpme}.
Take now a point $x_1$ with $|x_1|\le L< L_\infty$. It is clear that as  $m\to \infty$
$$
k^{1/(m-1)}\sim (c/m)^{1/m}\to 1, \qquad t^{-\alpha}, t^{-\beta}\to 1,
$$
the time convergence is uniform for any interval $I=[a,b]$, $0<a<t<b$.
Finally,  $L_1(m)^2 - t^{2\beta}\,|x_1|^2$ approaches $  L_\infty^2 - |x_1|^2$ which is bounded below by a multiple of $L_\infty-L>0$. Since $1/(m-1)\to 0$ we get
$$
\left( L_1(m)^2 - |x_1|^2 t^{-2\beta}\,\right)^{1/(m-1)}\to 1.
$$
This argument shows the uniform convergence of $B_m(x,t) $ to 1 on compact sets inside the ball of radius $L_\infty$. We leave to the reader the verification of the convergence to zero outside of the ball $B_{ L_\infty}(0)$.\qquad \nc \qed

\medskip

The new pattern (mesa pattern, see Figure \ref{fig:PMEmesa}) shows in particular that all values of $B_m(x,t;M)$ larger than 1 have disappeared because of the extremely strong propagation of the equation for very large $m$ when $u>1$, as reflected in the diffusivity coefficient $D(u)=mu^{m-1}$\nc,  that creates a sudden collapse at $t=0+$  (i.e., in a narrow strip just after $t=0$). It is interesting to see this new behaviour in semigroup terms as follows: in passing to the limit $m\to \infty$  we have lost the basic continuity property $C([0,T]: L^1(\ren))$ of the solutions constructed so far; such continuity does not hold for $B_\infty$ at $t=0$ (though is trivial for $t>0$!).

\noindent $\bullet$  The theory of the PME in the limit $m\to\infty$ has been developed by different authors for solutions with data in $L^1$,  see \cite{BBH81, BenIgb03, CaffFried87, Elliott86, FdmnHoll87, Sacks89},  and we will not develop it further here. Let us just say that by passing to the limit on the PME equations with finite $m$ and fixed initial data $ f(x)\ge 0$ we can define a kind of PME-$\infty$ flow that turns out to be constant for all times $t>0$.
$$
\lim_{m\to \infty} u_m(x,t)= u_\infty(x).
$$
The unique identification of $u_\infty(x)$ in terms of $f$ is then the concern of the theory. Here  a mesa appears again, now supported on  a set $\Omega$ obtained by solving  a variational inequality.
We will only mention here an illustrative case example taken from the mentioned literature.

\begin{theorem}\label{limit.f.minfty}
Let $u_m(x,t)$ be the solution of the PME  with exponent $m>1$ having initial data
$$
u_m(x,0)= f(x),
$$
where $f$ is a nonnegative and continuous function, that is radially symmetric ($f=f(|x|)$,  decreasing in $r=|x|$, and  such that $f(0)=H>0$. Then,

(i) If $H\le  1$ we have
$$
u_\infty(x):= \lim_{m\to \infty} u_m(x,t)=f(x) \quad \forall x\in\ren, \ t>0.
$$

(ii) If $H> 1$  there exists an $R>0$ such that for all $t>0$
\begin{equation*}
u_\infty(x):= \lim_{m\to \infty} u_m(x,t) =
\left\{
\begin{array}{l}  1  \qquad \quad   \forall x: |x|<R,  \\
 f(x) \qquad \forall x: |x|>R.
\end{array}
\right.
\end{equation*}
Here $R$ is determined by conservation of  mass, $\int f(x)\,dx=\int u_\infty(x)\,dx$. We may express it as
\begin{equation}\label{form.mesa.R}
 R=\sup_{r>0}\,\{ \int_{|x|\le r}f(x)\,dx \ge \mbox{\rm Vol}\,(B_r)=\omega_N r^N\}. \nc
\end{equation}
\end{theorem}

The next remarks discuss the issue of asymptotic attractors.

\noindent {\bf Remarks. } 1)
We see from this result that no collapse happens when $H\le 1$ while a very definite collapse with re-distribution of the mass happens when $H>1$. In the latter case the mesa is supported on $\Omega=B_R(0)$. We also see that in this case $u_\infty$  has a negative jump, $f(R)-1$, at $x=R$.

2)    Note that $f(R)=0$ may occur, for instance has a very large mountain shape near the origin, so that the collapse overruns the original tail of $f$. In that case $ u_\infty(x)$ is a clean mesa, which will be the asymptotic attractor for all those initial data. Il looks like the typical asymptotic behaviour of a Barenblatt solution.

3) The mesas are not the only asymptotic attractors. Indeed, for data $f(x)$ with  a positive tail that goes beyond $R$  we discover that the mesa is not the asymptotic attractor, in contradiction to the result  on unique profile of asymptotic behaviour stated in  \eqref{form.asb.bar} for the PME with finite $m$. We may conclude that the mesa is determined by the initial mass distribution, see formula \eqref{form.mesa.R}, and does not encapsulate (i.e., attract) data that are partly supported outside of this calculated support.

3) Insisting on the issue. For nonradial data $F(x)$ the set $\Omega$ is obtained as the noncoincidence set of a variational inequality and can present many different geometries. A simple example when the nonradial data $F(x)$ consists of two different copies of the initial data of Theorem \ref{limit.f.minfty} located at sufficiently different points,
$$
F(x)=f(x-x_1)+ f(x-x_2).
$$
The stabilization will now happen towards two disjoint ball-shaped mesas. Same construction with any number of disjoint mesas\nc. \qed

\medskip

 We address the reader to the above references for more general versions of this result of Theorem \ref{limit.f.minfty} for data in $L^1(\ren)$, as well as and detailed explanations. For an application of the idea $m\to \infty$ in the study of tumor growth in biology see  \cite{PerQV14}. The mesa problem with fractional diffusion has been studied by us in \cite{Vaz2015mesa}\nc.

\


\subsection{ Fast Diffusion. The critical exponent}\label{fde.cexp}
The self-similarity approach still works for a while when we pass to exponents $m<1$, the so-called {\sl fast diffusion range} of the
 PME equation $u_t=\Delta (u^m)$, which will be now called FDE. The name is due to the fact that the diffusion coefficient $D=mu^{m-1}$ is  very large \nc for small densities. Following up the derivation done just before and only taking into account that now $m-1<0$ we arrive at similar values for the similarity exponents and to a  self-similar, constant-mass  \nc solution $B(x,t)=B_m(x,t)$ given by
\begin{equation}\label{for.Bar2}
B(x,t)=t^{-\alpha}\left(C+\frac{\beta (1-m)}{2m}\frac{|x|^2}{t^{2\beta}}\right)^{-1/(1-m)}
\end{equation}
This formula  has a similar analytical form but a  very different shape from \eqref{for.Bar1} and this will greatly affect the general theory. Instead of compact support of the slow case, and the exponential decays of the Gaussian kernel, it has a power-like tail with decrease rate \
$B(x,t)\sim C(t)\,|x|^{-2/(1-m)}$ as $|x|\to \infty.$ This behaviour is called in stochastic processes a long tail or a heavy tail, and they have an importance in different  fields of application. The calculation is correct as long as $m\sim 1$, but when $m$ goes down we see that
both exponents
\begin{equation}\label{sim.exp.FDE}
\alpha=\frac{N}{2- N(1-m)}, \qquad \beta= \frac{1}{2- N(1-m)}
\end{equation}
grow, and in fact they diverge when $m$ goes down to the so-called \sl fast diffusion critical exponent\rm,  \cite{HerrPierre1984} \nc
$$
 m_c=(N-2)/N, \quad N\ge 2.
$$
\begin{figure}[ht!]
\centerline{\includegraphics[width=0.6\textwidth]{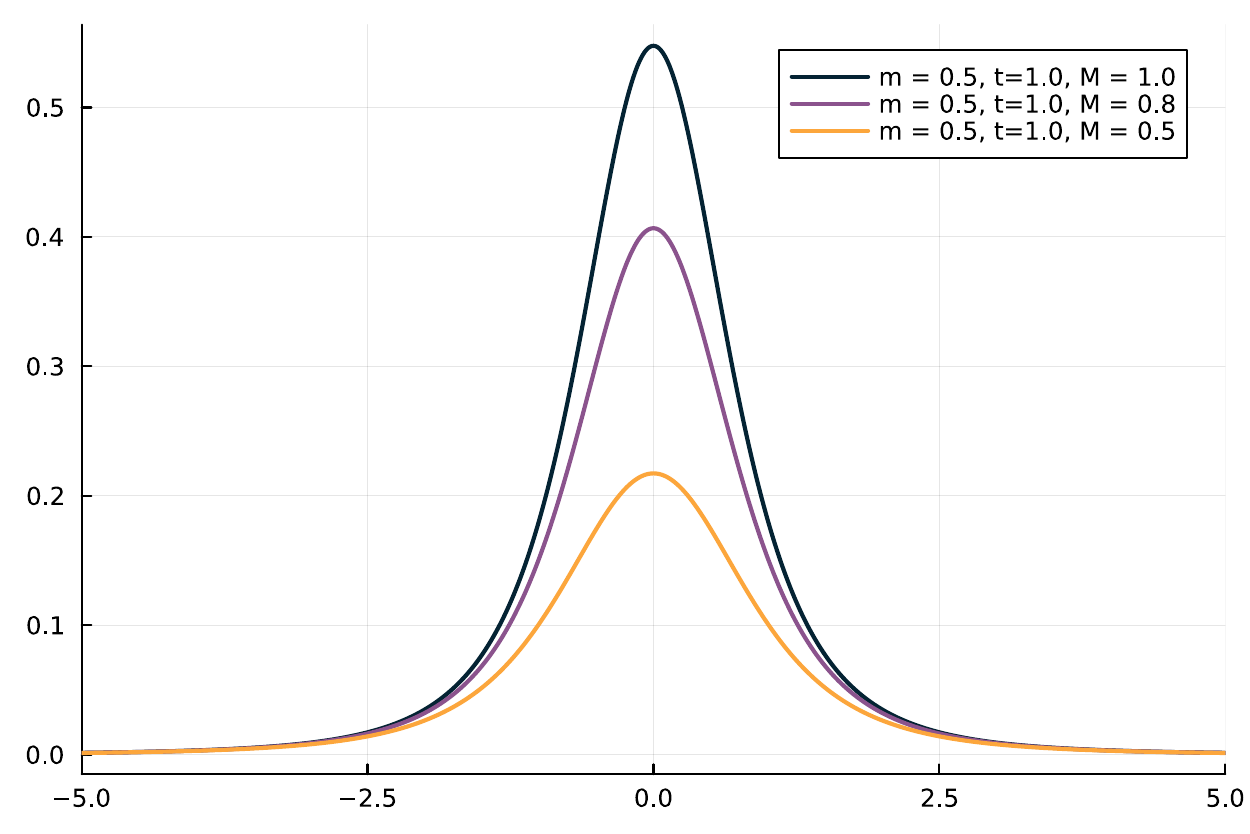}}
\vspace{-0.5cm}
\caption{Fast diffusion. Barenblatt profiles }
\label{fig:FDE}
\end{figure}

At the same time, the Barenblatt profile $B(\cdot,t)$ is integrable and conserves mass for $2/(1-m)>N$ which means precisely $m>m_c$. Therefore, we have found a fast diffusion Barenblatt solution in the \sl good fast diffusion range \rm $m_c<m<1$ but the calculation does not work for $m\le m_c$. In other words, the use of the Barenblatt solution as a solid foundation for the description of finite mass conservation breaks down at $m_c$.

\begin{theorem}\label{prop.mc.fde} Let $m\in (m_c,1)$, $m>0$.  Let $u(x,t)\ge 0$ be the semigroup solution of the FDE with initial data $u_0\in L^1(\ren)$, $u_0\ge 0$. Then, for every $t>0$ the mass is conserved:
\begin{equation}\label{mc.fdemc}
\int_{\ren} u(x,t)\,dx=\int_{\ren} u_0(x)\,dx.
\end{equation}
\end{theorem}

\noindent {\sl Proof.} We  show here a new  way to prove mass conservation for the solutions of FDE with $ L^1(\ren)$ data.
 Let $m_c<m<1$, and let us consider a  weak solution $u\ge 0$ of the FDE  defined in  $\re^N$ for a period of time $[0,T]$. We start from the equation, multiply it by a smooth test function \  $0\le \vp(x)\le 1$   with compact support in  $\re^N$, and integrate to get
$$
\frac{d}{dt}\int_{\ren} u\,\vp\,dx=\int_{\ren}  u^{m}\Delta\vp\,dx=\int_{\ren}  u^{m}\vp^{m}\frac{\Delta\vp}{\vp^{m}}\,dx
$$
Let us call $X(t)=\int_{\ren}  u\,\vp\,dx$  the weighted mass calculated at a certain time at $t$. Then, using H\"older's inequality we get
$$
 \left|\frac{dX}{dt}\right|\le \nc   X(t)^{m}\, Y^{1-m}, \quad
Y=\int_{\ren}  \left|\vp^{-m}\Delta\vp\right|^{1/(1-m)}\,dx.
$$
Note that  $Y(\vp)$ is a  positive  number that depends only on $\vp$, $m$ and $N$. We thus get an ODE inequality that will allow the mass to be controlled.

We need to estimate $Y$ in a more precise way by choosing suitable  test functions supported in large balls $B_R(0)$.  We choose a smooth $\vp=\vp(|x|) $  that  is  supported in $K=B_R(0)$, is positive in $K$, $\vp(x)=1$ in $B_{R/2}(0)$, and $\vp$ can be continued  in a smooth way as $\vp=0$ outside of $B_R(0)$. We may now improve the formula for $Y$ as follows. Recall that
$$
Y=\int_{R/2}^R \left|\vp^{-m}\Delta\vp\right|^{1/(1-m)}\,r^{N-1}dr.
$$
We assume that, as $r=|x|$ approaches the border of the supporting interval $r=R$, $\vp(r)$  approaches zero \nc in a very flat way that we will make finite the integral in the definition of the $Y$.  It is easy to see that when the value $\vp(r)$  looks like $O((R-r)^k)$ with $k>1$ large enough as $r\to R$, the integral $Y$ is finite. Then, $\Delta \vp$ is bounded and
$$
\Delta \vp \sim C(R-r)^{k-2}.
$$
Hence, $|\vp^{-m}\Delta \vp|\sim C(R-r)^{k-2-m k}$, so that the integral of $Y$ is a bounded number if $k(1-m)>1$.

We now take a first $R>0$ and some $\vp$ and define $\vp_n(x)=\vp(x/n)$, which is supported in $B_{nR}(0)$
Then if $X_n=X(\vp_n)$, $Y_n=Y(\vp_n)$,  we get
\begin{equation}\label{lme.yn}
Y_n= \int_{nR/2}^{nR} \left|\vp_n(r)^{-m}\Delta\vp_n(r)\right|^{1/(1-m)}\,r^{N-1}dr= C(R,\vp)\, n^{N-\frac2{1-m}}.
\end{equation}
where $N(1-m)-2<0$ precisely when $m>m_c$. Moreover,
\begin{equation}\label{lme.ode}
\left|\frac{dX_n(t)}{dt}\right|\le
 X_n(t)^{m} Y_n^{1-m}.
\end{equation}
The latter is an ODE that allows to control the change of mass  in any interval $0<s<t<T$:
\begin{equation}\label{lme.hp}
\left|X_n^{1-m}(s)-X_n^{1-m}(t)\right|\le (1-m)Y_n |t-s|\le C_1\,  n^{N-\frac2{1-m}} |t-s|.
\end{equation}
Letting $n\to\infty$ we get $|\int u(x,t)dx-\int u(s,t)dx|=0$ $s,t>0$, which means conservation of mass.  \qed

\noindent $\bullet$  This proof is due to \cite{HerrPierre1984}. The idea was extended to anisotropic FDE (see Subsection \ref{ssec.10.2}) in \cite{VazVSS-2023}.\nc

Here are some complements on the topic.

\noindent $\bullet$ {\bf Interesting connections.}
We have found the typical Barenblatt patterns for fast diffusion. These probability distributions are known in other parts of mathematics by other names and have other uses. The Barenblatt functions are for instance the optimal Aubin-Talenti functions for the Gagliardo-Nirenberg inequalities in functional analysis, see \cite{Aubin1976, Talenti76}. In fact there is a deep relation between Barenblatt profiles and Aubin-Talenti functions for $m_c < m < 1$. Note that we are dealing here with a parabolic-elliptic connection.

The connection continues with other profiles for other ranges of $m<m_c$ that we are going to explore next. The related topic involves entropy methods, improved functional inequalities and stability results, and there are many other  and more recent contributions on these topics, we refer for instance to \cite{Jungel2016} and  \cite{BonDGV-PNAS10} and follow-ups. These very interesting and active topics fall outside of the  scope of this survey, limited by choice to present in detail the close connection of our main topics.

\noindent $\bullet$  The semigroup approach to the existence theory works as before and we have the
Aronson-Bénilan inequality \cite{Aron1979} with no change, and the Bénilan-Crandall \cite{BC81b}, now with converse sign:
$$
(1-m)t\,u_t \le u.
$$

\noindent $\bullet$ Once the self-similar solution is found in the good range $m_c<m<1$, the rest of the chartered path can be traveled. We can construct a theory of existence and uniqueness of solutions with integrable initial data, though in this case there is infinite propagation and no free boundaries, there are positive power-like tails instead. Conservation of mass is proved.  The asymptotic behaviour with attraction to the Barenblatt is proved  for all $u_0\in L^1(\ren)$ with $\int u_0\,dx\ne 0$ just as in \eqref{form.asb.bar} and \eqref{asymp.pme.Li}.

\noindent {$\bullet$}  Also the $L^1$-$L^\infty$ smoothing effect holds in the range $m>m_c$ in the same  formulation
\eqref{form.L1Linf}, and optimality is achieved only by the Barenblatt solution and its space translations under the same conditions of Proposition \ref{prop.2.3}\nc.

\medskip

\subsubsection
{\bf The very singular solution, VSS}

\noindent When $m_c<m<1$ there is an a.e. finite limit of the family of Barenblatt solutions $B(x,t;m,M)$ as the mass $M$ goes to infinity. Indeed, this limit has a simpler algebraic expression and enjoys more invariance transformations:
\begin{equation}\label{for.fpme.vss}
V(x,t;m)=\lim_{M\to\infty}B(x,t;m,M)=a(m)\left({t}/{|x|^2}\right)^{1/(1-m)}\,
\end{equation}
where $a(m)=(\frac{2m}{\beta (1-m)})^{1/(1-m)}>0$. The VSS is a classical solution of the FDE for all $x\ne 0$ and $t>0$, and has a standing singularity at $x=0$ for all $t>0$. It is not a solution in any region that includes $x=0$. The total mass is of course infinite for all  $t>0$. The theory of solutions of the FDE with singular points or regions and possibly infinite mass is performed in \cite{ChaVa}.

Note that this kind of solution does not exist for $m\ge1$ since then
$$
\lim_{M\to\infty}B(x,t;m,M)=\infty \quad \mbox{everywhere for } \  x\in\ren, t>0.
$$
On the other end of the good FDE range, $m\to m_c$ with $N\ge 2$, we have from formula \eqref{sim.exp.FDE} $\beta\to\infty$, hence  $a(m)\to 0$, and the \nc VSS disappears everywhere away from $x=0$.

\medskip

\noindent {\bf Remarks.} 1) \bc There is an interesting connection between the  VSS of the FDE in the range of exponents $m_c<m<1$ and the  Bénilan-Crandall inequality (BCI)
$$
(1-m)t u_t\le u.
$$
Indeed, the  VSS  satisfies the BCI with equality sign at all points where the VSS is finite, i.e., for $x\ne 0$. In other words,  it is  an optimal function for the  BCI in the strong sense. On the other hand, if we consider the Barenblatt solutions we have agreement with the inequality only at the spatial tails
$$
\lim_{|x|\to\infty} \frac{tB_t(x,t;m,M)}{B(x,t;m,M)}=\frac1{1-m},
$$
but the inequality is far from holding true for every finite $x$. Actually, $B_t<0$ for $x\sim 0$\nc.

2) The PME/FDE with exponent $m \le 0$ makes perfect sense as a parabolic equation when written in the form $\partial_t u=\mbox{div}(u^{m-1}\nabla u)$, and this \sl very singular diffusion \rm range has been treated in detail in our monograph \cite{Vazquez2006} where full references are given.

3) It is to be noticed that for $N=1$ we have $m_c=-1$, and the extra range $m\in (-1,0]$ can be treated as part of the good fast diffusion range. We will not expand on this issue here since it is not essential in what follows and has been treated elsewhere, as in  \cite{ERV88, Vazquez2006, Vnonex92}.

\medskip

The previous study of our questions in the range $0<m<m_c$, $N\ge 3$, called \sl very fast diffusion range\rm, is a completely different story that we will address below. Let us advance the main idea: there is no integrable Barenblatt solution and conservation of mass does not hold, see Section \ref{sec.break}. The borderline case is very special and will attract our special attention, see Section \ref{sec.crit.exp}.

\newpage

\section{Mass conservation and self-similarity for the PLE}\label{sec.mc.ple}

We continue the analysis of the same issues in the nonlinear setting by another well researched example, the evolution $p$-Laplacian Equation, $\partial_t u=\mbox{div}\,(|\nabla u|^{p-2} \nabla u)$, usually abridged as $\partial_t u=\Delta_p(u)$, and called here PLE for short. The same approach can be put to work line by line
and we get the same form of the self-similar solution, now with exponents
\begin{equation}\label{form.simexpo.ple}
\alpha=\frac{N}{p+ N(p-2)}, \qquad \beta= \frac{1}{p+ N(p-2)},
 \end{equation}
both for the superlinear case $p>2$ and for the good sublinear case that ranges from  $p=2$ to the exponent for which $\alpha$ and $\beta$ diverge, which is
\begin{equation}\label{form.pc}
p_c=\frac{2N}{N+1},
 \end{equation}
\nc a value between 1 and 2 for all dimensions. The step by step calculation now leads to the formula for self-similar profile given by
\begin{equation}\label{form.plap.l2}
 F(y;p)=\left( C_0+\frac{2-p}{p} \beta^{\frac{1}{p-1}}\nc\,|y|^{\frac{p}{p-1}}\right)^{-\frac{p-1}{2-p}}
 \end{equation}
 for $p_c<p<2$, and
\begin{equation}\label{form.plap.g2}
 F(y;p)=\left( C_0- \frac{p-2}{p} \beta^{\frac{1}{p-1}} \nc |y|^{\frac{p}{p-1}}\right)_+^{\frac{p-1}{p-2}}
 \end{equation}
for $p>2$, with $N(p-2)+p=1/\beta$ and  $C_0>0$ is an arbitrary  constant such that can be determined in terms of  the initial mass $M$.  They are again called the Barenblatt solutions \cite{Bar56, BarentSSSIA}. We have  Mass Conservation, existence of fundamental solution in self-similar form and asymptotic convergence just as before, but now for $p>p_c$.

\noindent  {$\bullet$} In the PLE theory there is analogue of the Aronson-Bénilan estimate and it  was obtained by Esteban-Vázquez in 1990, see \cite{EV90}, and applies to all nonnegative solutions of the PLE defined in the whole space $\ren$. We first define a new ``pressure variable'' by the formula
\begin{equation}\label{form.plap.press}
v=\frac{p-1}{p-2}\, u^{(p-2)/(p-1)}
 \end{equation}
 when $p\ne 2$, $p>p_c$. The estimate reads
$$
\Delta_p v\ge -\frac{C}{t}, \quad C= \max\{(p-1),1\}\, \alpha,
$$
with $\alpha$ as in the Barenblatt solution. This estimate is optimal for the Barenblatt solutions when $p<2$. Using the equation for $v$ (see formula (6) of \cite{EV90}) we get
$$
\partial_t u \ge -\frac{C u}{t}.
$$
On the other hand,  the proof of the Bénilan-Crandall homogeneity inequality   \cite{BC81b} does not suffer any change, it reads
$$
(p-2)t\,u_t +u\ge 0.
$$

\noindent  {$\bullet$} Also the $L^1$-$L^\infty$ smoothing effect holds in the range $m>m_c$ in the above formulation
\begin{equation}\label{form.L1Linf.p}
\sup_x u(x,t)\le c(m,N)\, M^{p\,\beta}\,t^{-N\beta}\,,
 \end{equation}
and optimality is achieved by the Barenblatt solution and its space translations. The proof given in \cite{V82port} relies on symmetrization for the PLE and use of the Barenblatt solution. Note:   It is not known whether other solutions attain the optimal constant\nc.

\noindent {\bf  Asymptotic behaviour. Attraction to the Barenblatt.} \bc
Following the path used in the case of the PME,  we may prove that for any $u_0\in L^1(\ren)$ with $\int u_0\,dx\ne 0$ we have the asymptotic behaviour  as $t\to \infty$ of the informal form $u(\cdot,t)\sim B(\cdot,t, M)$, thus showing that the PLE Barenblatt solution also has a universal attraction power inside the solutions with the same mass. The precise formulation is formally similar to the one exhibited in Section \ref{ssec.pme} for the PME in. The $L^1$ version is
\begin{equation}\label{form.asb.bar.p}
\lim_{t\to \infty}\|u(\cdot,t)- B(\cdot,t;M)\|_1=0.
\end{equation}
where now $B(\cdot,t;M)$ is the Barenblatt solution of the PLE given in \eqref{form.plap.l2}, \eqref{form.plap.g2}.  The version of the asymptotic convergence result happens in $L^\infty$ norm and then a normalization factor is needed for optimal formulation:
\begin{equation}\label{asymp.pme.Li.p}
\lim_{t\to \infty} t^{\alpha}\|u(\cdot,t)- B(\cdot,t;M))\|_{\infty}=0,
\end{equation}
with $\alpha$ the PLE exponent found in \eqref{form.simexpo.ple}.
Convergence in all $L^q$ norms follows by interpolation for all $1\le q\le \infty$ with suitable time factors\nc.

\noindent  {$\bullet$} {\bf The 4 step method.} A main item in the proof that the Barenblatt fundamental solutions are attractors in the asymptotic theory is the uniqueness of such solutions. In the case of the PLE it was proved by Kamin and the author in \cite{KaminVaz1988} using in a strong way the explicit form of the existing self-similar candidate. Uniqueness is one of the key items in the 
 ``4-step method'' introduced in that paper \nc to study asymptotic behaviour of similar nonlinear parabolic equations. It proceeds as follows: \ after establishing the suitable existence theory and basic properties, we go through the following steps, formulated in \cite{KaminVaz1988}:

 (1) existence of a scaling group that transforms the solution $u(x,t)$ into new set of solutions, $u_n(x,t)$.  We need the initial data $u_n(0)$ to converge to a Dirac mass located at the origin.

 (2) finding a priori estimates that apply to the sequence $\{u_n\}$ of transformed solutions with uniform bounds,

 (3) finding enough compactness to pass to the limit $n\to\infty$ to get a limit solution $\overline{u}$ for $t>0$. The theory must show that the limit $\overline{u}$ is an acceptable solution. We check that it takes on a Dirac delta as initial datum

 (4) using a uniqueness result to identify $\overline{u}$ with  the source-type solution $B$ we have constructed for this model equation, which is a self-similar function.

 An additional step is needed: the convergence $u_n \to \overline{u}(t)=B(t)$ is easily interpreted as the asymptotic result $u(t)\to B$ as $t\to\infty$. This is just reinterpreting the scaling property of the equation.

In other words, the existence of suitable scaling laws together with good a priori estimates allows for a strong link from singular behaviour near $t=0$ to asymptotic behaviour as  $t=\infty$. This is quite impressive.\nc

\medskip

\noindent {\bf Notes.} A detailed study of self-similarity for the PLE was done in \cite{BidVer2006}.
A  theory of solutions with possibly growing initial data (as $|x|\to\infty)$)  is done in \cite{DiBH90}.

\medskip

\subsection{ The very singular solution of the Fast PLE}
Proceeding as in the Fast PME case, we pass to the limit $M\to\infty$ in the Barenblatt formula \eqref{form.plap.l2} to get
\begin{equation}\label{form.vss.plp}
V(x,t;p)=\lim_{M\to\infty}B(x,t;p,M)=a(p,N)\left({t}/{|x|^p}\right)^{1/(2-p)}\,
\end{equation}
with some $a(p,N)>0$ easy to calculate from formula \eqref{form.plap.l2}. As in the PME case, this solution does not exist for $p>2$, and on the other hand it tends to zero uniformly away from $x=0$ when $p$ tends to the critical value $p_c$ for $N\ge 2$.
\nc


\subsection{Limit $p\to \infty$ of the profiles. Sand piles}

\noindent {\bf Relation between the mass $M$ and the $C$ in \eqref{form.plap.g2}.}
For $p>2$ we have
$$
M=N\omega_N\int_0^\infty \left( C- k   r^{\frac{p}{p-1}}\right)_+^{\frac{p-1}{p-2}}r^{N-1}dr
$$
with  $k= \frac{p-2}{p} \beta^{\frac{1}{p-1}}$. Hence, with $r^{\frac{p}{p-1}}=(C/k)s$ we get
$$
M=N\omega_N \frac{p-1}{p} C^{\frac{p-1}{p-2}} (C/k)^{N \frac{p-1}{p}} D_p
$$
where
$$
D_p=\int_0^1 (1-s)_+^{\frac{p-1}{p-2}} s^{N\frac{p-1}{p}-1}ds= \mathcal B\left(\frac{2p-3}{p-2}, N\frac{p-1}{p}\right)
$$
where $\mathcal B(x,y)$ denotes Euler's Beta function.
This  determines  $C=C_p$ as a function of $M$.
In the limit $p\to\infty$ we get
$$
\lim_{p\to\infty} D_p= D_\infty=\int_0^1 (1-s) s^{N-1}ds=\frac1{N(N+1)},
$$
hence the mass formula becomes
$$
M=\frac{\omega_N}{N+1} C_\infty^{N+1},
$$
where $C_\infty= \lim_{p\to\infty} C_p$.

\begin{figure}[ht!]
\centerline{\includegraphics[width=0.6\textwidth]{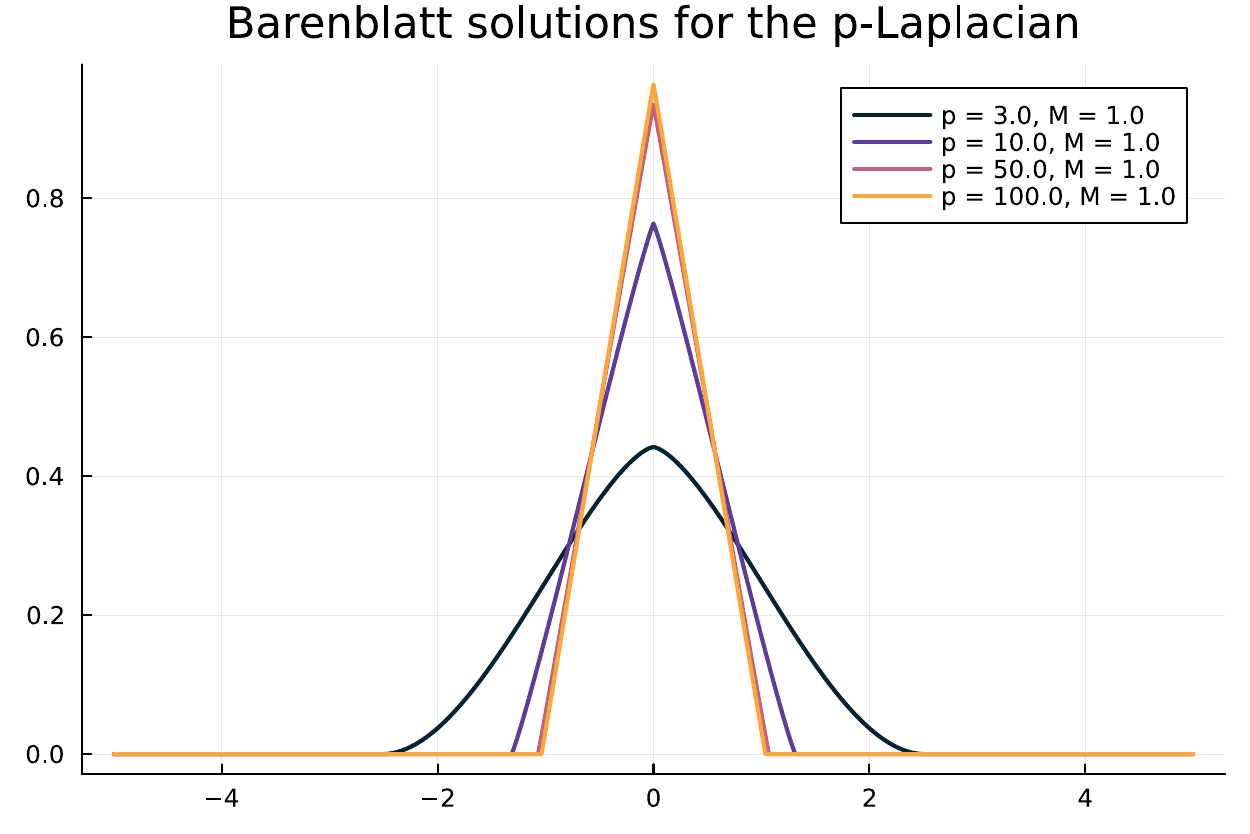}}
\vspace{-0.0cm}
\caption{\bc Formation of a conic sandpile in the limit $p\to\infty$}
\label{fig:plap}
\end{figure}

\noindent {\bf Limit.}
Once these preparations are ready, we may consider the problem with diffusion with large exponents already considered in Subsection \ref{ssec.pme.minfty} and take the limit as $p\to\infty$ of the Barenblatt solutions obtained for the PLE, see \eqref{form.plap.g2}. We are looking the patterns that replace in this case the mesas of the PME. A light effort shows that for the PLE The limit $p\to\infty$  is a sand pile function
\begin{equation}\label{form.plap.cone}
F(y;\infty) = \left( C_\infty -  |y|\right)_+, \qquad C_\infty^{N+1}= \frac{N+1}{\omega_N} M.
 \end{equation}

The graph \ref{fig:plap} shows the evolution of a Dirac delta by means of the PLE for different  values of $p>2$, some of them very large.
The initial collapse and sandpile formation  of the PLE as $p\to\infty$ is to be compared with the initial collapse into a mesa of the PME as $m\to\infty$.

\noindent {\bf Remarks.} 1) \bc
A complete analysis of the collapse of the initial value problem with $u_0\in L^1(\ren)$, $u_0\ge 0$, was done in  \cite{A-E-W, EFG97}. It is shown that for $t>0$ a collapse of the gradient happens towards compatibles values $|\nabla u| \le 1$. If the initial data $u_0$ already satisfy this gradient restriction, $|\nabla u_0| \le 1$, then the limit flow (i.e., the limit as $p\to \infty$ of the $p$-flow) reduces to a stationary state.

2) The fact that the sand pile profiles are not universal attractors for $p=\infty$ PLE solutions with initial data $u_0\ge 0$ having finite mass (even we add the condition  of boundedness) is proved with counterexamples similar to those shown in the case of the limit $m=\infty$ of the PME\nc.

\medskip

\noindent {\bf Note on the $p$-Laplacian systems.} Different applications suggest the interest in considering degenerate parabolic systems formed out of nonlinear diffusion equations of the above types. An important contribution in that direction is \cite{DBFr85} where DiBenedetto and Friedman consider a $p$-Laplacian system  of the form
$$
\partial_t u^j - \mbox{div}\,(|D u|^{p-2} D u^j)=F_j, \qquad 1\le j\le m.
$$
Note that here $u=(u^1,\cdots,u^m)$ and $Du$ denotes the whole gradient matrix. Existence and Hölder estimates are established for the gradient of local weak solutions of this system. Similar systems with PME nonlinearities can be considered, see
\cite{BerGHP, BerKam}. We leave to the interested reader to revise the literature on systems of those types, and consider for such systems the theory we have explained up to now.

\medskip

\noindent {\bf Comment on generic patterns in nonlinear diffusion.} The Gaussian pattern with its exponential expression is one of the most popular and relevant patterns is many areas of mathematics, not only in diffusion or heat propagation. Now, when we consider the Heat Equation as a particular case of the PME/FDE family we see that the patterns in all the nonlinear cases involve power functions and not exponentials. It follows that in some sense powers and not exponentials are the functions that we expect to find generically in our pattern formation.  Of course, the exponential Gaussian of the HE appears as a special limit of the power expressions of the Barenblatt profiles\nc. This comment applies very often in the rest of this work.

\newpage
\section{Extending the MC and the contraction laws}\label{ssec.emc}

The content of this subsection is new. There is an interesting version of the law of mass conservation that deals with the wider class of solutions that may have infinite mass when the difference of their initial data has finite integral. \bc In the final subsection we prove a new  and important result on lack of contraction in $L^p$  spaces for the PME, with $p>1$\nc.

\subsection{\bf Solutions with growing data. Optimal classes}\label{ssub.opt.th}
We start by examining larger classes of solutions. We want to extend the class of solutions to include non-integrable data, that may even grow as $|x|\to \infty$. We examine three cases.

 \noindent (I) in the case of the Heat Equation we recall the Widder theory,  \cite{Wid44}. According to it for every  Radon measure $\mu_0 $ (no sign restriction) such that the total variation $|\mu_0|$ satisfies
\begin{equation}\label{opt.exist.he}
\int_{\ren} e^{-c|x|^2}\, d|\mu_0| < \infty \qquad \mbox{for some }  \ c\in \re,
\end{equation}
there is a classical solution of the HE in a certain time interval $T=T(\mu_0)$ taking initial trace $\mu_0$, and the solution is representable by convolution with the Gaussian heat kernel. If $\mu_0 \ge 0$ the condition is optimal for existence and uniqueness of  solutions $u\ge 0$. The solutions with non-integrable data can be obtained as monotone limits of the $L^1$ solutions introduced as part of the $L^1$ semigroup. Condition \eqref{opt.exist.he} also ensures  existence of similar solutions  with changing sign but then this growth condition is not optimal.

\medskip

 \noindent (II) {\bf The PME.}  We start with the assumptions $u_0\ge 0$ and $u_0\in L^1_{loc}(\ren)$. The optimal class of solutions was established in the B\'enilan-Crandall-Pierre theory, \cite{BCP84}. It takes initial  data in the class of locally integrable data defined by the averaged growth condition
\begin{equation}\label{opt.exist.pme}
\sup_{R\ge 1} R^{-(N+ \frac2{m-1})} \int_{|x|\le R} u_0(x)\,dx <\infty\,.
\end{equation}
In that class there exists a unique weak solution $u: C([0,T) : L^1_{loc}(\ren))$ for some $T=T(u_0)>0$ that is obtained as monotone limit of the $L^1$ solutions introduced in Section \ref{ssec.pme} as the $L^1$ semigroup.  They are distributional solutions and satisfy the properties listed in the B\'enilan-Crandall-Pierre theory. In particular, the map: $t\mapsto u(t)$ is not only continuous from $[0,T)$ into $L^1_{loc}(\ren)$, but also  from $[0,T)$ into the weighted space $L^1(\ren, \rho_{\alpha}),$ with
$$
 \rho_{\alpha}(x)=(1+|x|^2)^{\alpha}, \qquad \alpha> (N/2)+ 1/(m-1).
$$
When $u_0$ is replaced by a bounded  Radon measure $\mu_0 $, the corresponding growth condition guarantees existence and uniqueness. We refer to Chapters 12 and 13 of \cite{VazquezPME2007} for complete details and notations.

\medskip

\noindent (III) {\bf  The FDE.} In this case the optimal class of initial data $u_0\ge 0$ consists of all $u_0\in L^1_{loc}$, cf. \cite{HerrPierre1984}\nc. There are no restrictions on the spatial growth of the solutions, but the regularity properties are worse when $m<m_c$.

\subsection{\bf Extending the MC}

 We recall  that  standard HE/PME/FDE theory   contains the property of the $L^1$ contraction stated in \eqref{ssec.pme}. It asserts that for any two general solutions $u$ and $v$ of the $L^1$ semigroup we have, as in formula \eqref{Licont},
$$
\int_\ren |u(x,t)-v(x,t)|\,dx\le \int_\ren |u_0(x)-v_0(x)|\,dx, \qquad t>0,
$$
and the equality sign need not hold. Same result for the HE, and the FDE\nc.
However, we will show next that equality holds if the solutions are ordered.

\begin{theorem}\label{thm.emc} Let $u$ and $v$ two solutions of the HE, the PME, or the FDE for $m>m_c$. Assume that their initial data satisfy $u_0\ge v_0\ge 0$ and $\int_{\ren} (u_0-v_0)\,dx=M_0$ is finite. Then for every $t>0$ for which both solutions exist we have $u\ge v$ and
\begin{equation}\label{signed.nass.4.1}
M(t):=\int_{\ren} (u(x,t)-v(x,t))\,dx=M_0.
\end{equation}
We will call this principle the conservation of the relative mass. We are not asking $u$ or $v$ to be integrable themselves, or asking for some behaviour as $|x| \to\infty$.\nc
\end{theorem}

\smallskip

The generality of the growth of the solutions is what makes the result interesting. The author has used a similar type of result for a very fast diffusion equation,  \cite{BBDGV}, and that  motivated work on this result, which is new in its generality for the nonlinear cases we consider. \nc

\noindent {\sl Proof.~}  We deal separately with the different equations. The HE case is easy.

\noindent (A) {\bf The HE.} It follows quickly from the  theory. By the analysis of point (I) of Subsection \ref{ssub.opt.th} the optimal class of nonnegative initial data for existence and uniqueness of solutions is given in formula \eqref{opt.exist.he}. This is the framework where we will work. Then, we  observe that the difference of two solutions in this class is another solution, $w=u-v$. Besides,  $w_0\ge 0$ implies that $w\ge 0$. For nonnegative solutions of the HE with finite $L^1$ data the mass $M(t)$ \nc is conserved as we have seen.

\medskip

\noindent (B) {\bf The PME.} We recall the assumptions $u_0\ge 0$ and $u_0\in L^1_{loc}(\ren)$. Once this is accepted, we work in the optimal class of solutions satisfying formula \eqref{opt.exist.pme}. We have seen above the properties of these solutions.
We will prove a conservation of mass similar to the result easily obtained for the Heat Equation\nc, but now the difference of solutions is almost never  a solution, hence we will use a completely different approach.

(i)  We take $T>0$ a time a bit less that the existence time $T(u_0)$ of the solution $u$,  and hence valid also for $v$ (by comparison)\nc. Let $w=u-v$ the difference of the solutions.  The  comparison result, see  \cite{BCP84}, makes sure that $w\ge 0$ for all $t\le T$\nc. We know that at $t=0$ the  mass of $w$ is finite $\int w(x,0)\,dx=M<  \infty$. We call the \sl relative mass \rm of $u$ and $v$ the mass of the difference of both solutions\nc. The $L^1$-contraction is a very well-known property that applies to differences of solutions of the PME (and other equations),  and this also applies to solutions in the optimal class (passing to the limit from the $L^1$ semigroup)\nc, so we derive the fact that $w(\cdot,t)\in L^1(\ren)$ for all $t\in (0,T)$ and
$$
\int w(x,t)\,dx\le M.
$$
It follows that we only need a lower bound, i.e., avoiding to lose part of the relative mass with time.

(ii) In order to get this more precise result we will use a suitable technique for controlling the evolution of the local version of the relative mass, defined as
\begin{equation}\label{form.lrm.pme}
Y(t,\varphi)=\int_{\ren} (u(t)-v(t))\,\varphi\,dx\,,
\end{equation}
where $\varphi(x)$ is  a properly chosen positive and smooth test function. We will assume that $\varphi$ is radially symmetric, $0\le \varphi\le 1$,     $\varphi(x)>0$ everywhere, $\varphi(x)=1$ for $|x|\le R_0$. Moreover, we will ask for a fast decay condition: $\varphi(x)=c_1\,|x|^{-k}$ for $|x|\ge 2R_0$ for some $R, c_1>0$. The exponent $k>0$ is fixed and large enough. So $\varphi$  decays like $\rho_{\alpha}$ with $k=2\alpha$ in the notation used a moment ago. Thus, if  \ $k>N+(2/(m-1))$ the function $Y(t,\varphi)$ is a well-defined
finite Lebesgue integral and continuous in time $t\in [0,T)$.

We observe that there is a constant $C_1>0$ such that\nc
$$
 |\Delta\varphi(x)| \le C_1\frac{\varphi(x)}{1 +|x|^2}\nc.
$$
This is important in what follows.
On the other hand, it is known from the already cited theory by Bénilan-Crandall-Pierre \cite{BCP84} that for every solution $u$ and for times $t$ in the existence period and away from zero we will have
$$
|u(x,t)|^{m-1}\le C(t)(1+|x|^2)\,.
$$
Here, $C(t)$ is a continuous function in $(0,T)$ that depends on the solution. This and the decay  $\varphi(x)\le c_1\,|x|^{-k}$  guarantee that the relative mass $Y(t,\varphi)$ is a well defined finite Lebesgue integral for $0<t<T$ \ if \ $k>N+(2/(m-1))$\nc.

Working in a formal way that can be justified as in \cite{BCP84} or \cite{VazquezPME2007}, see below, we get
$$
\frac{dY}{dt}= \int \frac{\partial}{ \partial t}(u(t)-v(t))\,\varphi\,dx=\int \Delta (u(t)^m-v(t)^m)\,\varphi\,dx=
\int  (u(t)^m-v(t)^m)\,\Delta\varphi\,dx,
$$
hence,
$$
|\frac{dY}{dt}|\le |\int (u^m-v^m)\,\Delta\varphi\,dx|\le m\int |u|^{m-1}(u-v)\,|\Delta\varphi|\,dx.
$$
 Hence, for a.e. $t$
\begin{equation}\label{form.emc.pme}
|\frac{dY}{dt}| \le  C_2(t) \int (u(t)-v(t))\,\varphi\,dx= C_2(t) \,Y(t),\nc
\end{equation}
 where $C_2(t)= mC_1 C(t)$. \nc This gives a control of possible growth of the local relative mass,  telling us that
$Y(t)$  is not only bounded but also Lipschitz continuous inside $(0,T)$.

To justify the calculation we better pass to the weak version using as a test function in space time $\zeta(t)\varphi(x)$,
where $\zeta(t)$ is a positive and smooth test function approximating the characteristic function of $[t_1,t_2]$ with $0< t_1,t_2<T$.
Indeed, instead of  \eqref{form.emc.pme} we may get
$$
|\int_{0}^{T}\zeta'(t)Y(t)\,dt| \le  \int_{0}^{T} C_2(t) \zeta(t)Y(t)dt=  \int_{0}^{T}dt\int C_2(t) \zeta(t)(u-v)\,\varphi\, dx.
$$
In the limit when we let $\zeta(t)$ tend to the characteristic function   of $[t_1,t_2]$  we get
\begin{equation}\label{form.emcw.pme}
|Y(t_2)-Y(t_1)| \le  \int_{t_1}^{t_2}\int C_2(t) (u(t)-v(t))\,\varphi\, dxdt.
\end{equation}
This is the useful integrated form\nc.


(iii) But we may sharpen \eqref{form.emc.pme} or \eqref{form.emcw.pme} and reduce the Lipschitz constant by using a further trick. We recall that $\Delta \varphi=0$ for $|x|\le R_0$, hence
$$
|\frac{dY}{dt}| \le C_2(t) \int_{|x|\ge  R_0} (u-v)\,\varphi\,dx
$$
Since the  integral defining the relative mass  $M(t)$ is bounded in time, then for any given $t=t_1$ and $\ve>0$ we actually have
$$
\int_{|x|\ge R} |u(t)-v(t)|\,dx\le \ve
$$
for $R(\ve,t)$ large enough. We want to control  $R(\ve,t)$ so it does oscillate with time. To avoid this problem we pass to the weak form \eqref{form.emcw.pme} and control the double integral
$$
\int_{t_1}^{t_2}\int_{|x|\ge R} |u(t)-v(t)|\,dx\le \ve.
$$
This allows to justify the calculations for $0<t_1<t<t_2<T.$

(iv) We may now  combine both facts and pass to the limit after  changing $\varphi$ into $\varphi_n(x)=\varphi(x/n)$ so that in practice the role of $R_0$ is played by  $nR_0$ and we have
$$
 |\Delta\varphi_n(x)| \le C_1\frac{\varphi_n(x)}{n^2 +|x|^2}\le C_1\frac{\varphi_n(x)}{1 +|x|^2}
$$
with same constant as before. Putting $Y_n=Y(t,\varphi_n)$ and taking $n$ large enough so that $nR_0> R(\ve)$  we get
$$
|Y_n(t_2)-Y_n(t_1)| \le \int \int_{|x|\ge nR_0} C_2(t)  (u-v)\,\varphi_n\,dxdt \le \overline C\ve.
$$
where $\overline C$ is the supremum of $C_2(t)$ in the interval $t_1\le t\le t_2.$ In the limit $n\to\infty$ we have $\varphi_n(x)\to 1$ and $Y_n(t)\to M(t)$ everywhere in a monotone way. We get
$$
|M(t_2)-M(t_1)| \le \overline C\ve.
$$
Since $\ve$ is arbitrary we conclude that $M(t_2)=M(t_1)$. Therefore, $M(t)$ is constant for $0<t<T$. It means that no mass is lost at any time $t>0$.

(v) Finally, we have to check the initial situation to exclude the possibility that $M(t)$ is a constant $M_1$ less than $M$. We argue as follows. First, given $\ve > 0$ we select a large ball such that $\int_{B_R} (u_0-v_0)\,dx\ge M-\ve$. Next, since both solutions $u$ and $v$ are continuous as functions $C([0,T): B_R)$, there exists a small time $t_1$ such
$$
\int_{B_R} (u(t_1)-v(t_1))\,dx\ge \int_{B_R} (u_0-v_0)\,dx-\ve\ge M-2\ve.
$$
But $\int_{B_R} (u(t_1)-v(t_1))\le M_1.$ This implies that $M_1=M$, \qed

\medskip

\noindent (C) {\bf  The FDE.} In this case the optimal class of initial data consists of all $u_0\in L^1_{loc}(\ren)$.
 We have concentrated first on the proof for nonnegative data and solutions (which are our main interest). We also refrain from considering measures as initial data treated in Chapter 13 of \cite{VazquezPME2007}\nc.

 We will consider the fast diffusion range in the good range $m_c<m<1$.
We will  prove the result by a different method that gives a much stronger result about the local evolution of the mass function in the FDE case. This is given by a beautiful result due to Herrero and Pierre \cite{HerrPierre1984}. The idea is as follows. With the same notations as before, we change the calculation of the local relative mass into
$$
|\frac{dY}{dt}|\le |\int (u^m-v^m)\,\Delta\varphi\,dx|\le C\int (u-v)^m\,|\Delta\varphi|\,dx.
$$
where we have used the algebraic inequality $u^m-v^m\le (u-v)^m$ valid for real numbers $0<v<u$ and $0<m<1$ (Note: it is also true for $v<u<0$.) Now we may use H\"older's inequality to get
$$
\int (u-v)^m\,|\Delta\varphi|\,dx\le \left(\int (u-v)\varphi\,dx\right)^m
\left(\int \psi(x)\,dx\right)^{1-m}
$$
where $\psi= |\Delta\varphi|^{1/(1-m)}\varphi^{m/(1-m)}$. It is not difficult to find examples of cut-off functions such that the last integral becomes finite, note that it does not depend on $t$ or $u,v$. Therefore, we get
$$
|\frac{dY}{dt}|\le K(\varphi)Y(t)^m,
$$
that upon integration produces the following local rule\nc

\begin{lemma}\label{lem.emc} Under the present conditions on $u_0$ and for  $m<1$ \nc we have
\begin{equation}\label{form,emc.fde}
|Y(t,\varphi)^{1-m}-Y(s,\varphi)^{1-m}| \le  K(\varphi)|t-s|\, \quad \mbox{for all } \ 0 \le s,t.
\end{equation}
\end{lemma}

The way to obtain from this local result the conservation of relative mass consists in taking
$\varphi_n(x)=\varphi_1(x/n)$ so that after an easy calculation we get $\varphi_n\to 1$ in $\ren$ and
$$
K(\varphi_n)=n^{-\gamma}K(\varphi_1),  \qquad \gamma=N(m-1)+2.
$$
Now observe that $\gamma>0$ whenever $m>m_c$. This shows the strong influence of the critical exponent in the qualitative and quantitative theory of the FDE. See further details of this type of calculations in \cite{VazVSS-2023}.

Passing to the limit $n\to \infty$ we get the conservation law since $n^{-\gamma}\to 0$. \quad \qed

\medskip

Here is an extension of the last result.

\begin{corollary}\label{cor.emc} The result of Theorem \ref{thm.emc} holds when  the data are not ordered, if the condition that $\int |u_0-v_0|\,dx<\infty$ holds.
\end{corollary}

\noindent {\sl Proof.~}  We accept solutions with the same growth conditions on the initial data so that the existence theory works\nc. First, note that there is no difficulty in the linear heat equation.

(i) Let us consider solutions of the PME or the FDE in the case where $u_0$ and $v_0$ are not ordered, so that
$$
\int_{\ren} (u_0-v_0)_+\,dx=M_1>0,  \ \int_{\ren} (v_0-u_0)_+\,dx=M_2>0,
$$
and then $M=M_1-M_2$. We now consider the  auxiliary function
$$
w_0(x)=\max\{u_0(x),v_0(x)\}=u_0(x) + (v_0-u_0)_+\,.
$$
This initial datum  produces a solution $w(x,t)$  (of the PME or  FDE) such that  $w(x,t)\ge u(x,t),v(x,t)$, and the initial data satisfy $\int (w_0(x)-u_0(x))\,dx=M_2$, $\int (w_0(x)-v_0(x))\,dx=M_1$. We conclude from the previous theorem that for all $t>0$
$$
\int_\ren (w(x,t)-u(x,t))\,dx=M_2,
\int_\ren (w(x,t)-v(x,t))\,dx=M_1,
$$
Hence, we arrive at the invariant quantity
$$
\int_\ren (u(x,t)-v(x,t))\,dx=M_1-M_2.
$$

(ii) The case where solutions have any signs. There is no problem for the heat equation. In the case of the PME we need to recall that $u^m$ must be replaced by $|u|^{m-1}u$. We recall that the Bénilan-Crandall-Pierre covers the case of signed solutions and the estimates that we need are similar. Moreover, we have inserted absolute values whenever they may help in the minor changes needed in the proof.

Similar arguments apply to the FDE proof with minor changes. Thus, when $v<0<u$ we need to use
the inequality $(u^m-v^m)= (u^m+|v|^m)\le C(m)(u+|v|)^m$. \qquad \qed

\medskip

\noindent {\bf Note.}  Mass conservation for the critical exponent will be investigated in Section \ref{sec.crit.exp}. MC will not hold for $m<m_c$ as we explore in Section \ref{sec.break}\nc.


\subsection{Lack of contraction in $L^p$ spaces for the PME}\label{ssec.lack}

We have mentioned that the PME generates a contraction semigroup with respect to the $L^1$ norm  as a basic property of the theory. We may wonder if this property may  be  extended to a similar contraction of the solutions of the PME with respect to \nc some or all the $L^p$ norms, $1<p\le \infty$. We recall that the property of contraction with respect to all  $L^p$ norms is true for the HE and it also well known in the case of the PLE $u_t=\Delta_{q}(u)$ for all $1\le p\le \infty$ and $1\le q\le \infty$.

We will show that the same conjecture for the PME is not justified in view of the examples that we present next.  This is new material to the author's knowledge and settles in part an old question\nc.

 \noindent $\bullet$ {\bf Example 1.} We will show that  contraction could be very far from happening in the case of the $L^\infty$ norm for the PME. We produce  an example in the class of solutions growing at spatial infinity. \nc We take the PME with exponent $m=2$ and consider the family of blow-up solutions
\begin{equation}\label{form.blowup.pme}
U(x,t; C,T)=\frac1{T-t}\left(C\,(T-t)t^{2\beta}+ k\,x^2\right)=k\frac{x^2}{T-t}+ C\,(T-t)^{-\frac{N}{N+2}},
\end{equation}
which is obtained by a known modification of the standard Barenblatt solutions \eqref{for.Bar1}. Remember that $\beta=1/(N(m-1)+2)=1/(N+2)$.   We readily observe an interesting fact: for $C_1, C_2\ge 0$ the difference of the two solutions $U_1(x,t)= U(x,t; C_1,T)$ and $U_2(x,t)= U(x,t; C_2,T)$ is constant in space and grows in time. More precisely, we have
$$
U_1(x,t)-U_2(x,t)=(C_1-C_2)\,(T-t)^{-\frac{N}{N+2}},
$$
and the difference blows up as $t\to T$ with a precise rate. Thus, if $C_1>C_2>0$ we have $\|U_1(0)-U_2(0)\|_\infty= C_1-C_2$ and
\begin{equation}\label{form.lack.lp}
\|U_1(t)-U_2(t)\|_\infty=  K(t)\,\|U_1(0)-U_2(0)\|_\infty, \quad K(t)= (1-(t/T)))^{-\frac{N}{N+2}},
\end{equation}
for $0<t<T.$ Therefore, the $L^\infty$ norm  of the difference not only grows in time, it even blows up at a finite time $T>0$ that can be chosen at will. The example implies the following consequence involving bounded solutions.

\begin{proposition}\label{prop.lackLp} For every $C, t_1>0$ there exist integrable  solutions of the PME $u_1, u_2$ defined in a time interval $0\le t\le  T$ with $T>t_1 $, and having bounded initial data $u_{1}(0)\ge u_2(0)\ge 0$, such that
\begin{equation}\label{form.lack.lp.v2}
\|u_1(t_1)-u_2(t_1)\|_\infty\ge C\|u_{1}(0)- u_{2}(0)\|_\infty.
\end{equation}
\end{proposition}

The justification proceeds by approximation of the growing initial data of the Example with a couple of increasing sequences of bounded functions, solving the PME for both data in the class of bounded solutions, $u_{1,n}\ge u_{2,n}\ge 0$  and examining that happens for $n\sim \infty$ with  formula  \eqref{form.lack.lp}. Actually we can use bounded $L^1$ functions for the approximations. We leave the details as an exercise\nc.

\medskip

 \noindent $\bullet$ {\bf Example 2.} A more elaborate counterexample can be found in some space $L^p(\ren)$, $p$ finite, for $m>2$. Using the same type of solutions  and using a similar idea, we may obtain differences that have a finite $L^p$ norm for some large $p$ at any time $0<t<T $, but the said norm blows up as $t\to T$.

For $m>2$ the formula of the pressure is the same, so that
$$
U_1-U_2= P_1^{\sigma}-P_2^{\sigma}\sim (P_1-P_2)/P^{1-\sigma},
$$
where $\sigma =1/(m-1)<1$, $1-\sigma=(m-2)/(m-1)\in (0,1)$, and $P$ is an interpolation between $P_1$ and $P_2$. Now we write
$$
P_i(x,t)=F_i(x,t)/(T-t), \qquad F_i(x,t)= C_iA(t)+kx^2,   \qquad A=(T-t)^{2\beta}.
$$
Therefore,
$$
\|U_1-U_2\|_p=c(C_1-C_2)(T-t)^{-\sigma+ 2\beta} Y(t),
$$
where
$$
Y(t)^p=\int (A +kx^2)^{-p(1-\sigma)}\,dx=c(p,N) A^{N/2-p(1-\sigma)},
$$
and $c(p,N)$ is finite if $p$ is large,  more precisely if $2p(1-\sigma)>N$, i.e., $p>p_*= N(m-1)/2(m-2)$,  $m>2$.

It follows that $\|U_1(t)-U_2(t)\|_p \approx (T-t)^{-\mu}$ where
$$
-\mu=-\sigma+ 2\beta +2\beta(N/2p)-2\beta(1-\sigma)=
$$
$$
-\sigma  +\beta(N/p) +2\beta\sigma=-N(m-1)\beta\sigma + \beta(N/p)=-N\beta(1-(1/p)),
$$
so that $\mu=N\beta(p-1)/p>0$. Therefore, $\|U_1(t)-U_2(t)\|_p $ is finite for $t<T$ and diverges as $t\to T$ if $p>p_*$ and $p>1$. In other  words, it means $p>N/2$ and $m-2>N/(2p-N)$. For example, we can obtain a growing estimate in $L^2(\ren)$ if $N=3$ and $m>5$, if $N=2$ and $m>3$, or if $N=1$ and $m>7/3$\nc.

\medskip

\begin{open} Of course, in the case of the heat equation the difference of solutions is a solution and all $L^p$ norms, $p>1$, decay in time. The question of finding a counterexample to the $L^p$ contraction with $p>1$ remains open for the PME in the range $1<m<2$, and for fast diffusion with $0<m<1$\nc.
\end{open}

\section{Breakdown of the MC for very fast diffusion}\label{sec.break}

The study of our set of questions for the  corresponding  very fast diffusion range of the FDE, or the fast diffusion range of the PLE, is a completely different story that we will address next. Let us advance two key ideas: there is no integrable Barenblatt solution and conservation of mass does not hold. Researchers  investigate what is the fate of that waning mass and what is the shape that those declining solutions mostly take, and we will address those issues here.
The borderline cases, $m=m_c$ and  $p=p_c$, are very special and will be dealt with later, see Sections \ref{sec.crit.exp} to \ref{sec.ncbl}.

An important idea behind this break-down is the extreme character  of the diffusion coefficient for very small densities \nc $u$. Indeed, since the PME can be written as a linearized heat equation \ $\partial_t u=\mbox{div}(D(x,t)\nabla u)$ \ with variable coefficient $D(u)=mu^{m-1}$, this coefficient goes to infinity when $u\to 0$ (and to zero as $u\to\infty$) whenever $m<1$. Moreover, when $m$ is far from 1 the variability of $D(u)$ causes a new phenomenon. Namely, it may happen that the diffusion  turns out to be too fast, for instance at the tails of integrable solutions, so that there is a catastrophe in the making: the  escape of some mass to infinity. We will examine next the quantitative aspects of these ideas.

\subsection{Extinction in the very fast diffusion range of the FDE}\label{ss.vfpme}

The theory of the PME/FDE posed in the whole Euclidean space departs strongly from the paradigm of mass conservation and corresponding self-similarity described in previous sections when we pass to the lower range \, $0<m<m_c$ (with $N\ge 3$). To start with, the existence and uniqueness of a class of integrable solutions for $u_0\in  L^1(\ren)$ can still be obtained by the Semigroup  Generation Theorem  as explained by B\'enilan and Crandall in  \cite{BC81}. We will call them the semigroup solutions. They are integrable but they need not be bounded for $t>0$ as happened in the PME. $L^1$ contraction still happens. For positive and bounded initial data they are $C^\infty$ smooth classical solutions for $t>0$, see for instance for details \cite{Vazquez2006}\nc.

However, there is a very striking result, called \sl extinction in finite time\rm, that departs from mass conservation as far as possible.  The basic quantitative result is due to \cite{BC81}. We rephrase it in our notation.

\begin{proposition}\label{thmBR.ext} Let $u(x,t)$ be a semigroup solution of the PME in the very fast diffusion range $0<m<m_c$ and let us assume that the initial data $u_0\in L^{p_*}(\ren)$ with $p_*=N(1-m)/2>1$. Then there is a finite time $T$ such that
$u(x,t)\equiv 0$ for all $t\ge T$. This is called the extinction time, and is has the quantitative estimate
$$
T\le c(m,N)\|u_0\|_{p_*}^{1-m}.
$$
\end{proposition}

\noindent {\sl Sketch  of proof.~} For a nonnegative solution the calculation is based on a simple estimate of the evolution of the $L^{p_*}(\ren)$ norm of $u(\cdot,t)$ that is shown to satisfy a nonlinear differential inequality: $X'+ CX^{(N-2)/N}\le 0$, where $X(t)=\|u(\cdot,t)\|_{p_*}^{p_*}$. This result was first proved by B\'enilan and Crandall \cite{BC81} and a proof is contained in the results about finite-time extinction  \nc of Chapter 5 of \cite{Vazquez2006}.  Just recall that the semigroup solution is continued as $u(t)\equiv 0$ for all $t\ge T(u_0)$\nc. We will explain the calculation below in the PLE case so that the reader will see the details in a similar situation. For that reason and to avoid repetitions, we skip them here.

 For a signed solution we write $-u_{0-}\le u_0\le u_{0+}$ with $u_{0+}=\max\{u_0,0\}$, $u_{0-}=\max\{-u_0,0\}$, and solve the FDE semigroup to find solutions $v(t)=S_t (u_{0-})$, $w(t)=S_t (u_{0+})$, both are nonnegative solutions to which the Proposition applies. Then use comparison to deduce that $-v(t)\le u(t)\le w(t)$. The result follows after observing that the extinction time is bounded above by the maximum of extinction times of $v$ and $w$. \nc\qed

This result is completed by a more general result where mass conservation is excluded for all nonnegative solutions. A proof of this result is given in \cite{VazFPL3-2021} applied to the fractional $p$-Laplacian equation. We will repeat the short proof there since it is easy and the application is quite general.

\begin{proposition}\label{prop.massto0} Let $0<m<m_c$ and let $u$ be a semigroup solution of the PME such that $u(x,t)\ge 0$ and it has finite initial mass, $\int_{\ren} u(x,0)\,dx=M>0$. Then,
$$
M_u(t):=\int_{\ren} u(x,t)\,dx  \to 0 \quad \mbox{as \ } \ t\to\infty.
$$
In other words, $u$ undergoes mass extinction in either finite or infinite time.
\end{proposition}

\noindent {\sl Proof.~}  Given $u_0(x)=u(x,0)\in L^1(\ren)$ we approximate it in $L^1$ by a sequence of bounded and integrable functions $u_{0,n}$. By the $L^1$ contraction this allows to get the  approximate solutions $u_n(x,t)$. By the contraction property, for every $t>0$ we have
$$
\|u(t)-u_n(t)\|_1\le \|u_{0}-u_{0,n}\|_1\le \ve_n\to 0.
$$\nc
According to Proposition \ref{thmBR.ext}  we know that the corresponding solution $u_n$ has a finite extinction $T_n$ time that depends on the norm of $u_{0,n}$ in $L^{p_*}(\ren)$.
Since $u_n(T_n)\equiv0$,
$$
\|u(T_n)\|_1 = \|u(T_n)-u_n(T_n)\|_1\le \ve_n.
$$
 Finally, we recall the basic property that the mass is a nonincreasing function of time to continue the bound
 $\|u(T_n)\|_1 \le \ve_n$ for all $t\ge T_n$.\qed

\medskip

\noindent {$\bullet$ \bf The fate of finite-mass solutions.}
One may wonder what is the role of self-similarity in this completely changed scenario. \bc We are expecting solutions with loss of mass and even extinction. Let us just mention the existence of the self-similar solutions called pseudo-Barenblatt solutions, see below.

The existence of self-similar solutions of the backwards-in-time type  and with finite mass \nc was first studied by King in his seminal paper \cite{Ki93}, where the main asymptotic results on fast diffusion were examined and surveyed by formal methods, leading in particular to the existence of a type of self-similar solutions with extinction, hence no mass conservation.
The similarity exponents $\alpha$ and $\beta$ are not determined from any conservation law, but from PDE considerations (existence of solution of a nonlinear eigenvalue problem), and it is called an {\sl anomalous exponent}, a favorite topic in Barenblatt's research and teaching, \cite{BarentSSSIA}. Barenblatt also uses in that respect the label \sl self-similarity of the second kind \rm or \sl type-II self-similarity\rm.

A rigorous analysis was then performed by  M. Peletier and Zhang \cite{PZ}, who showed the existence of the relevant asymptotic profiles as a consequence of a careful phase-plane analysis. A complete analysis of the issue is contained in Chapter 7 of our monograph \cite{Vazquez2006}, see specially Section 7.2. Summing up, there are plenty of self-similar solutions but there is no self-similar solution with the property of conserving a finite mass.

\medskip

\noindent {$\bullet$ \bf A geometrical flow.} There is one special exponent, $m_s=(N-2)/(N+2)$ (for $N\ge 3$) for which the fate of finite-mass solutions is described by a self-similar solution with simple exponents, $\alpha=1/(1-m)=(N+2)/4$, and $\beta=0$. We can write the solution as
\begin{equation}\label{form.sep.var.fde.0}
U(x,t)=(T-t)^{(N+2)/4}F(x).
 \end{equation}
Recall that $(T-t)$, instead of $t$, is the correct way of writing the time dependent terms when extinction will happen. Then, $F$ is found to be the explicit pattern
\begin{equation}\label{form.sep.var.fde}
F(x)=\kappa\,\left(\frac{a}{a^2+x^2}\right)^{(N+2)/2}
 \end{equation}
where $\kappa=\kappa(N)$ and $a>0$ is a free parameter.  This function is a case of the Aubin-Talenti maximizing functions for the Gagliardo-Nirenberg inequalities that are prominent in functional analysis, see \cite{Aubin1976, Talenti76}. \nc The special similarity with $\beta=0$ is called separation of variables. Separation of variables plays a big role in problems posed in bounded domains, see \cite{Vazquez2004}, but is not a frequent occurrence for problems in the whole space. We will return to a similar calculation when studying the PLE, and the fractional versions in Section  \ref{sec.cmass.fpl}.

In fact, this special case is important because the equation describes the evolution flow of the famous Yamabe problem.
The problem is explained in detail in \cite{Vazquez2006}, Section 7.5 and Appendix A.III. The fact that the profile is an attractor for all solutions with finite mass if $m=m_s$ was proved in \cite{dPS01}. This result selects this family of solutions as the one that matters in representing the whole class of solutions with finite mass (recall that they never have  constant mass in this scenario).

We point out that  formula \eqref{form.sep.var.fde} represents the stereographic projection of the ball. This geometrical connection  is important in the proof of asymptotic attraction. Geometric flows are important in connection with nonlinear diffusion, we will study below the same issue in dimension $N=2$  (where formally $m_s=0$) which represents Ricci flow in the plane and there we will find a similar separate-variables solution, see Section 8.1 (logarithmic diffusion). We will return again to geometric flows  in Section 13.7.

\noindent {$\bullet$}  We may also examine the limit of infinite mass of these special self-similar solutions to find a kind of very singular solution, hopefully given by a simple formula. We get
$$
\lim_{a\to 0}U(x,t;C,a)=\kappa\,\lim_{a\to 0}\frac{(C-a^2t)^{(N+2)/4}}{(a^2+x^2)^{(N+2)/2}}=\frac{K}{|x|^{N+2}}.
$$
This formula produces a solution everywhere but for $x=0$ where is has a standing singularity. Moreover,  it is stationary in time (we ask the reader to check this fact). Besides, it is the pointwise limit of standard smooth solutions that extinguish in finite time.

\medskip

\noindent {$\bullet$ \bf Pseudo-Barenblatt solutions in Very Fast Diffusion.} Finally, we may also wonder what happened to the Barenblatt solutions that played such an splendid role for $m>m_c$. In fact, they can be extended for $m<m_c$ by algebraic continuation, and we obtain the formulas
\begin{equation}\label{for.seudo.Bar2}
\widehat B(x,t)=(T-t)^{\alpha}\left(C+\frac{\beta (1-m)}{2m}{|x|^2}{(T-t)^{2\beta}}\right)^{-1/(1-m)}
\end{equation}
with $\alpha=1/(m_c-m)$, $\beta=\alpha/N$ (attention: some signs are opposite to the case $m>m_c$), and arbitrary parameter $T$ that represents the extinction time. $\widehat B$ is continued as 0 for $t\ge T$. We have used such solutions in \cite{BBDGV, BonGV-ARMA10, BonDGV-PNAS10} with the name of \sl pseudo-Barenblatt solutions\rm. They do not have the previous property of finite mass because of their typical ``fat tail'' at infinity, $B(x,t)\sim C(m,N,t)|x|^{-2/(1-m)}$\nc.  The corresponding  VSS is
\begin{equation}\label{for.fpme.vss.2}
\widehat V(x,t;m)=\lim_{M\to\infty} \widehat B(x,t;m,M)=
a(m) \left({(T-t)_+}/{|x|^2}\right)^{1/(1-m)}\,
\end{equation}
an example of separate-variables extinguishing solution. This function lives  in a borderline functional space as a function of $x$. Note that for $m_s=(N-2)/(N+2)$ the tail decay is $O( |x|^{2/(1-m)})=O( |x|^{(N+2)/2})$, which is half the decay of the separable-variables solution obtained in \eqref{form.sep.var.fde.0}.

\medskip

\noindent \bc {$\bullet$ {\bf Loss of mass via escape  of particles to infinity.}
In order to obtain a very clear explanation of  the phenomenon of mass loss in the very fast diffusion range $0<m<m_c$
we may make use of Lagrangian trajectories  to illustrate the loss by escape to
infinity. This is done in great detail in the book \cite{Vazquez2006}, Section 5.5.3, and Appendix  AII.
Let us recall the main facts of this approach. We consider the solution $u \ge 0$ as the  density of a substance composed of particles and we write the equations as an instance of the general mass conservation law of continuum mechanics,
\begin{equation}\label{la.cons.law.new}
u_t +\nabla (u\, {\bf v})=0.
\end{equation}
Then the pointwise velocity of the particles is  defined by
\begin{equation}\label{eq.veloc0}
{\bf v}=- u^{m-2}\nabla u=- \frac{1}{m-1}\nabla u^{m-1},
\end{equation}
as done in \eqref{pme.vel}.  The definition makes sense mathematically for all $m\ne 1$. Fore $m=1$ we have ${\bf v}= - \nabla \log( u)$.

We want apply the particle approach  to solutions having a specific asymptotic behaviour as $|x|\to \infty$.
We have to  integrate the kinematic equation for the trajectories $dx/dt = {\bf v}(x,t)$, where the vector field is known through the solution $u(x,t)$, to obtain the way in which particles that are initially  located far
away diverge to infinity with time. This needs precise information on the
behaviour of $\bf v$.

We refer to Section AII.2 of \cite{Vazquez2006} for the calculation of the behaviour for compactly supported solutions
of the PME, the HE, and also the FDE in the range $m_c < m < 1$. The conclusion is
that particles go to infinity with a power rate $x(t)\sim c t^{\beta}.$ Actually, we explain that
the result is tied to the curious linear dependence between velocity and position for
the particles: ${\bf v} = \beta x/t$, which characterizes the Barenblatt solutions.
However, when we apply the Lagrangian idea to solutions $u$ having finite extinction time of the FDE in the very fast range $0<m<m_c$, we find that the trajectories of blow up in finite time, blow up happens first or the trajectories that lie farther away from the origin, and at the extinction time $T$ all of them have reached infinity\nc.

This idea applies also the further models like the PLE that we study next\nc.

\subsection{The very fast diffusion range of the PLE}
There is a close relationship between the theories of the PME and the PLE that has been attested by numerous authors, see an early statement in \cite{Vtwo}, see also the equivalence transformation obtained by us in \cite{RazSanVaz2008} valid for all radial solutions\nc. Of course, there are also marked differences, specially in what regards regularity \cite{DiBenedetto93}, see also \cite{Junning1995}. In the present situation, parallelism wins. Thus we have a result about extinction in finite time.

\begin{proposition}\label{thm.ext.plp} Let $u(x,t)$ be a nonnegative solution of the PLE in the very fast diffusion range $1<p<p_c$ and let us assume that the initial data $u_0\in L^{q_*}(\ren)$ with $q_*=N(2-p)/p$. Then there is a finite time $T$ such that
$u(x,t)\equiv 0$ for all $t\ge T$.  The minimum such $T$ \nc is called the extinction time.
\end{proposition}

\noindent {\sl Proof. }
 Recall that $p_c=2N/(N+1)$ was calculated  in formula \eqref{form.pc}. We have $1<p_c<2$ for $N\ge 2 $ and $p_c=1$ for $N=1$. There is nothing to prove in 1D. We assume that $N\ge 2$.\nc.
We give the proof because it is simple and not readily available. The calculation is based on a simple estimate of the evolution of the $L^{q_*}(\ren)$ norm that we give in some more detail
\begin{equation*}
\begin{array}{c}
\displaystyle \frac{d}{dt}\int u^{q}\,dx=q\int u^{q-1}u_t\,dx=-q\int \nabla u^q |\nabla u|^{p-2}\nabla u\,dx=\\
\displaystyle -q(q-1)\int u^{q-2}|\nabla u|^p\,dx=-C(p,q)\int |\nabla (u^{(p+q-2)/p})|^p\,dx.
\end{array}
\end{equation*}
Using Sobolev's inequality we get
$$
\displaystyle \frac{d}{dt}\int u^{q}\,dx +C(p,q,N)\left(\int u^{\frac{(p+q-2)pN}{p(N-p)}}\right)^{(N-p)/N}\le 0.
$$
It follows that precisely for the value $q=q^*$ we obtain the nonlinear differential inequality $X'+ CX^{\,(N-p)/N}\le 0$ where $X(t)=\|u(\cdot,t)\|_{q_*}^{q_*}$. Extinction in finite time follows by integration of the differential  inequality. \qed

\medskip

Moreover, conservation of mass is excluded for the all finite-mass solutions as a consequence of Proposition \ref{prop.massto0} that applies without changes. Consequently, all finite mass solutions undergo mass extinction in finite or infinite time.

We may also wonder what happened to the Barenblatt solutions that played such an splendid role for $p>p_c$. \bc The results are analogous to those of the Very Fast FDE.  See \cite{Vazquez2006}. We refrain from giving the whole details to concentrate on other more informative topics\nc.

\medskip

\noindent {$\bullet$ \bf Special finite-mass solutions and separation of variables. } The issue of
the finite-mass solutions in the very fast range $1<p<p_c$ depends of an accurate analysis of the suitable
choices of self-similar solutions, an analysis that has not been fully done, but see \cite{BidVer2006}.
Now, there is a case in which we can try separation of variables again like we did in the PME, and put
$$
U(x,t)=(T-t)^{\alpha}F(x).
$$
Substituting into the PLE we get the two conditions: (i)  $\alpha-1=\alpha (p-1)$ in order to eliminate time, so that $\alpha=1/(2-p)$, and (ii) the Lane-Emden equation
$$
-\Delta_p F= \alpha\, F\,.
$$
We need a nonnegative, bounded, and bounded-mass solution $F$ defined in whole space. The literature is clear about this point, see \cite{BidVer1989, Sciun16}: such a solution exists for the precise value $p_s=2N/(N+2)$, i.e., the one for which the associated Sobolev exponent is $p_s^*=2$. Dimension is $N\ge 3$. There is even uniqueness of one radial solution up to translation an change of parameter\nc.
 In the particular case $p=p_s$ we get $\alpha= (N+2)/4$ and
\begin{equation}\label{form.sep.var.ple}
F(x)=\kappa(N)\,
\left(\frac{a^{N+2}}{a^{2N}+ |x|^{2N/(N-2)}}\right)^{N/2},
 \end{equation}
where and $a>0$ is a free parameter. The method applies in more generality, viz to the Doubly Nonlinear Equation \eqref{DNE}. We will work out that example in more detail in  Section \ref{sec.dnle} \nc  and recover the exponents and functions we have just stated.


\subsection{Positivity versus extinction in fast diffusion}\label{sec.masscon}

While for some exponent ranges it is known that nonnegative data produce global solutions that keep being positive for all time everywhere in $x$, in the range of what is called very fast diffusion (VFD) the semigroup solutions may vanish in finite time.
Here, we want to settle the dichotomy between positivity and extinction for  nonnegative solutions of the PME and PLE in the fast diffusion range.
For any solution $u(x,t)$ and  at any given time $t>0$, the following alternative holds for $u(\cdot,t)$ as a function of $x\in \ren$.

\begin{theorem}\label{thm.mc} Let $u$ be a nonnegative solution of the PME with $0<m<1$ with initial data $u_0\in L^1(\ren)$, $u_0\ge 0$. Then for every $t_1>0$  the space function $u(\cdot, t_1) $ is either strictly positive everywhere, or otherwise identically zero. If it is identically zero for $t_1$ it will keep being so for all $t>t_1$. The same result applies to the solutions of the fast PLE.
\end{theorem}

This property was proved in \cite{VazFPL3-2021} for the fractional $p$-Laplacian evolution equation \eqref{frplap.eq} and can be established here for the non-fractional models along the same lines. We ask the reader to check this statement if necessary.
The theorem  allows us to define the extinction time $T(u_0)$  as the first time $t_1$ where $u(x,t_1)$ is the trivial function, and equivalently, as the first time where the nonnegative solution is no more strictly positive and touches zero.

\subsection{Partial result for relative mass conservation}\label{ssec.prmc}

The law of relative mass conservation must be false for the FDE in the very fast diffusion cases since we have proved that even the relative mass between a positive solution $u$ and the zero solution goes to 0 as time increases. But it is well accepted that the disappearance of mass is due to the very fast speed of propagation of disturbances from zero when $u$ is very small. That leaves room to speculate.

Thus, if we consider the class of solutions that decay to zero as $|x|\to\infty$ with the Barenblatt rate $1/u \sim |x|^{2/(1-m)}$ (which means that $u$ is not integrable), then the conservation of the relative mass has been proved by Bonforte et al. in \cite{BBDGV}, Proposition 1, for all $m\le m_c$.
Again, we recall that these results are still to be proved in the corresponding range of the PLE.

\newpage

\section{Mass conservation for the critical exponents}\label{sec.crit.exp}

After the analysis we have just performed, there remains to analyze the situation for the FDE/PME with critical exponent $m_c=(N-2)/N$ for $N>2$, as well as the situation for the fast diffusion range of the PLE \nc with critical exponent $p_c=2N/(N+1)$ for $N>1$.

We have already seen that the construction of the fundamental solution as a self-similar solution that conserves mass must fail since the similarity exponents diverge. The remaining question we face is whether or not \nc  conservation of mass still holds. The answer is yes for both equations but the proof is not easy.

In the case of the PME such result was proved by Bénilan and Crandall already in 1981 in \cite{BC81} as part of more general results on continuity of solutions of porous medium semigroups with respect to nonlinearities. \bc We  give here a simpler proof as a variant of the conservation proof of Theorem \ref{prop.mc.fde}.

\begin{theorem}\label{thm.mc.fde} Let $m=m_c\in (0,1)$ and $N\ge 3$.  Let $u(x,t)\ge 0$ be the semigroup solution of the FDE with initial data $u_0\in L^1(\ren)$, $u_0\ge 0$. Then, for every $t>0$ the mass is conserved:
\begin{equation}\label{mc.fde.m.crit}
\int_{\ren} u(x,t)\,dx=\int_{\ren} u_0(x)\,dx.
\end{equation}
\end{theorem}

\noindent {\sl Proof.} (i) {\sc Preparation.}  Let $t_2>t_1>0$, and recall that $u(t)$ is a continuous function with values in $L^1(\ren)$. In particular we have $u(t_1)\in L^1(\ren)$. This means  that for every $\ve>0$ there exist $R(\ve)$ such that
$$
\int_{\{|x|\ge R(\ve)\}} u(x,t)\,dx \le \ve.
$$
By the Bénilan-Crandall inequality we have $(1-m)t u_t\le u$ so that for $t_2\ge t\ge t_1$ we have
$$
u(x,t)\le u(x,t_1) (t/t_1)^{1-m}\le u(x,t_1) (t_2/t_1)^{1-m}.
$$
This means that for  $t_2>t>t_1>0$ we have
$$
\int_{\{|x|\ge R(\ve)\}} u(x,t)\,dx \le C_1\ve, \quad C_1=(t_2/t_1)^{1-m}.
$$

\noindent (ii) {\sc Mass evolution.} We begin the argument as in Theorem \ref{prop.mc.fde} with similar notation. We consider a  weak solution $u\ge 0$ of the FDE  defined in  $\re^N$ for a period of time $[0,T]$. We start from the equation, multiply it by a smooth test function \  $0\le \vp(x)\le 1$   with compact support in  $\re^N$; indeed, we impose $0< \vp(x)$ for $|x|<2$,  $\vp(x)=1$ for $|x|\le 1$, and $\vp(x)=0$ for $|x|\ge 2$. We put $\vp_n(x)=\vp(x/n)$ for $n$ an  integer. As in  Theorem \ref{prop.mc.fde} we  integrate to get
$$
\frac{dX_n}{dt}= \frac{d}{dt}\int_{\ren} u\,\vp_n\,dx=\int_{\ren}  u^{m}\Delta\vp_n\,dx
$$
where $X_n(t)=\int_{\ren}  u\,\vp_n\,dx$  is the weighted mass calculated at a certain time at $t$. Using H\"older's inequality we get
$$
 \left|\frac{dX_n}{dt}\right| \le \left(\int_{\{|x|\ge n\}}  u\, dx\right)^{m}\, Y_n^{1-m},
$$
where $Y_n$ can be written as in Theorem \ref{prop.mc.fde}, \
$ Y_n=\int_{\ren}  \left|\vp^{-m}\Delta\vp\right|^{1/(1-m)}\,dx.$
This has been estimated as
$$
Y_n= \int_{n}^{2n} \left|\vp_n(r)^{-m}\Delta\vp_n(r)\right|^{1/(1-m)}\,r^{N-1}dr= C(\vp)\, n^{N-\frac2{1-m}}.
$$
Now, $N(1-m)-2=0$ precisely when $m=m_c$. We get $Y_n=Y_1$, a constant independent of $n$.

\noindent (iii) We gather both calculations and find that for $n\ge R(\ve)$ we get
$$
|\frac{d}{dt}\int_{\ren} u\,\vp_n\,dx|\le C_1\ve Y_1.
$$
Pass now to the limit $n\to \infty$, $\ve\to 0$, to get the conservation of mass form $t=t_1$ to $t=t_2$. Now we recall that these positive times are arbitrary.  \qed

\begin{corollary}\label{mc.fcm.sing} The  conservation result of Theorem \ref{thm.mc.fde} holds for signed solutions with $u_0\in L^1(\ren)$. The term mass should be replaced by first integral.
\end{corollary}

\noindent {\sl Proof.} In that case we consider the signed solution $u$ with  initial data $u_0$ and the nonnegative solutions $u_1$ and $u_2$ with initial data $u_1(x,0)=\max\{-u_0(x),0\}$ and $u_2(x,0)=\max\{u_0(x),0\}$ resp. By comparison we have the bounds
$$
-u_1(x,t)\le u(x,t)\le u_2(x,t).
$$
We can apply the previous analysis to $u_1$ and $u_2$ to conclude as in point (i)  above that the absolute integral $\int_{|x|\ge R} |u(x,t)|\,dx$  is small enough for every time interval $[t_1\le t\le t_2]$ with $0<t_1<t_2$ if integration extends outside of a large ball $B_R(0)$. The rest of the proof works.\qed
\nc
\

 We will discuss the fractional version of this equation in Section \ref{sec.cmass.fpme} and we will give there another  proof that shows how to deal with such borderline problems.

 \medskip

For the PLE the situation is a bit more complicated because we cannot integrate by parts twice to pass all the derivatives to the test function in the weak formulation. This makes the problem mathematically interesting. We will present a complete proof next. Note that the case $p>p_c$ has been settled in \cite{FinoFVespri2014} (for a larger class of equations), while for $1<p<p_c$ mass conservation does not hold, as proved in the previous section. Before addressing the proof we recall that $p_c=1$ for $N=1$ so we  assume that $N\ge 2$. The case $p_c=1$ turns out to be different and will be treated in Section \ref{sec.n1pc1}.

\begin{theorem}\label{thm.mc.cple} Let $p=p_c>1$ and $N\ge 2$.  Let $u(x,t)\ge 0$ be the semigroup solution of the fast PLE with initial data $u_0\in L^1(\ren)$, $u_0\ge 0$. Then, for every $t>0$ the mass is conserved:
\begin{equation}\label{mc.pc}
\int_{\ren} u(x,t)\,dx=\int_{\ren} u_0(x)\,dx.
\end{equation}
\end{theorem}

\noindent {\sl Proof.}  (i) {\sc Reduction step.} We first make a reduction in the class of data. We may always assume that $u_0\in L^1(\ren)\cap L^\infty(\ren)$ and that $u_0$ is compactly supported, say in the ball of radius $R_0$. If our form of mass conservation is proved under such assumptions, then it follows for all data $u_0\in L^1(\ren)$ as a consequence of the semigroup $L^1$ contraction property, which is a well-known property of this equation (as already discussed in Subsection \ref{ssec.pme})\nc. We recall that the mass $M(t)=\int_{\ren} u(x,t)\,dx$ is nonincreasing in time, so we are only worried about possible loss of mass during the evolution.

\medskip

\noindent(ii) {\sc Preparation step.} We now collect some needed delicate preliminary facts:

\noindent (a) If the initial support is contained in a ball, the solution is monotone in space outside a larger ball (along outward radial directions, as follows from  the Aleksandrov reflection principle that we explain in Appendix \ref{ssec.aleks}). We can check that it works for this equation and setting\nc.

\noindent (b) An important tool will be the condition of uniform small mass at infinity, which is derived from  the property of almost monotonicity in time: since for any $t>0$ $u(\cdot,t)$ is integrable there is a radius $R=R(\ve,u,t)>0$ such that
$$
\int_{|x|\ge R} u(x,t)\dx<\ve.
$$
We need the  bound which to be uniform in intervals of time. For that we may use the well known growth inequality $(2-p)tu_t\le u$, proved in \cite{BC81b} and applied to this equation in \cite{VazFPL2-2020}. Indeed, by integrating in time from $t_1>0$, we prove
 for any two times $0<t_1<t_2<T$  and every $\ve>0$ the mass of $u(t)$ outside a large ball of radius $R=R_\ve$  is less than $\ve$ for all $t\in (t_1,t_2)$, and  $R_\ve$ depends  on $\ve, t_1$ and $t_2$ (and $u$) \nc but nothing more.

\noindent (c)
Let us take a point $x$ such that $|x|=r\ge 4R_0$ and let us consider the annulus $A=\{r/4\le |x|\le r/2 \}$. Using the monotonicity mentioned  in (a), more precisely Lemma \ref{le.AleksAC}, we have for $y\in A$ and all $t>0$
$$
u(y,t)\ge u(x,t).
$$
After integrating w.r.t. $y$ in $A$ this gives
$$
|A|\, u(x)\le \int_{|y|\in  A} u(y,t)\,dy\le M(r/4),
$$
where $M(r)=\int_{|x|\ge r} u(y,t)\,dy$. If follows that
$$
u(x,t)\le c\,M(r/4)\,r^{-N},
$$
\nc Using this and the uniform integrability near infinity in the interval $t_1\le t \le t_2$, we conclude that
\begin{equation}\label{udecay}
u(x,t)|x|^N \le C\ve,
\end{equation}
where  $C=C(t_1,t_2)>0$ is universal and $\ve$ is the uniform small bound for the mass near infinity, i.e., it holds for $|x|\ge R_\ve$\nc.

\medskip

\noindent (iii) {\sc Main step.} The technical part of the proof starts in a rather standard way (as the similar proofs done recently in \cite{Vazquez2020, VazFPL2-2020} for fractional diffusion) and needs some delicate  technical treatment that we explain next.
We do a calculation for  the weighted mass. Taking a smooth and compactly supported test function $\varphi(x)\ge0 $, we have for $t_2>t_1>0$:
\begin{equation}\label{mass.calc1.plp}
\displaystyle \int u(t_1)\varphi\,dx-\int u(t_2)\varphi\,dx= \iint
|\nabla u|^{p-2}\nabla u \cdot \nabla \varphi \,dxdt.
\end{equation}
Space integrals are over $\ren$  and time integrals over $[t_1, t_2]$. In the proof we will use the sequence of test functions $\varphi_n(x)=\varphi(x/n)$ where $\varphi(x)$ is a cutoff function which equals 1 for $|x|\le 1$ and zero for $|x|\ge 2$. We take $n\ge 2R_\ve$, see point (b).  We will estimate the integral in \eqref{mass.calc1.plp}  in absolute value (by taking absolute value of the integrand). In the next paragraphs we will sometimes forget the time integrals momentarily for ease of writing.\nc

 We want to estimate the term in the right-hand side of \eqref{mass.calc1.plp} in absolute value, let us call such an integral $I_n$. We need only integration over $A_n=\{n\le |x|\le 2n\}$. By H\"older inequality we have
\begin{equation}\label{eq.inter}
\displaystyle I_n\le \left(\iint_{A_n} |\nabla u|^{p}\varphi_n\,dxdt \right)^{(p-1)/p}
\left(\iint_{A_n} |\nabla \varphi_n|^{p}\varphi_n^{-(p-1)}\,dxdt \right)^{1/p}.
\end{equation}
Now, when $p=p_c$ we have $p-1=(N-1)/(N+1)$, $(p-1)/p=(N-1)/2N$. Since the last integral is
$$
\iint |\nabla \varphi_n|^{p}\varphi_n^{-(p-1)}\,dxdt =c(p)\iint |\nabla (\varphi_n^{1/p})|^{p}\,dxdt,
$$
 we will ask $\zeta= \varphi^{1/p_c}$ to be smooth to make the integral convergent\nc. In that case the integral is bounded for each $n $ and is of the order of  $C\,n^{N-p}$ (times $t_2-t_1$). So we still have to prove that the local energy integral, $ \int \int_{A_n} |\nabla u|^{p}\varphi_n\,dxdt $, is small enough as $n\to\infty$. More precisely, if $p=p_c$ we need
\begin{equation}\label{est.energ}
\int \int_{A_n} |\nabla u|^{p}\varphi_n\,dxdt =o(n^{-N}).
\end{equation}
 Indeed, if this holds, going to \eqref{eq.inter} we find $I_n=o(n^{-q})$ with
$$
q= \frac{N(p-1)}{p}-\frac{N-p}{p},
$$
which vanishes for $p=p_c$, so that $I_n\to 0$\nc.

\smallskip

\noindent (iv) {\sc The local energy. } We derive the needed local estimate under the same assumptions on $u_0$.  Let us prove a first estimate.

\begin{lemma}\label{lem.energy1} Let $p=p_c$.  Let $u(x,t)\ge 0$ be the semigroup solution of the Fast PME with initial data $u_0\in L^1(\ren)\cap L^\infty(\ren)$, $u_0\ge 0$. Then for every $t_1>0$  we have the
local energy estimate
$$
\displaystyle \int_0^T \int_{3n\ge |x|\ge 2n} |\nabla (u(x,t)|^{p}\,dxdt\le   \int_{4n \ge |x|\ge n} u_0^2(x)\,dx+ C_1t_2 \,U_n^{p}n^{N-p}.
$$
\end{lemma}

\noindent {\sl Proof.} \nc We take a smooth cutoff function $\psi$ much as before, but now supported in an annulus $A=\{1\le |x|\le 4\}$ with $\psi(x)=1$ for $2\le |x|\le 3$,  and we get for $0\le t<t_1$ the identity
\begin{equation}\label{mass.calcnew}
\begin{array}{c}
 \displaystyle \frac12 \int u^2(x,0)\psi^2(x)\,dx-\frac12 \int u^2(x,t_1)\psi^2(x)\,dx=\\
 \displaystyle  \iint
 |\nabla u(x,t)|^{p-2} (\nabla u(x,t)\cdot \nc \nabla_x ( u(x,t)\psi^2(x))\,dxdt.
\end{array}
\end{equation}
We write the second term as $I+II$ where
\begin{equation*}\label{mass.calcnew2}
 \displaystyle I:= \iint
|\nabla u|^{p}\psi^2(x)\,dxdt.
\end{equation*}
This term is positive and will be the local energy. We also have
\begin{equation*}\label{mass.calcnew3}
\displaystyle  II:= 2 \iint u(x,t)\psi(x) \,|\nabla u(x,t)|^{p-2}\nabla u(x,t)\cdot\nabla \psi(x)\,dxdt.
\end{equation*}
Term $II$ does not have a role and has to be absorbed. We choose $\psi_n(x)=\psi(x/n)$ and the set $D=D_n$ where $\psi_n$ is not zero, a typical annulus at distance $n$ from the origin. Now if $u$ is bounded above by some $U(D_n)$ in $D_n$ we have (after dropping dependence on time for ease of notation)
\begin{equation*}\label{mass.calcnew4}
\displaystyle  II_n\le 2 \,U(D_n) \left(\iint_{D_n} |\nabla u(x)|^{p}\psi^2(x)\,dxdt \right)^{(p-1)/p} \, Y_n^{1/p}\,\nc.
\end{equation*}
By \eqref{udecay} we know that $U(D_n)\to 0$ with a rate  $U(D_n)\le U_n=C\ve n^{-N}$\nc. When $p=p_c$, the critical exponent, we get precisely
\begin{equation*}
 Y_n = \nc\iint_{D_n} |\nabla \psi_n|^{p} \psi_n^{2-p}\,dx\le C\,n^{N-p}(t_2-t_1),
\end{equation*}
where $C$ is independent of $n$. Then we get
$$
II_n\le C t_2^{1/p}U_n\, I_n^{(p-1)/p}n^{(N-p)/p}\,.
$$
Using Young's inequality:
$$
ab\le \frac1{p_1}a^{p_1}+ \frac1{q_1}b^{q_1}, \quad \mbox{with \ } \ a,b>0, \ p_1>1, \ q_1=p_1/(p_1-1),
$$
with $p_1=p$ we get
\begin{equation}
II_n\le \frac12 I_n+ C_1U_n^{p}t_2n^{N-p}.
\end{equation}
Going back to \eqref{mass.calcnew} we get
$$
\frac12 I_n\le \frac12 \int u^2(x,0)\psi_n^2(x)\,dx +C_1t_2 \,U_n^{p}n^{N-p}.
$$
This implies the estimate. \qed

\medskip

\noindent {\sc End of Proof.} We now go back to equation \eqref{eq.inter} and use formula \eqref{udecay} and Lemma \ref{lem.energy1} with initial time $t_1$  to conclude that
$$
\iint_{A_n} |\nabla u|^{p}\varphi_n\,dx\le U_n\int_{|x|\ge n}u(x,t_1)\,dx+ C_1t_2 \,U_n^{p}n^{N-p}\le
 C\ve \,n^{-N} + C\ve^p n^{-N(p-1)-p}\nc.
$$
Since $N(p-1)+p=N$ \nc precisely for $p=p_c$, this is the estimate we needed. \qed

\medskip

Though we are less interested in signed solutions we can extend the result to them.

\begin{corollary} The  conservation result of Theorem \ref{thm.mc.cple} holds for signed solutions with $u_0\in L^1(\ren)$. The term mass should be replaced by first integral.
\end{corollary}

\noindent {\sl Proof.} In that case we consider the signed solution $u$ with compactly supported initial data $u_0$ and the nonnegative solutions $u_1$ and $u_2$ with initial data $u_1(x,0)=\max\{-u_0(x),0\}$ and $u_2(x,0)=\max\{u_0(x),0\}$ resp. By comparison we have the bounds
$$
-u_1(x,t)\le u(x,t)\le u_2(x,t).
$$
We can apply the previous analysis to $u_1$ and $u_2$ to conclude as in points (a), (b), (c) above that the absolute integral $\int_{|x|\ge R} |u(x,t)|\,dx$  is small for every time interval $[t_1\le t\le t_2]$ with $0<t_1<t_2$ if integration extends outside of a large ball $B_R(0)$. Hence, the estimates of point (ii) of  the previous proof apply to signed solutions. The proof of points  (iii) and (iv) stays if we  take $u$ bounded  by $U_n$ in absolute value. This concludes the conservation formula \eqref{mc.pc} for signed $u_0$ with compact support. For general signed $u_0$ use the approximation argument of the Reduction Step.\qed
\nc

\newpage

\section{What happened to the fundamental solutions}\label{sec.fscexp}

We have said in Section \ref{sec.cm.conn}, right after formula \eqref{sim.exp.FDE}, \nc that the similarity exponents $\alpha$ and $\beta$ of the  fundamental solution blow up as $m$ goes down to $m_c=(N-2)/N$ in the PME, or as $p$ goes down to $p_c=2N/(N+1)$ in the PLE. So the question is: where has the standard form of the fundamental solution gone in that limit? Can we make sense of the limit in a qualitative and quantitative way?

In order to make sense of the results that follow, we recall that the PME/FDE can be written
as \ $\partial_t u=\mbox{div}(D(u)\nabla u)$ \ with variable coefficient $D(u)=mu^{m-1}$, and this coefficient goes to zero as $u\to\infty$ when $m<1$. Since such large densities are needed, at least for small times, in the case of the fundamental solution, there is a moment where diffusion breaks down and the equation is unable to spread around the initial point mass. This moment is reached in the borderline cases.

\subsection{The PME/FDE case } A careful inspection of the Barenblatt formula gives the answer to the disappearance of the fundamental solution as a solution of the equation. The following has been known for a number of years.\nc

\begin{theorem}\label{thm.dirac.pme} Let $N>2$ and let us denote by $B_m(x,t;M)$ the fundamental solution of the PME/FDE  with exponent $m>m_c$ and mass $M>0$. Then as $m\to m_c$ we have the singular limit
\begin{equation}\label{limit.bs.mc}
\lim_{m\to m_c} B_m(x,t;M)=M\delta(x) 
\end{equation}
in the sense of measures. In other words, the initial point mass remains ``completely frozen'' at its location for all times with the same mass. Actually, outside of any ball $B_R(0)$ around  $x=0$ the convergence to zero is  uniform thanks to the estimate
$$
B_m(x,t)\le C(N)\,\ve^{1/(1-m)}\,|x|^{-2/(1-m)}t^{1/(1-m)}, \quad x\ne0,
$$
valid if $\ve=m-m_c>0$.
\end{theorem}

\nc \noindent {\bf Remarks.} 1) The answer makes perfect sense if the limit is understood as a \sl measured-valued solution, \rm a well-known concept used in the theory of conservation laws  \cite{DiPerna, MeasureVS1996} and optimal transport. It has gained reputation in the study of other evolution equations thanks to the techniques of optimal mass transport  \cite{Villani1, Villani2, FigGlaudo, Santamb15} and the use of Wasserstein metrics. A complete theory does not exist for such type of solutions in the nonlinear diffusion setting, see some results in \cite{PorzioST13, PorzioST14, DGC23, VazquezMs2009}.

\noindent 2) The non-existence of a fundamental solution of the FDE in the critical case and below, due to the lack of diffusion, was first proved by Brezis and Friedman in \cite{BrezFried83}. They considered the Cauchy problem for the PME with initial Dirac delta and proved that then there exists no strong nonnegative solution if $0< m\le m_c$, $N\ge 3$. The approximations show that the Dirac mass does not diffuse into a nontrivial solution in the surrounding space. This is what we will see in the limit process.

\noindent 3) Quantitative information will be given in the proof. Indeed, for $\ve=m-m_c>0$ small we have the uniform estimate
$$
B_m(x,1)\le C\,\ve^{1/(1-m)}\,|x|^{2/(1-m)},
$$
on the vanishing of the approximations for $x\ne 0$. Note that the rate of decay $O(|x|^{-2/(1-m)})$ becomes in the limit $ O(|x|^{-N})$ which is not integrable, hence the need for a finer analysis\nc.

\medskip

\noindent {\sl Proof.} We will do the calculations for $t=1$; since the solution is self-similar, this is enough. Without loss of generality we may also take $M=1$ because of the scaling laws.

\noindent (i) First, we calculate precise asymptotic  estimates on  formula $B_m(x,t;M=1)$   when $m\to m_c$ at time $t=1$. Thus, we want  to  examine the limit $m\to m_c$ in the profile
$$
F_m(y)=\left(C_m+ k_m |y|^2\right)^{-\frac1{1-m}}.
$$
Now,  \nc  if $\ve=m-m_c>0$, we have
$$
\quad \alpha(m)=\frac{N}{2-N(1-m)}=\frac1{\ve}, \quad \beta(m)=\frac1{2-N(1-m)}=\frac1{N\ve}.
$$
We have
$$
k_m=\frac{(1-m)\alpha}{2mN}\sim  \frac1{k_1\ve}, \mbox{ with } \ k_1=N(N-2),
$$
This gives the first pointwise estimate
$$
F_m(x) \le C(N)\,\ve^{1/(1-m)}\,|x|^{-2/(1-m)}.
$$
To get an estimate for $t>0$ we write $B_m(x,t)=t^{-\alpha}F_m(|x|t^{-\beta})$, so that the previous estimate gives
$$
B_m(x,t)\le C(N)\,\ve^{1/(1-m)}\,|x|^{-2/(1-m)}\, t^{\mu},
$$
with $\mu>0$,  and more precisely
$$
\mu=\beta(\frac2{1-m}-N)=\frac{\beta \ve N}{1-m}= \frac{1}{1-m}.
$$

\noindent (ii) We now  examine the limit of $C_m$. It is determined by the mass condition
$$
\omega_N\int_0^\infty \left(C_m+ k_m r^2\right)^{-\frac1{1-m}}r^{N-1}dr=1.
$$
To work this out, we do as follows; we first put $r^2= s^2 C_m k_m^{-1}$ which leads to
$$
\omega_N \, C_m^{-\frac1{1-m}+\frac{N}2}k_m^{-N/2} \int_0^\infty \left(1+  s^2\right)^{-\frac1{1-m}}s^{N-1}ds=1.
$$
Note that
$$
\frac{N}2-\frac1{1-m}=\frac{N(1-m)-2}{2(1-m)}=\frac{-N\ve}{2(1-m)}\sim -\frac14 N^2\ve.
$$
We still  have to control the integral
$$
I_m= \int_0^\infty \left(1+  s^2 \right)^{-\frac1{1-m}} s^{N-1} ds,
$$
since it is convergent for all $m>m_c$ but divergent at $m_c$. In order to do that
we recall some formulas for Euler's Beta function:
$$
{\mathcal B}(p,q)=\int_0^1 t^{p-1}(1-t)^{q-1}\,dt=\int_0^\infty \frac{t^{p-1}\,dt}{(1+t)^{p+q}}=\frac{\Gamma(p)\Gamma(q)}
{\Gamma(p+q)}.
$$
Hence, if $s^2=t$ we get
$$
I_m= \int_0^\infty \left(1+  s^2 \right)^{-\frac1{1-m}} s^{N-1} ds=
\frac12 \int_0^\infty \left(1+  t\right)^{-\frac1{1-m}} t^{(N-2)/2}dt=\frac12 {\mathcal B}(N/2,q_m), \quad
$$
with
$$
q_m= \frac1{1-m}-\frac{N}2= \frac{N\ve}{2(1-m)}\sim \frac{N^2}4 \ve,
$$
Using this estimate for $q_m$ we get
$$
I_m=\frac12 {\mathcal B}(N/2,q_m)=\frac{\Gamma(N/2)\Gamma(q_m)}{2\Gamma(N/2+q_m)}\sim \frac{\Gamma(q_m)}{2}\sim \frac{2}{N^2\ve}.
$$
Summing, up we can estimate $C_m$ as
$$
C_m^{N^2\ve/4}\sim c_N \,(k_1\ve)^{N/2}{\mathcal B}(N/2,q_m)\sim c(N)\ve^{(N-2)/2}.
$$
It follows that $C_m\to 0$ very fast as $m\to m_c$ with the estimate
\begin{equation}\label{C_mbound1}
\log C_m \sim c_1 \frac{\log \ve}{\ve}, \quad c_1=\frac{2(N-2)}{N^2}.
\end{equation}

\noindent (iii) We now prove the limit of $F_m$. After the previous calculations we easily conclude that
$$
\lim_{m\to m_c} F_m(y)=0 \quad \mbox{ if} \ y\ne 0, \quad \lim_{m\to m_c} F_m(0)=\infty.
$$
Recall that for all $m\in (m_c,1)$ all $F_m$ are smooth, radially symmetric and monotone functions of $|y|$ and all have mass 1\nc.

 In order to conclude as stated in \eqref{limit.bs.mc} that the family $F_m$ converges to the unit delta function located at $y=0$, we still have to eliminate  the possible escape of mass to infinity. In other words, we must show that $F_M$ is a standard mollifying family when $\ve=m-m_c\to 0$. Here is the argument: we fix a radius $R$ and calculate the mass outside of the ball $B_R$, in $\{|y|\le R\}$. Putting again $r^2= s^2 C_m k_m^{-1}$  we get
$$
I_{R,m}=\omega_N\int_R^\infty \left(C_m+ k_m r^2\right)^{-\frac1{1-m}}r^{N-1}dr=
$$
$$
C_m^{-\frac1{1-m}+\frac{N}2}k_m^{-N/2} \int_{S_m(R)}^\infty \left(1+  s^2\right)^{-\frac1{1-m}}s^{N-1}ds
$$
where $S_m(R)=R (k_m/C_m)^{1/2}$ which is of the order  $S_m(R)\ge (1/\ve)^{c/\ve}$ which is super-exponentially far in the tail of the integral, so this integral is negligible with respect to the full integral:
$$
\frac{I_{R,m}}{I_{0,m}}=\frac{\int_{S_m(R)}^\infty \left(1+  s^2\right)^{-\frac1{1-m}}s^{N-1}ds}
{\int_{0}^\infty \left(1+  s^2\right)^{-\frac1{1-m}}s^{N-1}ds},
$$
 The approximate calculation of $I_{R,m}$ is easy for $\ve\sim 0$ so that $S_m(R)$ is very large, and it gives the desired fast relative decay, $I_{R,m}/I_{0,m}\sim S_m(R)^{-\ve N^2/2} \sim c_1 \ve^{(N-2)/2}$. \qed

\medskip

For further information about the limit equation, we direct the reader to
 Section 5.6 of \cite{Vazquez2006}. Moreover, when $N=2$ the limit exponent is $m_c=0$. In that case the limit is carefully worked out in \cite{Vazquez2006}. Curiously, mass concentration still happens but conservation of mass is broken\nc! Actually, there is extinction in finite time, see Section \ref{sec.ricci}.

\medskip

\noindent {\bf Blowup of the constant of the smoothing effect.} It is well-known that solutions of the PME with $m>m_c$ and integrable initial data say of unit mass, are uniformly bounded for all fixed $t>0$. Moreover, there is a \sl best constant \rm of this functional inequality
$$
K(m,N)=\sup\left\{u(x,1): \ x\in \ren,  u_0\ge 0, \int u_0\,dx=1\right\},
$$
where it is understood that $u(x,t)$ is the solution of the PME with initial data $u_0\ge 0$.
Schwarz symmetrization allows to prove that this best constant is achieved by the Barenblatt fundamental solution, \cite{V82sym, Vazsym05}.  The calculation in (i) shows how the uniform boundedness of the solutions disappears in the limit in the form
$$
K(m,N)=\sup_{x\in \ren} B_m(x,1)=\sup_{y\in \ren} F_m(y)= C_m^{-\frac1{1-m}}\sim C_m^{-\frac{N}2}\
$$
which grows super-exponentially as $m\to m_c$
\begin{equation}\label{limit.smeff}
\log K(m,N)\sim \frac{N-2}{N} \, \cdot \frac{\log(1/(m-m_c))}{m-m_c}.
\end{equation}

It is also interesting to give information about the size of the region around $0$ where the solution $F_m$ has concentrated. This can be done by just estimating the distance $d_m$ at which $F_m$ takes a half of  the maximum value $F_m(0)$. An easy calculation gives
$$
d_m^2\sim  c C_mk_m,
$$
which gives in rough approximation $\log d_m\sim \frac12 \log C_m$, so that a rough version of conservation of the unit integral is obtained, $\log(K_md^N_m)\sim 0$ in some weak sense (with a small error compared to both factors).

\subsection{Large-time behaviour of general finite-mass solutions}\label{sec.expdecay}

The lack of a fundamental solution for the FDE with $m=m_c$ raises the question of whether there exists some typical behaviour for suitable initial data or at least some generic bounds. Something is known. The question is touched upon in Section 5.6 of \cite{Vazquez2006} where the following interesting result about exponential decay is given.

\begin{theorem} Let $m=m_c$, $N\ge 3$. Solutions of the PME with data in $L^1(\ren)\cap
L^\infty(\ren)$ decay in time according to the rate
\begin{equation}\label{cep.expo}
\log (\|u(t)\|_\infty ) \sim - C(N)M^{-2/(N-2)} t^{N/(N-2)}\nc
\end{equation}
 as $t\to \infty$, with $M=\|u_0\|_1$. This rate is sharp.
\end{theorem}

This very fast and delicate decay rate was rigorously  established by Galaktionov, Peletier and the author in \cite{GPV00}, based on previous formal analysis by King \cite{Ki93}.  Sharp rate means that it cannot be improved as an upper bound\nc.

\smallskip

\begin{open} A similar analysis for the PLE in the critical exponent $p_c$ is not known to the author.
\end{open}


\subsection{The limit of the fundamental solution for the fast PLE}

 In this subsection  ${B}_p(x,t;M)$ will denote the fundamental solution of the PLE with exponent $p>p_c$ and mass $M>0$. \nc  We leave it to the reader to repeat the above procedure and get the precise asymptotic results as follows

\begin{theorem}\label{delta.pc} As  $p\to p_c$ we have   for all $t>0$ the singular limit $B_p(x,t;M)$
\begin{equation}\label{limit.bs.pc}
\lim_{p\to p_c} { B}_p(x,t:M)=M\delta(x)
\end{equation}
in the sense of measures in $\ren$\nc. In other words, the initial point mass remains frozen at its location for all times with the same mass. Moreover, for $\ve=p-p_c>0$ small we have the uniform estimate
$$
B_p(x,1)\le C\left(\ve\,|x|^{-p}\right)^{1/(2-p)}\,.
$$
that implies uniform convergence to zero for $x$ away from the origin as $p\to p_c$\nc.
\end{theorem}

\noindent {\sl Proof.} Here are some details. If $\ve=p-p_c>0$, we have
$$
\quad \alpha(p)=\frac{N}{p-N(2-p)}=\frac{p_c}{2\ve}, \quad \beta(p)=\frac1{p-N(2-p)}=\frac1{(N+1)\ve}.
$$
Also,
$$
F_p(y)=\left(C_p+ k_p |y|^{\frac{p}{p-1}}\right)^{-\frac{p-1}{2-p}},
$$
with
$$
k_p=\frac{2-p}{p} \beta^{\frac{1}{p-1}}\sim  k_1(N,p)\,\ve^{-\frac{1}{p-1}}.
$$
 We only need to fix the constant $C_p$ to have unit mass. Finally, $C_p$ is determined by the mass condition
$$
\omega_N\int_0^\infty \left(C_p+ k_p r^{p/(p-1)}\right)^{-\frac{p-1}{2-p}}r^{N-1}dr=1,
$$
To work this out, we first put $r= s (C_p/ k_p)^{(p-1)/p}$ which leads to
$$
\omega_N \, C_p^{-\frac{p-1}{2-p}+\frac{N(p-1)}{p}}k_p^{-\frac{N(p-1)}{p}} \int_0^\infty \left(1+  s^{p/(p-1)}\right)^{-\frac{p-1}{2-p}}s^{N-1}ds=1.
$$
Note that
$$
(p-1)(\frac{N}{p}-\frac1{2-p})=(p-1)\frac{N(2-p)-p}{p(2-p)}=\frac{-(p-1)(N+1)\ve}{p(2-p)}\sim -c(N)\ve\,, \quad
$$
with $c(N)=(N+1)^2(N-1)/2N$. We still  have to control the integral
$$
I_p= \int_0^\infty \left(1+  s^{p/(p-1)}\right)^{-\frac{p-1}{2-p}}s^{N-1}ds,
$$
since it is convergent for all $p>p_c$ but divergent at $p_c$. If  we put $s^{p/(p-1)}=t$ and recall some formulas for Euler's Beta function ${\mathcal B}(p,q)$: we get
$$
I_p=  c \,\int_0^\infty \left(1+  t\right)^{-\frac{p-1}{2-p}} t^{\frac{N(p-1)-p}{p}}dt=c\, {\mathcal B}(\frac{N(p-1)}{p},q_p), \quad
$$
with
$$
q_p= \frac{p-1}{2-p}-\frac{N(p-1)}{p} \sim c(N) \ve,
$$
Therefore,
$$
I_p=c \, {\mathcal B}(\frac{N(p-1)}{p},c(N) \ve)\sim\frac{\Gamma((N-1)/2)\Gamma(c(N) \ve)}{\Gamma(c(N) \ve+ (N-1)/2)}\sim \frac{c_1}{\ve}.
$$
Summing, up we can estimate $C_p$ as
$$
C_p^{\,c_1\ve}\sim c_N \,(\ve)^{N/p}{\mathcal B}(N/2,q_m)\sim c(N)\ve^{(N-p)/p}.
$$
It follows that $C_p\to 0$ very fast as $p\to p_c$ with the estimate
\begin{equation}\label{C_mbound2}
\log C_p \sim c_2 \frac{\log \ve}{\ve}, \quad c_2=c_2(N).
\end{equation}
This ends the proof with a result quite similar to the PME one. \qed

\medskip

The smoothing effect constant blows up as $p\to p_c$ in a similar way as seen above for the PME.
We will explain in Section \ref{sec.n1pc1} what happens for $N=1$, $p_c=1$.


\section{Non-conservation in some borderline cases}\label{sec.ncbl}

There are only two extreme cases of equations with critical exponents where mass conservation does not hold, namely, the fast PME with $m=0$ and the fast PLE with $p=1$. The candidates to be Barenblatt solutions, i.e., those obtained as limits as the parameters converge $m\downarrow 0$ ( resp. $p\downarrow 1$) still remain concentrated at $t=0$, but a new curious phenomenon appears: a very precise mass loss is obtained at the same time in such cases until the ``almost frozen solutions'' disappear. The rate of mass decay is constant, a characteristic of this type of diffusion. Therefore, we cannot speak of conservation of mass for the solutions of these equations, but instead of a {\sl universal law of mass loss} which involves all finite solutions in a certain class and applies at all times with the same rate.

\subsection{Logarithmic diffusion in $N=2$}\label{sec.ricci}
In order to pass to the limit $m\to 0$ in the PME it is important to take a look the equation
$$
\partial _t u= \Delta(u^m)= \mbox{div\,}\big(mu^{m-1}\nabla u\big)
$$
and notice the extra factor $m$ that would make the equation too degenerate in the limit. It is easy to avoid that factor by rescaling time, $t'=mt$, so that we write
$$
\partial_{t'} u= \mbox{div\,}\big(u^{m-1}\nabla u\big)
$$
The solutions of such form of the PME are then proved to converge  when $m\to 0$  to solutions of the limit equation
\begin{equation}\label{logdiff.eqn}
\partial_{t'} u= \mbox{div\,}\big(u^{-1}\nabla u\big)=\Delta (\log \, u)\,,
\end{equation}
well known as the \sl logarithmic diffusion equation\rm, of course for suitable initial data. The equation and its relation with the PME are studied in detail in Chapter 8 of \cite{Vazquez2006}, based on previous work by Daskalopoulos and Del Pino \cite{DP95}, DiBenedetto and Diller \cite{DBD96} and Esteban, Rodriguez and V\'azquez \cite{VER96, REV97}, where the existence of solutions with mass loss was detected and studied. In the sequel we use the rescaled version of the PME with notation $t$ instead of $t'$. The connections with 2D Ricci flow have made this equation a popular item after \cite{Ham88}. See also \cite{ChenLi1991}.

Concerning the limit of the fundamental solution of the PME as $m\to 0$ the following {\sl concentration plus extinction} result  is proved in page 91 of  \cite{Vazquez2006}.

\begin{theorem}\label{delta.pme0} Let $N=2$ and let $B_m(x,t;M)$ the fundamental solution of the PME with exponent $m>0$ and mass $M>0$. Then as $m\to 0$ we have the singular limit
\begin{equation}\label{limit.bs.mc2}
\lim_{m\to 0} B_m(x,t;M)=(M-4\pi t)_+\delta(x).
\end{equation}
In other words, a part of the initial point mass remains frozen at its location but another part escapes to infinity, so that the whole solution vanishes in finite time $T=M/4\pi$. Moreover, for $m>0$ small we have the uniform estimate
$$
B_m(x,1)\le (4m \, t/|x|^2)^{\frac1{1-m}},
$$
which is independent of $M$.
\end{theorem}

\begin{figure}[ht!]
\centerline{
\includegraphics[width=0.5\textwidth]{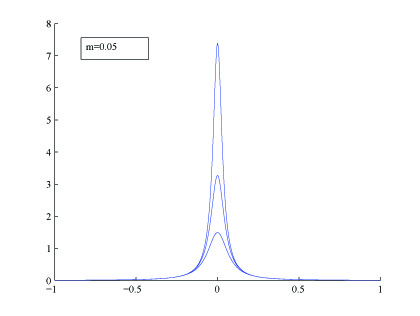}
\quad \includegraphics[width=0.5\textwidth]{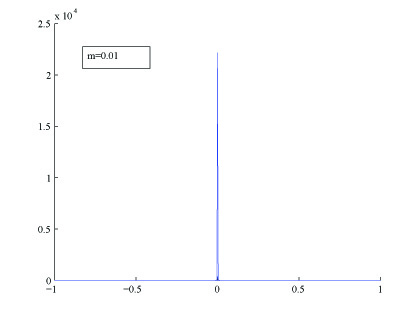}
}
\vspace{-0.5cm}
\caption{Left: Fast diffusion close to Dirac delta  for $m> 0$. Right: The needle pattern formed as  $m\to 0$}
\label{fig:FDE.m.0}
\end{figure}

The phenomenon of mass loss for the planar logarithmic diffusion when different nonnegative initial data $u_0$ are prescribed is studied in detail in \cite{Vazquez2006}, and the ``magical number'' $4\pi$ appears as the mass loss rate (per unit time) of the so-called maximal solutions. The universal law of mass loss holds for all maximal solutions with finite mass\footnote{ Since the initial value problem does not have uniqueness of solutions, a special class of uniquely defined maximal solutions is chosen.}, and the universal rate is
$$
\frac{dM}{dt}=-4\pi,
$$
as long as $u$ is not the null function.  We point that  there is no uniqueness of solutions for the Cauchy problem and higher mass loss rates are also accepted by other solutions. The  class of maximal solutions is not the only interesting one on geometrical grounds. In fact, $8\pi$ is the preferred rate on geometrical grounds for the application to model the Ricci flow in 2D \footnote{ Such solutions do not belong to the maximal class of course.}. This rate class contains a very special solution in explicit form
\begin{equation}\label{form.Ricci2}
U(x,t;T)=\frac{8(T-t)}{(1+ |x|^2)^2}.
\end{equation}
that we should add to our collection of highly representative self-similar solutions of diffusive-convective processes. Moreover, if we introduce an admissible parameter $a>0$ we get the family of solutions
$$
U(x,t;a,T)=\frac{8a^2(T-t)}{(a^2+ |x|^2)^2},
$$
that upon letting $a\to 0$ gives the measure-valued solution
$$
\widehat U(x,t;T)=8\pi(T-t)_+\delta(x).
$$
It clearly loses $8\pi$ units of mass per unit of time. Note that in this case the difference of two
solutions has a form that is constant in time for all $L^p$ norms
$$
U(x,t;T_1)-U(x,t;T_2)=8(T_1-T_2)\frac{1}{(1+ |x|^2)^2}
$$
This holds until one of them finds its extinction.

\medskip

\noindent {\bf Remark.} Let us recall what happens to the very singular solutions. Taking the limit $m\to0$
in formula  \eqref{for.fpme.vss}  we see that when $N=2$ we have $a(m)\sim 4m^2\to 0$ so that the limit is identically zero for $x\ne 0$, $t\to 0$, \nc and $\widehat U(x,t;m=0) $ is just an infinite Dirac mass located at $x=0$.

\noindent {\bf Remark.} The non-existence of fundamental solutions of the equation in dimension $N\ge 3$ was discussed by Hui in  \cite{Hui94}. Note that in those dimensions it is not a critical problem, since $m_c>0$.

\noindent {\bf Remark.} If we extend the scope of our investigation by allowing exponents $m<0$ then immediate extinction (i.e., in zero time) was proved by the author in \cite{Vnonex92} for all finite  mass solutions for all exponents $m<0$ if $N\ge 2$. This strong discontinuity at the initial time $t=0+$ may remind us of similar discontinuities for $ t=0+$  for the PME semigroup when $m=\int$ of the PLE when $p=\infty$.

\medskip


\subsection{The PLE limit with $p=1$ in 1D. Total Variation Flow}\label{sec.n1pc1}

This is a novel investigation of this paper. Surprisingly, it has a very similar flavor.

\begin{theorem}\label{delta.pc1} Let $N=1$ and let $B_p(x,t;M)$ the fundamental solution of the PLE with exponent $p>1 $ and mass $M>0$. Then as $p\to 1$ we have the singular limit
\begin{equation}\label{limit.bs.pc1}
\lim_{p\to 1} B_p(x,t)=(M-2t)_+\delta(x) \quad \forall x\in    \ren, \ t>0.
\end{equation}
In other words, a part of the initial point mass remains frozen at its location but another part escapes to infinity, so that the whole solution vanishes in finite time, $T=M/2$. Moreover, for $\ve=p-1>0$ small we have the uniform estimate
$$
B_p(x,1)\le (2\ve  t/|x|^p)^{\frac1{2-p}},
$$
which is independent of $M$.
\end{theorem}

\noindent {\sl Proof.} Here are the details that show how the fundamental solution in this case loses mass as time goes and disappears completely in finite time\nc, contrary to what happens in Theorem \ref{delta.pc}.
If $N=1$ we have $p_c=1$ and we put $\ve=p-1>0$, we follow the proof of the mentioned theorem and get
\ $\alpha(p)=\beta(p)=1/(2p-2)=1/(2\ve).$ Also,
$$
F_p(y)=\left(C_p+ k_p |y|^{\frac{p}{p-1}}\right)^{-\frac{p-1}{2-p}}, \quad
k_p=\frac{2-p}{p} \beta^{\frac{1}{p-1}}\sim  (2\ve)^{-1/(p-1)}=\big(\frac1{2\ve}\big)^{1/\ve}.
$$
The free constant $C_p$ is used to fix the mass to the value $M>0$, i.e.,
$$
2\int_0^\infty \left(C_p+ k_p r^{p/(p-1)}\right)^{-\frac{p-1}{2-p}}\,dr=M,
$$
Putting  $r= s (C_p/ k_p)^{(p-1)/p}$ we get
$$
2 \, C_p^{-\frac{p-1}{2-p}+\frac{p-1}{p}}k_p^{-\frac{p-1}{p}} \int_0^\infty \left(1+  s^{p/(p-1)}\right)^{-\frac{p-1}{2-p}}\,ds=M.
$$
Note that
$$
(p-1)(\frac{1}{p}-\frac1{2-p})=(p-1)\frac{2-2p}{p(2-p)}v \sim -2\ve^2.
$$
We now  have to control the integral
$$
I_p= \int_0^\infty \left(1+  s^{p/(p-1)}\right)^{-\frac{p-1}{2-p}}ds,
$$
We put $s=t^{(p-1)/p}$ and recall some formulas for Euler's Beta function to get
$$
I_p=  \frac{p-1}{p} \,\int_0^\infty \left(1+  t\right)^{-\frac{p-1}{2-p}} t^{\frac{p-1}{p}-1}dt= \ve \,{\mathcal B}(\frac{p-1}{p},q_p), \quad
$$
with
$$
q_p= \frac{p-1}{2-p}-\frac{p-1}{p} \sim 2\ve^2.
$$
Therefore,
$$
I_p=\ve \, {\mathcal B}(\frac{\ve}{p},2\ve^2)\sim \ve\frac{\Gamma(\ve)\Gamma(2\ve^2)}{\Gamma(\ve+ 2\ve^2)}\sim \frac{1}{2\ve}.
$$
Summing up, we can estimate $C_p$ as
$$
C_p^{\,2\ve^2}\sim 2 M^{-1}\,(2\ve)\, \frac{1}{2\ve} \sim \frac{2}{M}.
$$

$\bullet $ Going back to the expression of the solution, we have
$$
B_M(x,t)=t^{-\alpha}F_p(x\,t^{-\beta};M)=t^{-\alpha}\left(C_p+ k_p\, |y|^{\frac{p}{p-1}}\right)^{-\frac{p-1}{2-p}}
 $$
with $y=x\,t^{-\beta}$. Since $\alpha=\beta= 1/(2(p-1))$
we get
$$
 B=\left( C_p t^{(2-p)\beta/(p-1)} + k_p \,t^{-1/(p-1)}|x|^{p/(p-1)}\right)^{-\frac{p-1}{2-p}},
$$
Inserting the value of $C_p$ we get
$$
 B_M(x,t) \sim\left( \left(\frac{2t^{2-p}}{M}\right)^\frac1{2(p-1)^2} + k_p \left(\frac{|x|^p}{t}\right)^{1/(p-1)}\right)^{-\frac{p-1}{2-p}}
$$
It immediately follows that for every $x,t$ we have the bound
$$
 B(x,t)\le k_p^{-\frac{p-1}{2-p}}\frac{t^{\frac{1}{2-p}}}{|x|^{-\frac{p}{2-p}}}\sim (2\ve  t/|x|^p)^{\frac1{2-p}}
$$
which tends to zero uniformly as $p\to 1$ away from $x=0$. Actually the bound
$$
B(x,t)\le \left(\frac{2t^{2-p}}{M}\right)^\frac{-1}{2(p-1)(2-p)}
$$
implies that when $t^{2-p}>M/2$ we have
$$
\lim_{p\to 1} B_M(x,t)=0 \quad \mbox{uniformly in } x\in\ren,
$$
which means extinction in finite time $T=M/2$. Moreover, for all fix $x>0$ we have for $p\sim 1$
$$
-B_x=\frac1{2-p}\frac{(2\ve t)^{\frac{1}{2-p}}}{|x|^{\frac{2}{2-p}}}\sim \frac{2\ve t}{|x|^2}.
$$
This means that mass flowing  out from the interval $(-R,R)$ is
$$
\frac{dM_R(t)}{dt }= 2\left.|B_x|^{p-1}\right|_{x=R}\sim 2,
$$
the same for all large $|x|$. This means that $M(t)=M-2t$ as long as $M(t)>0$, hence for $0<t<T=M/2$. \qed

\medskip

\noindent {\bf Remark.} The sup norm of the solution behaves as $p\to 1$ like
$$
 B_p(0;M)= \left(\frac{2t^{2-p}}{M}\right)^{-\frac{1}{2(p-1)(2-p)}},
$$
which tends to  infinity for $0<t<T$ and to zero for $t>T$. This reflects the described Dirac delta formation.

\subsection{Behaviour of general solutions}\label{ssec.behgs}

The limit equation of the PLE when $p\to 1$ is the so-called Total Variation Flow $\partial_t u=\mbox{div}(\nabla u/|\nabla u|)$, see \cite{ACM2004} for a general reference. In 1D it reads
\begin{equation}\label{eq.tvf}
\partial_t u=\left(\frac{u_x}{|u_x|}\right)_x.
\end{equation}
The behaviour of the solutions in 1D with different initial data \nc is described by Bonforte and Figalli in \cite{BonFig2012}, Proposition 2.10, as follows.

\begin{proposition} Let $u (x, t )$ be the solution to the Cauchy problem for the TVF posed in $\re$ starting from a nonnegative and compactly supported initial datum $u_0 \in L^1 (\re)$ with mass $M=\int_\re u_0(x ) \,dx$. Then there is extinction in finite time and the following
formula holds
\begin{equation}\label{massform.tvf}
\int_\re u( t , x ) \,dx = \int_\re u_0(x ) \,dx- 2t
\end{equation}
for $0<t<T$, where $T = T (u_0) = M/2$ is the extinction time.
\end{proposition}

This means that the universal mass loss is established for that class of solutions with rate $dM/dt=-2$. The authors go on to establish the rate and form of extinction for suitable solutions. A  very simple proof of inequality \eqref{massform.tvf} goes as follows.

\noindent {\sl Proof. }  The concept of solution to equation of \eqref{eq.tvf} comes from the semigroup construction as used before in the previous sections on PME and PLE. Being in 1D allows to use the following very important transformation: if we integrate equation \eqref{eq.tvf} we obtain solutions of the so-called dual PME equation: $v_t=\varphi(v_{xx})$ with nonlinearity $\varphi(s)=\mbox{sign}(s)$. This new equation admits a theory of viscosity solutions that satisfy the maximum principle. This known theory justifies the calculations that follow. We observe that an easier way to justify the comparison of the integrated solutions is using the Shifting Comparison Principle introduced by the author in \cite[Lemmas 2.1, 2.2]{Vazquez83}.


Now we apply it to a solution $u$  of \eqref{eq.tvf} with mass $M$ and support in the interval $I=[a,b]$ and the fundamental solutions $u_1(x,t)= M\,\delta(x-a)$ and $u_2(x,t)= M\,\delta(x-b)$,
and then integrate to find the corresponding $v$ solutions we find the inequality
$$
v_1(x,t)\ge v(x,t) \ge v_2(x,t)\,
$$
because this inequality is certainly true at $t=0$. But both solutions $v_1$ and $v_2$ are explicit step functions:
$$
v_1(x,t)=0 \mbox{ for } x< a, \quad v_1(x,t)=(M-2t)_+ \mbox{ for } x\ge a,
$$
and similarly $v_2$ with jump at $x=b$. We conclude that for all $0<t<M/2$ the function $u(x,t)$ is supported in $[a,b]$ and
and $v(\infty,t)=M-2t$, This proves \eqref{massform.tvf}. \quad \qed

\medskip

\noindent {\bf Relation with previous Theorem \ref{delta.pc1}}. The equation is invariant under a group of transformations that preserves mass so that if $u$ is a solution, also
$$
u_k(x,t)=ku(kx,t), \quad k>0
$$
is a solution. Using the results of \cite{BonFig2012} one proves that as $k\to\infty$
$$
\lim_{k\to\infty} u_k(x,t)= (M-2t)_+\delta(x).
$$
More precisely, in Theorem 2.11 of the same paper the following family of explicit solutions of the problem is constructed with separate-variables form:
\begin{equation}\label{form.step}
U(x,t)=\frac1{a}(T-t)_+ \quad \mbox{for } \ |x| \le a, \ a>0,
 \end{equation}
and $U=0$ otherwise. This function is a compactly supported step function for every $0<t<T$, and the whole solution loses mass at the rate equal to 2 per unit time as predicted, $M(t)=2(T-t)$. As $a\to0 $ we get the degenerate fundamental solution \eqref{limit.bs.pc1}.

\medskip

\noindent {\bf Other examples.}
1) The construction of explicit extinction solutions for the TVF was addressed by Bellettini et al. in \cite{BellCN2002}. In the planar case  the results are linked to the denoising problem in the context of image reconstruction. They find the explicit solution
$$
u(x,y):=(1-\lambda_\Omega t)_+\chi_\Omega(x)
$$ as the unique entropy solution of the TVF in $\re^2$ with initial datum $ u_0(x)=\chi_\Omega(x)$, if $\Omega$ is a convex set satisfying some additional conditions.   See also \cite{BellCN2005}.

2) Let us recall what happens to the very singular solutions. Taking the limit $p\to1$
in formula  \eqref{for.fpme.vss}  we see that when $N=1$ we have $\beta=(2(p-1))^{-1}$, hence $a(p)\sim 2(p-1)\ \to 0$ so that the limit is identically zero for $x\ne 0$, $t>0$, and $V(x,t;p=1) $ is just an infinite Dirac mass sitting at $x=0$.

\medskip

We want to stress at this point the remarkable parallelism of results between TVF in 1D (with mass loss given by precise formula $M(t)= M-2t$ and the Ricci flow in 2D (with mass loss given by  $M(t)= M-4\pi t$.

\medskip

 \noindent {\bf 1D relation PME with the PLE.}
There is however another correspondence between PME and PLE worth mentioning. It is well-known that in 1D we can obtain a transformation of the usual class of integrable solutions $u(x,t)$ of the PME into bounded solutions $v(x,t)$ of the PLE by just integrating in space:
\begin{equation}\label{fn.trans}
v(x,t)=\int_a^x u(y,t)\,dy.
\end{equation}
The base point $a$ in the definition can be chosen at convenience, normally is taken as $a=-\infty$ or $a=0$.
The correspondence works in both directions, using $u=v_x$, the exponents transform as  $p=m+1$, and no problems arise for $m>0$. See this and a more general  equivalence transformation in \cite{RazSanVaz2008}.

In what follows we refer to the solutions of the previous part  of \ref{ssec.behgs} with the letter $v$, not $u$ since $u$ is reserved for $v_x$. So, when we consider nonnegative solutions $v$ of the PLE and pass to the limit $p\to 1$, as done in this section, we go to solutions $u$ of the PME/FDE
in the limit $m\to 0$. Now, we recall that $m=0$ was not identified as a critical exponent for the PME/PLE in Section  \ref{sec.fscexp} while $p=1$ is here. There is no real contradiction since the study of Section  \ref{sec.fscexp} is done for positive solutions, while here the correspondence maps those solutions $v$ into a signed solutions $u$ of the FDE, the limit equation as $m\to 0$ being
$$
u_t=(u/|u|)_{xx}=(\mbox{sign}(u))_{xx}.
$$
We are not interested in pursuing further the many  novelties of signed solutions. The correspondence \eqref{fn.trans} will be revisited in the study of conservation laws in Subsection \ref{ssec.HJ}. For the sake of curiosity, we ask the reader to determine which solution of the sign-PME corresponds to the solution \eqref{form.step} of
Equation \eqref{eq.tvf}.   \hfill \qed

\medskip

 \noindent {\bf Comment.} For further information we point out that in the paper \cite{BonFig2012} the explicit dynamic of solutions is obtained (explicit evolution
of simple functions), not only an explicit loss of mass/extinction time formula, together with the
sharp asymptotic behaviour (convergence to the separate variable solutions  with the ``right support”), with
counterexamples about possible rates of convergence. And also the limit $m \to 0+$ in Fast Diffusion is analyzed
and gives two possibilities: either the logarithmic diffusion or the sign diffusion (the latter is not
mentioned here).\nc

\newpage

\section{The Doubly Nonlinear Equation}\label{sec.dnle}

 As an extension of the study of the previous nonlinear diffusion models, we present in this section a model that combines the two nonlinearities that appeared in preceding sections. The equation is called Doubly Nonlinear Equation (DNLE). The main features will be examined and the main results stated (rather  informally)\nc, but we will not delve on the theory with the same amount of detail because of the reasons of space and goal given earlier. \nc

By the name Doubly Nonlinear Equation we mean the equation DNLE
\begin{equation}\label{DNE}
\partial_t u=\mbox{div}(|u|^{a}|\nabla u|^{p-2}\nabla u).
\end{equation}
with parameters $a$ and $p$. It can be written in the equivalent form
\begin{equation}\label{DNE2}
\partial_t u=c\,\mbox{div}(|\nabla (u^m)|^{p-2}\nabla (u^m)),
\end{equation}
where $a=(m-1)(p-1)$ and $c>0$ is an innocent constant that we may forget in this section.  The equation is posed in $\re^N$ with $N\ge 1$.

In case $p=2$ and $a=m-1$ we recover the standard PME, while for $a=0$ we get the PLE. There is an alternative way of writing de DNLE  as
\begin{equation}\label{DNE3}
\partial_t (|v|^{\mu-1}v)=c\,\mbox{div}(|\nabla v|^{p-2}\nabla v),
\end{equation}
by just putting $v=|u|^{m-1}u $ in \eqref{DNE2}, so that $u =|v|^{\mu-1}v $ with $\mu=1/m>0$.

\medskip

By choosing the above combination of nonlinearities we hope that many of the techniques and results of Sections \ref{sec.cm.conn}-\ref{sec.ncbl} can be recovered. Thus, a theory of existence and uniqueness has been produced using the Crandall-Liggett generation theorem  \cite{CrandallLiggett71} based on implicit time discretization, using the fact that the doubly nonlinear elliptic operator is $m$-accretive in $L^1(\ren)$. We will use next the notation of \eqref{DNE2}.

\bc Doubly nonlinear equations have a wide range of applications in various physical contexts, including the flow of non-homogeneous non-Newtonian fluids and the modeling of  underground water. The equation has been proposed by Leibenzon in modeling turbulent fluid problems,  \cite{Leib}, in 1947.  Is also appears in glaciology,  see \cite{C5}. This equation has been investigated mathematically  by different authors, we should mention  Kalashnikov's well-known survey paper \cite{Kalsurvey} of 1987.   Earlier work is due to \cite{Bamb77, JLL69, Rav70}.  See the book \cite{ADSh2002} and  \cite{Vazquez2006}, Chapter 11, for further studies. See more in Subsection \ref{otherissues} below\nc.

\subsection{Self-similarity}\label{dnl.sss} The fundamental solution is known to exist in the appropriate range, that turns out to be precisely: $m(p-1) + (p/N)t>1$. As expected, the formula is
\begin{equation}\label{form.sssol.dne}
U(x,t)=t^{-\alpha}F(x\,t^{-\beta}),
 \end{equation}
Its similarity exponents are now
\begin{equation}\label{self.expo.dnl}
\alpha=\frac{1}{m(p-1)- 1 + (p/N)}, \qquad \beta=\frac{\alpha}{N},
\end{equation}
In order to  get finite mass diffusion  solutions we need $\alpha, \beta>0$, which implies $ m(p-1)+ (p/N)>1$.
 The space profile has different forms in three different ranges:

\noindent {\bf \bc (1)} The range of exponents that corresponds to slow diffusion and free boundaries is given by \ $m(p-1)>1$ and then the profile of the fundamental solution is
\begin{equation}\label{ss.dnle2s}
F(\xi)=\left(C-
k\,|\xi|^{\frac{p}{p-1}}\right)_+^{\frac{p-1}{m(p-1)-1}}.
\end{equation}
with any $C>0$ and \bc some $ k= k(m,p,N)>0$\nc. In the notation of \eqref{DNE3} the exponent range condition is $\mu+1< p$.

\noindent {\bf \bc (2)} The space profile is given in the good fast diffusion range, $1-(p/N)< m(p-1)<1$, by the formula
\begin{equation}\label{ss.dnle2}
F(\xi)=\left(C+
k\,|\xi|^{\frac{p}{p-1}}\right)^{-\frac{p-1}{1-m(p-1)}}.
\end{equation}
with any $C>0$ and some $ k= k(m,p,N)>0$. We leave it to the reader to check that these are weak solutions of the DNLE.
In the notation of \eqref{DNE3} this range is
$$
p-1<\mu<  N(p-1)/(N-p).
$$

\begin{figure}[ht!]
\centerline{\includegraphics[width=0.45\textwidth]{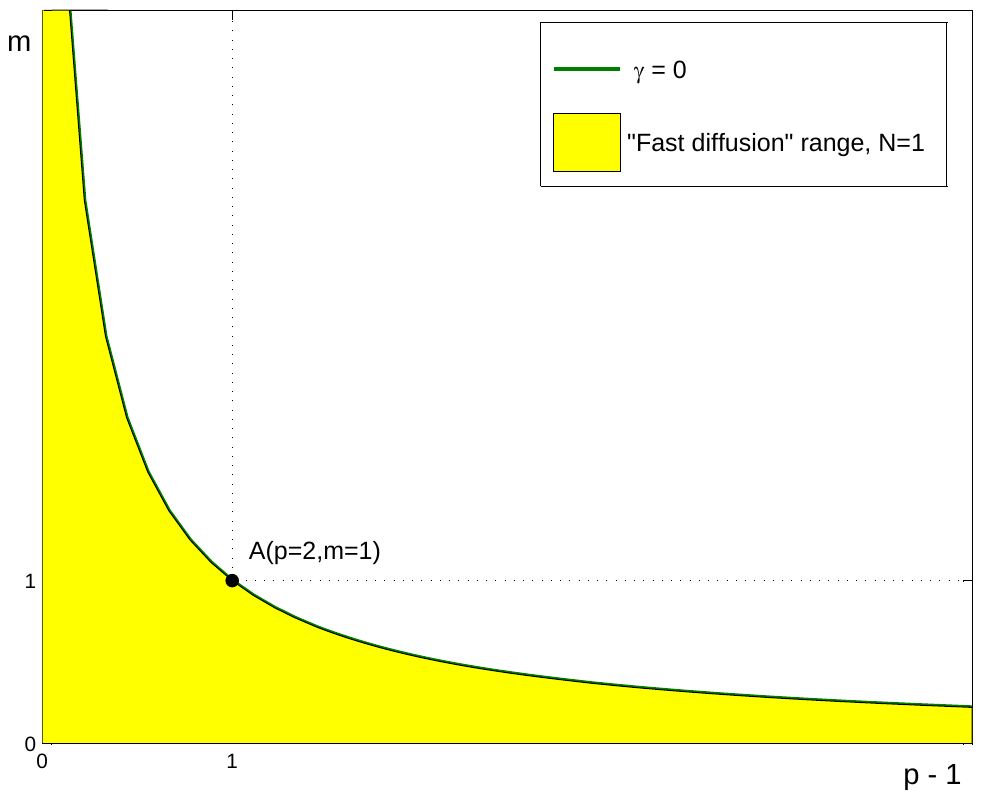} \quad
\includegraphics[width=0.45\textwidth]{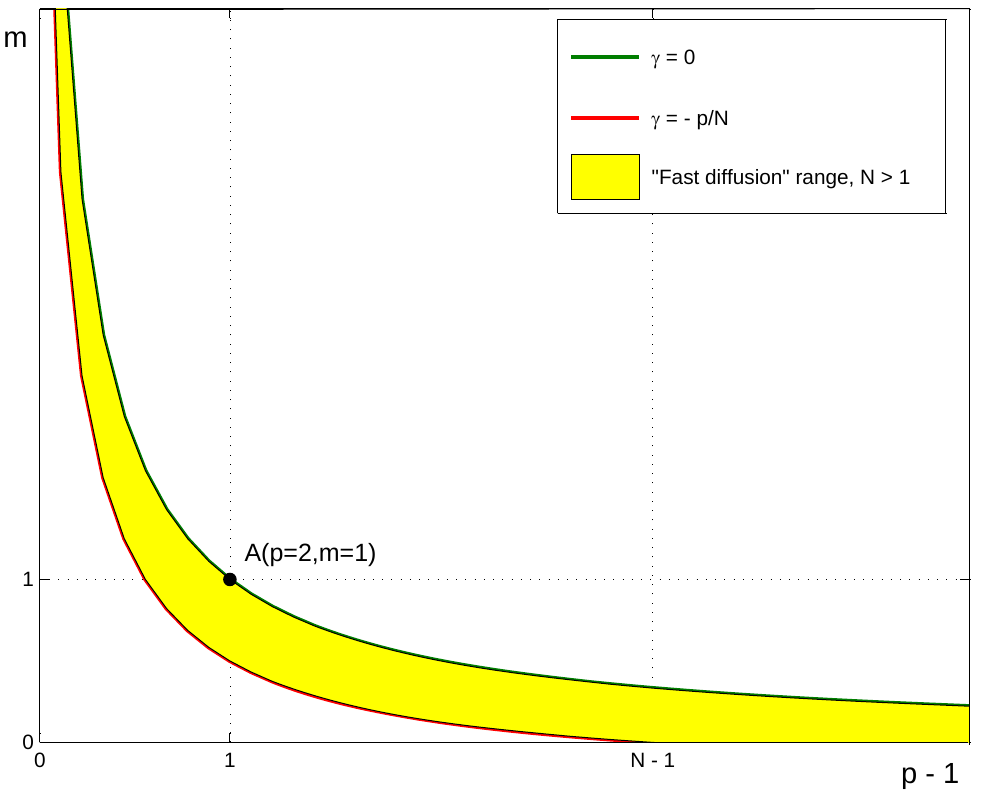}}
\vspace{-0.0cm}
\caption{Regions of slow and fast diffusion. Here, $\gamma=m(p-1)-1$.\\
Below the red line on the right-hand figure there is no conservation of mass}
\label{fig:DNLrange}
\end{figure}

\noindent $\bullet$ In this range we may let the mass of the solution $M$ tend to infinity and get the Very Singular Solution with profile
\begin{equation}\label{ss.dnl.vss}
F_\infty(\xi)=k_1(m,p,N) |\xi|^{-\frac{p}{1-m(p-1)}},
\end{equation}
and space-time form \bc
\begin{equation}\label{ss.dnl.vss2}
U(x,t)=k_1(m,p,N) \,t^{\frac{1}{1-m(p-1)}}|x|^{-\frac{p}{1-m(p-1)}}. \nc
\end{equation}

\noindent {\bf (3)} \bc  There is a line of ``quasi linear diffusion'' in the sense of having a power homogeneity given by $m(p-1)=1$, i.e., $\mu=p-1$.

\begin{proposition}\label{prop.9.dnlin} Let $N\ge 1$, $m>0$, $p>1$ and  $m(p-1)=1$. The fundamental solutions have exponents $\alpha=N/p$ and $\beta=1/p$ and the profile is also explicitly calculated for all $m>0$ as
\begin{equation}\label{ss.dnle3}
F(\xi)=C_1\,\exp(-k  |\xi|^{(m+1)/m}),
\end{equation}
where $C>0$ is  a constant that may be chosen at will and \ $ k=\beta^{m}/(m+1)$.
\end{proposition}

Let us calculate the profiles  for $m>0$, $m(p-1)=1$. We have
$$
-\beta \, (r^N F)'=-(r^{N-1}|(F^m)'|^{p-1})'
$$
where $F(r)$ must be  positive, radial (i.e., $F=F(r)$, $r=|\xi|$) and radially nonincreasing. Then, arguing as before we get
$$
\beta \, r F =|(F^m)'|^{p-1}=m^{p-1}F^{(m-1)(p-1)}|F'|^{p-1}.
$$
Now we use $m(p-1)=1$ to get the simpler equation
$$
|F'|^{p-1}= \beta m^{1-p} \,r\, F^{p-1}, \qquad F'(r(= -\frac{\beta^m}{m}\,r^{m}\, F(r).
$$
Integrating we get $\log(F)=C-k r^{(m+1)/m}$,  $ k=\beta^{m}/(m+1).$ \qed

\medskip

\noindent {\bf Remarks.} i) For $m=1$, $p=2$ we obtain  again the typical Gaussian function.
Notice that for values $m\ne 1$, this DNLE is \bc not a linear  parabolic equation, but it has  however self-similarity with Gaussian-like profiles.

ii) See \cite{StVaz2013} for the problem in a bounded domain\nc.

\medskip

\noindent {\bf \bc (4)} \nc Finally, the critical situation on the lower part of the good fast range happens now not for a single value but for a whole critical line
$$
m(p-1) + (p/N)=1.
$$
\bc Below that line there will extinction in finite time. We will revise the topic in Subsection \ref{dnl.more}
 We will not  enter into the details of the limit line here\nc.

\subsection{Inequalities  and  asymptotics}\label{dnl.prop}

\noindent $\bullet$ The Crandall-Liggett homogeneity inequality is obtained as in the PME case and reads
$$
(m(p-1)-1)\,tu_t+u \ge0
$$
which  is valid for all nonnegative solution of the Cauchy Problem as long sa $m(p-1)\ne 1$.

\noindent  {$\bullet$} The $L^1$-$L^\infty$ smoothing effect holds in the range of the Barenblatt solutions in the usual formulation

\begin{proposition}\label{prop.dnl.LLse} Let  $m(p-1)+(p(N)>1$. For every nonnegative  semigroup finite-mass  solution $u(x,t)$ of the DNLE we have
\begin{equation}\label{form.L1Linf.p1}
\sup_x u(x,t)\le B_M(0,t)= c(m,p,N)\, M^{p\,\beta}\,t^{-\alpha}\,,
 \end{equation}
and optimality of the constants is achieved by the Barenblatt solution and its space translations.
\end{proposition}

We recall that the proofs given in \cite{V82sym, V82port, Vazsym05} rely on symmetrization for the PME, and PLE equations and uses of the Barenblatt solution. Note:  It is not known whether other solutions attain the optimal constant\nc.

\noindent $\bullet$ \bc Asymptotic behaviour\nc. The large time behaviour shows that total mass in the only parameter the solution retains in the large-time asymptotics in first approximation. What is surprising is that the result is completely formally analogous to the one for
the heat equation, for the PME and the PLE.

\begin{proposition}\label{prop.dnl.ab}  Let  $m(p-1)+(p(N)>1$. For every nonnegative  semigroup  solution $u(x,t)$ of the DNLE with initial datum $u_0\in L^1(\ren)$ and every $t>0$ we have
\begin{equation}\label{form.asb.bar.dnl}
\lim_{t\to \infty}\|u(\cdot,t)- B_M(\cdot,t)\|_1=0.
\end{equation}
where $M=\int_{\ren} u_0(x)\,dx$\nc.
\end{proposition}

The proofs of all three nonlinear models, PME, PLE and DNLE, are quite similar.


\subsection{Finite-mass solutions without MC and separation of variables}\label{dnl.more}
We want to expand on the contents of  Section \ref{ss.vfpme} by revisiting the issue of the finite-mass solutions in the very fast range of the DNLE equation that now corresponds to the region of pairs $(m,p)$ such that $m>0$, $p>1$ and $m(p-1) +(p/N)<1 $.  \bc We are expecting solutions with loss of mass and even extinction. Let us just mention the existence of the pseudo-Barenblatt solutions
\begin{equation}\label{dnl.pseudo}
U(x,t)= (T-t)_+^{-|\alpha|}F(x\,(T-t)^{-|\beta|}),
\end{equation}
with $\alpha$ and $\beta$ as in  \eqref{self.expo.dnl}, and the profile $F(\xi)$ is as in \eqref{ss.dnle2}. It must be observed that the very fast range condition $m(p-1)-1+(p/N)<0$ implies that the decay rate of the spatial profile of the solutions is
$$
\frac{p}{1-m(p-1)}<N,
$$
hence we are not dealing any more with solutions with finite mass, they are larger. However, they extinguish in finite time as expected from our  detailed analysis in Section  \ref{ss.vfpme}.

We  also have the pseudo VSS solution
\begin{equation}\label{ss.dnl.vss3}
U(x,t)=k_1(m,p,N) \,(T-t)^{\frac{1}{1-m(p-1)}}|x|^{-\frac{p}{1-m(p-1)}}. \nc
\end{equation}
with arbitrary $T>0$. This special function has an integrable singularity at $x=0$ but it is integrable for large $x$.

The rest of results od Section \ref{ss.vfpme} can be studied and many are quite similar. We think this is a good exercise for the reader.
Detailed results will depend of an accurate analysis of the suitable
choices of self-similar solutions, an analysis that has not been fully done.

\medskip

\noindent  {\bf Separation of variables.} This is similar to what happened in the PME but it has an extra difficulty. \nc We try the form
$$
U(x,t)=(T-t)^{\alpha}F(x).
$$
Substituting into the DNLE we get the two conditions:

(i) we put $\alpha-1=\alpha m(p-1)$ \  in order to eliminate time from the equation, so that $\alpha=1/(1-m(p-1))$, and

(ii) the Lane-Emden equation
$ -\Delta_p F^m= \alpha\, F\,,$ that is better recognized as
$$
-\Delta_p G= \alpha\, G^q\,,
$$
with $G=F^m$ and $q=1/m$. We need a nonnegative, bounded, and finite-mass solution $F$ defined in whole space. The literature is clear about this point, see \cite{BidVer1989, Sciun16}: such a solution exists in the range $1<p<N$ for the precise value $q=p^* -1, $ where $p^*=Np/(N-p)$ is the Sobolev exponent associated to $p$. There is even uniqueness of one radial solution up to space translation and re-parametrization:
\begin{equation}\label{form.sep.var.dnle}
G(x)=  F^m(x)=
\kappa(N,p)\,\left(\frac{a}{a^{p}+ |x|^{p/(p-1)}}\right)^{\frac{N-p}{p}},
 \end{equation}
where and $a>0$ is a free parameter. We need to state the relation between $p$ and $q$ ($q=p^*-1$) in terms of  $m=1/q$ and $p$. We get in this way the line of solutions with separation of variables:
$$
m(p-1) + (m+1)\frac{p}{N}=1.
$$
As for the profiles, it is pointed out in \cite{Sciun16} that these are the Aubin-Talenti optimizers of the Sobolev inequalities
from $W_0^{1,p}(\ren)$ into $L^{p^*}(\ren)$. Here is a clear relationship between asymptotic profiles in fast diffusion and minimizers in the Calculus of Variations.

\medskip

\noindent  {\bf Note.} When $p=2$ we get $m=(N-2)/(N+2)$, the value called $m_s$ in Section \ref{ss.vfpme} when studying the Very Fast PME, see pattern  function in \eqref{form.sep.var.fde}. On the other hand, for $m=1$ (case of the PLE), we get $p_s= 2N/(N+2)$.

\subsection{Other issues and equations}\label{otherissues}
Many of the issues addressed in this paper for the PME and the PLE still hold in this combined setting, but a number of problems are still open for the DNLE, and it should be workable to answer them. Concerning the theory and basic estimates there are many results by a large number of authors, see the early papers  \cite{EV86, Vespri92}. The first proof of asymptotic behavior was done (in one space dimension)  in \cite{EV88}. See also \cite{AkStef2010, BogDMS18, Ivanov97, StVaz2013}. For a formal asymptotic analysis see \cite{King16}. For very recent work on the equation see \cite{Vestb2024}. Doubly nonlinear equations involving the $p$-Laplacian for the value $p=1$ are rather special, see a recent work in \cite{MMT23}.  The questions of mass conservation seem to be only partially known.

\medskip

\noindent {\bf Other types of related equations.} (1) As a variant of the doubly nonlinear equations we can mention the flux-limited models, related to the so-called relativistic heat equation
$$
\partial_t u=\mbox{div\,}\left( \frac{u\,\nabla u}{\sqrt{u^2+ \kappa |\nabla u|^2} } \right).
$$
Such  equations have been  studied in \cite{KR97, Brenier03, Cher2003, ACM2006, ACMS2012, CarrCM2013} and other references.
We will  see somewhat similar equation forms below under the heading ``geometric flows''.

\medskip

\noindent (2)  There is a variant of the $p$-Laplacian equation, often called the {\sl normalized $p$-Laplacian equation}  that, in its evolution version, takes the form
\begin{equation}\label{eq.rpl}
\partial_t u= |\nabla u|^{2-p}\Delta_p (u)=\Delta u + (p-2)\frac1{|\nabla u|^2}\sum_{i,j} \partial_i u\,\partial_j u\,\partial^2_{ij} u
\end{equation}
The normalized $p$-Laplacian has non-divergence form and appears in stochastic tug-of-war games, hence the name game-theoretical $p$-Laplacian, \cite{Peres08}. There is a viscosity solution theory for such equations and no conservation of mass holds if $p\ne 2$. This model has been studied by many authors in recent times, see for instance \cite{Does11, ImbJinP19, JuuKaw06, JuuLindMan01, ParvVaz20, PorVaz12}. The limit $p\to \infty$, written as
\begin{equation}\label{eq.rpl.inf}
\partial_t u= \frac1{|\nabla u|^2}\sum_{i,j} \partial_i u\,\partial_j u\,\partial^2_{ij} u
\end{equation}
is also interesting. As a mathematical problem it is related to the theory of absolutely minimizing functions, cf. \cite{AronssCJ04}. This limit case may reminds us of the mesa problem for the PME but is not at all similar in the mathematical content.  Though this equation is heavily nonlinear in several dimensions, in $N=1$ it just becomes the Heat Equation, cf. \cite{PorVaz12, PorVaz13}\nc.

\medskip

\noindent (3) There are many studies on Doubly Nonlinear  equations of a more general form like
$$
\partial_t b(u)=\mbox{div}(D_\xi F(x,\nabla u)),
$$
 where $b(u)$ is a monotone scalar function and $F(x,\xi)$ a suitable function with  convexity and coercivity properties, cf. \cite{BogDMS18}\nc.

\medskip

\noindent (4) Another direction for expanding the scope of doubly nonlinear evolution is to replace the $p$-Laplacian operator by a nonlocal version, as done in \cite{HyndLind16} where the following doubly nonlinear and nonlocal model is studied
$$
|v_t|^{p-2}v_t+(-\Delta_p)^s v=0,
$$
where $p\in (1,\infty),$ $ s\in (0,1)$, and $(-\Delta_p)^s$ is the fractional $p$-Laplacian, see next sections for nonlocal operators.

In summary, the questions of mass conservation seem to be only very partially known in these more general models. Reasonable conjectures can be made for cases where mass conservation holds.


\vskip3.5cm

\newpage

{\bf \huge Part 2}

\section{Diffusion with fractional operators}\label{sec.nonlinfrac}

The  considerations made in previous sections apply to similar problems arising in the more novel topic of  diffusion involving nonlocal integro-differential operators, in particular fractional Laplacian operators. These equations have attracted in recent years an increasing amount of attention due to mathematical interest in different mathematical disciplines and many ongoing applications in Science.

\bc As mentioned in \cite{VazTNLD2016}, replacing the Laplacian operator by fractional Laplacians is motivated by the need to represent processes involving anomalous diffusion. In probabilistic terms, it features long-distance interactions instead of the next-neighbour interaction of random walks and the short-distance interactions of their continuous limit, the Brownian motion. The main mathematical models used to describe such processes are the fractional Laplacian operators, since they have special symmetry and invariance properties that makes for a richer theory.  These operators are generators of stable L\'evy processes that include jumps and long-distance interactions. They  reasonably account for observed anomalous diffusion where particle movement deviates essentially from standard Brownian motion (and found in materials and biological systems). There are  applications in continuum mechanics (elasticity, crystal dislocation, geostrophic flows,...), phase transition phenomena, population dynamics, optimal control, image processing, game theory, finance, and others. There are also applications in Mathematical Physics,  modeling complex systems with memory effects and multi-scale behaviors. We refer to the abundant literature, like  \cite{Apple2009, Bert96, ContTankov, GiOsh08, MMM, MK2000, VIKH, Woy2001}, see also Section 1.2 of \cite{Vaz14}. Let us also point out that this type of equations  can be studied via very efficient numerical approximations\nc.

Summing up, fractional diffusion and other non-local diffusion equations offer a useful alternative for describing systems where the standard diffusion equation fails to accurately represent the underlying physical processes. They are particularly useful in situations involving complex media, anomalous transport, and non-Gaussian statistics. This relevance has spurred the already high interest of the topic from the analysis point of view.


\subsection{The fractional Laplacian}\label{ssec.fracLap}

 The most studied nonlocal integro-differential operator is the fractional Laplacian operator, defined for suitable functions $f(x)$, $x\in \ren$, as
	\begin{equation}\label{eq.frlap}
	(-\Delta)^s f(x) = c(N,s)  PV \int_{\ren} \frac{f(x) - f(y)}{|x-y|^{N+2s}} \,dy,
	\end{equation}
	where $s\in(0,1)$ is the fractional parameter, mainly $s=1/2$ but all values $0<s<1$ are now considered. The standard value of the constant ia
$$
c(N,s) := \frac{2^{2s}s \,\Gamma ( \frac{N}{2} + s)}{ \pi^{\frac{N}{2}} \Gamma(1 - s)},
$$
and $PV$ denotes that the integral is taken in the principal value sense, cf. the classical references \cite{Landkof66, Stein1970}, as well as recent ones like \cite{BucurBook16, CaffarelliSilv2007, DiNPV12, GarofaloMR3916700, Kwasnicki2017}.

The corresponding evolution equation we will consider is the fractional heat equation (FrHE)
\begin{equation}\label{eq.frachteq}
\partial_t u + (-\Delta)^s u=0
\end{equation}
posed in the Euclidean space $x\in \ren$, $ N\ge 1$, for $t>0$.  We supplement the equation with an initial datum
\begin{equation}\label{init.data}
\lim_{t\to 0} u(x,t)= u_0(x),
\end{equation}
with suitable initial conditions, typically $u_0\in L^q(\ren)$ with $1\le q\le \infty.$ We have shown that existence, qualitative behaviour, mass conservation and asymptotic behaviour can be established much as it was done for the heat equation, see \cite{BPSVal2014, BonSirVaz2017, VazquezFHE2018}.  When comparing the fractional heat equation with the usual HE, this one is often referred as  standard diffusion equation, or equation with the standard Laplacian\nc.  Elliptic equations with nonlocal integro-differential operators are now a fast developing parallel field, see for instance \cite{BucurBook16, FRRO24, RosOton, RosOton}.

\subsection{Mass conservation and its connections  for the FrHE}\label{ssec.fhe}

The analysis we have performed in Subsection \ref{ssec.cmheq} for the Heat Equation can be adapted to the Fractional Heat Equation, $\partial_t u + (-\Delta)^s u=0$, for all $s\in(0,1)$.  In this case there are no major surprises. Under the assumptions of  self-similarity and mass conservation  we arrive at the conjectured self-similar form
$$
U(x,t)=t^{-\alpha}F(x\,t^{-\beta}),
$$
where now $\alpha=N/(2s)$ and $\beta=1/(2s)$,  which are larger quantities than in the standard case $s=1$. These values are again a consequence of dimensional considerations. The Brownian scaling is now replaced by $|x|\sim t^{\beta}$. The profile $F$ is not so easy to find for all $s$. There is an explicit case for the value $s=1/2$ and then the profile is explicit
\begin{equation}\label{fs.fpme}
F(y)=C\, (1+|y|^2)^{-(N+1)/2}.
 \end{equation}
The reader will recognize the Poisson kernel for the Laplacian in the upper half-plane in $\re^{N+1}$. For other values of $s\in (0,1)$ an explicit form is not known, but Blumental and Getoor  \cite{BlumGet1960} used Fourier analysis   to find a sharp description of the behaviour, and this is all we need. Indeed, for all $s\in (0,1)$ the profile $F$ is a $C^\infty$ smooth function with radially symmetric form, it is decreasing in $r=|y|$ and behaves as $r\to\infty$ like the precise power tail
$$
F(y)\sim C(N,s) \,|y|^{-(N+2s)}\,.
$$
We have already found these fat tails in the Barenblatt profiles for the fast PME and PLE, two nonlinear models we have studied above involving the standard Laplace operator. They reappear here in a linear model with nonlocal diffusion.

Once this is settled, the representation formula holds, and CM holds for all solutions with $u_0\in L^1(\ren)$.  Nonnegative solutions are indeed positive with a tail at infinity equal or larger than the one predicted by the fundamental solution, \cite{BPSVal2014, BonSirVaz2017}. Finally, as $t\to\infty$ we have convergence of all finite-mass solutions to the fundamental solution with the same mass, just as in \eqref{he.asymp}, \eqref{he.asymp.Li} with $G$ replaced by $U$. \bc The detailed analysis of this asymptotic  convergence is rigorously proved in our paper \cite{VazquezFHE2018}. \bc Furthermore, a spectral theory that enables to get an asymptotic development has been recently announced by Chan et al. in \cite{ChanFMM24}\nc.

\medskip

\subsection{  Nonlinear evolution equations with fractional diffusion}
These equations  arise in many contexts: in the quasi-geostrophic flow model, in boundary control problems, in surface flame propagation, in statistical mechanics, in financial mathematics, and many others.

\medskip

\noindent {\bf (NL1)} A popular nonlinear evolution model is the so-called fractional porous medium equation (FrPME)
\begin{equation}\label{eq.fracpmeq}
 \partial_t u + (-\Delta)^s (|u|^{m-1}u)=0
\end{equation}
posed in the Euclidean space $x\in \ren$, $ N\ge 1$, for $t>0$.  We supplement the equation with an initial datum \eqref{init.data} as before,
where now the correct initial space is $u_0\in L^1(\ren)$.  It was extensively studied in the last decade, the basic theory is covered in \cite{dP11, dP12, dP17}. The Bénilan-Crandall theory allows to find a mild solution for every $u_0$ in that class and the set of solutions forms a continuous semigroup of $L^1$ contractions. These solutions are unique, they are bounded for positive $t$ and satisfy the equation in the usual weal sense. The self-similar solution with mass conservation has been found by rather standard methods under the assumption that $m> m_{cs}:=(N-2s)_+/N$.  The important question of uniqueness is solved in \cite{VazqJEMS2014}

\begin{theorem} \label{thm.Bs} For every choice of parameters $s\in(0,1)$ and $m>m_{cs}$, and every mass $M>0$, equation \eqref{eq.fracpmeq} admits a unique fundamental solution; it is a nonnegative and continuous weak solution for \ $t>0$ and takes a Dirac initial data $M\,\delta$ as a trace in the sense of Radon measures. Such solution has the self-similar form
\begin{equation}\label{sssf}
u_M^*(x,t)=t^{-\alpha} F(|x|\,t^{-\beta})
 \end{equation}
for exponents $\alpha$ and $\beta$ that can be calculated in terms of $N$ and $s$ in a dimensional way, precisely
\begin{equation}\label{scale.expo}
\alpha=\frac{N}{N(m-1)+2s}, \qquad \beta=\frac{1}{N(m-1)+2s}\,.
\end{equation}
The profile function  $F(r)$, $r\ge 0$, is a  bounded, positive and  smooth function, it is radially monotone and  goes to zero at infinity.
\end{theorem}

The precise rate of decay of $F$ is delicate and can be found in \cite{VazqJEMS2014}.  The following attraction theorem completes the basic picture

\begin{theorem}\label{thm.exlimit} Let $m>m_{cs}$, let $u_0=\mu \in {\cal M}_+(\ren)$ a Radon measure, and let $M=\mu(\ren)$. If $u^* _M$ be the self-similar Barenblatt solution with mass $M$, then  we have
\begin{equation}\label{conver.express}
\lim_{t\to\infty} t^{\alpha}\,|u(x,t)-u^*_M(x,t;M)|=0\,,
\end{equation}
and the convergence is uniform in $\ren$.
\end{theorem}

Conservation of mass was proved for the semigroup solutions with initially finite mass if $m>m_{cs}$, and is false for $m<m_{cs}$,
Results about existence of VSS for a certain range of $m>m_{cs}$ are examined in the paper, as well as extinction in finite time for $m<m_{cs}$.
See more about the topic in the introduction to Section \ref{sec.cmass.fpme}.

\medskip

\noindent {\bf (NL2)} A nonlinear version of the fractional Laplacian is given by the fractional $p$-Laplace operator that may be introduced as follows.
The nonlocal energy functional for functions $f(x)$ defined in $\ren$ (or in an open subdomain)
\begin{equation}\label{Jsp1}
{\mathcal J}_{p,s}(f)= \frac1{p}\nc \int_{\ren }\int_{\ren } \frac{|f(x)-f(y)|^p}{|x-y|^{N+sp}}\,dxdy\,.
\end{equation}
is a power-like functional with nonlocal kernel of the $s$-Laplacian type that has attracted a great deal of attention in recent years. It is just the $p$-power of the Gagliardo seminorm, used in the definition of the $W^{s,p}$ spaces (fractional Sobolev, Slobodeckii or Gagliardo spaces) with seminorm and norm given by
 $$
 \quad [f]_{s,p}^p=p{\mathcal J}_{p,s}(f), \qquad \|f\|_{s,p}^p=\int |f|^p\,dx + p{\mathcal J}_{p,s}(f),
 $$
 cf. \cite{AdFou, DiNPV12}. We  consider the functional  ${\mathcal J}_{p,s}$ for   exponents $0<s<1$ and $1<p<\infty$ in dimensions $N\ge 1$. Its  subdifferential (or Euler-Lagrange operator)
${\mathcal L}_{s,p}$  is the nonlinear operator defined  a.e. by the formula
\begin{equation}\label{frplap.op}
\displaystyle \qquad {\mathcal L}_{s,p}(f):= \, P.V.\int_{\ren}\frac{\Phi(f(x)-f(y))}{|x-y|^{N+sp}}\,dy\,,
\end{equation}
where  we write $\Phi(z)=|z|^{p-2}z,$ and $PV$ means principal value. It is a usually called the $s$-fractional $p$-Laplacian operator. It is then well-known from general theory,  see Brezis \cite{BrezisBkOMM73}, that ${\mathcal L}_{s,p}$ is a maximal monotone operator in $L^2(\ren)$ with dense domain.  A multiplicative constant $c(s,p,N)>0$ is usually inserted into the  definition \eqref{frplap.op}. This constant will not affect our evolution results, as shown below, in any essential way. In any case, the choice of a suitable constant is discussed in Section 5 of \cite{dTGCV21}, mainly in connection with the limit $s\to 1$. \nc

We will study the corresponding gradient flow, i.e., the evolution equation
\begin{equation}\label{frplap.eq}
\partial_t u + {\mathcal L}_{s,p} u=0,
\end{equation}
posed in the Euclidean space $x\in \ren$, $ N\ge 1$, for $t>0$. We will often refer to it as the EFPL equation (evolution fractional $p$-Laplacian equation). We take initial data as in \eqref{init.data},
where in principle $u_0\in L^2(\ren)$. However, the theory shows that Equation \eqref{frplap.eq} generates a continuous nonlinear semigroup in any $L^q(\ren)$ space, $1\le q<\infty$, in fact it is a non-expansive semigroup for every $s,p$ and $q$ as specified.

%
%
Ir was studied in our previous papers \cite{Vazquez2020, VazFPL2-2020, VazFPL3-2021}. We refer to it as the {\sl fractional $p$-Laplacian evolution equation}, FrPLE for short.  Motivation and related equations for this model can be seen in the  \cite{Vazquez2020} and its references. There, the superlinear case $p>2$ was studied. The case $1<p< 2$, usually called fast diffusion, has been later treated in \cite{VazFPL2-2020}, mainly for the class of solutions that are $L^1$-integrable with respect to the space variable.

\medskip

\noindent {\bf Remark.} In the limit as $s\to 1$ we obtain ${\mathcal L}_{s,p}(u)\to \Delta_p(u)$ but  only if the formula \eqref{frplap.op} is normalized with a suitable factor $c(s,p,N)$. Likewise $(-\Delta)^s f\to \Delta f$\nc.

\noindent {\bf Reminder.} In this part of the survey we assume the values $0<s<1$, $0<m<\infty$ and $1<p<\infty$.


\noindent{\bf Outline of recent and new results.} Next, we will explore the connection between self-similarity, fundamental solutions and mass conservation in some examples of  nonlinear fractional diffusion. \bc Our attention is not so much focused on the general theory that advanced readers can easily find in the cited literature, but attention to novelties and open directions regarding mass conservation\nc. In Section \ref{sec.cmass.fpme} we examine nonlinear fractional diffusion, where the main progress has been done in recent years. Here we find  two further new contributions of the present paper, namely the proofs of conservation of mass in the critical fractional nonlinear cases. The critical case of the FrPME is treated in Theorem \ref{thm.mc.fpme}. Mass conservation in the critical case of the FrPLE occupies Section \ref{sec.cmass.fpl}, see Theorem \ref{thm.mc.fple}.
%

 \newpage
\section{MC in Fractional Nonlinear Diffusion}\label{sec.cmass.fpme}

 It is known that the law of mass conservation holds for strong solutions of the Fractional PME (FrPME for short) in the range \ $m>m_c=(N-2s)/N$, $m>0$, when posed in the whole space with initial data in $L^1(\ren)$. This was proved by our team in \cite{dP12} as part of a rather complete semigroup theory for this equation, we will often call them the semigroup solutions. In the range $m>m_c$ the existence of self-similar Barenblatt solutions is established and boundedness inequalities  of the $L^1$-to-$L^\infty$ type hold. There is no property of finite propagation, due to the influence of the nonlocal operator, and the spatial tails are described. The role of the fractional Barenblatt solutions as asymptotic   attractors is proved in \cite{Vaz2015}\nc.

On the other hand, for every $0<m<m_c$ we can find plenty of  semigroup solutions (also constructed in \cite{dP12}) that are shown to have finite extinction time  \cite{VazqJEMS2014, VazVol2014, VazVol2015}, and the results of Section \ref{ss.vfpme} apply. Hence, we are only left with the open question of MC for the critical case $m=m_c$.

Let us mention in passing that there is a theory of existence and uniqueness for ``large data'' (typically solutions $u(x,t)$ whose data  grow as $|x|\to\infty$) valid or the fractional diffusion of Fast Diffusion type, see \cite{BonVa2014}. The topic is not dealt in this survey, but  extinction in time estimates and lower bounds are interesting.
On the other hand, the case $ m < 0$ has been studied in \cite{BonSeVaz2016} proving the non-existence results in the fractional case, analogous to the ones obtained by the author in 1997, mentioned in the ``standard (non-fractional) part” of this survey.\nc

In the case of the fractional $p$-Laplacian equation, FrPLE, the situation is similar, see Section  \ref{sec.cmass.fpl} below.

\subsection{Fractional  FDE with critical exponent}\label{sec.cmass.fpme.cr}

 The case of the fractional PME/FDE \nc was settled in \cite{dP12} by a complicated method that we could not extend to other critical cases. We present here a new proof that can be extended to the $p$-Laplacian after several modifications.  We assume that $m_{cs}>0$, i.e., $N\ge 2$, or $N=1$ and $1>2s$.

\begin{theorem}\label{thm.mc.fpme} Let  $m=m_c>0$.  Let $u(x,t)\ge 0$ be the semigroup solution of the FrPME in the whole space $\ren$ with initial data $u_0\in L^1(\ren)$, $u_0\ge 0$. Then for every $t>0$ the total mass is conserved
$$
\int_{\ren} u(x,t)\,dx=\int_{\ren} u_0(x)\,dx.
$$
\end{theorem}

\noindent {\sl Proof.}
(i) {\sc Reduction Step.} Before we proceed with the proof we make some reductions following the ideas of the proof of Theorem \ref{thm.mc.cple}:   we may always assume that $u_0\in L^1(\ren)\cap L^\infty(\ren)$ and $u_0$ is compactly supported. If our form of mass conservation is proved under these assumptions, then it follows for all data $u_0\in L^1(\ren)$ by the semigroup contraction property with respect to the $L^1$-norm, a property that holds for the FrPME, as well as the PME and the PLE.

(ii) {\sc Mass control.} An important tool in our proof will be the uniform control of the mass located near infinity, reflected in the following lemma.

\begin{lemma} For any two times $0<t_1<t_2<T$ before the possible extinction, the mass of $u(t)$ outside a large ball of radius $R_\ve$  is less than $\ve$ for all $t\in (t_1,t_2)$, and  $R_\ve$ depends also on $t_1$ and $t_2$ but nothing more.    \rm
\end{lemma}

\noindent $R_\ve$ depends on the solution of course. The proof uses two ingredients: the fact that the solutions are integrable in space for every $t$ and the property of weighted monotonicity in time, that comes from the B\'enilan-Crandall estimate  $(1-m)\,t\,u_t\le u$. See more details in the proof of Theorems \ref{thm.mc.cple} and  \ref{thm.mc.fple}.

\medskip

(iii) The technical part of the proof starts in the style of the mass conservation proofs in \cite{Vazquez2020, VazFPL2-2020}. We do a calculation for  the weighted mass and later we will pass to the limit in a sequence of test functions. Taking a smooth and compactly supported test function $\varphi(x)\ge0 $, the weighted mass loss for the interval $t_2>t_1>0$ is written as
\begin{equation}\label{mass.calc1.pme}
\displaystyle \int u(t_1)\varphi\,dx-\int u(t_2)\varphi\,dx=
\iint u^m(x)\mathcal (-\Delta)^s \varphi(x)\,dx=I(t_1,t_2,\var),
\end{equation}
with space integral over $\ren$  and time integral over $[t_1, t_2]$. We start with  $\varphi(x)$ being a smooth cutoff function which equals 1 for $|x|\le 2$ and zero for $|x|\ge 3$. \bc It is easy to see that in that  case $|(-\Delta)^s \varphi(x)|\le C |x|^{-(N+2s)}$; this large decay will be  important. Then we  introduce the sequence of test functions $\varphi_n(x)=\varphi(x/n)$. We will take $n\ge 2R$, $R\ge R_\ve$. Let us note that \nc
$$
\mathcal (-\Delta)^s \varphi_n(x)=n^{-2s}\,\mathcal (-\Delta)^s \varphi(x/n).
$$


(iv) Our goal is to prove that $I(t_1,t_2,\var_n)\to 0$ as $n\to \infty$. In the next paragraphs we will forget the time integral momentarily in $I(t_1,t_2,\var)$ for ease of writing. It is convenient to split the space integral in two regions for the calculations involving the mass loss.  We have an outer region $A_R=\{x: |x|\ge R\}$ and an inner region $B_R=\{x : |x|\le R\}$, both fixed at any time $t\in (t_1,t_2)$. First, we have to estimate
\begin{equation*}
 \displaystyle I_n(A_R):= \int_{A_R}  u^m(x,t)\mathcal (-\Delta)^s \varphi_n(x)\,dx, \quad I_n(B_R):= \int_{B_R}  u^m(x)\mathcal (-\Delta)^s \varphi_n(x)\,dx.
 \end{equation*}
We estimate the contribution of outer region as
$$
I_n(A_R) \le \left(\int_{A_R}  \,u(x,t)dx\right)^m
\left(\int_{A_R}\big|(-\Delta)^s \varphi_n(x)\big|^{1/(1-m)}\,dx\right)^{1-m}
$$
Now, the first integral is less than $\ve$ uniformly in $t$ by what was said above. For the second factor we recall that $1-m_c=2s/N$, so that putting $x/n=y$ we get
 $$
 \int_{A_R} \big|(-\Delta)^s\var_n(x)\big|^{N/2s}\,dx\le \int_{\ren} \big|(-\Delta)^s\var(y)\big|^{N/2s}\,dy=C,
 $$
 and this constant does not depend on $n$ or $t$. We conclude that $I_n (A_R)\le C \ve $ where $C$ does not depend on $\ve$ or $n$.

In the inner region $B_R$  the argument bears on the other factor. If
$$
I_n(B_R):= \int_{B_R}  u^m(x)\mathcal (-\Delta)^s \varphi_n(x) \,dx,
$$
we use the fact that  when $x\in B_R$ we have \ $|(-\Delta)^s \varphi_n(x)|\le C\,n^{-2s}$\
with a uniform constant if $ R\ll n$. This estimate can derived from the formula for $(-\Delta)^s \varphi_n(x)$ since $\varphi_n(x)-\varphi_n(y)=0$ for $|y|\le 2n$ if $x\in B_R$. Hence,
$$
\left|I_n(B_R)\right|\le \left(\int_{B_R}  \,u(x,t)\,dx\right)^m
\left( \int_{B_R}\big|(-\Delta)^s \varphi_n(x)\big|^{1/(1-m)}\,dx\right)^{1-m}\le
Cn^{-sp} M^m
$$
that tends also to zero as $n\to\infty$ ($M$ is the $L^1$-norm of the initial data). The conclusion is that \ $I(t_1,t_2,\var)\to 0$ as $n\to\infty$, hence
$$
\displaystyle \int_{\ren} u(x,t_1)\,dx= \int_{\ren}  u(x,t_2)\,dx.
$$
Letting $t_1\to 0$ and using the fact that $u(t)\in C([0,\infty);L^1(\ren))$,  this ends the proof for the special initial data. \qed

\medskip

\begin{corollary} The  conservation result of Theorem \ref{thm.mc.fpme} holds for signed solutions of the FrPME with $u_0\in L^1(\ren)$. The term mass should be replaced by first integral.
\end{corollary}\nc

\noindent {\sl Proof.}  In that case we begin with the reduction step and consider the signed solution $u$ with compactly supported initial data $u_0$ and the solutions $u_1$ and $u_2$ with initial data $u_1(x,0)=\max\{-u_0(x),0\}$ and $u_2(x,0)=\max\{u_0(x),0\}$ resp. By comparison we have the bounds
$$
-u_1(x,t)\le u(x,t)\le u_2(x,t).
$$
We can apply the previous analysis to $u_1$ and $u_2$ to conclude as in points (a), (b), (c) above that the absolute integral $\int_{|x|\ge R} |u(x,t)|\,dx$  is small for every time interval $[t_1\le t\le t_2]$ with $0<t_1<t_2$ if integrated outside of a large ball $B_R(0)$. So the Lemma of previous proof applies to signed solutions. The proof of point (iii) needs to replace $u^m$ by $|u|^{m-1}u$, and in point (iv) we should use estimates in absolute value. \hfill \qed
\nc

\noindent {\bf Application to the standard PME.} The above proof works also for the limit case $s=1$ where we recover the classical Porous Medium Equation with critical exponent $m_c=(N-2)/N$ if $N>2$. It represents a new, hopefully simpler proof.

\medskip

\noindent  {\bf Remark.} The proof of mass conservation is not true for the planar log-diffusion equation, $\partial_t u=\mbox{div}(u^{-1}\nabla u)$, which is the critical case for $N=2$, $s=1$, since we know that it offers a well-known phenomenon of extinction in  finite time  that we have described in Section \ref{sec.ricci}. A similar extinction phenomenon happens for $N=1$, $s=1/2$. We supply the following explicit separate-variable solution 
\begin{equation}\label{explicit.s12}
U(x,t)=\frac{2(T-t)}{1+|x|^2} \,\quad\, x\in\re, \ 0<t<T.
\end{equation}

\subsection{MC for Fractional PLE with critical exponent}\label{sec.cmass.fpl}

We have already shown that conservation of mass holds for strong semigroup solutions of the fractional PLE with  initial data in $L^1(\ren)$ in the range $p>p_c=2N/(N+s)$. This was  proved in \cite{Vazquez2020} for $p>2$ and in \cite{VazFPL2-2020} for $p_c<p<2$. It is precisely the range where finite-mass self-similar solutions exist
and the CM paradigm of Section \ref{sec.cm.diff} applies. It is not an easy proof in both cases. On the other hand, for every $1<p<p_c$ (and $N\ge 2$) we have proved in \cite{VazFPL2-2020} that there are plenty of semigroup solutions with finite extinction time,  and moreover all finite-mass initial data give rise to solutions that lose all their mass with passing time.  Thus, we are only left with the case $p=p_c$. Note that $1<p_c<2$ for all $N\ge 1$.

\begin{theorem}\label{thm.mc.fple} Let $p=p_c$.  Let $u(x,t)\ge 0$ be the semigroup solution of the Cauchy Problem for the FrPLE with initial datum $u_0\in L^1(\ren)$, $u_0\ge 0$.
Under the further assumption that $sp_c<1$, we prove that for every $t>0$ the mass is conserved
$$
\int_{\ren} u(x,t)\,dx=\int_{\ren} u_0(x)\,dx.
$$
\end{theorem}

\noindent {\sl Proof.}    The proof we give is based on the basic idea of the proof just shown for the FrPME. It is adapted after substantial changes for the FrPLE, and it also needs in principle the extra assumption $sp_c<1$ to ensure that suitable duality methods work. Such an assumption was introduced in \cite{VazFPL2-2020}. It is always true for $N=1$, or for $N\ge 2$ and $s\le 1/2$. Note that  $sp_c=N(2-p_c)$.

Here is the full detail of the proof.

(i)  {\sc Reduction.} As in the proof for the FrPME, we can make a reduction on the class of data. We may always assume that $u_0\in L^1(\ren)\cap L^\infty(\ren)$ and $u_0$ is compactly supported, say in the ball of radius $R_0$. If our form of mass conservation is proved under these assumptions, then it follows for all data $u_0\in L^1(\ren)$ by the semigroup contraction property.

\medskip

(ii) {\sc First mass control.} Under those extra assumptions, we can derive the condition of uniform small mass near infinity. We first recall the almost monotonicity in time, $(2-p)tu_t\le u$ due to scaling. Using this monotonicity and the integrability in $t_1>0$ we conclude that for any two fixed times $0<t_1<t_2<T$  the mass of $u(t)$ outside a large ball of radius $R=R_\ve$  is less than $\ve$ for all $t\in (t_1,t_2)$, and  $R_\ve$ depends also on $t_1$ and $t_2$ but nothing more.

We also recall another known fact:  if the initial support is contained in a ball, the solution is monotone in space (along outward cones, as follows from the Aleksandrov argument that we have used for instance in \cite{VazFPL3-2021}), see Appendix B. Using this and the uniform integrability, we conclude that for every constant $c>0$ there \bc is an $R_0$ such that \nc for large $r=|x|\ge 2R_0$ we have $u(x,t)\le cMr^{-N}$ for all $t>0$. Moreover,  we also know that the solution is uniformly positive before a possible extinction time  (a fact that we will exclude with our proof),  see argument in \cite{VazFPL3-2021}.

\medskip

(iii) {\sc Integrations of weighted mass.} The technical part of the proof starts as the similar proofs in \cite{Vazquez2020, VazFPL2-2020} and takes into account ideas of Theorem \ref{thm.mc.fpme} but it needs  a much more difficult technical treatment. First, we  recall that  $sp_c<N$, something that will  happen in all dimensions.
We do a calculation for  the weighted mass. Taking a smooth and compactly supported test function $\varphi(x)\ge0 $, we have for $t_2>t_1>0$:
\begin{equation}\label{mass.calc1}
\begin{array}{l}
\displaystyle \left|\int u(t_1)\varphi\,dx-u(t_2)\varphi\,dx\right|\le \iiint
\left|\frac{\Phi(u(y,t)-u(x,t))(\varphi(y)-\varphi(x)}{|x-y|^{N+sp}}\right|\,dydxdt\\[10pt]
\le \displaystyle \left(\iiint  |u(y,t)-u(x,t)|^p\,d\mu(x,y)dt\right)^{\frac{p-1}{p}}
\left(\iiint |\varphi(y)-\varphi(x)|^p\,d\mu(x,y)dt\right)^{\frac{1}{p}}\,,
\end{array}
\end{equation}
where $d\mu=|x-y|^{-N-sp}dxdy$. Space integrals are over $\ren$  and time integrals over $[t_1, t_2]$. We use the sequence of smooth test functions $\varphi_n(x)=\varphi(x/n)$ where $\varphi(x)$ is a cutoff function which equals 1 for $|x|\le 2$ and zero for $|x|\ge 3$. We take $n\ge 2R$. Then, we have to consider different regions for the efficient calculation of the multiple integrals. Note that  in the right-hand side of \nc \eqref{mass.calc1} we estimate integrals in absolute value (by taking absolute value of the integrand). In the next paragraphs we will forget the time integrals momentarily for ease of writing.

\medskip

\noindent $\bullet$ We first deal with the outer-outer region $A_R=\{(x,y): |x|,|y|\ge R\}$. We have to estimate
\begin{equation*}
 \displaystyle I(A_R):= \iint_{A_n} \frac{|\Phi(u(y,t)-u(x,t))|\,|\varphi_n(y)-\varphi_n(x)|}{|x-y|^{N+sp}}\,dydx\,dt.
 \end{equation*}
\bc We will need to use  $\mathcal M$,  the operator already introduced in Section 2 of \cite{VazFPL3-2021} by the formula
\begin{equation}\label{frplap.opM}
({\mathcal M}\var)(x):= P.V.\int_{\ren}\frac{|\var(x,t)-\var(y,t)|}{|x-y|^{N+sp}}\,dy\,.
\end{equation}
It has properties similar to $\mathcal L_{s,p}$ when $sp<1$. In using this operator the assumption $sp<1$ is needed\nc.

  We can split the absolute integral $I(A_R)$ in two and calculate
$$
 \displaystyle I(A_R;1):= \iint_{A_R} u(x,t)^{p-1} \frac{\,|\varphi_n(y)-\varphi_n(x)|}{|x-y|^{N+sp}}\,dydx\,dt
$$
$$
\le \int_{A_R}  \,u(x,t)^{p-1}\, \mathcal M(\var_n)(x)\,dx\le \big(\int_{A_R}  \,u(x,t)dx\big)^{p-1}.
\big(\int_{A_R} \mathcal M(\var_n)(x)^{1/(2-p)}\,dx\big)^{2-p}
$$
Now, the first integral is less than $\ve$ uniformly in $t$ by what was said above, while the second is bounded by
 $$
 \int_{A_R} (\mathcal M\var_n)(x)^{1/(2-p)}\,dx\le \int_{\ren} n^{-sp/(2-p)}(\mathcal M\var)(x/n)\,dx=
 \int_{\ren} (\mathcal M\var)(y)\,dy\le C.
 $$

 We proceed symmetrically with the integral $I(A_R;2)$ that uses  $ u(y,t)^{p-1}$ instead of $ u(x,t)^{p-1}$. Putting both together, we conclude that $I(A_R)\le C \ve $ where $C$ does not depend on $\ve$ or $n$.

\noindent $\bullet$  The argument in the inner region $B_n=\{(x,y): |x|,|y|\le 2n\}$ is simple. Since
$\varphi_n(x)-\varphi_n(y)=0$, we see that the contribution to the integral \eqref{mass.calc1} is zero.
Recall that we are using $n\ge 2R$.

 \noindent (iv) {\sc Integrations in cross regions}. We still have to make the analysis in  other regions in order to cover the whole domain $x,y\in \ren$. An option is to consider the two cross regions
$C_n=\{(x,y): |x|\ge 2n ,|y|\le R\}$ and $D_n=\{(x,y): |x|\le R ,|y|\ge 2n\}$. Both are similar so we will  only work out the contribution  in $D_n$. The idea is that we have an extra estimate: $|x-y|>n$ that avoids the singularity in the weight of the integrand. Forgetting again the time integrals   for the moment, we have $\var(x)=1$, $\var(y)=0$, $u(y,t)\le cMr^{-N}$ so that
$$
I(D_n)\le C(p)\,(I_1(D_n)+ I_2(D_n)),
$$
where
 \begin{equation*}
 \begin{array}{c}
  \displaystyle I_1(D_n)\le\iint_{D_n} |u(x,t)|^{p-1}\,d\mu(x,y)\le \\
\displaystyle  \big(\int_{B_R} \,|u(x,t)|^{p-1} dx
  \big(\int_{|x-y|>n} |x-y|^{-N-sp} dy\big)   \big)\\
 \displaystyle  \le Cn^{-sp} \int_{B_R} |u(x,t)|^{p-1}\,dx\,.
 \end{array}
\end{equation*}
Since $0<p-1<1$, we have
$$
\int_{B_R} |u(x,t)|^{p-1}dx\le R^{N(2-p)}\left(\int_{B_R} |u(x,t)|\,dx\right)^{p-1}
\le R^{N(2-p)}\|u(x,t_1)\|^{p-1}.
$$
Therefore, $ I_1(D_n)\le CR^{N(2-p)} n^{-sp} $ which tends to zero as $n\to\infty$ with a power rate.
As for $I_2(D_n)$, using the known space decay of the solution as $|y|\to \infty$ we have:
$$
 \displaystyle I_2(D_n)\le\iint_{D_n} |u(y,t)|^{p-1}\,d\mu(x,y)\le
  (cM)^{p-1} R^N\int_{|y|\ge 2n} |y|^{-N-sp}y^{-N(p-1)}\,dy
  $$
 where we have integrated first the $x$ variable. Then $I_2(D_n)\le Cn^{-sp-N(p-1)}R^N\,,$
which tends to zero as $n\to\infty$.

Same analysis for $I(C_n)$. Note that these regions overlap but that poses no problem. This concludes the estimate of the right-hand side of \eqref{mass.calc1}\nc.

\medskip

(v) Going back to the formula \eqref{mass.calc1}, we conclude in the limit $n\to\infty$ that the mass is conserved for the whole time interval
$$
\int_\ren u(x,t)\,dx= \mbox{constant} \qquad \mbox{ for all } \ t_1\le t\le t_2.
$$
Since $t_1<t_2$ are arbitrary we conclude that mass is conserved,
$$
\int_\ren u(x,t)\,dx= \int_\ren u_0(x)\,dx \qquad \mbox{ for all } \  t\ge  0.
$$
In particular, there will be no extinction in finite time.
\qed


\medskip

\noindent {\bf Conservation of integral for signed solutions.} The above proof implies the conservation in time of the integral $\int u(x,t)\,dx$ for signed solutions of the FrPLE, with the same arguments done there.

\medskip

\noindent {\bf Remark for the standard PLE.} The above proof works also for the limit case $s=1$ where we recover the classical $p$-Laplacian with critical exponent $p_c=2N/(N+1)$,  \cite{Vazquez2006}.

\begin{open} The problem of proving mass conservation for the solutions of the fractional PLE, treated in this section, is still open  if the restriction $sp<1$ does not hold. One might wonder if the proof available for $p>p_c$ can be adapted by some kind of continuity. We have been unable to do it, the uniform control at infinity used there (due to available super-solutions) is missing and a new idea is needed. It could use the delicate local energy inequalities of Section \ref{sec.crit.exp}, that must be worked out.
\end{open}


\section{Loss of mass in fractional nonlinear flows }\label{sec.loss.vf.fpme}

We have addressed in Section \ref{sec.break} the phenomenon of loss of mass for the very fast diffusion ranges of the PME and the PLE. A very similar situation happens for the fractional models that we have under scrutiny at this moment. In the fractional PME we will consider the analogous lower range $0< m< (N-2s)/N$; in the fractional PLE in the range $1<p<p_c=2N/(N+s)$, cf. the references  \cite{VazqJEMS2014} and  \cite{VazFPL2-2020} resp.
It was proved that there is no integrable Barenblatt solution in those ranges, moreover  conservation of mass does not hold and finite mass goes to zero in infinite or finite time. An important problem is to investigate the detail of the evolution of that waning mass, and what is the shape that those declining solutions mostly take.

In that respect, there are still many open gaps. Thus, the analysis of the existence of self-similar solutions with possibly anomalous exponents has not been done (to our knowledge).  The general considerations about the extreme character  of the diffusion for large and small intensities of $u$ made in the local flows can still be applied here. The equivalents of Propositions \ref{thmBR.ext} and \ref{prop.massto0} are still true.

We will add next relevant information about the issue of mass decay and extinction in time complementing what is known from the mentioned references. Our task is constructing suitable  self-similar solutions that will show the typical behaviour to be expected. Our contribution consists of two kinds of special solutions that take explicit or semi-explicit forms.

\subsection{Very singular solutions in fractional very fast diffusion }

\noindent {\bf FrPME.} It is proved in \cite{VazqJEMS2014} that in the above-critical range $m_c=(N-2s)/N<m<m_1=N/(N+2s)$ there exist singular  VSS solutions, with the expected scaling of the form
$$
V(x)=C\,t^{1/(1-m)}\,|x|^{-2s/(1-m)}\,,
$$
where the constant $C$ is calculated in a very explicit way as a function by $N,m,s$, see full discussion in Section 13 of \cite{VazqJEMS2014}. In this case the VSS is the limit of the standard finite-mass Barenblatt solutions. On the other hand, in the very fast range $0<m<(N-2s)/N$ there is a \bc singular function with a similar form\nc. It was calculated in \cite{VazVol2015}, section 6.1, and has the form
$$
V(x)=C_1\,(T-t)^{1/(1-m)}\,|x|^{-2s/(1-m)}\,,
$$
where $C_1>0$ is calculated much as $C$ before. Now the VSS ceases to be integrable at infinite as a space distribution, and it also undergoes extinction in finite time.  The extinction rate is $O((T-t)^{1/(1-m)})$, valid for every $m<m_c$, but the appropriate norm is  none of the $L^p$ spaces, but in the $L^{p,\infty}(\ren)$ (Marcinkiewicz space) with $p=N(1-m)/2s$.

\medskip

\noindent {\bf FrPLE.} The problem of finite time extinction for the FrPLE was studied in \cite{VazFPL2-2020}. The VSS for $1<p<p_c=2N/(N+1)$ takes the form
$$
V(x)=C_1\,(T-t)^{1/(2-p)}\,|x|^{-sp/(2-p)}\,,
$$
cf. Theorem 16.2 of \cite{VazFPL2-2020}.


\subsection{Separate-variables solutions in fractional nonlinear equations}
 We begin the study with the Fractional PME. We may use the method used Section \ref{ss.vfpme}, when studying the geometrical flow, to obtain a separate-variables solution of the form
$$
U(x,t)=(T-t)^{\alpha}F(x),
$$
that implies that $\alpha=1/(1-m)$ and $(-\Delta)^s F^m =\alpha F$. Assume that $N\ge2$ or $N>2s$ if $N=1$. Writing $G=F^m$ we get the fractional elliptic semilinear problem
$$
(-\Delta)^s G =\alpha G^q, \qquad q=1/m>1.
$$
This equation was studied in  \cite{CotTav04} in the wider study of the best constants for fractional Sobolev inequalities, see also
\cite{ChenLi1991, Lieb}. More precisely, it is  proved that the elliptic equation has a positive, finite-mass solution if and only $q+1$ is the fractional Sobolev exponent, $2N/(N-2s)$, hence $q=(N+2s)/(N-2s)$. Moreover,  the explicit solution is unique up to translations and dilations and is given by
$$
G(x)=c\left( \frac{a}{a^2 +|x|^2}\right)^{\frac{N-2s}{2}}.
$$
Summing up, we get for the precise exponent $m_s=(N-2s)/(N+2s)$ the smooth finite-mass solution
$$
U(x,t)=\kappa\,(T-t)^{\frac{N+2s}{4s}}\left( \frac{a}{a^2 +|x|^2}\right)^{\frac{N+2s}{2}},
$$
where $\kappa(N,s)>0$. This solution exhibits the extinction rate $U(\cdot,t)=O((T-t)^{\frac{N+2s}{4s}})$ in all $L^p$ spaces.

\medskip

\noindent {\bf Separate-variables solutions for the general case.} We try separate-variables solution of the form
$$
U(x,t)=(T-t)^{\alpha}F(x),
$$
on the fractional DNLE equation: $\partial_t u + (-\Delta)_p^s u^m=0$. We were originally interested in the case $m=1$ but the algebraic calculations are not more difficult and the relation with the Calculus of variations more clear. We assume $N>2s$ and the fast diffusion condition $m(p-1)<1$.
As usual, the Ansatz of separation of variables implies first that $\alpha=1/(1-m(p-1))>0$, and then
$$
(-\Delta)_p^s (F^m) =\alpha F \quad \mbox{or better} \quad (-\Delta)_p^s G =\alpha G^{q},
$$
where $q=1/m$. Existence of a radially symmetric and positive solution of this fractional elliptic solution is proved  by Brasco et al. in \nc \cite{Brasco2016} in the form of extremal functions for the corresponding fractional critical Sobolev inequality, and this analysis leads to the necessary condition $q+1=Np/(N-sp)$. This in turn means
$$
m(p-1) + (m+1)\frac{sp}{N}=1.
$$
When this exponent condition is satisfied a solution is obtained but it is not explicit (to our knowledge).
The authors conjectured a possible explicit form that has not been found. They prove some convenient properties, in particular the expected decay at infinity. Theorem 1.1 of \cite{Brasco2016} says:  {\sl Let $G$ be any minimizer for the fractional elliptic equation in a suitable energy space. Then $G$ is bounded, has constant sign, it is a radially symmetric and monotone function with
$$
\lim_{|x|\to \infty} |x|^{\frac{N-s p}{p-1}} G(x) = K.
$$
for some constant $K>0$.} In the fractional PLE case, we put $m=1$  and the separation formula holds with $2=Np/(N-sp)$, hence $p_s=2N/(N+2s)$ and $F=G$. On the other hand, we get time decay with exponent $\alpha=1/(2-p_s)=(N+2s)/4s$.

\vskip 1.5cm

\newpage

{\bf \huge Part 3}

\section{ On some nonlinear hyperbolic conservation laws}\label{ssec.ncl}

The study of first-order equations of the type called Conservation Laws interested Riemann around 1860 \cite{Riemann} as part of the study of gas dynamics (more precisely, the propagation of sound waves). He discovered and examined in that framework the curious phenomenon of shock formation, by which we mean the fact that the physical process tends to produce solutions that exhibit a finite jump at certain moving locations. In this way he pushed forward the mathematical topic of ``differential equations with discontinuous solutions'', mainly undeveloped by then. To be precise, we refer to evolution processes such that, even when the data are smooth and bounded, the physical solution becomes essentially discontinuous after some time.

Since the equations of gas dynamics form a system  characterized by its mathematical difficulty\nc, see \cite{Sm83},
much work  concentrated on scalar conservation laws of the possibly nonlinear form
\begin{equation}\label{cleq.nd}
\partial_t u + \sum_i \partial_{x_i} F_i(u)=0,
\end{equation}
where the $F_i(x,t,u)$ are real-valued functions defined for $x  \in \ren$, $t>0$ and $u\in \re$,
often under the restriction $u\ge0$. When the functions $F_i$ are smooth, these first-order equations can be solved by the method of characteristics at least for small times, and the occurrence of discontinuities is explained as the existence of points at which the characteristics that propagate different values of $u$ come to  intersect. This was studied in works like  \cite{Lax73} and \cite{Dafermos73}. In 1970 \cite{Kruzhkov70} Kruzhkov identified the correct concept of solution to obtain existence and uniqueness of a physical acceptable generalized solution  for scalar laws in $N$ space dimensions \nc with possible jump discontinuities. This is called entropy solution or Kruzhkov solution. The acceptable initial datum in his seminal paper was bounded, $u_0\in L^\infty(\ren)$.

The multi-dimensional law \eqref{cleq.nd} was treated by Crandall \cite{Crandall72} via the theory  of nonlinear semigroups \nc acting on convex subsets of Banach spaces. Here the appropriate  space is $L^1(\ren)$ where the semigroup acts as a contraction, see also B\'enilan's thesis, \cite{BeTh72}. The  Kruzhkov concept of solution has to be extended to cover unbounded data, which is a delicate step. In other words, at this time it was realized  that the mathematical theory of nonlinear conservation laws shared many properties in common with other evolution equations generating nonlinear semigroups of contractions in $L^1$, like the heat equation, the PME and the PLE, i.e., the ones we have studied before.

Following the lead of the early works by B\'enilan  and Crandall, we are interested in finding  connections of our previous investigation on mass conservation and self-similarity with the theory of hyperbolic conservation laws, and we will find it in the case one-dimensional scalar conservation laws of the form \eqref{cleq.1d} with power-like nonlinearities $\varphi$. Since our objectives in this section are limited we will only consider nonnegative solutions. There are many difficult topics of the general theory that we do not touch at all and many topics in the proposed  connection to be explored.

\medskip


\subsection{One-dimensional scalar conservation laws}

We explore next the theory of the scalar conservation
\begin{equation}\label{cleq.1d}
 u_t + \varphi(u)_x = 0
\end{equation}
where the mathematical study has been undertaken and is essentially complete, cf. \cite{Bressan00, Dafermos16, Serre99}. We will focus our attention in some strong analogies with our study of nonlinear diffusion models done in Section \ref{sec.masscon}, plus the case of failure of mass conservation that we will find as a limit case.

According to existing theory we may assume that  $\varphi$ is a continuous real function with $\varphi(0)=0$ and $\varphi(\re)=\re$. The assumption that $\varphi$ is monotone nondecreasing is not necessary but we will make it here, as well as $u\ge0$, for the sake of getting very precise conclusions in a rather moderate space.

There is a standard  nonlinear semigroup theory that allows to prove existence and uniqueness of bounded entropy solutions that enjoy three basic properties, valid for all $u_0, v_0 \in L^\infty(\re)$:

(i)   Comparison:   $u_0\le v_0$  a.e. implies $u(\cdot,t)\le v(\cdot,t)$ for a.e. $x$ and all $t>0$.

(ii)  Contraction: $u_0 - v_0\in L^1(\re) $  implies that $u(\cdot,t), v(\cdot,t)\in L^1(\re)$ and
$$
{ \|u(x,t)- v(x,t)\|_1\le \|u_0(x)- v_0(x)\|_1. \nc}
$$

(iii) Conservation (of mass): $u_0- v_0\in L^1(\re) $  implies that $u(\cdot,t), v(\cdot,t)\in L^1(\re)$ and
$$
{ \int_{\re} (u(x,t)- v(x,t))\,dx=\int_{\re} (u_0(x)- v_0(x))\,dx\,.\nc}
$$
These results are found in \cite{Serre04} in greater generality.  The reader will  see in our presentation below \nc
the extent of the analogy with the results for the PME proved in Section \ref{ssec.emc}.

\medskip

\noindent {$\bullet$ \bf Two propagation speeds.}  An interesting consideration deals with the concept of speed of propagation associated to an equation of the form \eqref{cleq.1d}. It is well-known that away from the shocks the equation can be integrated by the method of characteristics
$$
\frac{dx}{ \varphi'(u)}=dt,
$$
along which the value of $u(x,t)$ is kept constant. This gives the notion of characteristic propagation speed, ${\bf v}_c=\varphi'(u)$. On the other hand, \eqref{cleq.1d} can be seen as a mass conservation law of the form \eqref{cont.eq}, $ u_t + ({\bf v}\,u)_x= 0$, and then we must identify the mass speed as ${\bf v}_m=\varphi(u)/u$. This second version will hold at the shocks we will find below, and is known as the Rankine-Hugoniot condition.

We will focus in particular on conservation of mass and existence of fundamental solutions for all powers $\varphi(u)=u^m$, $m>0$.
The conclusion fails in the limit $m\to 0$ where the fundamental solution degenerates and we will find extinction in finite time with precise mass rates. We will study in detail this limit case that reminds us of the analogous cases in nonlinear diffusion.


\subsection{Convex nonlinearities. Slow transport. Shocks}

$\bullet$ { \bf Burgers equation. } As a starting point, we recall that the simplest of these laws is the Burgers equation,
\begin{equation}\label{eq.burgers}
\partial_t u + 2u \partial_x u = 0,
\end{equation}
corresponding to the choice  $\varphi(u)=u^2$. If $u$ has two signs, there are two options for the nonlinearity, $\varphi(u)=u^2$ and $\varphi(u)=u|u|$, with different theories. Here, we will keep the assumption $u\ge0$. The theory of this equation was studied by Burgers as a model for turbulence
in 1948 \cite{Burgers48}. He used the approximation of the above mentioned equation by a viscous version
\begin{equation}\label{eq.burgers.vis}
u_t+2u \,u_x=\mu\, u_{xx}
\end{equation}
by just adding a viscosity term $\mu\, u_{xx}$ with small viscosity parameter $\mu>0$. The resulting parabolic equation  can be solved by a linearization trick based on integration  and  Hopf-Cole transformation, \cite{Hopf1950}. This is a well-known computation in the literature, First we introduce $v(x,t)=\int_{-\infty}^x u(y,t)\,dy$ that solves
$$
v_t+ v_x^2=\mu\, v_{xx},
$$
and has  a multiple of the Heaviside function as initial datum: $v(x,0)=0$ for $x<0$;  $v(x,0)=M$ for $x>0$. Then,  we put $z=e^v$ and we see that $z(x,t)$   solves the heat equation with initial a step function from $z=1$ for $x<0$ to $z=e^M$ for $x>0$. Solving this equation explicitly we obtain an expression for $u$ with $\mu>0$\nc.
Passing to the inviscid limit $\mu\to 0$, one obtains the fundamental solution (with a Dirac delta as initial datum)
$U(x,t) $ with the explicit expression of the original Burgers' equation:
$$
U(x,t)=\frac{x}{2t} \qquad \mbox{ for } \ 0\le x\le  2 t^{1/2},
$$
and $U=0$ otherwise. $U(\cdot,t)$ has triangular shape with a vertical side on the right\footnote{To be compared with all the diffusive shapes we have seen before.}. See a detailed presentation in \cite{Lax73}\nc.

Looking for more detail, we see that the jump discontinuity  is located at the ``shock line'' $X(t)=2\,t^{1/2}$ with a height $h(t)=u(X(t)-)-0=2\, t^{-1/2}$. This line plays the role of the free boundaries in nonlinear diffusion flows  in the sense that it separates the positivity region where the equation holds in the classical sense from the empty region where $U=0$. Note that now the acceptable solution (the entropic solution) is discontinuous along the shock line\nc.

\medskip

\begin{figure}[ht!]
\centerline{\includegraphics[width=0.45\textwidth]{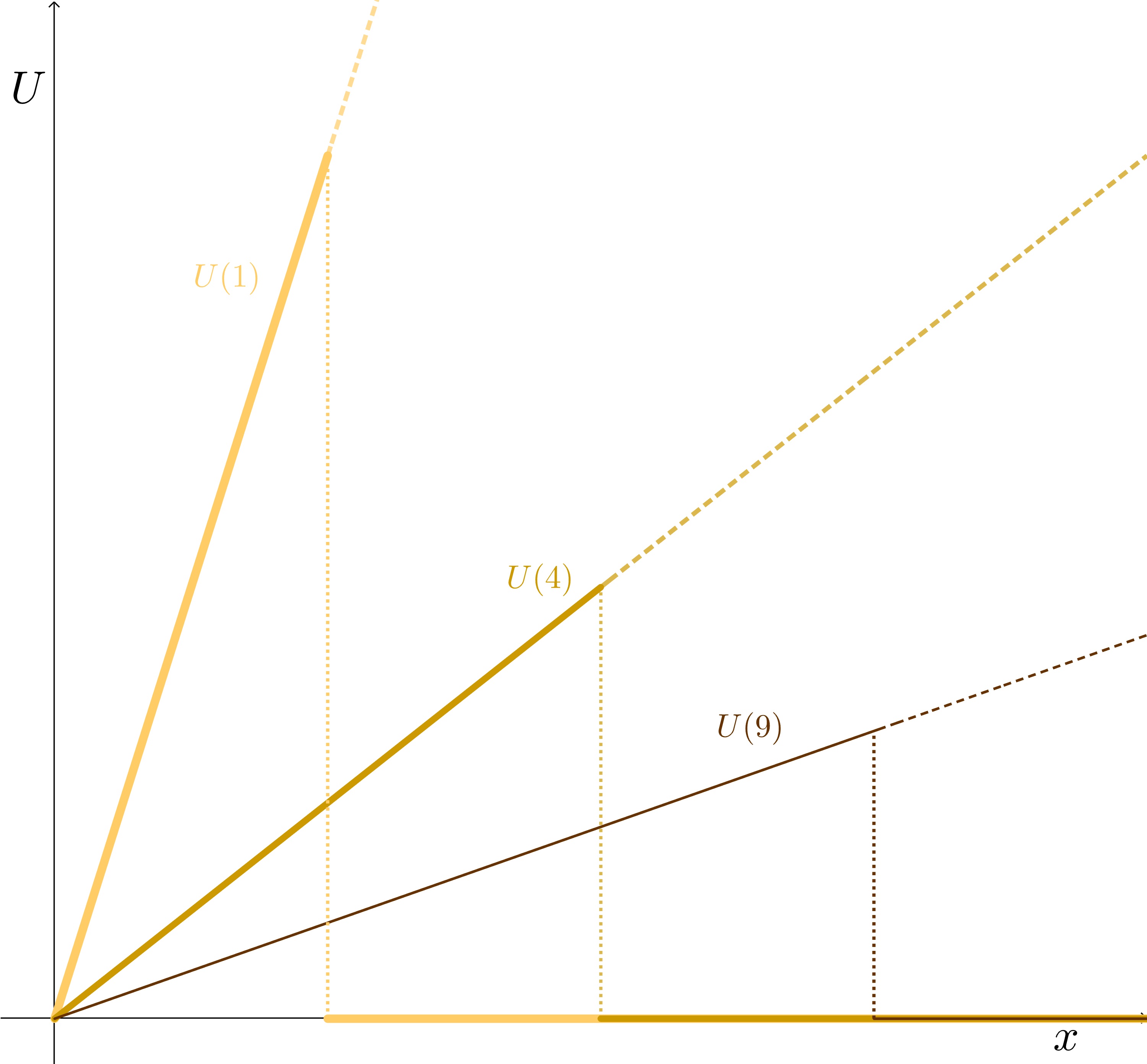} \qquad
\includegraphics[width=0.45\textwidth]{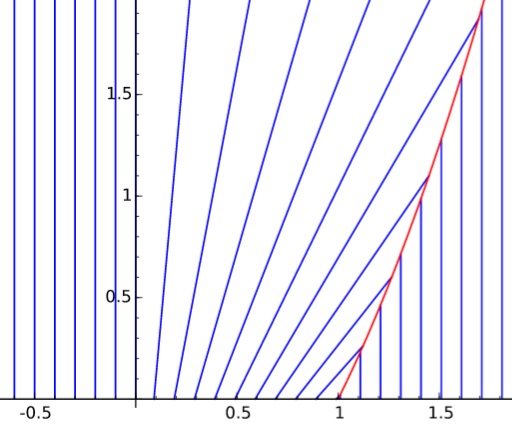}}
\vspace{-0.5cm}
\caption{Burgers equation. To the right: Solution at different times.  For $t=0$ it is a Dirac delta\nc. To the left: Graphics with  characteristics and  shock line}
\label{fig:Burgers}
\end{figure}

$\bullet$ { \bf Convex powers. } The above computation can be generalized to the study of nonnegative solutions of the conservation law \eqref{cleq.1d} with a convex power nonlinearity, $\varphi(u)=u^m$, $m>1$. Using the approach of Subsection \eqref{ssec.pme} we look for self-similar solutions of the form
$$
U(x,t)=t^{-\alpha}F(x\,t^{-\beta})\,.
$$
After substitution into the  conservation law, algebraic calculation, and identification of the powers of time that appear, and use of the mass conservation rule, we arrive at  the exact values of the self-similar exponents: $\alpha=\beta=1/m$. Then, the profile $F$ satisfies
$$
(F^m)'=\frac1{m} (yF)'\,.
$$
Integration and the assumption of finite mass leads to \ $F^{m-1}=y/m$ whenever $ F\ne 0$.  This in turn leads to the explicit solution with formula
\begin{equation}\label{cleq.Nw}
U^{m-1}(x,t)=\frac{x}{mt} \qquad \mbox{ for } \ 0\le x\le X_m(t)= C t^{1/m},
\end{equation}
and $U=0$ otherwise. Now $U^{m-1}$ has a triangular shape. This self-similar function is called an \em $N$-wave\rm.
 Here, $C>0$ is a free constant that reflects the total mass $M$  in the form $M=a(m) \,C^{\,m/(m-1)}$. Then, we can write
\begin{equation}\label{cl.shock}
X_m(t)=\xi(m)\, (M^{(m-1)}t)^{1/m}, \qquad \xi(m)= m\,(m-1)^{-(m-1)/m}.
\end{equation}
The reader will easily check that the smooth part of the solution is composed of straight characteristic lines with speed ${\bf v}_c=mu^{m-1}$, while the shock speed is
${\bf v}_m=X'(t)=u^m/u=u^{m-1}$, which is slower than the incoming characteristics from the left. This is the geometric requirement for a shock to be admissible in an entropy solution.

The fact that $U(x,t;M)$ takes as initial data the  $M$ times the Dirac delta  is immediate.
There is a theory of entropy solutions with initial data a nonnegative finite measure and the uniqueness of such solutions comes from the work of Liu and Pierre \cite{LiuPierre84}, and it implies in particular the uniqueness of the (entropy) fundamental solution.  It is also immediate by the comparison principle that bounded solutions with compactly supported data produce solutions with compact support in space for all $t>0$, and this leads easily to the law of mass conservation.  Furthermore, they also prove  that the solution \eqref{cleq.Nw}  provides a description of the  asymptotic behavior of the solutions initiated by an arbitrary nonnegative integrable function.

\medskip

 $\bullet$  {\bf Remark on $N$-waves.} We borrow the name from the form of the entropy solutions of the Riemann problem where the initial data are constant for $x<0$ (left state) and $x>0$ (right state). Then the entropy solution has the self-similar form $U(x,t)=F(x/t),$ and the profile $F$ has a piecewise continuous form between two shock lines that is called $N$-wave. In our case the restriction $u_0$ implies that $F$ has the form of half the typical Riemann $N$-wave. \nc
\medskip

$\bullet$ { \bf Remark on the change of variables. } If $u$ is a solution of the conservation law $u_t+(u^m)_x=0 $ with $m>1$, it is easy to see that the change of variables $w=(m/2)u^{m-1}$ produces a solution of the Burgers' equation, the one with $m=2$, as long as we are justified in applying the chain rule. Surprisingly, this fact is not justified and has important consequences.

By direct inspection of the formulas we see that this is exactly true for the fundamental solutions on the left of the shock, where the solutions are smooth. However, the equivalence if far from complete, since the argument does not include what happens at the shock. Actually, while she shock moves as $X(t)=c\, t^{1/m}$ for the fundamental solution $U$ with exponent $m\ge 2$, it is $X(t)=c \,t^{1/2}$ for Burgers's equation. This may look like a contradiction, a mere change of variables changes the location of the shocks\nc.

Actually, the location of the shock is determined by the law of mass conservation which applies to different functions: in the first case we conserve $\int_\re U(x,t)\,dx$, in the Burgers case  we conserve the mass of $w$, i.e., $\int_\re U^{m-1}(x,t)\,dx$.

\medskip

$\bullet$ { \bf The pointwise inequalities. } In order to stress the analogy with the diffusive models of section \ref{sec.cm.conn},  we point out the entropy solutions form a semigroup of $L^1$ contractions with initial data in $X=L^1(\re)_+$, Besides, we have the Bénilan-Crandall estimates \cite{BC81b}
$$
(m-1)t \,u_t\ge -u, \quad (u^{m-1})_x\le \frac1{mt}\,.
$$
Using them, is easy to prove the $L^1$-$L^\infty$ smoothing effect in the  form
\begin{equation}\label{scl.sm.eff}
\sup_x u(x,t)\le \sup_x U(x,t;M)= U(X(t),t)= c(m)\,M^{1/m}\,t^{-1/m},
\end{equation}
with  $c(m)=(m-1)^{-1/m}$. This formula holds with $\le$ inequality for all entropy solutions $u$ with finite mass $M$. The optimal constant is achieved by the $N$-wave solution with the value that we have just indicated. Moreover, in this case it can be easily shown that equality with the optimal constant is achieved only if $u$ is the $N$-wave or a translate thereof. Compare with Proposition \ref{prop.2.1} for the PME.



\subsection{Concave nonlinearities. Rarefaction waves with fat tails }

 We now consider equation $u_t+(u^m)_x=0 $  with  $0<m<1$, the so-called concave case. We still work with nonnegative solutions. Note that in this case  $\varphi(u)=u^m$ is not smooth at $u=0$, and $\varphi'(0)=\infty$.
The self-similarity analysis can be applied as before, but now the characteristic lines of the self-similar solution that we construct tend to diverge so that we lose the shock discontinuity on the right-hand side (i.e., for large $x$), while it appears on the left-hand side. The self-similar solution $U=U(x,t;m)$ is
\begin{equation}\label{cleq.Nw2}
U^{1-m}(x,t;m)=\frac{mt}{x} \qquad \mbox{ for } \  X_m(t)= C t^{1/m} \le x<\infty,
\end{equation}
and $U=0$ on the left of the shock line  $ X_m(t)$. The smooth solution on the right-hand side of the  shock line forms the so-called \sl rarefaction wave \rm and is composed of characteristic lines. The conditions for admissible shocks are met. We also get $M=a(m) \,C^{-m/(1-m)}$
and

\begin{equation}\label{shock.mlss1}
X_m(t;M)=\xi(m)\,M^{-(1-m)/m}t^{1/m}, \qquad \xi(m)=m(1-m)^{(1-m)/m}.
\end{equation}
These formulas are algebraically very similar to the convex case, but the geometric interpretion is quite different.
Summing up, while the convex case offers solutions with compact support, with leading shocks and finite propagation, the concave case offers solutions with trailing shocks, moreover they exhibit a positive rarefaction wave ahead of the shock, hence infinite propagation of perturbations from the zero level and a definite fat tail at infinite of the order $u(x,t)=O(x^{-1/(1-m)})$ as $x\to +\infty$.

\medskip

\noindent $\bullet$ The entropy solutions still form a semigroup of $L^1$ contractions with initial data in $X=L^1(\re)_+$. Besides, we have the Bénilan-Crandall estimates in the opposite
$$
(1-m)t \,u_t\le u, \quad (u^{m-1})_x\le \frac1{mt}\,.
$$
Note that $m<1$ so that the latter inequality is a lower bound for $u_x$. The $L^1$-$L^\infty$ smoothing effect still holds in the same form \eqref{cl.shock}. Only difference is that the constant becomes $c(m)=(1-m)^{-1/m}$. We point out that this constant  goes from $c(0)=e$ to $c(1)=\infty$ when $m$ goes from 0 to 1. Note the curious fact that $c(m)$ has a singularity only for $m=1$, i.e., for the equation of linear transport, that we will examine below.

\medskip

\begin{figure}[ht!]
\centerline{\includegraphics[width=0.45\textwidth]{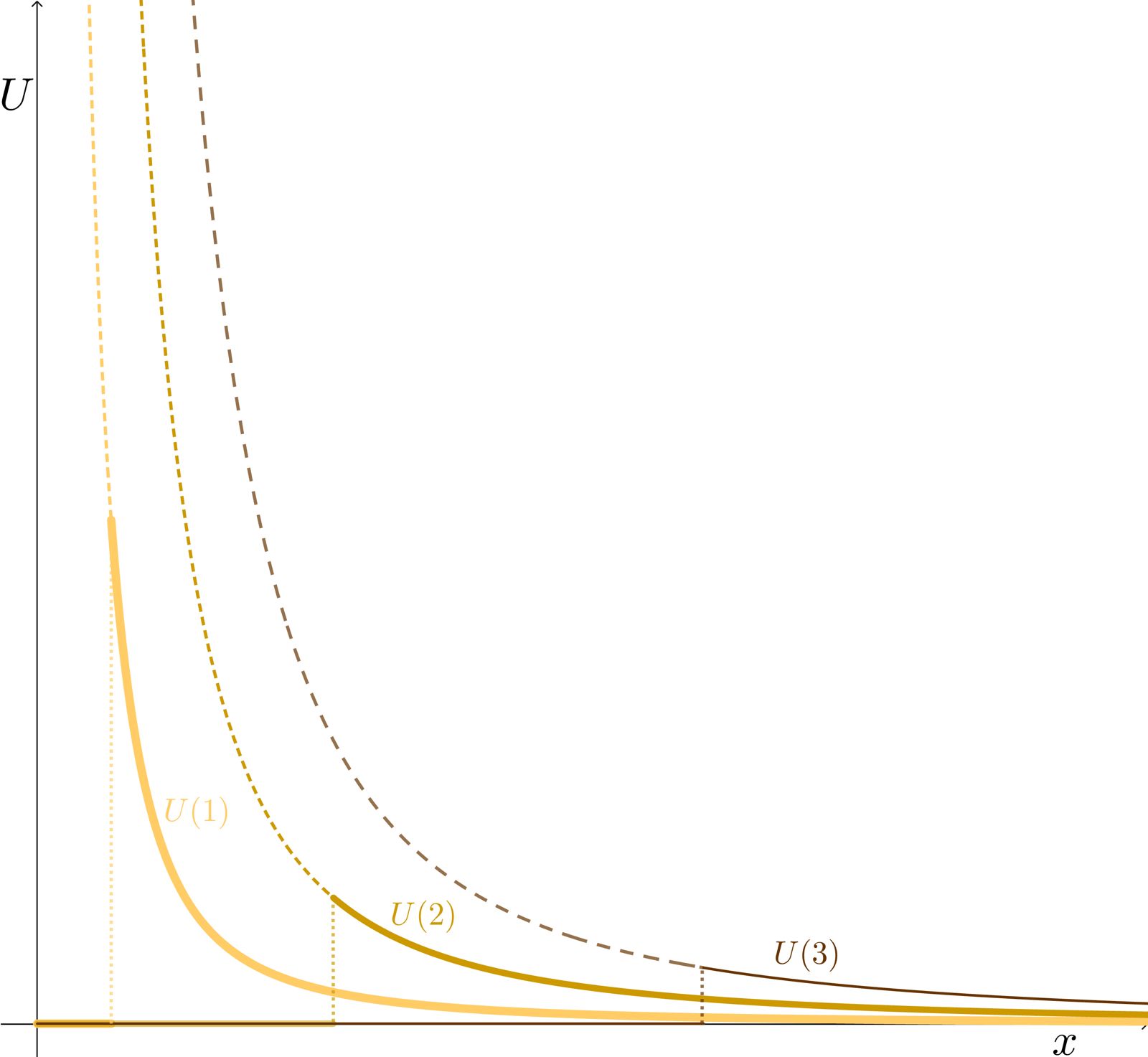} \quad
\includegraphics[width=0.45\textwidth]{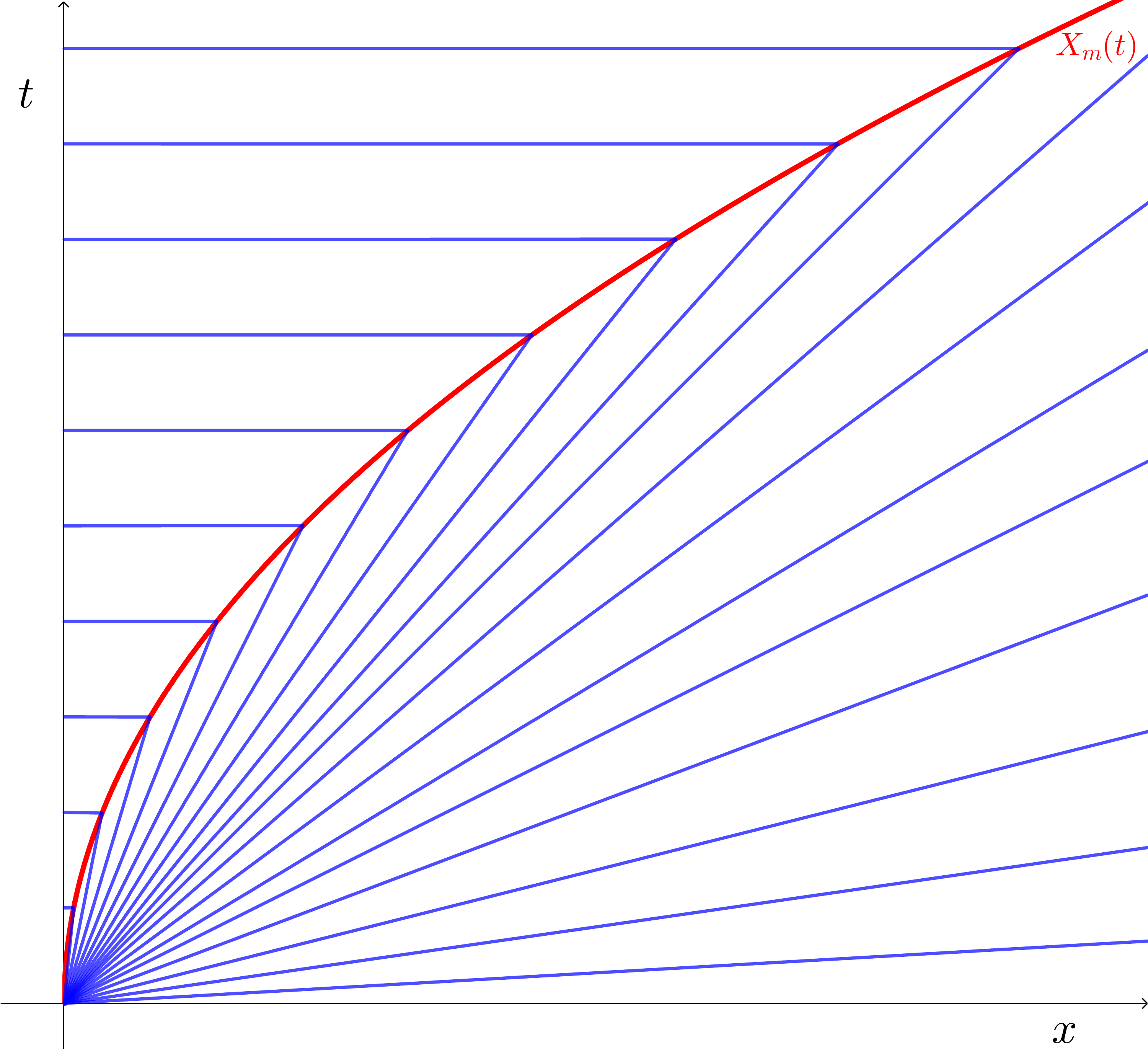}}
\vspace{-0.0cm}
\caption{Solution at different times; shock and characteristics, for  $m=0.5$}
\label{fig:Shock.m}
\end{figure}

\smallskip

The justification of the validity of  solution \eqref{cleq.Nw2} and its role in the description of the  asymptotic behavior of the solutions initiated by an arbitrary nonnegative integrable function can be found in the literature, cf. \cite{BeTh72}. To begin with, Kruzhkov theory is  extended with Crandall’s $m$-accretive semigroup approach. Conservation of mass holds for all $0<m<1$. We refrain from giving full details for fear of excessively extending this survey\nc. Again, the similarity with the fast diffusion case of the PME and PLE is clear, with minor  differences.

\medskip

$\bullet$ { \bf Remark on the change of variables. } \nc In this case we can still use the transformation $w\sim u^{m-1}$ to obtain an equivalence of the smooth parts of the solutions, but this idea does not produce more useful information. Instead, we may resort to a new idea, i.e., exchanging the roles of $x$ and $t$ as follows. Put $x'=t$, $t'=x$ and $w(x',t')=v^m (x,t)$. Then $w$ satisfies the equation
\begin{equation}
\partial_{t'}w + \partial_{x'}w^{m'}=0, \quad \mbox{with} \ m'=1/m,
\end{equation}
which is an equation of the convex type, well studied before. The validity of this transformation can be justified for entropy solutions defined in the positive quadrant $\{x>0, t>0\}$ with zero data on the axes. Actually, we can easily check that the shock locations are correctly predicted. It has to be noted that $U$ conserves in time the mass $\int_{ X(t)}^\infty U(x,t)\,dx$, while $w$ conserves in space the integral
$$
\mbox{Mass}(W,t')=\int_0^{X'(t')} W(x',t')\,dx'=\int_0^{T(x)} U^m(x,t)\,dt.
$$

\noindent  {\bf Remark.} The reader may check \nc by direct inspection of formulas \eqref{cl.shock} and \eqref{shock.mlss1}
that the exchange of the $x$ and $t$ variables also exchanges  the shock locations exactly, in the sense that not only the exponents but also the leading factors agree, and finally the value of the respective masses is the same.
\medskip

\noindent {\bf Viscous conservation laws.} We may consider different diffusions as regularization terms, starting with the Laplacian viscosity approximation term, $\mu \Delta u$, natural in fluid mechanics, to porous medium and $p$-versions, and we may use fractional Laplacian operators and other integro-differential versions. Examples are abundant in the literature, see for instance \cite{Alib2014}. When we use models with viscosity terms the limit of zero viscosity is usually called the {\sl inviscid limit}. Proving that the obtained limit function is admissible in the sense of entropy solutions is a basic job.

\medskip

\noindent {\bf Remark. } As an interesting connection, in \cite{GPV00} the authors study the large-time behavior of the solutions of the initial-value problem for the nonlinear diffusion equation $ u_t=\Delta (u^m)$ posed in the whole space $\ren$  in dimensions $N\ge 3$ with nonnegative initial data $u(x,0)\in L^1(\ren)$  when the exponent takes on the critical value $m=(N-2)/N<1$. When the solutions are radial a clever change of variables allows to pass to a viscous version of the conservation law \eqref{cleq.1d} with the same exponent.

\medskip


\subsection{Comparison with linear transport}

The simpler linear transport equation: \ $u_t+cu_x=0$  is very illustrative. This equation can be integrated along characteristics that are straight parallel lines: $X(x_0,t)=x_0+ct$, along which the value of $u$ is kept constant so that
$$
u(x,t)=u_0(x-ct)\,.
$$
Note that there is no intersection of characteristic lines so that no shock occurs, and no new singularity appears if it is not contained in the initial data. Note that the construction of the fundamental solution by any reasonable approximation method gives the measure-valued solution
$$
U(x,t)=M\,\delta(x-ct)\,.
$$
Note that the limit of the fundamental solutions obtained in the convex and concave cases gives precisely this result as $m\to 1$. In particular we may check that $\xi(m)\to 1$ as $m\to 1$,
and the coefficient $c(m)$ in the $L^1$-$L^\infty$ smoothing result tends to infinity.

\noindent {\bf Nonlinear interaction.} Summing up, the interaction of different wave speeds for the nonlinear conservation laws produces a dispersion that regularizes the initial point mass, something that does not occur in the linear transport equation that produces parallel transport with no interaction, hence no regularization.

\subsection{ The question of very singular solutions}

This question was introduced in Section \ref{sec.cm.conn} in the study of the PME/FDE equation where they existed only for $m<1$. If we pass to the limit in the formulas for fundamental solutions studied in this section we see that the $M\to \infty$ limit of formula \eqref{cleq.Nw} for $m>1$ gives the unbounded classical solution
\begin{equation}\label{cleq.Nw.Minf}
U_\infty^{m-1}(x,t)=\frac{x}{mt} \qquad \mbox{ for } \ 0\le x<\infty, \ 0<t<\infty
\end{equation}
which is a  very simple rarefaction fan. On the other hand, for $m<1$ we pass to the limit in formula \eqref{cleq.Nw2} and get
\begin{equation}\label{cleq.Nw2.Minf}
U_\infty^{1-m}(x,t)=\frac{mt}{x} \qquad \mbox{ for } \ \mbox{ for } \ 0\le x<\infty, \ 0<t<\infty
\end{equation}
which is again a rarefaction fan in the first quadrant. The actual finite mass solutions are portions of these surfaces $u=u(x,t)$ bordered by suitable shocks.


\subsection{Connection with Hamilton-Jacobi equations. Mass analysis}

This is an old connection that implies a better understanding of the two related equations. We work in 1D.
When we integrate in $x$ a solution $u(x,t)$ of the scalar conservation law $u_t+\varphi(u)_x=0$, we obtain a function
\begin{equation}\label{fn.hj}
v(x,t)=\int_a^x u(y,t)\,dy
\end{equation}
that is called the mass function or distribution function of $u(x,t)$, and accumulates the mass contained between $y=a$ and $y=x$ at time $t$. The base point $a$ in the definition can be chosen at convenience, normally is taken as $-\infty$, 0, or $\infty$. With any choice of $a$, function $v$
solves, at least formally, the associated equation
\begin{equation}\label{eq.hj}
v_t+ \varphi(v_x)=0,
\end{equation}
that can be classified as a Hamilton-Jacobi equation. This connection does not exist in several variables.
Similarly, in the viscous case we start from the equation  $u_t+\varphi(u)_x=\mu \, u_{xx}$, and end up with
\begin{equation}\label{eq.hj.visc}
v_t+ \varphi(v_x)=\mu \, v_{xx}.
\end{equation}
Note: It is common in the literature to find the equation in the form \bc $v_t+\varphi(v_x)=\mu \, v_{xx}$ \nc as the result of using the alternative definition of $v$ as
$$
v(x,t)=\int_x^{\infty} u(y,t)\,dy\,,
$$
that amounts to taking $a=\infty$ and changing the sign of $v$. So the versions are equivalent.
In the case of the Burgers equation we get the well-known Hamilton-Jacobi equation
$$
v_t+ (v_x)^2=\mu \, v_{xx},
$$
that can be integrated by the Hopf-Cole transformation.  This made possible the first analytical progress in this area.

\medskip

\noindent $\bullet$ {\bf Theory.}
While the theory of scalar conservation laws is well-posed in the context of entropy solutions for finite-mass,  in the integrated equation \eqref{eq.hj} we are looking for a theory that applies to bounded solutions, possibly increasing, with definite values at infinity. This is supplied by the concept of viscosity solutions of  Hamilton-Jacobi equations that can be studied in books like  \cite{Barles94, Lions82}. There is an equivalence between finite-mass entropy  solutions $u$ of equation \eqref{cleq.1d} and viscosity solutions $v$ of \eqref{eq.hj} obtained by integration if $\varphi$ is regular, this result was folklore and Caselles gave a clear proof in \cite{Case92}.

\medskip

\noindent $\bullet$ {\bf Self-similar solutions for the HJ equation.} If we perform the corresponding integrations in equation \eqref{cleq.1d}   when $\varphi(u)=u^m$, we obtain for $m>1$ the explicit self-similar solution $V=V(x,t;m,M)$ given as
\begin{equation}\label{eq.hj.mge1}
V=k(m)\frac{x^{m/(m-1)}}{t^{1/(m-1)}}, \quad k(m)=\frac{|m-1|}{m^{m/(1-m)}} \quad \mbox{ for } \ 0\le x \le X_m(t),
\end{equation}
while $V=0$ for $x\le 0$ and $V=M$ for $x\ge X_m(t)=\xi(m)\, (M^{(m-1)}t)^{1/m})$. $V(X(t),t)=M$  means that $k(m)=\xi(m)^{-m/(m-1)}$.  On the other hand, for $0<m<1$ we get
\begin{equation}\label{eq.hj.mle1}
V=M- k(m)\frac{t^{1/(1-m)}}{x^{m/(1-m)}}  \quad \mbox{ for } \ X_m(t)\le x <\infty,
\end{equation}
while $V=0$ for all $x\le X_m(t)$. Precisely, the coefficient $k(m)$ is computed so that $V(X_m(t),t)=0$.

In both cases, the distribution function $V$ increases from 0 to $M$ as expected.


\subsection{The critical limit $m=0$.  Mass analysis}\label{ssec.HJ}
In the study of fast diffusion equations we have found the limit case $m=0$ for the PME/FDE in dimension $N=2$, and the limit case $p=1$ for the fast PLE in dimension $N=1$, both with special limits for the fundamental solutions. We recall the existence of standard fundamental solutions was closely tied to mass conservation. In fact, in those limit cases mass loss did occur.

We will make use of the integrated distribution  in the study of the limit $m\to 0$ of the family of solutions we have discussed for equation \begin{equation}\label{cleq.1d.m}
 u_t + (u^m)_x = 0.
\end{equation}
with a concave nonlinearity. Then we also have a delicate limit (in 1D) as $m\to 0$. First of all, direct inspection of the fundamental solutions shows that
$$
\lim_{m\to 0} U(x,t;m,M)=0,
$$
and this limit is uniform for $t\ge T$ and $|x|\ge c>0$, with $T,c$ arbitrary positive constants. Therefore, the mass distributions $U(\cdot,t;m,M) $ disappear from any bounded set as $m\to0$. There are two alternatives: either the mass $M$  concentrates at $x=0$ or escapes to infinity. We will show next that both things happen at the same time for a certain time until the whole mass is depleted.

\begin{theorem}\label{thm.claw.mloss}  As $m\to 0 $ we have the following limit
\begin{equation}\label{cleq.lim.m0}
\lim_{m\to 0} U(x,t;m,M)=(M-t)_+\,\delta(x)
\end{equation}
in the sense of measures.  The limit  ``solution'' undergoes total extinction at the time $T=M$.
\end{theorem}

 \noindent  {\sl Proof.}  The reason for extinction in finite time in the limit $m=0$ is that mass escapes to infinity at a rate of 1 unit of mass per unit of time as we will show.

The proof is based on a precise analysis of the amount of mass contained in exterior sets when $m>0$. The idea may have a general use in similar models but the computation is not difficult in the case of concave conservation laws since the amount of mass contained in different regions is exactly controlled by the distribution function $V(x,t)$  defined by formula \eqref{fn.hj} applied to $U$, integrating from $a=-\infty$\nc. Thus, for every $t>0$, the support of the solution is contained in the region
$(X(t),\infty)$ and we have in the limit
$$
\lim_{m\to 0} X(t)= 0 \quad \mbox{if} \ t/M< 1\,
$$
$$
\lim_{m\to 0} X(t)= \infty \quad \mbox{if} \ t/M> 1 \,
$$
If you look at the distribution function, which is a monotone function in $x$, we conclude that $V(x,t)=0$ vanishes identically for $t> M$.

Now, we look at the intermediate levels. The amount of mass contained in an appropriate exterior region $(R(t),\infty)$ is  $M'<M$, and the mass inside the left-hand side region is $M-M'<M$ iff
$$
V(R(t),t;m,M)=M-M',
$$
i.e.,
$$
R^{-m/(1-m)}= \frac{M'}{k(m)t^{1/(1-m)}}=\frac{M'}{M}\,X(t)^{-m/(1-m)},
$$
$$
R(t)= \xi_m \,(t/(M')^{1-m})^{1/m}\,.
$$
We see  that $R(t)$ behaves like $X(t)$ with $M'$ instead of $M$. Looking at the convergence of the distribution functions $V(x,t; M)$, we see that for $t\ge M'$ an amount  $M-M'$ has escaped to infinity, while the  rest
is confined in a small neighbourhood of $x=0$. Arguing for all $0<M'<M$ the proof is done. \qquad \qed

\begin{theorem}\label{thm.claw.mloss.gen} Let us consider initial data $u_0\in L^1(\re) $  supported in a finite interval $[a,b]$\nc. It we solve the scalar conservation law for $m>0$, pass to the limit $m\to 0 $ and  obtain a limit solution $u(x,t) $,  the following mass loss holds:
\begin{equation}\label{cleq.lim.m0.mass}
\int_\re u(x,t) =(M-t)_+
\end{equation}
and the support of $u(x,t)$ does not leave $[a,b]$ with time\nc.  The limit  solution undergoes total extinction at the time $T=M$.
\end{theorem}

\noindent {\sl Proof. } (i) If we integrate the  equations $u_t +(u^m)_x=0$  with $m>0$ we obtain solutions of the corresponding HJ equation $v_t=(v_{x})^m$ with nonlinearity $\varphi(s)=(s)^m$. This equation satisfies the maximum principle and so does the limit $m\to 0$ that corresponds to nonlinearity $\varphi(s)=\mbox{sign}(s)$. Now we apply it to a $v$, solution of \eqref{eq.tvf} with mass $M$ and support in the interval $I=[a,b]$, and compare it to the fundamental solutions $u_1(x,t)= M\,\delta(x-a)$ and $u_2(x,t)= M\,\delta(x-b)$. If we  integrate all of them, find the corresponding $v$ solutions we find the inequality
$$
v_1(x,t)\ge v(x,t) \ge v_2(x,t)\,
$$
because this inequality is certainly true at $t=0$. But both solutions $v_1$ and $v_2$ are explicit step functions:
$$
v_1(x,t)=0 \mbox{ for } x< a, \quad v_1(x,t)=(M-t)_+ \mbox{ for } x\ge a,
$$
and similarly $v_2$ with jump at $x=b$. We conclude that for all $0<t<M/2$ the function $u(x,t)$ is supported in $[a,b]$ and $v(\infty,t)=M-t$, for $t<M$. This proves \eqref{cleq.lim.m0.mass}. \qed


 Recall that for any given $t>0$ the total mass of $U(x,t;m,M)$ is $M$ whenever $m>0$.
 (i) Let us fix $M=1$ for simplicity. We have just seen that the mass in any region  not near 0 nor infinity goes uniformly to 0 as $m\to 0$, We recall that $\xi(m)\sim m/e$ when $m$ is close to 0 , hence $X(t)\sim mt^{1/m}/e$. This means that $X(t)\to 0$ as $m\to 0$ if $t\ge 1$, and to $\infty$ if $t>1$.

 (ii) To examine the possible escape of some mass to infinity we need some  precise calculations.
Thus,  for every fixed $R>0$ we have for all $m<1$ small
$$
I_R(t)=\int_{|x|\ge R} U(x,t)\,dx= m^{1/(1-m)}t^{1/(1-m)}\int_{R_1}^{\infty} x^{-1/(1-m)}\,dx.
$$
where  $R_1=\max\{R, X(t)\} $. Hence,
$$
I_R(t)= (1-m)\,(m^m t)^{1/(1-m)} R_1^{-m/(1-m)}.
$$
When $0<t\le 1$ we have $R\ge X(t)\sim mt^{1/m}/e$ for all small $m$, and then we get
$$
I_R(t)= (1-m)\,(m^m t)^{1/(1-m)} R^{-m/(1-m)}
$$
and this integral tends to $t$ as $m\to 0$. So in the limit a fraction of the mass has escaped to infinity. On the other hand, when $R\le X(t)$
$$
I_R(t)= I_{X(t)}=M=1\,,
$$
hence all the mass has escaped to infinity. We obtain the expected result when $M=1$.

(ii) For $M\ne 1$ we use a time scaling the equation. If $u(x,t)$ is a solution with unit initial mass $M>0$, then for every $M>0$
$$
\widehat u(x,t)=M \, u(x, M^{m-1}t)
$$
is a solution with initial mass $M$. The claim follows. \qquad \qed

\medskip

\noindent{\bf Examples of solutions involving step functions}

\noindent {\bf Example E1.} We first examine the case of a simple step function as initial data for the scalar conservation law with $1>m>0$. To be precise, we take initial data a  function of the form
$$
u_0(x)=1 \qquad  \mbox{for } -M\le x\le 0,
$$
and zero otherwise. By using the technique of characteristics and admissible shocks, we easily check that for small times the solution consists of narrower step function with same height 1 in the interval $-M+t\le x\le mt$ plus a standard rarefaction fan on the right-hand side. Prove that at time $T=M/(1-m)$ this solution becomes identical with the fundamental solution we described in  \eqref{cleq.Nw2}. This is a remarkable case of complete coincidence of two different solutions in finite time that is typical of nonlinear conservation laws.

\medskip

\noindent {\bf Example E2.} We now pass to the limit  $m\to 0 $. We consider as initial data the step function
$$
u_0(x)=1 \quad \mbox{for } \ x\in (-a,0),  \ u_0(x)=0  \quad \mbox{otherwise. }
$$
We may construct the solution starting with these initial data for the equation with $m>0$ and pass to the limit $m\to 0$. Applying the method of characteristics the top value $u=1$ propagates with speed $m\,u^{m-1}\to 0$. An admissible shock arises on the left with speed $f(u)/u=u^{m-1}=1$. In the limit we get the ``solution''
$$
u(x,t)=1 \quad \mbox{for } \ -a+t\le x\le 0,
$$
if $0<t<T=a$ and $u_0(x)=0 $ otherwise. Thus, we have $M(t)=a-t$ and we get confirmation of the mass loss in the case of this particular solution. In physical terms, the evolution can be explained as  follows: there is shock on the left-hand side of the profile that eats into the remaining step function at a rate of one unit mass per unit time until the whole mass is depleted.
By a simple rescaling of the equation we can obtain the solution with the same mass law when initial data a step function different space support
$$
\widehat u(x,t)= C \quad \mbox{for } \ x\in (-a+t/C ,0], \ \mbox{when } t<aC,
$$
and $\widehat u(x,t)=0$ otherwise. Here $C$ and $a$ are free positive parameters. Then we have $ \widehat M_0=aC$ and $\widehat M(t)=\widehat M_0-t$. Any space translation also qualifies as solution, with the same properties.

\medskip

\noindent {\bf Remarks.} This result has to be compared with the Total Variation Flow in 1D considered in Subsection \ref{ssec.behgs}, another typical of borderline case in a theory. There the mass loss was also found to be constant, but the rate was $2$ instead of 1. There, the evolution of a step function was considered and a different explicit solution found, formula \eqref{form.step}. See also \cite{BonFig2012}, mentioned in Section \ref{sec.ncbl}. We recommend the reader to compare both evolutions, which modify the initial mass distribution in a different form: indeed, the rate of mass loss  here is $1$, while it is 2 in the 1D TVF.

\smallskip

\begin{open}
 We also propose to solve the problem with initial data composed of two step functions. Will they propagate as two step functions of heights as in the original display? This produces more difficult and interesting examples\nc.
 \end{open}


\subsection{Limit $m=0$. The log version }

There is another version of the equation obtained after a simple rescaling to allows to find a different limit. Thus, we introduce into the equation an innocent coefficient $a>0$ and write the equation as
$$
u_t+ a\, (u^m)_x=0.
$$
All the previous theory applies after replacing the old time, let us call it $t'$, by $at$.
The choice that we need here is $a=1/m$. Note that this means that the process acts faster than before, in fact very fast as $m\to 0$. We note that whenever $u(x,t)$ is a solution of the original equation, then
$$
\widehat u(x,t)= u(x,t/m)
$$
is a solution of $ \widehat u_t+ \widehat u^{m-1} \widehat u_x=0. $
The fundamental solution is then
$$
\widehat U^{1-m}(x,t;M)=U^{1-m}(x,t/m;M)=\frac{t}{x}
\quad \mbox{ for } \  \widehat X_m(t)= X_m(t/m) \le x<\infty.
$$
In the limit $m\to 0$ we should find a kind of generalized solution of the  equation $\widehat u_t+ \widehat u^{-1}\widehat u_x=0$, i.e.,
$$
\widehat u_t+ (\log(\widehat u))_x=0.
$$
Regarding the mass extinction, we have increasingly faster movements, so the conclusion of mass escape to infinity is still valid, but now extinction happens for all positive times.

\subsection{Nonlocal models. Conservation laws with long-range effects}
More recently, there arose an interest in nonlocal versions of the conservation laws we have been studying, they come from  different applications.

\noindent $\bullet$ In one example, Biler et al. \cite{BilerGM} study the  initial value problem for the nonlinear and
nonlocal equation
$$
\partial_t u + |u_x|\,((-\partial_{xx})^s u)=0, \qquad x\in \re, \ t>0,
$$
with  $0<s<1$. It arises as a 1D model in Head's theory of dislocations, the dislocation  density being represented by $u_x$. Upon differentiation we get a 1D version of the  fractional PME model of Subsection \ref{PME.nonlocalpress}. The exponent in this case is $s=1/2$.

\medskip

\noindent $\bullet$ On the other hand, equations of the form
$$
\partial_t u + \partial_x\left(f(u)\, V(u\ast K)\right)=0
$$
are proposed in \cite{Betan11}  as models for sedimentation
of small equal-sized spherical solid particles dispersed in a viscous fluid, where the local solids
volume fraction $u = u(x, t)$ is sought as a function of depth $x$ and time $t$. Here, $f$ is a suitable function depending on the application and $K$ is a smooth convolution kernel. The nonlinearity $V$ is called a hindered settling factor.  The authors show that the \nc theory of (Kruzhkov) entropy solutions  can be applied.

\medskip


\noindent{\bf Final general remarks.} A large literature on this topic studies the above questions for equation \eqref{cleq.1d} with quite general nonlinearities $\varphi$. It must remarked that the concave case has not been studied in the same detail as the convex case.

We will not enter into the more general study of systems of conservation laws in several space variables, that is a world in itself. Relevant information can be found in the abundant literature, see  for instance \cite{Bressan00, Bressan23, Dafermos16, Glimm70, GodRav91, Holden2015, LeFloch02, Liu74, Liu76, Serre99, Zheng2001}.

\vskip 1.5cm

\newpage

{\bf \huge Part 4}

\section{ Some related diffusion  problems}\label{sec.comments1}

In the sequel we present a list of selected topics that have a close relationship with the main contents of the survey and hopefully offer useful additional information. The choice reflects personal  expertise and interests of the author. Before proceeding further, let us recall that our work deals with problems in the whole space, while work on similar problems in bounded domains is also very extensive and takes a very different direction, in particular in what regards mass conservation. There will be no special reference to bounded domains in what follows. We will also avoid discussing here equations of the type \ $\partial_t u=A(u)+f$, \ where the forcing term $f$ has an important influence, but see Section \ref{sec.models}.

In what follows next, we describe selected models where the concepts, tools and results of the preceding sections have a direct application, and some that depart in a more significant way. A number of other models are briefly commented at the end of the section. They are less central to our presentation.


\subsection{Fractional PME with nonlocal pressure}\label{PME.nonlocalpress}

It is a nonlocal and nonlinear diffusive equation that uses the Darcy law plus a nonlocal pressure given by a Riesz potential to arrive at the formulation
\begin{equation}\label{CV.eq}
\partial_t u=\mbox{div}(u \nabla p), \qquad   p=(-\Delta)^{-s}u\nc,
\end{equation}
with $0<s<1$. \bc Apart from its mathematical interest the problem had a number of physically interesting motivation.  Indeed, equations of the more general form $u_t=\nabla \cdot (\sigma(u)\nabla {\cal L}u)$ have appeared recently in a number of applications in particle physics. Thus, Giacomin and Lebowitz consider in   \cite{GL1} a lattice gas with general short-range interactions
and a Kac potential $ J_\gamma(r)$ of range $\gamma^{-1}$, $\gamma>0$ . Scaling
spacelike with $\gamma^{-1}$ and timelike with $\gamma^{-2}$, and passing
to the limit $\gamma\to 0$, the macroscopic density profile
$\rho(r,t)$ satisfies the equation
\begin{equation}
\qquad \frac{\partial \rho}{\partial t}=\nabla   \cdot\left[\sigma_s(\rho)\nabla
\frac{\delta F(\rho)}{\delta \rho }\right]
\end{equation}
Here  $F(\rho)=\int f_s(\rho(r))\,dr- (1/2)\iint J(r-r')\rho(r)\rho(r')\,drdr'$,
where $f_s(\rho)$ is the (strictly convex) free energy density of the
reference system, and $\sigma_s(\rho)$ is the mobility of the system with only short-range
interactions. The use of the model in the study of phase segregation  is treated in the book \cite{GLM2000}.

On the other hand, the equation with $s=1/2$ has been proposed by Head \cite{Head72} as the equation of motion of a dislocation continuum, and then $u$ is the dislocation density and the space dimension is $N=1$. The mathematical investigation of this case was performed by Biler et al. in \cite{BilerGM}, though it is done in terms of the integrated equation $v_t+|v_x|\Lambda(v)=0$, where $\Lambda$ is the L\'evy operator of order 1, which is equivalent to $(-\partial^2_{xx})^{1/2}$. \nc

\medskip

\noindent $\bullet$ {\bf The Cauchy problem.} In our context initial data are given $u_0\in L^1(\ren)$, $u_0\ge 0$. This Cauchy problem was studied  by Caffarelli et the author \cite{CaffVa2011a}, and previously in 1D by Biler et al. in \cite{BilerGM}.
 In this case there is existence and uniqueness of some type of viscosity solution,  obtained after integrating the equation in $x$, yielding an equation of Hamilton-Jacobi type\nc. Mass is conserved and the fundamental solution is explicit, quite similar to the classical Barenblatt profiles that we have met in Section \ref{sec.cm.conn}. We have
$$
B_s(x,t)=t^{-\alpha}F_s(x\,t^{-\beta})
$$
with \ $\beta=(N+2-2s)^{-1}, \ \alpha=N\beta $, and
$$
\quad F_s(y;s)=(C-k|y|^2)_+^{1-s}
$$
The constant $C>0$ allows to determine the mass of the solution. The explicit form in 1D was first found in \cite{BilerGM}\nc. The fundamental solution is also explicit  for $N>1$, and conserves mass, all the campactly supported solution also do. The uniqueness of the self-similar solution has been proved in \cite{Due2016} (in the sense of weak-strong uniqueness, see Corollary 4.2) and in \cite{Serfaty2020} where it appears in the context of the mean field limit for Coulomb-type flows.

\medskip

\bc \noindent $\bullet$ {\bf Boundedness.}
 We state here the result about passage from initial $L^1$ to $L^\infty$ for positive times that is proved in \cite{CaffSorVaz}.

\begin{proposition} Let $u$ be a weak solution of Cauchy Problem with $u_0 \in L^1(\ren)$
and u0 decreases exponentially as $|x|\to \infty$. Then there exists a positive
constant $C_1$ such that for every $t > 0$
$$
\sup_{x\in\ren} |u(x, t)| \le C_1 \,t^{-\alpha}\|u_0\|_{L^1(\ren)}^{\gamma}
$$
with $\alpha = N/(N + 2 - 2s)$, $\gamma = (2- 2s)/(N + 2 - 2s)$. The constant $C$ depends only on $N, s$.
\end{proposition}
Furthermore, H\"older continuity of the solutions is proved in \cite{CaffSorVaz, CaffVaz15}.

\medskip

\noindent $\bullet$ {\bf  Asymptotic behaviour.} The behaviour of finite-mass solutions for large time  has been established in \cite{CaffVa2011a} by entropy methods, see also \cite{CarrHSV}.

\begin{proposition} Let $u(x,t)\ge  0$ be a weak solution of the Cauchy Problem with bounded
and integrable initial data such that $u_0\ge 0$ has finite entropy. Then, as $t\to\infty$ we have
$L^1(\ren)$ convergence
$$
\|u(x,t)-B_s(x,t;C)\|_1 \to  0,
$$
and the uniform convergence
$$
t^\alpha \|u(x,t)-B_s(x,t;C)\|_{L^\infty(\ren)} \to 0.
$$
Recall that $\alpha=N/(N+2-2s)^{-1}$. The constant $C$ is determined by the rule of mass equality:
$\int u(x,t)\,dx=\int U_c(x,t)\,dx$. \nc
\end{proposition}

\begin{open}
However, the problem still lacks a full uniqueness and comparison theorem in several dimensions, so the mathematical knowledge is only partial. \end{open}

\bc \noindent $\bullet$ {\bf More general equations.}  There were extensions of these results  involving general power nonlinearities in equation \eqref{CV.eq}, like
$$
\partial_t u=\mbox{div}(u \nabla (-\Delta)^{-s}u^{m-1}),
$$
studied in \cite{BilerIK}, or
$$\partial_t u=\mbox{div}(u^a \nabla (-\Delta)^{-s}u),
$$
studied in  \cite{StdTVA2018, StdTVA2019} \nc with $a>0$ and $0<s<1$. For a survey of such models of fractional diffusion see \cite{VazAbel10}. Conservation of mass for these models is proved in \cite{StdTVA2019}, see Theorem 5.2.\nc

\subsubsection{The non-diffusive limit}
When we  let the fractional exponent $s\to 1$ in the equation \eqref{CV.eq} we get a ``mean field'' equation arising in superconductivity and superfluidity,
$$
\partial_t u=\mbox{div}(u \,\nabla N(u)),
$$
 where $N(u)$ denotes the Newtonian potential. This is a non-diffusive model that has vortices as typical solutions and has been studied by a number of authors, like \cite{AmSerf08, E1994, ChRS96, LinZh00, SerfVa2014}. Conservation of mass holds. Recently this vortex model has been studied for different mobilities $m(u)$ in \cite{CarrGCV22, CarrGCV22b}, in the form \
$$
\partial_t u=\mbox{div}(m(u) \nabla (-\Delta)^{-1}u).
$$
\bc Conservation of mass and existence of fundamental solutions has been studied and settled for such models\nc.

\newpage

\subsection{Diffusion on manifolds}\label{ssec.manif}

The main operator in the diffusion theory, the Laplacian, can be defined on any many-dimensional Riemannian manifold $(M,g)$, where $g=(g_{ij})$ denotes the Riemannian metric.  It is called the \sl Laplace–Beltrami operator\rm, and acts on a smooth scalar function $f=f(x)$ as follows:
\begin{equation}
\Delta_g(f)=\frac1{\sqrt{|g|}}\sum_{i,j=1}^N  \partial_i\left(\sqrt{|g|}\,g^{ij}\,\partial_j f\right),
\end{equation}
where $|g|= det\, g$. Recall that the volume element in the manifold is
$$
dv=\sqrt{det\, g}\, dx_1\cdots dx_N.
$$
The heat equation is then written as
\begin{equation}
\partial_t u =\Delta_g u
\end{equation}
for functions $u(x,t)$. The study of diffusion on manifolds is strongly tied to the study of the solutions generated by this equation. It is a classical topic in differential geometry \cite{Gri2009, RosenbBk97, Yau78}.

The basic example of manifold is the Euclidean space $\ren$ where the action of previous sections was played.
The heat flow on manifolds may behave quite differently from the Euclidean heat flow and this depends on the geometrical properties of the manifold.
Thus, conservation of mass is a common occurrence but it may fail in some strange manifolds that are called stochastically incomplete.
Stochastic incompleteness of a Riemannian manifold $M$ amounts to non-conservation of total mass for the heat semigroup on $M$.
Nonlinear characterizations of stochastic completeness have been recently studied in \cite{Grillo20}.

The most typical examples of curved manifolds are the sphere \ ${\SP^{N}}$ and the hyperbolic space $\Hn$. Since the sphere is a compact manifold, conservation of mass for the heat flow is immediate, so we will not insist in this example.

\medskip

\noindent{\bf Flows in the hyperbolic space.} We may parameterize in different ways the hyperbolic manifold $\Hn$ in any dimension $N\ge 2$. One way that is convenient for analysis  is to take as point  space for the hyperbolic space the Euclidean set $\ren$, take  a pole point $o\in \ren$, and take the hyperbolic metric given by the line element
$$
(ds)^2= (dr)^2+ ({\Sh \,r})^{2} (d\theta)^2
$$
where $x=(r,\theta) $ is the expression in polar coordinates. Thus, $r$ will be the geodesic distance from $x$ to $o$, and $\theta\in {\SP}^{n-1}$ is the angular coordinate. In this way $\Hn$ happens to be a homogeneous Cartan-Hadamard manifold.

\noindent $\bullet$ Then the Laplace-Beltrami operator reads
$$
\Delta_g(f)= ({\Sh}\,r)^{1-N}\, \partial_r\left(({\Sh}\,r)^{ N-1}\,\partial_r f\right) + ({\Sh}\,r)^{-2} \Delta_{\theta} f.
$$
where $\Sh \, r$ denotes the hyperbolic sine function, and $\Delta_{\theta} $ is the Laplace-Beltrami operator on $\SP^{N-1}$.

The heat kernel on the hyperbolic space was studied in detail by Grigorian in \cite{Gri98}. Though this kernel, or fundamental solution, has an increased expansion rate with respect to the typical expansion rate of the Euclidean kernel, which is called \sl ballistic effect\rm,  it does conserve mass, and so do all finite-mass solutions. The ballistic effect and the fine asymptotic convergence of general finite-mass solutions towards the fundamental solutions is studied in our recent work \cite{VazHyp22}.

\noindent $\bullet$ The porous medium equation on the hyperbolic space  $\Hn$ is
\begin{equation}
\partial_t u =\Delta_g (u^m),
\end{equation}
with $m>1$. It was first studied by us in \cite{Vaz2015}. On the other hand, the fast diffusion  case $0<m<1$ was studied by Grillo and Muratori in \cite{GrilloM14}. This case is interesting since we already know that in the Euclidean setting there is a critical exponent $m_c<1$  below which the mass is gradually lost, see Section \ref{sec.break}. Paper \cite{GrilloM14} deals with nonnegative solutions of the fast diffusion equation posed in $\Hn$ under two conditions:  \\
(i) $u$ is radial, i.e., $u(r,t)$, and \\
(ii) $m$ lies in the exponent range $(m_s,1)$ where $m_s=(N-2)/(N+2))$.\\
As the authors say, the situation on negatively curved manifolds is different from the Euclidean one since vanishing of solutions in finite time often occurs not only for $m$ close to 0, but for all $m < 1$, cf. \cite{BonGrVaz2008}. In fact, this is  somewhat similar to what happens in the case of the homogeneous Dirichlet
problem on bounded Euclidean subdomains of $\ren$. We conclude that conservation of mass happens only for $m\ge 1$.

\noindent $\bullet$ A number of the results proved for flows on $\Hn$ can be generalized to more general manifolds with negative curvature, especially in the Cartan-Hadamard class, see our work in \cite{GrilloMV17, GrilloMV19}.

There is also more recent work on fractional heat equations on $\Hn$. The fractional porous medium equation on the hyperbolic space was studied in \cite{Berchio20}. For solutions with large data see \cite{GrilloMP23}.

The is a large literature on nonlinear equations of elliptic and parabolic type on  such manifolds, that we will not refer here to avoid distractions from our main topics.

\newpage


\subsection{Anisotropic diffusion} \label{ssec.10.2}

A number of equations of PME type or PLE type have been studied in anisotropic or inhomogeneous settings by choosing exponents that depend on the space variable $ x $ and sometimes also on $t$.
We will study here two model cases with pure power nonlinearities.

\subsubsection{Anisotropic porous medium diffusion}
There are models of nonlinear evolution that combine a mechanism of nonlinear diffusion with the occurrence of strong anisotropy.  The combination of nonlinear diffusion of Porous Medium  type with strong anisotropy was recently considered by Feo, Volzone and the author  in \cite{FVV23}
by studying the {anisotropic porous medium equation}
\begin{equation}\label{APM}
u_t=\sum_{i=1}^N(u^{m_i})_{x_i x_i}, \quad \tag{APME}
\end{equation}
with different exponents $m_i>0$ for $i=1,...,N$,  $N\geq2$. The problem is posed for $t>0$ and $x\in\ren$. The study was done in the fast diffusion range: $0<m_i<1$ for all $i$.
We consider solutions to the Cauchy problem to \eqref{APM} with nonnegative initial data
$ u(x,0)=u_0(x),$ where $u_0\in L^1(\mathbb{R}^N)$, $u_0\ge 0$. We put $M:=\int_{\mathbb{R}^N} u_0(x)\,dx$, so-called total mass.
The analysis of the problem shows marked differences between having linear exponents $m_i=1$, fast exponents $m_i<1$, or slow exponents $m_1>1$, so that we have decide to put the algebra of self-similar solutions first. The work follows precedents by \cite{S01, S01bis, JS06, JS05, H11} and other authors.


\noindent{\bf Formal analysis of self-similarity.} The self-similar solutions that we look for have the curious anisotropic form
\begin{equation}\label{sss.an}
U(x,t)=t^{-\alpha}F(t^{-a_1}x_1,..,t^{-a_N}x_N)
\end{equation}
with constants $\alpha>0$, $a_1,..,a_n\ge 0$ to be suitably chosen below. We substitute this formula into equation \eqref{APM}.
After writing $y=(y_i)$ and $y_i=x_i \,t^{-a_i}$, we see that time is eliminated as a factor in the resulting equation on the condition that:
$$
\alpha(m_i-1)+2a_i=1 \quad \mbox{ for all } i=1,2,\cdots, N.
$$
We also want integrable solutions that will enjoy the mass conservation property, which
implies $\alpha=\sum_{i=1}^{N}a_i$. Then, \eqref{APM} and \eqref{sss.an} become an equation for the profile function $F(y)$, that must satisfy the following nonlinear anisotropic stationary equation in $\mathbb{R}^N$:
\begin{equation}\label{StatEq}\sum_{i=1}^{N}\left[(F^{m_i})_{y_iy_i}+\alpha \sigma_i\left( y_i F\right)_{y_i}
\right]=0.\end{equation}
Conservation of mass must also hold : $\int U(x,t)\, dx =\int F(y)\, dy=M<\infty $ \ {\rm for} $t>0$.

The conditions on the exponents are worked out. Putting $a_i=\sigma_i \alpha,$ we get unique values for $\alpha$ and $\sigma_i$:
\begin{equation}\label{alfa.aniso}
\alpha=\frac{N}{N(\overline{m}-1)+2},
\end{equation}
where $\overline{m}:=\frac1{N}{\sum_{i=1}m_i}$ is the additive average, and
\begin{equation}\label{ai.aniso} \ 
 \sigma_i= \frac{1}{N}+ \frac{\overline{m}-m_i}{2}.
\end{equation}
 We see that \ $\sum_{i=1}^{N}\sigma_i=1$, and in particular $\sigma_1=1/N$ in the isotropic case.

\noindent{\bf Consequences: } (1) For  a diffusion process spreading in space along all the space variables we want $\sigma_i>0$ for all $i$, which implies that
\begin{equation}\label{aniso.h2}
m_i<\overline{m}+2/N \qquad \forall i,
\end{equation}
which is a condition to avoid large separation between the exponents.

(2) Moreover, we expect the diffusion process to  obey the maximum principle  so that the self-similar solution  decay in time in maximum value like a power of $t$, and this implies the condition $\alpha>0$, which in turn means the following condition on the exponents $N(\overline{m}-1)+2$, i.e.,
\begin{equation}\label{aniso.h3}
\overline{m}> 1-(2/N),
\end{equation}
or in another form, $\sum_{i=1}^{N}\sigma_i> N-2$. This condition is important for the case of fast exponents as we will see.

\subsubsection{\bf Anisotropic fast PME}
We will focus here on the  fast diffusion range (fast diffusion in all directions, $0<m_i<1$ for all $i$)  that has some interesting features, in particular regarding mass conservation. In the study we find a critical exponent
\begin{equation*}\label{m bar}
     m_c =1-\frac{2}{N}.
\end{equation*}
We recall that  $m_{c}$ was also critical  for the isotropic fast diffusion equation (\emph{i.e.} the PME/FDE equation with $m_{1}=m_{2}=...=m_{N}=m<1$).
In this anisotropic setting will always assume the conditions (H1f)
\begin{equation}\label{H1f}
\noindent   \qquad 0< m_i< 1,  \qquad  \mbox{for all} \ i=1,...,N,
\end{equation}
as well (H2)
\begin{equation}\label{H2}
  \noindent  \qquad \overline{m}:=\frac1{N}{\sum_{i=1}m_i} >m_c, \quad \mbox{i.\,e.,  } \sum_{i=1}m_i>N-2.
\end{equation}
Condition (H2) is essential for our results, being a necessary and sufficient condition in the isotropic case for the existence of fundamental solutions that conserve mass. A phenomenon that happens again here. On the other hand, (H1)  is a condition of {\sl fast diffusion in all directions} that is made here for convenience, it is only necessary for the analysis that we do. Note that conditions (H1) and (H2) imply the small separation condition (H3), i.e., \eqref{aniso.h2}.

Paper \cite{FVV23} contains the construction of a theory of existence and uniqueness of continuous weak solutions, the properties of the ensuing evolution semigroup, for instance it is a semigroup of $L^1$-contractions with bounded orbits; mass conservation holds for these solutions.

The point of interest for us here is the existence of a self-similar analysis from which we derive the existence and uniqueness of a self-similar fundamental solution for every mass value and its main properties, like positivity, boundedness and decay; and finally the asymptotic behaviour as $t\to \infty$ of all nonnegative solutions to the Cauchy Problem with the above conditions on data and exponents.
For more information we refer to  \cite{FVV23}

\begin{theorem}\label{afde.fundamental solution} Under the conditions {\rm (H1f)} and {\rm (H3)}, for any mass $M>0$ there is a unique self-similar fundamental solution $U_M(x,t)\ge 0$ of equation \eqref{APM} with mass $M$. The profile $F_M$  of such a solution is a  positive, bounded and smooth function that is monotone decreasing in all outgoing coordinate directions. Moreover,  function $F_M(0,x_i,0)$  decays like $O(|x_i|^{-2/(1-m_i)})$ along the $X_i$ axis, for each index $i$.
For all $M>0$  $F_M$ is just a scaled version of $F_1$ with parameter $M$.
\end{theorem}

It is proved that there exists a suitable solution of the elliptic equation \eqref{StatEq}, which is the anisotropic version of the equation of the Barenblatt profiles in the standard PME/FDE, as was seen before Section \ref{sec.cm.conn}.
 The solution is indeed explicit in the isotropic case, see \eqref{for.Bar2}.  We do not get any explicit formula for $F$ in the anisotropic case.

\noindent{ \bf Remark.}
We recall that by solution we understand a weak energy solution, see the theory in Section  2 of this paper. In the end, the solution of this class is proved to be smooth, so the weak energy solution is indeed a classical solution of the equation.

\begin{theorem}\label{afde.thmasympto}
Let $u(x,t)$ be the unique solution  of the Cauchy problem for equation \eqref{APM} with nonnegative initial data $u_{0}\in L^{1}(\re^{N})$
 under the restrictions {\rm(H1f)} and {\rm (H3)}. Let $U_{M}$ be the unique self-similar fundamental   solution with the same mass as $u_{0}$. Then,
\begin{equation}\label{L1conv}
\lim_{t\rightarrow\infty}\|u(t)-U_{M}(t)\|_{1}=0.
\end{equation}
The convergence holds in the $L^{p}$ norms, $1\le p <\infty$,  in the proper scale
\begin{equation}\label{Lpconv}
\lim_{t\rightarrow\infty}t^{\frac{(p-1)\alpha}{p}}\|u(t)-U_{M}(t)\|_{p}=0,
\end{equation}
\end{theorem}
 Exponent $\alpha$ is given in \eqref{alfa.aniso}.

\noindent $\bullet$ In a further work  \cite{VazVSS-2023} the very singular solutions are constructed for the same range of exponents.
In a recent work \cite{CH25} it is proved that for the range $ \sum_{i=1}m_i<N-2$ extinction in finite time is a common occurrence so that is no Conservation of Mass.

\subsubsection{Fast anisotropic flow. Mass conservation via  local mass estimates}

We assume that $0<m_i<1$ and we are interested on conservation of mass in the critical exponent case
\begin{equation}\label{H2c}\tag{H2c}
\sum_{i=1}^N m_i= N-2,
\end{equation}
i.e., \ $ \overline{m}=m_c=(N-2)/N.$ We want to prove conservation of mass for the solutions of the Cauchy Problem.
Such conservation is already known for the supercritical case $ \overline{m}>m_c$ but we will give a proof that can be extended to the critical case. On the other hand, in the subcritical case $ \overline{m}<m_c$ we expect loss of mass and extinction in finite time.

\begin{theorem}\label{P uniq  comp supp}
Let us assume $ \overline{m}>m_c=(N-2)/N$ holds.
Let  $u$  be a nonnegative constructed weak solution   of the Cauchy problem with a nonnegative bounded and $L^1$ datum $u_0$. Then for every $t>0$ we get
 \begin{equation}\label{masscons}
\int_{\mathbb{R}^N} u(x,t)\, dx= \int_{\mathbb{R}^N} u_0(x)\, dx=M_0.
 \end{equation}
\end{theorem}
\begin{proof} (i) Let  $M(t):=\int_{\mathbb{R}^N} u(x,t)\,  dx$. By the property of monotonicity of the $L^p$ norms in time we know that  $M(t)$ is a monotone function and $M(t)\le M_0$ for all $t>0$. So we need to prove that $M(t)\ge M_0$. We know from the definition of solution that $u\in C([0,\infty):  L^1(\ren))$. Given a  cutoff function $\vp$, we introduce a notation  for the local mass of the  solution
 $$
 X(t;\vp)=\int_{\mathbb{R}^N} u(x,t)\,\vp( x)\,  dx.
 $$
 We assume that the test function   satisfies $0\le  \vp(x)\le 1, $  is nonnegative, smooth and  decays at infinity. It will further specified during the proof. Note that $ X(t;\vp)\to M(t)$ as $\vp\to 1$. We will sometimes omit the  $\vp$ and just write $X(t)$\nc.

(ii) This uses a variation of the  ideas  of \cite{VazVSS-2023}. We consider the  equation for $u$,  multiply it by $\vp$ and integrate by parts to obtain
 \begin{equation}\label{bb1}
 \frac{dX(t)}{dt} = \frac{d}{dt}\int_{\mathbb{R}^N}u\,\vp(x)\, dx=
\sum_{i=1}^N \int_{\mathbb{R}^N}u^{m_i}\,\partial^2_{x_ix_i} \vp(x)\, dx.
 \end{equation}
  Now we choose $\vp(x)=\vp_1(x_1)\vp_2(x_2)\cdots \vp_N( x_N)$. We get an integrated expression for the evolution of $X(t)$ of the form
 \begin{equation}\label{ev.diff}
\displaystyle  X(t_2)-X(t_1)
=\sum_{i} \int_{t_1}^{t_2}\int_{\mathbb{R}^N} u^{m_i}(x,t)\, \partial^2_{x_ix_i}\vp_i(x)\, dx_i \, \Pi_{j\ne i} (\vp_j(x_j)\, dx_j)dt.
 \end{equation}
 Note that the $N$ summands in the left hand result of applying  each differentiation term $\partial^2_{x_ix_i}$ on the corresponding $\vp_i(x_i)$  factor.

 Applying H\"{o}lder's inequality  we obtain
\begin{equation}\label{ev.intd}
 | X(t_2)-X(t_1)| \le \sum_i \left(\int_{t_1}^{t_2} \int_{\ren} u \phi\,dx\right)^{m_i} \,W_i^{1-m_i},
\end{equation}
where
$$
W_i = \int_{t_1}^{t_2} \int_{\mathbb{R}^N}
|\vp_i^{-m_i}\partial^2_{x_ix_i}\vp_i(x_i)|^{1/(1-m_i)}dx_i
 \Pi_{j\ne i} (\vp_j(x_j)\, dx_j) dt.
$$
so that $W_i = (t_2-t_1)Y_i$ with
\begin{equation}\label{y_i}
 Y_i = \int_{\mathbb{R}^N} |\vp_i^{-m_i}\partial^2_{x_ix_i}\vp_i(x_i)|^{1/(1-m_i)}dx_i\,
 \Pi_{j\ne i} (\vp_j(x_j)\, dx_j).
\end{equation}
Note that the terms $Y_i$ are purely geometrical, independent of $u$. Their contribution is important and turns out to be delicate.

 (iii)  We must be  more precise with the choice of the test functions. We will choose a $\vp$ supported in a compact subset $\mathcal K$ of $\mathbb{R}^{N}$. More precisely, we take as $\mathcal K$ the anisotropic box $\mathcal K=\prod_{i=1}^N[-L_i,L_i]$ with $L_i>0$. The lengths $L_i$ will grow in different ways depending on the anisotropy, so it will be convenient to use a common normalization. In that sense we choose a smooth $\phi(y)$, $y\in \re$, such that such that $0\le \phi(y)\le 1$, $\phi(y)> 0$ in the interval $[-1,1]$,  $\phi(y)= 1$ in $|y|<1/2$, and  $\phi(y)=0$ for $|y|\ge 1$.
 Now we put
 $$
 \vp(x)=\phi(\frac{x_1}{L_1})\,\phi(\frac{x_2}{L_2})\cdots \phi(\frac{x_N}{L_N})=\Pi_{i=1}^N \,\phi(\frac{x_i}{L_i}).
 $$
 We will also need a very fast convergence $\phi(y)\to 0$ as $|y|\to 1$ so that the integrals that appear below will be convergent.
 We may calculate the first element in  \eqref{y_i} and evaluate the rest by permutation of the indices. Since
$$
\partial^2_{11}\vp_i(x)= L_1^{-2}\partial^2_{11}(\phi(y_1))
$$
with $y_1=x_1/L_1$, inserting it into the formula for $Y_1$ and integrating separately we get $Y_1= C_1 Z_1 L_2\dots L_N $
 with $C_1=(\int_{-1}^1 \phi(y)\,dy)^{N-1}$ and.
$$
Z_1 =\int_{-L_1}^{L_1} \left|\vp_1^{-m_1}\partial^2_{11}\vp_1\right|^{1/(1-m_1)}\,dx_1=
L_1^{1-(2/(1-m_1)}\int_{-1}^1\left|\phi(y)^{-m_1}\partial^2_{11}\phi(y)\right|^{1/(1-m_1)}\,dy,
$$
At this point, we need to choose $\phi$ so that  $Z_1$ is finite. The last integral is only delicate near  the boundary points $y=\pm 1$ of the interval $|y|\le 1$.  It is easy to see that the integral  is finite if  $\phi(y)$  approaches $y= 1$ like $C(1-y)^k$, and $k>2$ is large, and the same rate holds as $y$ approaches $y=-1$ like $C(1+y)^k$. Indeed, near $y=1$ $\partial_{1}\phi(y)$ decays like   $ck(1-y)^{k-1}$ and $\partial^2_{11}\phi$ is bounded
$$
\partial^2_{11}\phi(y) \sim ck(k-1)(1-y)^{k-2}.
$$
 Hence, $|\phi^{-m_1}\partial^2_{11}\phi|\sim c_1(1-y)^{k-2-m_1 k}$, so that the last integral of $Z_1$ is a bounded number if $k(1-m_1)>1$ and the  same rate of approach to zero happens at $y=-1$. A similar calculation is done for the other $i$, then we need $k(1-m_i)>1$ for all $i$.

Summing up, we get  $Z_i= C_i'' L_i^{1-(2/(1-m_i)}$. All together,
$$
Y_i=C'''_i L_1L_2\cdots L_N L_i^{-2/(1-m_i)},
$$
where the constants $C'''_i$ depend only on  the $m_i$ and $N$. We recall that the integrals  $Y_i$ are geometric factors that depend only on the size and form of the common test $\phi$.

(iv) We now choose the sizes of the intervals $[-L_i,L_i]$ in a convenient anisotropic way.  It works as follows. We choose the factor $L_i^{2/(1-m_i)}$ as constant and equal to $R_{\star}>0$, $L_i=R_{\star}^{(1-m_i)/2}$.
We get for all $i$
$$
{Y}_i= D_i \,R_{\star}^{-\mu}, \quad  \mu= -\sum_{i} \frac{1-m_i}{2}+1
=\frac{N}2(\overline{m}-m_c)\ge0.
$$
We see that under our basic FDE assumption $\overline{m}>m_c$ we  get $\mu>0$ and all the ${Y}_i$ go to zero if $R_{\star}\to \infty$ with the  same rate. Moreover, all the $L_i$ go to infinity in a precise anisotropic way as $R_{\star}\to \infty$, $L_i=R_{\star}^{(1-m_i)/2}$. We use the notation $\mathcal K (R_{\star})$ for the corresponding anisotropic support box. Going back to \eqref{ev.intd} we get
 \begin{equation}\label{afd.estim}
 |X(t_2) - X(t_1)| \le \sum_i D_i' M_0^{m_i}(t_2-t_1)R_{\star}^{-\mu(1-m_i)}.
 \end{equation}
Letting  $R_{\star}\to \infty$ this estimate ends the proof of conservation of the total mass in the case $\overline{m}>m_c$ .
\medskip

(v) {\bf Mass conservation for the critical case.}  On the other hand, when we are in the  critical case $\overline{m}=m_c$ and we repeat the previous analysis we get $\mu=0$ and the   anisotropic $Y_i$ (i the way that we have chosen befor) are constant independent of the basic measure $R_{\star}$, a fact that will make the rest of the analysis a bit more difficult. Actually, the simpler analysis in the critical exponent of the anisotropic case meets the same difficulty.

In order to solve it, we go back to the evolution estimate \eqref{ev.intd} and realize that when we apply H\"older inequality to  formula \eqref{ev.diff} we may use the  fact that the integrals that we later split into two factors had non-zero integrand only on compact sets. Let's take a closer look at this. Actually, choosing a fixed $i$ then for every $t$ the expression
$$
u^{m_i}(x,t)\, \partial^2_{x_ix_i}\vp_i(x) \, \Pi_{j\ne i} (\vp_j(x_j)
$$
is supported inside the  subset
$$
\mathcal A_i({L_1,\cdots, L_N})=\{x: L_i/2<|x_i|<L_i; \quad \ |x_j|<L_j \ \forall j\ne i\}.
$$
Thus, we  may write
 \begin{equation}\label{ev.intd2}
 | X(t_2,\vp)-X(t_1,\vp)| \le  \sum_i C_i (t_2-t_1)^{1-m_i} \left(\int_{t_1}^{t_2} \int_{\mathcal A_i} u \,dx\right)^{m_i}.
\end{equation}
It is convenient to observe that $\mathcal A_i({L_1,\cdots, L_N})$ is included in
$$
\mathcal B_i(L)=\{x: L/2<|x_i|:  \quad \ x_j\in\re \ \forall j\ne i\},
$$
if $L\le L_i$.  We may replace ${\mathcal A_i}$ by $\mathcal B_i(L)$  in the $i$-th integrand of \eqref{ev.intd2} if $ L\le L_i$. Note that for every fixed $i$ the sets $\{\mathcal B_i(L): L>0\}$ are a monotone decreasing family with  increasing $L$, and  \ $\cap_{L> 0} \mathcal B_i(L)$ is the empty set. This allows to bound every term in the right-hand side of \eqref{ev.intd2} by a small constant.
Let us do it for the first term,
$$
\int_{t_1}^{t_2}\int_{\mathcal A_1} u \,dx\le \int_{t_1}^{t_2}\int_{\mathcal B_1(L)} u \,dx.
$$
The integrability of $u$ in $\ren\times (t_1,t_2)$ in the form $u\in C([0,\infty):  L^1(\ren))$ implies that given a small $\ve>0$
we have $\int_{t_1}^{t_2} \int_{\mathcal B_1(L)} u \,dx\le \ve$ if $L$ is large enough, $L\ge l(\ve,1)$. The corresponding term in \eqref{ev.intd2}
is then bounded above by $C_1(t_2-t_1)^{1-m_1}\ve^{m_1}$ whenever $L_1=R_{\star}^{(1-m_1)/2}\ge l(\ve,1)$.
In the same way we estimate the rest of the terms and get
$$
| X(t_2,\vp)-X(t_1,\vp)| \le  \sum_i C_i (t_2-t_1)^{1-m_i}  \ve^{m_i}.
$$
if $R_{\star}\ge l(\ve,i)^{2/(1-m_i)}$ for all $i=1,2\cdots,N$. Consequently for given $t_2>t_1>0$ we get
$$
\lim_{R_{\star}\to \infty} | X(t_2,\vp)-X(t_1,\vp)| =0.
$$
In other words, $M(t_2)=M(t_1)$. The proof is done. \end{proof}

\subsubsection{\bf  Anisotropic slow PME}

 \bc The case where there are slow diffusion exponents $m_i>1$ is considered in a subsequent work \cite{FVV25}, The restrictions on the exponents to obtain constant-mass self-similarity solutions are now two: first, (HIs)
\begin{equation}\label{H1s}
\noindent   \qquad m_i> 1,  \qquad  \mbox{for all} \ i=1,...,N.
\end{equation}
In this anisotropic slow setting will also assume the conditions (H3)
\begin{equation}\label{H3}
m_i<\overline{m}+2/N \qquad \forall i,
\end{equation}

Here is the result about self-similar solutions
\begin{theorem}\label{afde.fundamental solution.s} Under the conditions {\rm (H1s)} and {\rm (H3)}, for any mass $M>0$ there is a unique self-similar fundamental solution $U_M(x,t)\ge 0$ of equation \eqref{APM} with mass $M$. The profile $F_M$  of such a solution is a  positive, bounded and continuous function that is monotone decreasing in all outgoing coordinate directions. Moreover, the profile  $F_M(x)$  is continuous and compactly  supported. For all $M>0$  $F_M$ is just a scaled version of $F_1$ with parameter $M$.
\end{theorem}

Moreover, Theorem \ref{afde.thmasympto} holds in the new setting, with the same formulas for convergence and same $\alpha$, \eqref{alfa.aniso}.

\medskip

\noindent $ \bullet$ {\bf Anisotropic convergence of the supports}.  Let us describe how the spatial support of a nontrivial solution $u(x,t)\ge0$ with compactly supported and bounded initial data evolves in time. In fact, it does in the same anisotropic way as the fundamental solution as we will explain.

Let us recall some notations for the sets that will deal with. For a continuous function $f=f(x)$ defined in a set $A\subset \ren$  we denote the positivity set by
$$
\Omega(f)=\{x\in A:  f(x)>0\}.
$$
When $f=F_M$, the fundamental profile, we write $\Omega(F_M)$ for its positivity set.  For the time-dependent functions $U_M(x,t)$ and $u(x,t)$ we write
$$
\Omega(U_M,t)=\{x\in \ren:  U_M(x,t)>0\},
$$
We observe that when $x=(x_1,\cdots,x_N)\in \Omega(U_M,t)$ and we put
$$
y_i=t^{-\alpha \sigma_i t}x_i  \qquad  \mbox{for all } i_1,\cdots,N,
$$
then $y=(y_1,\cdots,y_N)\in \Omega(F_M).$
This anisotropic contraction of the support will be written here as $\mathcal{C}(t)$, i.e.,
$$
\Omega(F_M)=\mathcal{C}(t)(\Omega(U_M,t)).
$$
The same contraction will be applied to the support of the solution $u(x,t)$ to estimate how its support approaches  $\Omega(F_M)$, the contracted support of  $U_M(x,t)$.

\begin{theorem}\label{exp.cs}  Under the above conditions  we have the  set convergence
\begin{equation}\label{supp.estim.sharp}
\lim_{t\rightarrow\infty} d_H(\mathcal{C}(t)\Omega(u,t), \Omega(F_M))=0.
\end{equation}
Here $d_H$ is the standard Hausdorff distance between sets in $\ren$.
\end{theorem}

See proof in Section of \cite{FVV25}. In particular the support of $\Omega(u,t)$ grows in every coordinate direction $\vec{e}_i$ pointing along the $X_i $ axis like
$$
\overline {x}_i(t)\sim c_i\,t^{\alpha \sigma_i t} \quad \mbox{ for every $i=1,\cdots,N$}.
$$

\noindent{\bf Remark.} Much more accurate approximations have been obtained in the isotropic case (the PME), see \cite{KKochV18, Seis14}, in the form of asymptotic series. These approximations they are not known here.

\subsubsection{\bf  Anisotropic PME/FDE of mixed type}

The study of anisotropic diffusion containing exponents in fast range $m_i<1$, the linear case $m_i=1$ and the slow case $m_i>1$ is very interesting and is work in progress, \cite{FVV25M}\nc.

\subsubsection{Anisotropy in $p$-Laplacian equations}
 Nonlinear diffusion coupled with spatial anisotropy
was studied by Feo, Volzone and the author in \cite{FVV21} in the model called APLE, ``anisotropic  $p$-Laplacian equation'':
 \begin{equation}\label{APL}
 u_t=\sum_{i=1}^N(|u_{x_i}|^{p_i-2}u_{x_i})_{x_i}\quad \mbox{posed in
 } \ \quad Q:=\mathbb{R}^N\times(0,+\infty),
 \end{equation}
We take $N\geq2$ and $p_i>1$  for $i=1,...,N$. As in the case of the anisotropic PME, there is a proper choice of the parameters $p_i$
to ensure the existence of self-similar fundamental solutions with fixed mass, the conservation of mass for finite-mass solutions. The case where all $p_i$ are the same is called orthotropic diffusion. Then, the operator in the right-hand side of \eqref{APL} is $p$-homogeneous.

In this case the first condition on the exponents is taken  to be in the range of fast diffusion in all directions
\begin{equation}\label{H1p}
1<p_i<2 \qquad \mbox{for all } \ i = 1,\cdots,N.  \tag{H1p}
\end{equation}
We recall that in the orthotropic  fast diffusion equation (i.e., equation \eqref{APL} with $p_1 = p_2 =\cdots=p_N = p < 2$, hence $p$-homogeneous), there is a critical exponent,
\begin{equation}\label{p crit}
p_c(N):=\frac{2N}{N+1}
\end{equation}
such that $p>p_c$ is a necessary and sufficient condition for the existence of fundamental
solutions, cf. \cite{Vazquez2006}. Note also that $1<p_c(N)<2$ for $N\ge 2$.

The second condition on the exponents in the general non-orthotropic case, is
\begin{equation}\label{p_isum}
\sum_{i=1}^{N}\frac{1}{p_i}<\frac{N+1}{2},\tag{H2p}
\end{equation}
that is crucial in the study. We we may also write it in terms of $p_c$  as: $\bar p>p_c$, where $\bar p$ is the inverse-average
\begin{equation}\label{p bar}
\frac{1}{\overline{p}}=\frac{1}{N}\sum_{i=1}^{N}\frac{1}{p_i}.
\end{equation}
We point out that (H2p)  excludes the presence of (many) small exponents $1<p_i< p_c$ close to 1.
We ask the reader to compare the different algebra of conditions \eqref{p_isum} and (H2)=\eqref{H2}.

The study of self-similarity follows the lines of the similar study that we have sketched for the anisotropic porous medium equation, we refer to \cite{FVV21} for full details. To note that we do not get any explicit formula for the profile $F$ in the general anisotropic case
but the profile is explicit for the orthotropic case.

\medskip

\noindent {\bf Problem.} Prove conservation of mass as  in Theorem \ref{thm.mc.cple} in the  orthotropic case.

\medskip

\noindent{\bf Other forms of anisotropy or inhomogeneity.} A number of other equations have been studied in anisotropic or inhomogeneous settings. See \cite{AntShm2015}, where the basic equation is
\begin{equation}
\partial_t u=\mbox{div} (|u|^{\gamma(x,t)} \nabla u).
\end{equation}
Many variants are considered, among them the so-called $p(x,t)$-Laplacian. The equations are called ``PDEs with nonstandard growth conditions''. The  monograph offers a rather complete exposition of parabolic PDEs with nonstandard growth conditions, and a few results on the hyperbolic case. There is also a large literature of elliptic equations and calculus of variations functionals with nonstandard growth.  Conjectures say that when $\gamma(x,t)$ is keeps bounded \nc away from zero at sufficient rate then conservation of mass should hold, a good proof is needed.


\medskip

\subsection{Diffusion with convection. Fokker-Planck equations}
They generalize the plain diffusion equations. The linear version takes the form
\begin{equation}\label{FPS}
\displaystyle {\partial_ t}u +\mbox{div\,} ( u\,\mathbf {v})=
\mbox{div\,} (D\,\mathbf {\nabla }u )+ F\,,
\end{equation}
where $D$ is the (possibly variable) diffusivity matrix, $\bf v$ is the convective velocity, and $F$ is the forcing term, and accounts for sources or sinks. They can be seen as  diffusive counterparts of the continuity equation \eqref{cont.eq}. In the physics literature they are known as Fokker-Planck equations or Smoluchowski equations. In many  cases $\bf v$ derives from a potential $U$, i.e., ${\bf v}=-\nabla U$. In many  applications to fluid mechanics $\bf v$ is incompressible, i.e., \ $\nabla \cdot \bf v = 0.$ There is an extensive literature on the subject from the viewpoints of physics and mathematics, from the early  papers like \cite{Smo1915} to monographs like \cite{Risk1984, BKRS}.
When $F=0$ conservation of mass is natural, but the presence of sources or sinks change the law into the form
$$
\frac{dM}{dt}=\int F\,dx+ BC,
$$
where $BC$ accounts for the possible boundary terms that may also cause loss or gain of mass, as mentioned above.

We may pass to the nonlinear diffusion versions by putting $D=D(u,\nabla u)$. $F$ may include terms of reaction or absorption in the form $F(u,\nabla u)$.

\medskip

\noindent{\bf Fokker-Planck reduction.}
 There is a typical Fokker-Planck equation with quadratic potential $U=|x|^2$ that is reducible to plain diffusion by a change of variables, and it plays a role in the analysis of heat equation flows. Indeed, the change
\begin{equation}\label{trans.fp1}
y=x\,(t+1)^{-1/2}, \quad \tau=\log(1+t), \quad u(x,t)=(t+1)^{-N/2}v(y,\tau)
\end{equation}
allows to pass from $   \partial_t u=\Delta u$ to the Fokker-Planck equation for $v(y,\tau)$
$$
\partial_\tau v = \Delta_y v + \frac 12 \mbox{div\,}(y\,v)=\mbox{div\,}\left(v \, \nabla(\log v+ \frac12 y^2)\right).
$$

Transformation \eqref{trans.fp1} can be also applied to the solutions of the PME and PLE to obtain Fokker-Planck versions of those equations (called nonlinear FkPE for short), see \cite{VazquezPME2007}, Section 18.4. Note that the exponents $1/2$ and $N/2$ have to be changed into $\alpha$ and $\beta$ as given in Section \ref{ssec.pme} above. It has been amply observed that the asymptotic analysis is much more convenient  when done in terms of the Fokker-Planck versions, because it takes the form of convergence to equilibrium. Entropy dissipation methods have been quite successful in this analysis, see e.g. \cite{CarrTosc00, CJMTU, BBDGV, Jungel2016}. Wasserstein distances are used in that context, see \cite{AmGiSav05, CarrMcV03, CarrMcV06, Otto2001}.

The formation of singularities in finite or infinite time is a hot topic in current research, in particular for diffusion-aggregation problems, see a relevant example in \cite{CarrGCV2021}.

\medskip

\noindent{\bf Convection with fractional diffusion.} As is by now well known, nonlinear evolution equations with fractional diffusion arise in many contexts, in particular combined with convection or drift terms. A celebrated paper by Caffarelli and Vasseur \cite{CaffVass2010} studies drift-diffusion equations with fractional diffusion for a quasi-geostrophic model. It has the diffusion convection form \eqref{FPS} with an incompressible velocity and the right-hand side of the equation is $-(-\Delta)^{1/2}u$, a fractional diffusion term.

\medskip
\subsection{Equations and systems of aggregation-diffusion }
\bc There is an extensive literature covering the aggregation-diffusion models which, in a first general formulation can be written as
\begin{equation}\label{EQ.agdiff}
\partial_t u = \mbox{\rm div}\left(\mathfrak{m}(u)\nabla (U'(u) + V + W \star u)\right).
\end{equation}

In this notation $\mathfrak{m}$ is called mobility, $U$ is a non-linear function dealing with diffusion, $V$ is
a confining potential, and $ W $ is a suitable interaction term called aggregation potential. This family of
equations usually comes as a generalisation of the Keller-Segel model together with nonlinear PME-type diffusion. It is a very rich topic with many directions being explored nowadays.  The standard Keller-Segel system in chemotaxis \cite{JaLuck92, KellerS71, Bln2006} is one of most popular systems coming from Biology, where a standard diffusion process runs parallel with a strong aggregation force of chemical origin that leads to an equation of the above form where $W$ is the Newtonian potential. The standard mobility is linear $\mathfrak{m}(u)=cu$
and then the evolution \eqref{EQ.agdiff} i a 2-Wasserstein gradient flow \cite{AmGiSav05, Santamb15, Santamb2017}; the solutions are probability measures. \bc A doubly parabolic Keller-Segel system is
\begin{gather}\label{EQ.agdiff.dp}
\partial_t u = \Delta u+ \mbox{\rm div}(u\nabla v),\nonumber\\
\tau \partial_t v = \Delta v + u- au. \nonumber
\end{gather}
where $u\ge 0$ is a population density and $v\ge 0$ is the concentration of an agent substance, \cite{BL13}.
The ensuing competition \nc between diffusion and aggregation  may concentrate the substance in one or several points, a blow-up phenomenon that depends critically from the amount of mass at hand. In two space dimension the critical mass value is $8\pi$,  \cite{Bln2006, BKLN, BKLN2}. Bear this in mind to compare with other critical mass numbers we have found.

\bc Aggregation-diffusions generalizing Keller-Segel models have been studied in a series of papers, like \cite{CalvCarr06, CalvCarrH17, CarrHVY19, CarrHMV18}. There is a long discussion in the literature about evolution asymptotics and steady states for this equation, and the analysis usually is based on distinguishing in terms of the mass. At the PDE level, conservation of mass
in this setting usually requires, besides the usual PME assumptions of $\Phi'(u) = mU''(u)$,
some behaviour of $V, W$ at infinity\nc.

There is currently much activity in the field by many authors that we will not report here due to its extension. Let us mention a work of us showing mass concentration in infinite time see \cite{CarrGCV2021}. \bc This happens with nonlinear fast diffusion terms. Further work on the issue is reported in \cite{Bln2008, CarrGCY23, CarrFJDG24}\nc. Many systems of nonlinear diffusion equations, possibly combined with other types of equations, exist in the literature, but a systematic approach seems to be pending. See some examples in
\cite{CaffGual20, Desv15, DBFr85, Jungel2016, KimKiahm2020, MarkRS1990}. An  interesting effect of nonlinear mobility in nonlinear diffusion is {\sl saturation} while conserving mass, which leads to possibly double free boundaries and interesting effects not present in other
situations. This has been analysed recently in \cite{BailoCarHu23, CarrFJDGSiam25} in bounded domains.
Interesting  survey papers are  \cite{CarCraYao19, DGC23}.

\subsection{Diffusion in Mathematical Biology} Mathematics plays a big role in  Biology in many respects. The books \cite{EpMaz2013} and \cite{Perth15} cover many applications of nonlinear diffusion models and methods. Special attention concerns Hele-Shaw mechanical models of tumor growth that involve a reaction diffusion for tumor density and nutrients. Relevant work is found  in \cite{PerQV14, MelPert17} and other references. These systems use many of the above tools but do not conserve mass. On the other hand, the kinetic formulation of conservation laws is studied in \cite{Perth02}, and transport equations in biology are studied in \cite{Perth07}.

\medskip

\subsection{The infinity fractional  Laplacian}

\bc
We explain  next a limit case of the  fractional Laplacian equation.
In 2012 Bjorland, Caffarelli and Figalli introduced  in \cite{BjCaFi12} equations involving the so-called \sl infinity fractional  Laplacian \rm as  models for a nonlocal version of the ``tug-of-war'' game. Following their explanation, instead of flipping a coin at every step,  every player chooses a direction and it is an $s$-stable L\'evy process that chooses both the active player and the distance to travel. The  infinity fractional Laplacian, with symbol $\ifl$, is a nonlinear integro-differential operator. The corresponding equation in with standard derivatives is the {\sl normalized $p$-Laplacian equation} in the limit $p=\infty$ listed in \eqref{eq.rpl.inf}.

Among the slightly different definitions, we consider here the  definition introduced in \cite{BjCaFi12} (see also \cite{DTEnLe22} given by
\begin{equation}\label{eq:def1}
\ifl\phi(x):= C_s \sup_{|y|=1}\inf_{|\tilde{y}|=1} \int_{0}^\infty\left(\phi(x+\eta y)+\phi(x-\eta \tilde{y})-2\phi(x)\right) \frac{ \eta }{\eta ^{1+2s}} \qquad
\end{equation}
with $s\in(1/2,1)$. The constant is usually taken as  $C_s=(4^ss\Gamma(\frac{1}{2}+s))/(\pi^{\frac{1}{2}}\Gamma(1-s))$ but the value is irrelevant for the  results we have in mind. In their paper \cite{BjCaFi12} the authors study two stationary problems involving the infinity fractional  Laplacian elliptic problems posed in bounded space domains, namely, a Dirichlet problem and a double-obstacle problem.

 In paper \cite{dTEJ23} we consider the evolution problem
\begin{equation}\label{ifl.ev}
\partial_t u(x,t)= \ifl u(x,t),  x \in \ren , t > 0,
\end{equation}
with  $s\in(1/2,1)$ and $N\ge 2$ and initial data $u(x,0)=u_0(x)$, $x\in \ren$. Actually, when $N=1$ the operator $- \ifl$ is just the usual linear fractional
Laplacian operator $(-\Delta)^s$ of order $s$, and equation \eqref{ifl.ev} is just the well-known fractional heat equation, see also what is said in Subsections \ref{ssec.fracLap} and \ref{ssec.fhe}.

However, for $N\ge 2$ the operator is nonlinear so a new theory is needed. We develop in \cite{dTEJ23} an existence theory of suitable viscosity solutions for the above parabolic problem based on
approximation with monotone schemes. We show that the obtained class
of solutions enjoys a number of good properties.  As in the elliptic
case \cite{BjCaFi12}, we lack a uniqueness result in the context of
viscosity solutions. However, we are able to prove an important
comparison theorem relating two types of solutions, classical and
viscosity ones. As a counterpart, we also obtain uniqueness of classical solutions.
Moreover, we show that for smooth, radially symmetric
functions that are nonincreasing along the radius in $\ren$ with $N\geq2$, the operator $-\ifl$
reduces to the classical fractional Laplacian $(-\Delta)^s$ in
dimension $N=1$. This is a crucial property. The following consequence is important for our topic of mass conservation.

\begin{theorem}\label{ifl.nmc} Let $N\ge 2$, and let $u_0(x)=u_0(r)\ge 0$ be a smooth and  radially symmetric
function, nonincreasing along the radius $r=|x|$. Let $u(x,t)$ the classical viscosity solution of Problem \eqref{ifl.ev}.
Then conservation of mass does not hold.
\end{theorem}

 Since the equation for such data reduces to an 1D equation that conserves the mass defined in 1D,  it cannot conserve mass in N-dimension because of the Jacobian factor $r^{N-1}$.\nc

\subsection{Geometric flows}
Geometric flows are a very important class of nonlinear diffusion problems. There are many options. We mention some of them.

\noindent {\sc Flows in image processing.} We have already seen the total variation flow. On the other hand,  Alvarez et al. \cite{AlvDD16} discuss the modeling issue as it is important in  image processing, and they arrive at the equation
\begin{equation}
\partial_t u=\beta\,\left(\mbox{div\,}(\frac{\nabla u}{|\nabla u|}) \right)\,|\nabla u|\,,
\end{equation}
where $\beta(s)$ is a nondecreasing function given normally by a power of $s$. They work in dimension 2, but the results admit generalizations to higher dimensions. Generation of a solution is done by the implicit time discretization method of \cite{CrandallLiggett71} acting on the space $BUC(\ren)$ of bounded uniformly continuous functions. There is no natural conservation law at hand.

\medskip

\noindent {\sc The mean curvature flow.}
According to \cite{CMP15}, the mean curvature flow (MCF) is the negative gradient flow for area. This is a nonlinear partial differential equation for the evolving hypersurface that is formally similar to the ordinary heat equation, with some important differences. MCF behaves like the heat equation for a short time with the solution becoming smoother and small-scale variations averaging out. However, after more time, the nonlinearities dominate and the solution becomes singular. To understand the flow, one must understand the singularities it goes through. This is a difficult task that has been studied for  in the last century to model things such as cell, grain, and bubble growth. End of quote.

\medskip

An important case is that of graph surfaces. Then, for
$\Omega\subset \ren$, the graph of $u : \Omega\times [0, T ] \mapsto \re$ is a hypersurface in $\re^{N+1}$ moving by mean curvature flow. One may compute that the non-parametric mean curvature flow of a graph can be written as
\begin{equation}
\partial_t u=\sqrt{1+ |Du|^2}\, \partial_i
\left(\frac{\partial_i u} {\sqrt{1+ |Du|^2} }
\right)=\Delta u-\frac1 {\sqrt{1+ |Du|^2} } \, D^2 u\,(Du, Du),
\end{equation}
notations taken from \cite{BW}. \bc See also \cite{Mant11, RitSin2010} \nc for further information.

\medskip

\noindent {\sc The curve-shortening flow.} It is a process that modifies a smooth curve in the Euclidean plane $\re^2 $ by moving its points perpendicularly to the curve at a speed proportional to the curvature.
The curve-shortening flow is an example of a geometric flow, indeed it is the one-dimensional case of the mean curvature flow. As the points of any smooth simple closed curve move in this way, the curve remains simple and smooth. We take this information from the detailed article in \bc Wikipedia, \cite{CSF Wikipedia},  \nc where plenty of references are given, and the parabolic equation is written as
$$
\frac{\partial C}{\partial t}=\frac{\partial^2 C}{\partial s^2}=\kappa \,n,
$$
a kind of geometric interpretation of the heat equation. Let us point out that closed curves evolving by this flow  lose the enclosed area at a constant rate, which amounts to  $2\pi$ units of area per unit of time. This is the closest equivalent we can find to some of our previous nonlinear diffusion models with critical exponents, where mass is also lost with a universal constant rate.
The curve-shortening flow is applied in image analysis and other applied topics.

\medskip

\subsection{Other equations in non-divergence form. Anomalous exponents}\label{ss12.6}
When the divergence form is lost, then mass conservation is not a natural occurrence. However, parts of the programme can be saved in modified form. We will give two examples.

\noindent {\bf \ref{ss12.6}.1.} The following heat equation with two diffusivities depending on the sign of $u_t$
\begin{equation}\label{elpleq}
\partial_t u +\gamma |\partial_t u |=\Delta u, \    0<|\gamma|<1,
\end{equation}
was introduced by Barenblatt and Krylov (1955) as a model for \sl elasto-plastic filtration \rm (when charging and discharging are not equally efficient for the flow), see \cite{BarenSiva69,
Kamin69}. The mathematical theory was  studied in \cite{KaminPelVaz1991} where existence and uniqueness for the Cauchy problem with continuous initial data vanishing at infinity is established and the evident failure of the MC law was verified (it is natural from the extra term in the equation). A main goal of the paper was to construct a self-similar solution $B(x,t)\ge 0$ of the form
\begin{equation}\label{elpleq.sss}
 B(x,t)=t^{-\alpha/2}F(y), \qquad y=|x|/t^{1/2},
\end{equation}
for some special $\alpha>0$ that depends on $\gamma\in(-1,1)$. The function $\alpha(\gamma)$ is shown to be strictly increasing and continuous, tending to infinity as $\gamma\to 1 $ and to $\max\{N-2,0\}$ as $\gamma\to -1$. Conservation of  mass holds only for $\gamma=0$, $\alpha(0)=N$.

The exponent $\alpha$ is not determined from a conservation law like mass conservation, but from PDE considerations (existence of solution of an eigenvalue problem), and it is called an {\sl anomalous exponent}, a favorite topic in Barenblatt's research and teaching, \cite{BarentSSSIA}. Barenblatt also uses in that respect the label \sl self-similarity of the second kind\rm.   A detailed analysis of the dependence of the anomalous decay exponent $\alpha(\gamma)$ with varying $\gamma$ is made in \cite{Pel93}.

As a bonus, paper \cite{KaminPelVaz1991} also shows (after a delicate analysis) that the above self-similar solutions of the second kind give the asymptotic behaviour of a large class of nonnegative solutions of the elasto-plastic heat equation with integrable initial data.  \bc This is much in the spirit of the standard version of the HE and also the PME and PLE. The analysis of problems with anomalous exponents is a quite delicate subject, \cite{BarentSSSIA, BarenSiva69, ArVa95, EggFont}  \nc.

As an extension of this theory the authors of \cite{MenOr19} identify the sequence of self-similar profiles $F=F(y;\alpha_k)$, with possibly changing sign, by closed representation in terms of confluent hypergeometric functions. Employing specific properties of these special functions, oscillatory and asymptotic aspects of the profiles $F$ are obtained. They prove that there is an increasing and unbounded sequence of exponents $\alpha_0<\alpha_1<\cdots$, as in other diffusion equations. Here $\alpha_0$ corresponds to the exponent $\alpha$ of the nonnegative solution mentioned in \eqref{elpleq.sss}\nc.

\medskip

\noindent {\bf  \ref{ss12.6}.2.} The so-called \sl dual porous medium equation\rm,
\begin{equation}\label{eqn.n}
\partial_t u =|\Delta u|^{m-1}\Delta u, \    m>0,
\end{equation}
is another formal relative of the PME and the PLE. Its non-divergence form comes again in the way, and conservation of mass fails unless $m=1$. This poses the problem  of whether a theory of self-similar solutions with anomalous exponents can be developed. A positive answer was found in $N=1$ in the paper \cite{BerHuVa93}. There the equation
\begin{equation}\label{eqn.n1}
\partial_t u =|u_{xx}|^{m-1}\, u_{xx},
\end{equation}
was differentiated once and twice to obtain equations for $v=u_x$ and $w=u_{xx}$:
$$
\partial_t v=(|v_{x}|^{m-1}\, v_{x})_x, \quad \partial_t w=(|w|^{m-1}\, w)_{xx},
$$
so the nonnegative self-similar solution we were looking for is just a convenient changing-sign self-similar solution of the PME (twice integrated in space) with fast decay at infinity (or even compact support). These solutions of the PME have been classified by Hulshof in 1991 \cite{Hulshof91}. In this way, the ``anomalous exponent'' was identified. This is crucial to determine the asymptotic decay rate of the bounded, nonnegative solutions of  \eqref{eqn.n1}.

Let us examine a bit further the nonnegative self-similar solution constructed in \cite{BerHuVa93} for $m>1$ in 1D. It has the usual form
$$
Z(x,t)=t^{-\alpha}\,F(x\,t^{-\delta})
$$
with some precise exponents $\alpha$ and $\delta>0$. These exponents are uniquely determined. The resulting profile $F$ is known to be compactly supported, continuous and symmetric. It follows that the solution has finite mass for every $t>0$.

Close inspection of the exponents is interesting. Firstly, we have the algebraic condition $\alpha(m-1)+ 2m\delta=1$, needed in order to satisfy the equation. However, conservation of mass is not true! Indeed, we have shown in \cite{BerHuVa93} that for all $m>1$ we have $\alpha>\delta$ so that
$$
M(t)=\int_{\re}Z(x,t)\,dx= C\,t^{\delta-\alpha}.
$$
We see that the total  mass decreases from the value $M(0)=+\infty$ (i.e., a very singular initial mass) to $M(+\infty)=0$ with a definite decay. Let us finally point out that this another case of anomalous exponents, something that frequently occurs when needed conserved quantities fail. It is not clear whether we may call this self-similar function a fundamental solution in the sense of the Barenblatt solutions of the PME theory. In any case, further analysis of this equation would be welcome.

\medskip

\noindent \bf Note\rm. Besides the mathematical interest as a simple example of
a fully nonlinear and degenerate parabolic equation, equations like \eqref{eqn.n1} appear in some
problems in elasticity with damping \cite{DuLi72} and problems of Bellman-Dirichlet type. The existence and uniqueness of solutions of equations of the general type
\begin{equation}
\partial_t u =\beta(\Delta u)
\end{equation}
with monotone \bc nondecreasing $\beta$\nc, has been studied by a number of authors, like \cite{Strauss66, Konishi73, BenHa75}, and is also an important tool in the studies of the PME, see \cite{CarrVaz03}. \bc The monotonicity of $\beta$ ensures the degenerate ellipticity of the equation. viscosity solutions are designed to handle these equations. Comparison principle is then expected for such equations\nc.

\newpage

\medskip

\subsection{Other nonlinear diffusion equations}
 We collect here a number of equations that have been related to our main models in different ways.

\noindent{\bf Flows in heterogeneous media.} Porous media equations in heterogeneous media,
\begin{equation}
\rho(x)\,\partial_t u = \Delta u^m,
\end{equation}
were proposed by Kamin and Rosenau \cite{RosenKam82} as a model in plasma physics. The density $\rho(x)>0$ measures the inhomogeneity of the medium and is assumed to go to zero at infinity (meaning almost vacuum in the far field). The main point in the mathematical analysis that we are performing is that conservation of mass takes the form
$$
\int \rho(x) u(x,t)\,dx=\mbox{constant}.
$$
This modified formulation has been proved for a  wide range of densities $\rho$ and exponents $m$. Self-similar solutions of Barenblatt type are constructed and asymptotic convergence is proved, \cite{GKKV04, KaminReyVaz2010}. But cases of loss of mass are reported when $\rho$ decays very rapidly at infinity, and this phenomenon has been partially studied,  \cite{GHP97}. The critical exponents in this study are $\rho(x)\sim |x|^{-2}$ and $\rho(x)\sim |x|^{-N}$.

The radial solutions of this problem have been shown in \cite{Vaz2015} to be related to heat flows in negatively curved manifolds. Therefore, the  result  on non-conservation of mess has a counterpart in the  section about flows in manifolds.

    \medskip

\noindent{\bf Reducible equations.}  We mean types of equations where the conservation law
is not apparent, but it can be identified and used after convenient transformations.

\noindent $\bullet$ An example is
\begin{equation}
\partial_t \rho = \rho^{\alpha}\Delta \rho,\qquad -\infty<\alpha<1,
\end{equation}
that can be reduced to the standard PME via a change of variables.  Indeed, if $u=\rho^{1-\alpha}$ then it satisfies $\partial_t u= c \Delta u^m$, with $m=1/(1-\alpha)$. Thus, $0<\alpha<1$ leads to the slow PME ($m>1$) with conservation of mass and free boundaries. The trick is that we must take as mass the integral the expression
$$
M_\alpha(t)=\int_\ren \rho^{1-\alpha}(x,t)\,dx,
$$
and then it is conserved in time. We leave to the reader the adaptation of the cases where $\alpha<0$ and the identification of the critical exponent $\alpha_c=-2/(N-2)$. Finally, in case $\alpha>1$ we find that $u$ is an inverse power of $\rho$, so the interpretation may be far-fetched.

\smallskip

\noindent $\bullet$ A similar situation is posed by the so-called pressure equations of the type
\begin{equation}
p_t=\sigma(p)\Delta p + \gamma |\nabla p|^2.
\end{equation}
They have been studied by numerous authors like \cite{BeDPassoUg, BeUg, ChaVa, CaVa1999, BranVa2005}.
When $\sigma (u)=c\,u$ the transformation into a PME works by the formula $p=k\,u^{m-1}$ with suitable constant  $k$. A viscosity theory is available for those equations. \bc There is no mass conservation law for those equations, unless $\sigma(p)=p$, in which case we fall into the PME with $m=2$\nc.

\medskip

\noindent {\bf   Equations with more general structure conditions.} One may wonder if the main results of the programme can be obtained under more relaxed conditions on the form of the equation, like
\begin{equation}
\partial_t H(u )=\mbox{\,div\,} {\bf A}(x,t,u,\nabla u).
\end{equation}
Enormous progress has been done in specific cases. Particular cases related to our presentation above
are the generalized filtration equations of the form $\partial_t u=\Delta \Phi(u)$ with $\Phi$ monotone nondecreasing \cite{DaskaBk} or the corresponding gradient diffusion flows, $\partial_t u= \mbox{div} (D(|\nabla u|)\nabla u)$, like the Perona-Malik  equation  of interest in image processing, \cite{PMalik90}. These equations are known in the image processing literature as anisotropic diffusion, cf. \cite{Weickert98}.  The name ``anisotropic'' here refers to  a concept that has nothing to do with our studies in Subsection \ref{ssec.10.2}. Indeed, the underlying idea of the technique is to reduce image noise without removing significant parts of the image content, typically edges, lines or other details that are important for the interpretation of the image.
For an application  to image contour enhancement see \cite{BarVa}.

\subsection{Brief reference to other related models}\label{sec.models}

Analyses of self-similarity and mass conservation exist for a number of other related topics involving diffusion and other effects of physical or mathematical interest. We will mention some of them without aiming at a minimum coverage because of lack of space and/or expertise.

\noindent $\bullet$ Nonlocal diffusion problems with integrable kernels. This a topic with much  activity that can be followed for instance in the book \cite{AMRT2010}. Also, this topic is treated in a  recent publication, \cite{Rossi2020}.

    \medskip


 \smallskip

\noindent $\bullet$ There is a surge of activity on evolution processes driven by fractional in time derivatives (of the Caputo or Riemann-Liouville types) combined with linear and nonlinear space operators. Since this is worth treating in itself and much is still going on, let us only mention a book reference \cite{GalWarm20}, and also the paper \cite{ChanGCV22}, where we have worked on the issue.

    \medskip

\noindent $\bullet$  Work on higher-order equations proceeds maybe a minor pace, but the applications are interesting. Thus, fourth order equations appear thin film modeling and analysis. See \cite{Smyth88, BerFr90, CarrTosc02}. For examples of higher-order fractional diffusion see \cite{Imbert11, SegVaz2020}.

    \smallskip

We continue with two large topics where diffusion plays a partial but very important role.

\medskip


\noindent {\bf Reaction-Diffusion Equations.} This is a huge topic with specific problems, among them blow-up. The basic equation is Fujita's model \cite{Fujita66}
\begin{equation}\label{eq.diff.reac}
\partial_t u = \Delta u + f(u)
\end{equation}
with the power-like nonlinearity $f(u)=u^p$, $p>1$, as the most typical case. The main issue is the possibility of  blow-up in finite or infinite time, and the mode of explosion if blow-up happens. There is a huge literature on this topic, that we will not discuss here to avoid distractions from our main topics. For  convenient monographs see \cite{HarJend15, QS2009}, the Russian classic \cite{GKMS}, and the survey paper \cite{GalVaz2002}. There is no conservation of mass in such models. In fact, blow-up in reaction diffusion is a kind of opposite to mass loss in fast diffusion. Indeed, in the latter case mass wanders to infinite where it is lost, while in reaction diffusion (a part of the mass) concentrates into an infinite singularity at some point or many points of space.

A different direction is taken when the reaction term is negative, in that case it is called absorption. The typical equation is
\begin{equation}\label{eq.diff.abs}
\partial_t u = \Delta u - f(u)
\end{equation}
where $f(u)$ can be $u^p$ with any real exponent $p$, even negative. See the monograph \cite{GalVazBook}. Here mass if lost by the mechanism of annihilation via absorption into/by the medium. In these models the study of mass conservation for finite-mass solutions is replaced by the study of mass evolution, with phenomena like blow-up and extinction that demand very careful analytical tools. Much is known on those issues.

\medskip

A still different reaction type deals with the so-called Fisher-KPP models that give rise to the appearance of travelling waves, a very important object in many applications. The most typical reaction term is $f(u)=u(1-u)$ for $0<u<1$. See the original papers \cite{Fisher} and  \cite{KPP}, and works like \cite{AronWein1975, VolVol94}.

An alternative to the latter model is the so-called ``combustion model'', where $f_\ve(u)$ is nonnegative, smooth, has support in a right neighbourhood of $u=$:  $[0,\ve]$, and for every $\ve>0$ we impose the normalization  $\int f_\ve(u)du=1/2 $ (or any other fixed constant). Such approach comes from the derivation of the model for flame propagation at what is called the ``high activation energy limit'' under a number of simplifying assumptions,
as proposed by Zeldovich et al. \cite{ZFK38}. In this respect, we must consider  $\ve\sim 0$. The limit problem when $\ve\to 0$ was studied by a number of authors, cf.   \cite{BerCaffNi90, CaffVaz95, GalHuVaz97}.

\newpage

\subsection{Mention of other equations and systems}\label{sec.other}

There are other important topics in PDEs and applications that have shown interesting connection with the issues we have discussed in nonlinear diffusion. Let us mention some of the most important.\rm

\noindent {-First order Hamilton-Jacobi equations.

 \noindent {-Other geometrical equations.  Ricci flows,  mean curvature flow, ... }

\noindent {-Diffusion in Image Processing.}

\noindent {-Diffusion and Optimal transportation.}

\noindent {-Schr\"odinger and wave equations.}

\noindent {-Particle methods, numerical diffusion, computational challenges.}

\noindent {-Euler inviscid flows. Viscous fluid flows.}

\noindent {-Self-similar  turbulence diffusion.}

All of these topics have their own rich theory with many peculiarities and should be discussed elsewhere.
Some of them have been briefly mentioned above.

\newpage

{\bf \huge Appendices}
\appendix
\section{List of main used abbreviations}

A number of abbreviations for key concepts are used in the text. They are listed  below.

\begin{table}[h]
\caption{ }
\centering
\begin{tabular}{clccl}
  & & &  \\
  & CE & & &Continuity equation \\
  & DNLE & & &Doubly nonlinear equation  \\
  & FDE  & & &Fast diffusion equation  \\
  & FrHE  & & &Fractional heat equation   \\
  & FkPE  & & &Fokker-Planck equation \\
  & FrPLE  & & &Fractional  p-Laplacian equation  \\
  & FrPME  & & & Fractional porous medium equation \\
  & HE  & & &Heat equation \\
  & HJ  & & & Hamilton-Jacobi\\
  & IMC  & & & Integrated Mass conservation \\
  & L1C & & & $L^1$ contraction property\\
  & MC  & & &Mass conservation \\
  & NLD  & & & Nonlinear diffusion \\
  & PDE  & &  &Partial differential equation  \\
  & PLE  & & & p-Laplacian equation \\
  & PME  & & & Porous medium equation\\
  & SGT  & & & Semigroup Generation Theorem\\
  & TVF  & &  &Total variation flow \\
  & VFD &  & & Very fast diffusion
\end{tabular}
\label{Table}
\end{table}

%
%
%
%
%
%
%
%
%
%
%
%
%


\

\section{Aleksandrov's Reflection Principle}
\label{ssec.aleks}

 The Reflection Principle is an important tool in the analysis of the geometrical aspects of some elliptic and parabolic equations that are invariant under symmetries. It was introduced by A.D. Aleksandrov, \cite{Alek60}. It exploits the invariance of the equation
under reflection with respect to hyperplanes to derive comparison along  directions. We will need it in the sequel as a technical tool to ensure monotonicity properties of the solutions of the PME posed in $\ren$. Since the argument is quite general it applies to other equations and settings. We gather here the pertinent results and brief explanations taken from the monograph \cite{VazquezPME2007} for easier reference.

We need some notation. Any $H$,  hyperplane of $\ren$, divides $\ren$ into two
half-spaces $\Omega_1(H)$ and $\Omega_2(H)$. We denote by
$\pi=\pi_H$ the specular symmetry that maps a point $x\in\Omega_1$
into its symmetric image with respect to $H$, $\pi(x)\in\Omega_2$.
\index{Aleksandrov's Reflection Principle}

\begin{lemma} [Lemma 9.17 of \cite{VazquezPME2007}] \label{lem.ap} Let $u \ge 0$ be the unique solution of the Cauchy
problem for the HE or the PME with initial data $u_0\in L^1(\ren)$ and
assume that for a given hyperplane $H$ we have
\begin{equation}
u_0(\pi_H(x))\le u_0(x) \quad \mbox{for all } x\in \Omega_1(H).
 \end{equation}
Then, for all times
\begin{equation}
u(\pi_H(x),t)\le u(x,t), \quad \mbox{for all } x\in \Omega_1(H).
 \end{equation}
\end{lemma}

\noindent {\sl Proof.} \rm We recall that both equations are invariant when applying spatial rotations and translations. Performing these operations  we may assume that $H=\{x_1=0\}$, so that
 $$
\pi(x_1,..., x_N)= (-x_1,..., x_N)
 $$
and  $\Omega_1=\{x_1>0\}$. By approximation we may assume that the solutions
are continuous and even smooth, even at $t=0$. We consider in $\widehat Q= \Omega_1\times (0,\infty)$ the solution  $u_1=u$ and a second solution
$$
u_2(x,t)=u(\pi(x),t).
$$
By the symmetry invariance, $u_2$ is also a solution of the HE/PME; it
has initial values $u_2(x,0)\le u_1(x,0)$ in $\Omega_1$ by assumption.
The boundary values on $\Sigma=H\times (0,T)$ are the same.
If we are able to justify the Maximum Principle for these
solutions, then
$$
u_2(x,t)\le u_1(x,t)
$$
in $\widehat Q$, which proves the result. \qed

As the reader may have observed, there is nothing very particular
of the HE or PME in the proof. Indeed, the Reflection Principle holds
for the typical parabolic equations with space independent
coefficients, like the Heat Equation, the Fast Diffusion Equation,
the $p$-Laplacian equation, the Stefan problem, and many others. The fractional equivalents need specific proofs,
 as the one provided in \cite{VazqJEMS2014}\nc.

 This principle is quite  useful when applied to the study
of the monotonicity properties of the solutions with compactly
supported data.  Here is the corresponding {\sl
Monotonicity Lemma} along outgoing directions, as stated in Lemma 14.25 of \cite{VazquezPME2007}:

\begin{lemma} \label{lem.ap0} Let $u\ge 0$ be a solution of the Cauchy
problem for the HE or the PME with initial data
supported in the ball $B_R(0)$, $R>0$. Then for every $x_0\in
\ren$ such that $|x_0|>R$ and every $t>0$, $u(x,t)$ is monotone
nonincreasing along the ray $l(x_0)=\{x= sx_0: s\ge 1\}$ in the
sense that
\begin{equation}\label{eq.ad1}
u(s_2x_0,t)\le u(s_1x_0,t) \quad \mbox{if} \,\,\, s_2\ge s_1\ge 1.
\end{equation}
\end{lemma}

\noindent {\em Proof.} The application of Alexandrov's Reflection
Principle proceeds as follows: we draw the hyperplane $H$ which is
mediatrix between the points $x=s_2x_0$ and $y=s_1x_0$ in the above
situation. It is easy to see that $H$ divides the space $\re^N$ into
two half-spaces, one $\Omega_1$ which contains $y=s_1x_0$ and the support of
$u_0$, and another one, $\Omega_2$, which contains $x=s_2x_0$ and where
$u_0=0$. We consider now the initial and boundary-value problem in
$\widehat Q= \Omega_1\times (0,T)$. Two particular solutions of
this problem are compared: one of them is $u_1$, the restriction of
$u$ to $\widehat Q$, another one is
 $$
u_2(z,t)= u(\pi(z),t), \quad z\in \Omega_1,
 $$
where $\pi$ is the specular symmetry with respect to the
hyperplane $H$. By the Reflection Principle,
 $$
u_1(z,t)\ge u_2(z,t) \quad \mbox{for} \,\,\, z\in \Omega_1, \ t>0.
 $$
Putting $z=y$ we have $\pi(z)=x$ so that $u(y,t)\ge u(x,t)$ as
desired.\qed

\medskip

Actually,  a stronger version  of the result was given in the early paper by Aronson-Caffarelli, \cite{AC83}, Proposition 2.1, that was very influential in the sequel. We give an adapted statement

\begin{lemma}\label{le.AleksAC} Let $u$ a nonnegative solution of the Cauchy problem for the HE or the PME with initial data supported in the ball $B_R(0)$. Then,
$$
\inf_{\{|x|\ge r\}} u(x,t)> \max_{\{|x|=r+2R\}} u(x, t)
$$
for all $r,t>0$.
\end{lemma}

These results will be used in the delicate uniqueness proofs performed in Section \ref{sec.crit.exp} and following, starting with Theorem
\ref{thm.mc.cple}\nc.

\vskip1.5cm

\noindent {\textbf{\large \sc Acknowledgments.}} Author partially funded by Projects  PGC2018-098440-B-I00 and PID2021-127105NB-I00 from MICINN (Spain). Partially performed as an Honorary Professor at Univ. Complutense de Madrid. The author thanks Matteo Bonforte, David G\'omez Castro,  Félix del Teso,  and other colleagues for their interest in this work and useful suggestions. David and Félix  supplied graphic material. \bc The current version has been significantly improved thanks to the referees' attention to detail and some very insightful suggestions.\nc

\vskip0.7cm

\newpage


\

\


\noindent {\sc Address:}

\noindent Juan Luis V\'azquez. Departamento de Matem\'{a}ticas, \\
Universidad Aut\'{o}noma de Madrid,\\ Campus de Cantoblanco, 28049 Madrid, Spain.  \\
e-mail address:~\texttt{juanluis.vazquez@uam.es}\\
webpage:  https://verso.mat.uam.es/\~juanluis.vazquez/

\

\end{document}